\documentclass[11pt,reqno]{amsart}
\usepackage{amsmath}
\usepackage{amsfonts}
\usepackage{amssymb, amsthm, epsfig}
\usepackage{mathrsfs}
\usepackage{hyperref}
 \usepackage{yhmath}
\usepackage{cite}
\usepackage[utf8]{inputenc}
\usepackage[T1]{fontenc}
\usepackage{geometry}
\usepackage{graphicx}
\usepackage{verbatim}
\usepackage{color}

\theoremstyle{plain}
\newtheorem{thm}{Theorem}[section]

\newtheorem{lem}[thm]{Lemma}
\newtheorem{prop}[thm]{Proposition}

\theoremstyle{definition}
\newtheorem{defn}{Definition}[section]
\theoremstyle{remark}
\newtheorem{rem}{Remark}[section]

\numberwithin{equation}{section}%{section}

\allowdisplaybreaks

\DeclareMathOperator{\arcctg}{arccot}

\DeclareSymbolFont{lettersA}{U}{pxmia}{m}{it}
\DeclareMathSymbol{\piup}{\mathord}{lettersA}{"19}

\makeatletter

\newcommand{\Rmnum}[1]{\expandafter\@slowromancap\romannumeral#1@}
\makeatother
\begin{document}
\title[]
{Wave interactions for the two-dimensional self-similar Euler equations of nonconvex gases}

\author{Geng Lai}

\address{Department of Mathematics, Shanghai University,
Shanghai, 200444, P.R. China;  Newtouch Center for Mathematics, Shanghai University, Shanghai, 20044, P.R. China.}
\email{\tt laigeng@shu.edu.cn}

\author{Wancheng Sheng}
\address{Department of Mathematics, Shanghai University,
Shanghai, 200444, P.R. China;  Newtouch Center for Mathematics, Shanghai University, Shanghai, 20044, P.R. China.}
\email{\tt wcsheng@shu.edu.cn}

%\subjclass[]{35L65, 35L60, 35L67.}
\date{\today}

\begin{abstract}
Gases with nonconvex equations of state may differ significantly from polytropic ideal gases; for example, rarefaction and sonic shocks may occur. This paper studies a 2D Riemann initial-boundary value problem for the potential flow equations with a nonconvex equation of state. We construct global-in-time self-similar solutions by analyzing eight distinct types of wave interactions, including composite waves consisting of a sonic/double-sonic shock coupled with one or two fan waves. Three main difficulties arise in studying these wave interactions: the self-similar potential flow system is of mixed type, with its type at each point determined by the local pseudo-Mach number; the sonic or transonic nature of the shocks is not known a priori; and sonic singularities arise on sonic shocks. To overcome these difficulties, we proceed as follows. First, we introduce a density-dependent invariant region and apply the maximum principle to establish the hyperbolicity principle of the 2D self-similar Euler equations. Second, we establish a priori estimates for the shock types by differentiating the Rankine–Hugoniot conditions along the shocks, comparing the shock and acoustic speeds, and using the Oleinik-Liu entropy condition. Finally, we employ the hodograph transformation together with regularization to handle the various singularities on sonic shocks. These wave interactions serve as building blocks for more complicated 2D Riemann problems and as benchmarks for numerical algorithms for nonconvex conservation laws. We also propose several open problems.
\end{abstract}

\maketitle
\tableofcontents

\section{Introduction}
The global existence of entropy solutions to the Cauchy problem for the multidimensional compressible Euler equations remains a challenging open problem. It is therefore beneficial to study special cases such as two-dimensional (2D) Riemann problems, which are Cauchy problems whose initial data are constant along each ray emanating from the origin. 2D Riemann problems are important in theoretical analysis, scientific computation, and engineering applications. Solving them can elucidate the physical structure of waves and their interactions, and can also help improve numerical algorithms; see \cite{Ben,CCY1,CCY2,Glimm,LaxLiu,LZY,MP,Rinne1,Rinne2,ZhangandZheng,Zheng1} and references therein.

Several classes of two-dimensional (2D) Riemann initial value and initial-boundary value problems for the compressible Euler equations and related models have been studied. Sheng and Zhang \cite{Sheng3} constructed global piecewise smooth solutions with delta waves for 2D Riemann problems of the transportation equations. Dai and Zhang \cite{Dai} studied gas expansion into a vacuum for the pressure gradient system and established the existence of a global classical solution. Li \cite{LJ1} applied the hodograph transformation to obtain a solution in the hodograph plane for the expansion-into-vacuum problem of the Euler equations. Subsequently, Li and Zheng \cite{LJ4} established a global-in-time solution in the physical plane via the characteristic decomposition method proposed by Li, Zhang and Zheng \cite{LJ3} for 2D pseudo-steady simple waves. Li and Zheng \cite{LJ5} also obtained a global solution to a 2D Riemann problem with four bi-symmetric rarefaction waves. Further results on rarefaction wave interactions in 2D Riemann problems can be found in \cite{CZheng,Hu,LJ2,Zhao,Lai4,Lai2,Lai5}.

In recent decades, many important results have been obtained for 2D Riemann problems involving transonic shocks. Canic, Keyfitz and Kim \cite{Canic} established the existence of a solution to the regular reflection of a planar shock off a solid wall for the unsteady transonic small-disturbance model. Zheng \cite{Zheng} studied the regular reflection problem for the pressure gradient system, transforming it into a free boundary problem with an oblique derivative boundary condition for a degenerate elliptic equation, and established the existence of a global solution. Chen and Feldman \cite{ChenF1} established a global solution to the regular reflection problem for the potential flow equation. Bae, Chen and Feldman \cite{BCF} obtained the optimal regularity of the solution across the sonic circle, while Chen, Feldman and Xiang \cite{ChenF3} established the convexity of the transonic shock. Elling and Liu \cite{Elling1} proved an ellipticity principle guaranteeing that if the 2D self-similar potential flow equation is elliptic on the boundary of a bounded region, then it is uniformly elliptic in the interior. They subsequently used this principle and topological degree theory to construct a global solution to supersonic flow onto a solid wedge \cite{Elling}. Serre \cite{Serre} investigated several configurations of the 2D Riemann problem for a Chaplygin gas and constructed the corresponding global piecewise smooth solutions. For further advances on 2D Riemann problems for the potential flow equations, we refer the reader to \cite{BCF,Canic1,CGQ,CCHLQ,ChenF2,CX,ChenF3,ChenF4,CQ,Elling11,Kim} and references therein.

As far as we know, while many results are available for 2D Riemann problems of polytropic ideal gases, much less is known for nonconvex gases. This paper is concerned with a 2D Riemann problem for the potential flow equations of nonconvex gases.
The 2D potential flow equations consist of the conservation law of mass and Bernoulli's law for the density and velocity potential $(\rho, \Phi)$:
\begin{equation}\label{Potential}
\left\{
            \begin{array}{ll}
            \rho_t+(\rho \Phi_x)_x+(\rho \Phi_y)_y=0, \\[4pt]
         \displaystyle \Phi_t+\frac{1}{2}\big(\Phi_x^2+\Phi_y^2\big)+h(\tau)=0;
            \end{array}
       \right.
\end{equation}
see, e.g., Chen and Feldman \cite{ChenF1}.
Here, $\rho$ is the density;
$\Phi$ is the velocity potential, i.e., $(\Phi_x, \Phi_y)=(u, v)$; $(u, v)$ is the flow velocity; $h$ is the enthalpy, which satisfies
\begin{equation}\label{han}
h'(\tau)=\tau p'(\tau);
\end{equation}
$p$ is the pressure; $\tau$ is the specific volume. By the second equation of (\ref{Potential}) and equation (\ref{han}), the density $\rho$ can be expressed as a function of $\Phi_t$ and $\nabla\Phi$, i.e.,
 \begin{equation}\label{rho}
\rho=\rho(\Phi_t, |\nabla\Phi|^2).
\end{equation}

\begin{figure}[htbp]
\begin{center}
\includegraphics[scale=0.55]{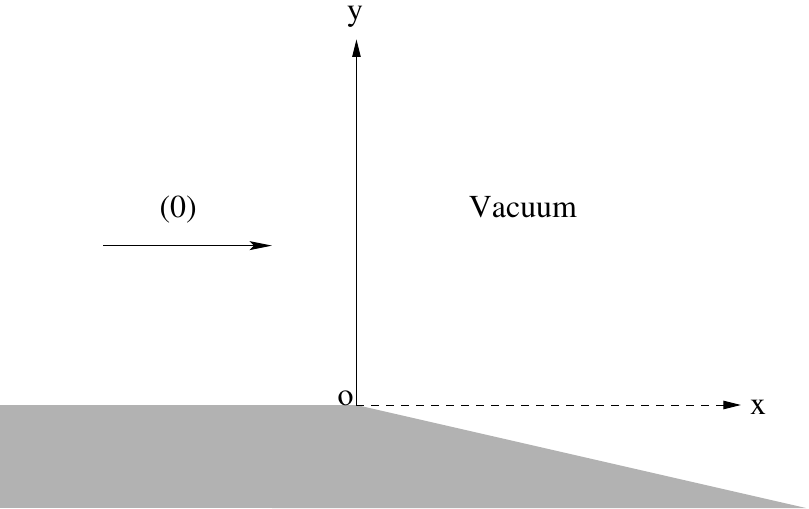}
\caption{\footnotesize Riemann initial-boundary value problem.}
\label{Initial-data}
\end{center}
\end{figure}

In this paper, we study a
 2D pseudo-steady supersonic flow past a rarefactive ramp for the 2D potential flow equations (\ref{Potential}) for nonconvex gases. 
Two-dimensional supersonic flows past compressible and rarefactive ramps
are fundamental in gas dynamics, as they give rise to global or local flow structures in various settings; see, e.g., Courant and Friedrichs \cite{CF}.
Assume that there is an infinitely long wedge with a straight horizontal wall, a sharp corner located at the origin $\mathrm{O}$, and a straight ramp; see Figure \ref{Initial-data}. The wall of the wedge is given by
 $$
 \partial \mathcal{W}=\big\{(x, y)\mid y=0, x<0\big\}\cup\big\{(x, y)\mid (x, y)=r(\cos(\theta_w-\pi), \sin(\theta_w-\pi)), r>0\big\},
 $$
 where $\theta_w\in (0, \pi)$.
At the initial time, there is a uniform oncoming supersonic flow with constant velocity $(u_0, 0)$ and constant density $\rho_0$ in the second quadrant, and the remaining region above the wedge is a vacuum\footnote{The method in this paper also applies to the case where the right initial state is a constant state free of vacuum; see Section 12.}. The oncoming flow will then ‌expand‌ into the vacuum and ‌be deflected‌ at the corner $\mathrm{O}$.
To investigate the flow structure near the corner‌,
 we consider (\ref{Potential}) with the initial data
 \begin{equation}\label{initial}
(\rho, \Phi)(0,x,y)=\left\{
              \begin{array}{ll}
                (\rho_0, u_0x), & \hbox{$x<0$, $y>0$;} \\[2pt]
               \text{vacuum}, & \hbox{$x>0$, $y\geq 0$;} \\[2pt]
  \text{vacuum}, & \hbox{$y<0$, $x>y\cot\theta_w$,}
              \end{array}
            \right.
 \end{equation}
and the slip condition along the wedge boundary $\partial \mathcal{W}$:
 \begin{equation}\label{boundary}
(\rho\Phi_x, \rho\Phi_y)\cdot {\bf n}\big|_{\mathbb{R}^{+}\times \partial\mathcal{W}}=0,
 \end{equation}
 where ${\bf n}$ is the exterior unit normal to $\partial \mathcal{W}$.
The problem (\ref{Potential}), (\ref{initial})--(\ref{boundary}) is in fact a 2D Riemann initial-boundary value problem and also captures  several important types of elementary wave interactions.

 %======================
%==========

\begin{figure}[htbp]
\begin{center}
\qquad\includegraphics[scale=0.4]{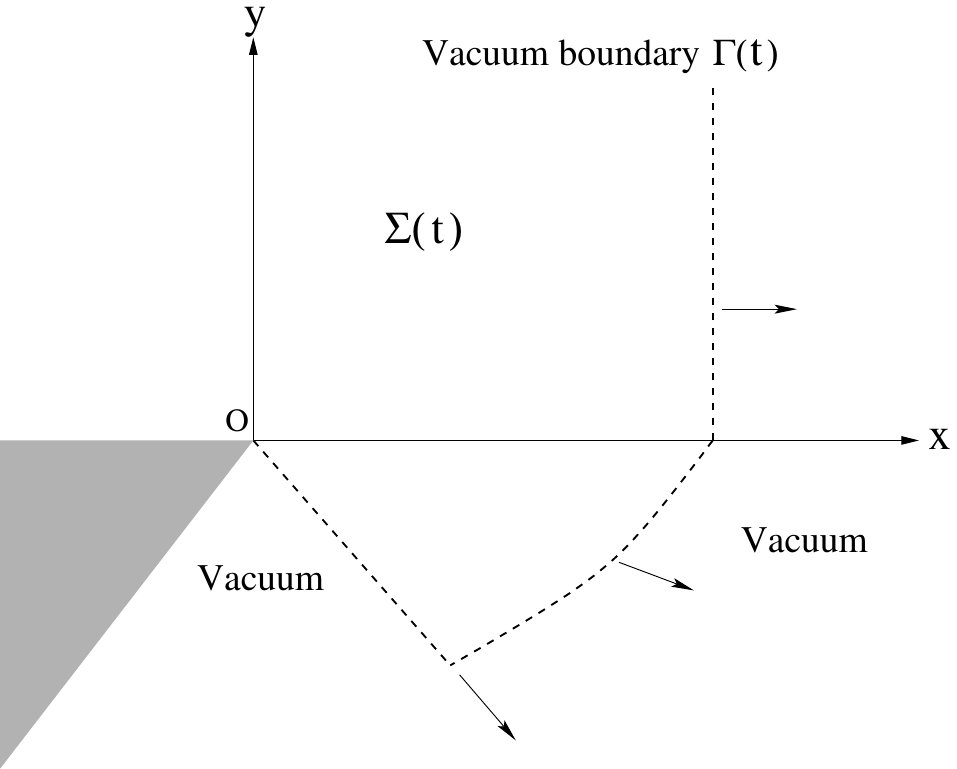}\qquad \quad\includegraphics[scale=0.4]{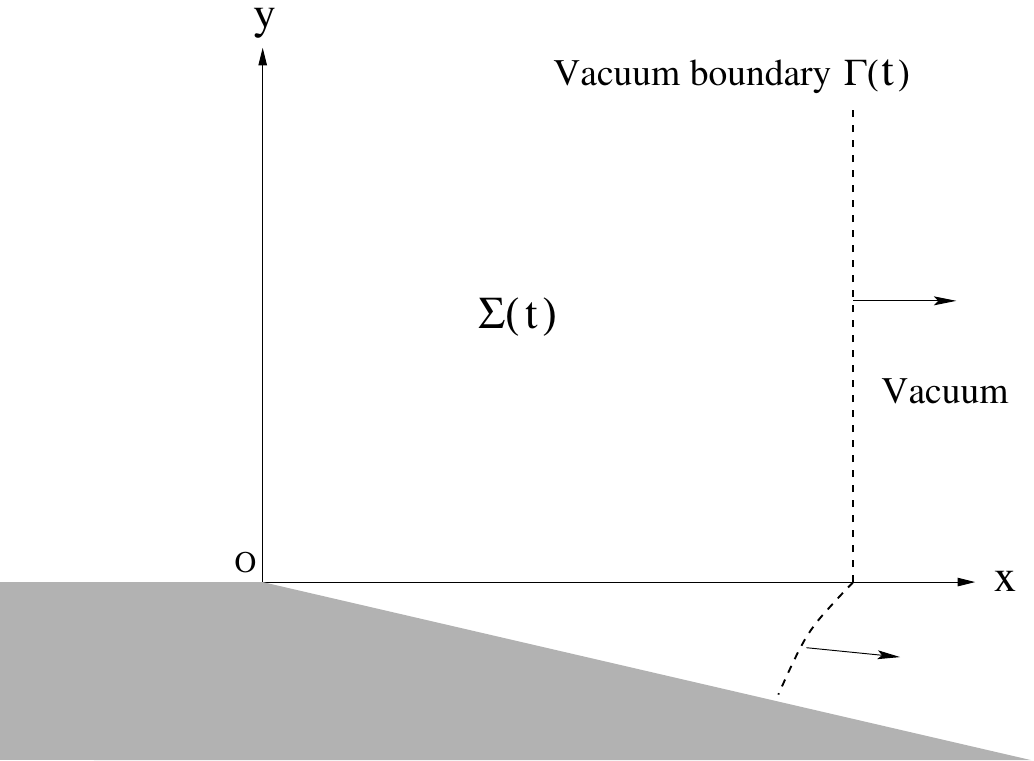}
\caption{\footnotesize Domain $\Sigma(t)$ and the gas-vacuum interface $\Gamma(t)$.}
\label{Vcummboundary}
\end{center}
\end{figure}

For the 2D Riemann problem (\ref{Potential}), (\ref{initial})--(\ref{boundary}), the flow will expand into the vacuum and then form a time-dependent vacuum boundary $\mathbb{R}_{+}\times\Gamma(t)$. Let $\Sigma(t)$ be the domain enclosed by $\partial \mathcal{W}(t)$ and $\Gamma(t)$, where $\partial \mathcal{W}(t)$  is a part of $\partial \mathcal{W}$; see Figure \ref{Vcummboundary}.
For any $t>0$, the flow is a vacuum outside the region $\Sigma(t)$.
The solution to the problem (\ref{Potential}), (\ref{initial})--(\ref{boundary}) may be discontinuous. Therefore, we need to consider weak solutions, possibly piecewise smooth. We first give the definition of a weak solution to the problem (\ref{Potential}), (\ref{initial})--(\ref{boundary}).
\begin{defn}
(Weak solution of  Problem (\ref{Potential}), (\ref{initial})--(\ref{boundary})) A function $\Phi\in W_{loc}^{1,1}(\mathbb{R}_{+}\times \Sigma(t))$
is called a weak solution to the problem (\ref{Potential}), (\ref{initial})--(\ref{boundary}) if $\Phi$ satisfies the following properties:
\begin{itemize}\label{def1}
 \item $-\Phi_t-\frac{1}{2}\big(\Phi_x^2+\Phi_y^2\big)\geq h_{\infty}$ a.e. in $\mathbb{R}_{+}\times \Sigma(t)$, where $h_{\infty}=\lim\limits_{\tau\rightarrow+\infty}h(\tau)$;
  \item For $\rho(\Phi_t, |\nabla\Phi|^2)$ defined in (\ref{rho}),
$$
\Big(\rho(\Phi_t, |\nabla\Phi|^2), \rho(\Phi_t, |\nabla_{(x, y)}\Phi|^2)|\nabla\Phi|\Big)\in \Big(L_{loc}^1([0,+\infty)\times \overline{\Sigma(t)})\Big)^2;
$$
\item The density $\rho(\Phi_t, |\nabla\Phi|^2)$ and $\nabla\Phi$ are continuous a.e. in $\mathbb{R}_{+}\times \Gamma(t)$, and
$$
\rho=0\quad \text{and}\quad (\Phi_x, \Phi_y, 1)\cdot {\overrightarrow{n}}=0 \quad\text{a.e. in} \quad \mathbb{R}_{+}\times \Gamma(t),
$$
where $\overrightarrow{n}$ represents the normal vector of the vacuum surface $\mathbb{R}_{+}\times \Gamma(t)$;
  \item For every $\Psi\in C_{c}^{\infty}([0,+\infty)\times \mathbb{R}^2)$,
$$
\iiint\limits_{\mathbb{R}_{+}\times \Sigma(t)} \rho(\Phi_t, |\nabla\Phi|^2)\Psi_t+
\rho(\Phi_t, |\nabla\Phi|^2)\nabla\Phi\cdot\nabla\Psi{\rm d}x{\rm d}y{\rm d}t+
\iint\limits_{\Sigma(0)}\rho_0\Psi(0, x,y){\rm d}x{\rm d}y=0.
$$
\end{itemize}
\end{defn}

 %======================
%==========

We intend to construct a piecewise smooth weak solution to  (\ref{Potential}) with jumps for $(\Phi_t, \nabla\Phi)$ across a surface $\mathcal{S}$ with unit normal $\overrightarrow{n}=(n_t, n_x, n_y)$. The requirement of the weak solution of
(\ref{Potential}) in the sense of Definition \ref{def1} yields the Rankine-Hugoniot jump condition across $\mathcal{S}$:
\begin{equation}
[\Phi]=0\quad \text{and}\quad
[\rho]n_t+[\rho\Phi_x]n_x+[\rho\Phi_y]n_y=0,
\end{equation}
where the square brackets denote the jump in the corresponding quantity across $\mathcal{S}$. The discontinuity surface $\mathcal{S}$ is called an admissible shock if the specific volume $\tau$ across $\mathcal{S}$ satisfies the following Oleinik-Liu entropy condition:
$$
-\frac{2h(\tau_f)-2h(\tau_b)}{\tau_f^2-\tau_b^2}\geq -\frac{2h(\tau_f)-2h(\tau)}{\tau_f^2-\tau^2}\quad \text{for all}~ \tau\in\big(\min\{\tau_f, \tau_b\}, \max\{\tau_f, \tau_b\}\big),
$$
where $\tau_f$ and $\tau_b$ denote the specific volumes upstream and downstream of $\mathcal{S}$, respectively.

\vskip 4pt

The 2D Riemann problem (\ref{Potential}), (\ref{initial})--(\ref{boundary}) admits self-similar solutions of the following form:
\begin{equation}
\rho=\rho(\xi, \eta)\quad \text{and}\quad  \Phi=t\vartheta(\xi, \eta),
\end{equation}
where $(\xi, \eta)=(x/t, y/t)$.
Then, the pseudo-potential function $\varphi=\vartheta-\frac{1}{2}(\xi^2+\eta^2)$ satisfies the following system:
\begin{equation}\label{42501}
\left\{
  \begin{array}{ll}
   (\rho \varphi_\xi)_\xi+(\rho \varphi_\eta)_\eta+2\rho=0,\\[4pt]
  \displaystyle\frac{1}{2}\big(\varphi_\xi^2+\varphi_\eta^2\big)+h(\tau)+\varphi =0.
  \end{array}
\right.
\end{equation}

Let $(U, V)=(u-\xi, v-\eta)$ be the pseudo-velocity. Then 
$$
\varphi_\xi=U\quad \text{and} \quad \varphi_\eta=V.
$$
The 2D pseudo-steady Euler equations for potential flow read
\begin{equation}\label{PEC}
\left\{
  \begin{array}{ll}
   (\rho U)_\xi+(\rho V)_\eta+2\rho=0,\\[4pt]
  \displaystyle  U_\eta-V_\xi=0
  \end{array}
\right.
\end{equation}
supplemented by Bernoulli's law in the $(\xi, \eta)$-plane
\begin{equation}\label{5802}
\frac{1}{2}q^2+h(\tau)+\varphi=0,
\end{equation}
where $q=\sqrt{U^2+V^2}$.

Using (\ref{5802}), the first equation of (\ref{PEC}) can be written as
\begin{equation}
  (c^{2}-U^{2})u_{\xi}-UV(u_{\eta}+v_{\xi})+(c^{2}-V^{2})v_{\eta}=0,\label{mass}
\end{equation}
where $c=\sqrt{-\tau^2 p'(\tau)}$ is the speed of sound.
Therefore, for smooth flow, system (\ref{PEC}) can be written as the following matrix form:
\begin{equation} \label{matrix1}
\begin{array}{rcl}
\left(
 \begin{array}{cc}
c^{2}-U^{2} & -UV \\
  0 & -1\\
  \end{array}
  \right)\left(
           \begin{array}{c}
             u \\
             v \\
           \end{array}
         \right)_{\xi}+\left(
                         \begin{array}{ccccc}
                          -UV &  c^{2}-V^{2}\\
                           1 & 0 \\
                         \end{array}
                       \right)\left(
                                \begin{array}{c}
                                  u \\
                                  v \\
                                \end{array}
                              \right)_{\eta}=\left(
                                               \begin{array}{c}
                                                 0 \\
                                                0 \\
                                               \end{array}
                                             \right).
                                             \end{array}
\end{equation}

The eigenvalues of (\ref{matrix1}) are determined
by
\begin{equation}\label{21}
(V-\lambda U)^{2}-c^{2}(1+\lambda^{2})=0,
\end{equation}
which gives
\begin{equation}\label{23}
\lambda=\lambda_{\pm}=\frac{UV\pm
c\sqrt{U^{2}+V^{2}-c^{2}}}{U^{2}-c^{2}}.
\end{equation}
Therefore, system (\ref{PEC}) is a mixed-type system and the type at each point is determined by the local pseudo-Mach number $q/c$. It is hyperbolic if and only if $q/c>1$ (supersonic) and elliptic if and only if $q/c<1$ (subsonic).
For supersonic flow, system (\ref{PEC}) has two families of wave characteristic curves defined as the integral curves of
\begin{equation}\label{24}
C_{\pm}:~\frac{{\rm d}\eta}{{\rm d}\xi}=\lambda_\pm.
\end{equation}

%\begin{defn}
%(Weak solution of {\bf Problem} (\ref{Potential}), (\ref{initial})--(\ref{boundary})) A function $\varphi\in W_{loc}^{1,1}(\Sigma)$
%is called a weak solution of the expansion in vacuum problem (\ref{Potential}), (\ref{initial})--(\ref{boundary})) if $\varphi$ satisfies the following properties:
%\begin{description}\label{def2}
 % \item[(i)] $-\varphi_t-\frac{1}{2}\big(\varphi_{\xi}^2+\varphi_{\eta}^2\big)\geq h_{\infty}$ a. e. in $\Sigma$, where $h_{\infty}=\lim\limits_{\tau\rightarrow+\infty}h(\tau)$;
 % \item[(ii)] For $\rho(\varphi, |\nabla\varphi|^2)$ defined in,
%$$
%\Big(\rho(\varphi, |\nabla\varphi|^2), \rho(\varphi, |\nabla\varphi|^2)|\nabla\varphi|\Big)\in \Big(L_{loc}^1(\overline{\Sigma}\Big)^2;
%$$
%\item[(iii)]  The density $\rho(\varphi, |\nabla\varphi|^2)$ and $\nabla\varphi$ are continuous a.e. in $\Sigma$, and
%$$
%\rho=0\quad \text{and}\quad (\varphi_\xi, \varphi_\eta)=0\quad\text{a. e. in}\quad \Gamma;
%$$
%where $\bf{n}$ represents the normal vector of the vacuum surface $R_{+}\times \Gamma(t)$;
 % \item[(iv)] For every $\psi\in C_{c}^{\infty}(R^2)$,
%$$
%\iiint\limits_{\Sigma}
%\rho\nabla\varphi\cdot\nabla\psi-2\rho\psi{\rm d}\xi{\rm d}\eta=0.
%$$
%\end{description}
%\end{defn}

 %======================
 \begin{figure}[htbp]
\begin{center}
\qquad\includegraphics[scale=0.38]{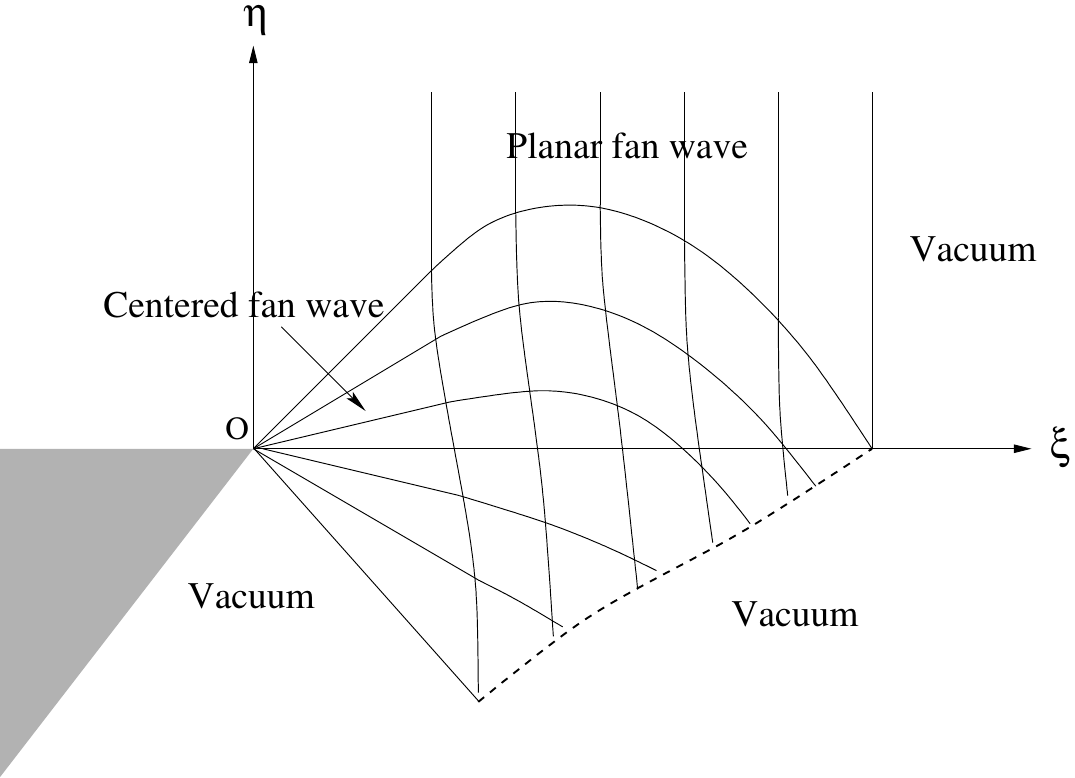}\qquad \qquad \includegraphics[scale=0.38]{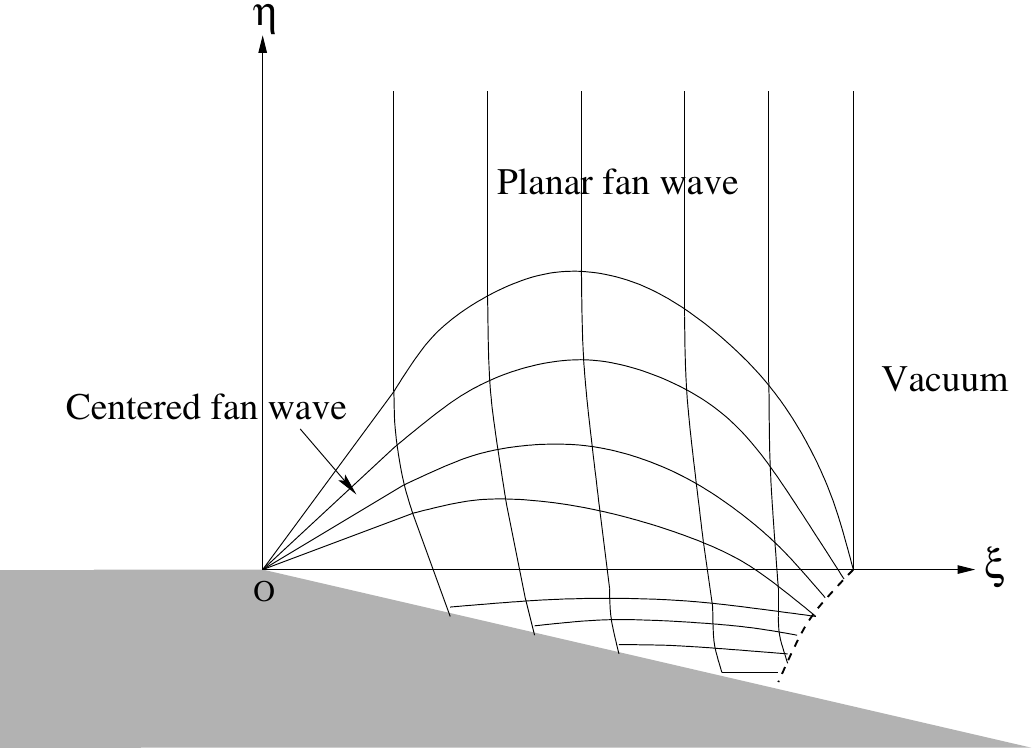}
\caption{\footnotesize Interaction of fan waves. Left: $\theta_w<\theta_v$; right: $\theta_w>\theta_v$.}
\label{Polytropic}
\end{center}
\end{figure}
%==========

When the oncoming fluid is a polytropic ideal gas, away from the corner the oncoming flow and vacuum are connected by a planar fan wave, and the oncoming flow is deflected by a centered simple wave to either a constant state or vacuum at the corner. Then the problem can be reduced to the interaction between the planar fan wave and the centered simple wave, together with the reflection of a simple wave off the solid ramp. The wave configuration in the $(\xi, \eta)$-plane is illustrated in Figure \ref{Polytropic}.
Sheng  et al. \cite{Sheng1,Sheng2} established the existence of a global classical solution  provided that $u_0>c_0$ and the wedge angle $\theta_w$ lies in a certain range, where $c_0=\sqrt{-\tau_0^2p'(\tau_0)}$.
Subsequently, Lai and Sheng \cite{Lai6}
 established the existence of a global classical solution to the problem (\ref{Potential}), (\ref{initial})--(\ref{boundary}) in the case $u_0\leq c_0$, in particular for $u_0=0$.

For ideal gases, the isentropes are always convex downward in the $(\tau, p)$-plane.
A natural question is whether the above results
can be extended to more general gases, for example, gases with nonconvex isentropes or non-ideal gases.
The possibility that isentropes are nonconvex was first explored independently by Bethe \cite{Bethe} and Zel'dovich \cite{ZE} for van der Waals gases, and by Thompson and Lambrakis \cite{Thompson2} for gases with high molecular complexity. Thompson \cite{Thompson1}
introduced a dimensionless thermodynamic quantity
$$\Gamma=\frac{\tau^3}{2c^2}\Big(\frac{\partial^2 p}{\partial \tau^2}\Big)_{s}$$
called the fundamental derivative of gasdynamics to
delineate the dynamic behavior of compressible fluids, where $p$ is the pressure, $\tau$ the specific volume, $s$ the specific entropy and $c$ the speed of sound.
Fluids with $\Gamma<0$ within a certain thermodynamic region are now commonly referred to as Bethe-Zel'dovich-Thompson (BZT) fluids.
This type of fluid may differ significantly from polytropic ideal gases. For example, rarefaction shocks and sonic shocks may occur; see e.g., \cite{BBKN,Cra,NSMGC,Thompson2,Weyl,ZGC}. For more results on BZT fluids, we refer the reader to the survey papers \cite{GCPS,K}.

%=============
\begin{figure}[htbp]
\begin{center}
\includegraphics[scale=0.39]{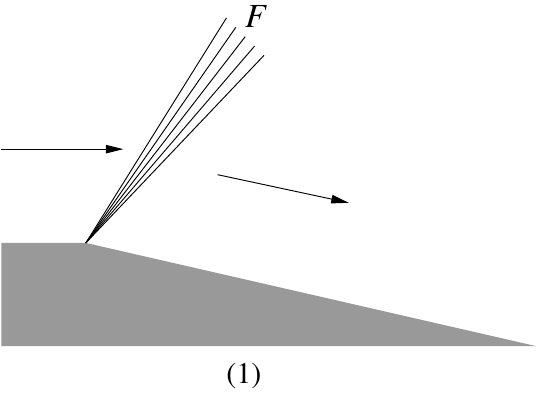}\qquad\includegraphics[scale=0.39]{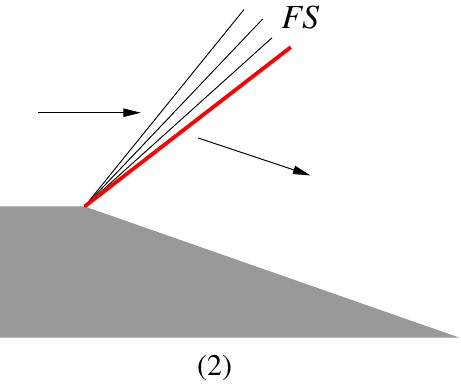}\qquad \includegraphics[scale=0.39]{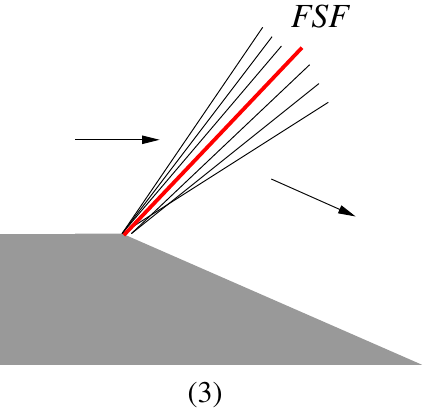} \\
\includegraphics[scale=0.39]{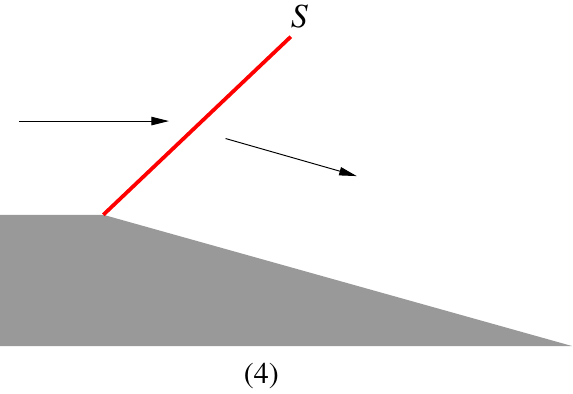}\qquad\includegraphics[scale=0.39]{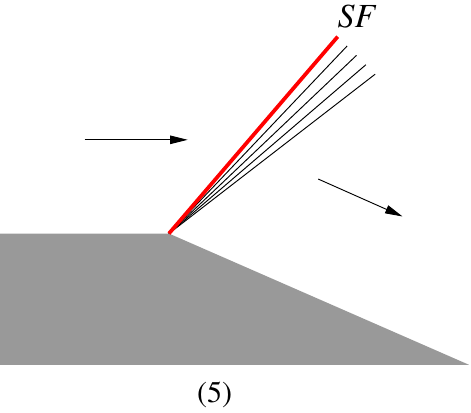}
\caption{\footnotesize 2D steady supersonic flows of BZT fluids past a rarefactive ramp.}
\label{Supersonicramp}
\end{center}
\end{figure}
%==========

It is well known that the one-dimensional (1D) Riemann problem plays a key role in understanding the wave structure of the compressible Euler equations.
Oleinik \cite{Oleinik} introduced an entropy condition that ensures the uniqueness of discontinuous (weak) solutions to scalar nonlinear conservation laws, even when the flux function is non-convex.
Chang and Hsiao \cite{ZH} studied the 1D Riemann problem for the scalar conservation law with a nonconvex flux function.
Wendroff \cite{Wen1,Wen2} studied the 1D Riemann problem for the Euler equations with a nonconvex equation of state and constructed composite wave solutions composed of shock waves and simple waves.
 Liu \cite{Liu1,Liu2} introduced an extended entropy condition and established the existence and uniqueness of entropy solutions to the 1D Riemann problem.
Godlewski and Seguin \cite{GS} studied the 1D Riemann problem for an isothermal p-system with a van der Waals equation of state complemented by the Maxwell's law.
We also refer the reader to the survey paper by Menikoff and Plohr
 \cite{MP} and  the  monograph by LeFloch and Thanh \cite{LeFloch2} for the 1D Riemann problem of more general gases.
 LeFloch \cite{LeFloch1} also developed  the mathematical theory of classical and nonclassical shock waves for the Cauchy problem in hyperbolic conservation laws.

Vimercati, Kluwick and Guardone \cite{VKG1} constructed self-similar solutions to 2D steady supersonic flows of BZT fluids past compressible and rarefactive ramps, a problem treated as a 1D Riemann problem.
For the BZT fluids, many nonclassical wave configurations can occur in 2D steady supersonic ramp flows.
 At an rarefactive ramp, the wave configuration  may be a fan (${\it F}$) wave, a fan–shock (${\it FS}$) composite wave, a fan–shock–fan (${\it FSF}$) composite wave, a shock (${\it S}$) wave, or a shock–fan (${\it SF}$) composite wave; see Figure \ref{Supersonicramp}. 
Vimercati, Kluwick and Guardone \cite{VKG2} also studied the morphology of nodes formed by shock interactions (including crossing, overtaking, splitting and refraction) in 2D steady flows of BZT fluids.
Li and Sheng \cite{Li-Sheng1,Li-Sheng2} and Wang {et al.} \cite{Wang1,Wang2} constructed self-similar solutions to the supersonic ramp problem for van der Waals gases.
More recently, Lai \cite{Laip2} studied the stabilities of the oblique self-similar fan-shock-fan and shock-fan-shock composite waves in 2D steady supersonic flows of BZT fluids.

\begin{figure}[htbp]
\begin{center}
 \includegraphics[scale=0.48]{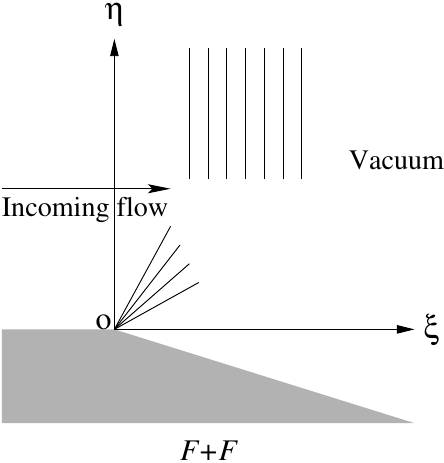}\qquad\qquad \includegraphics[scale=0.48]{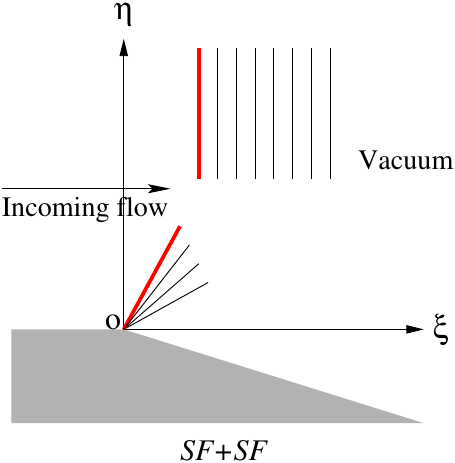}\qquad \qquad \includegraphics[scale=0.48]{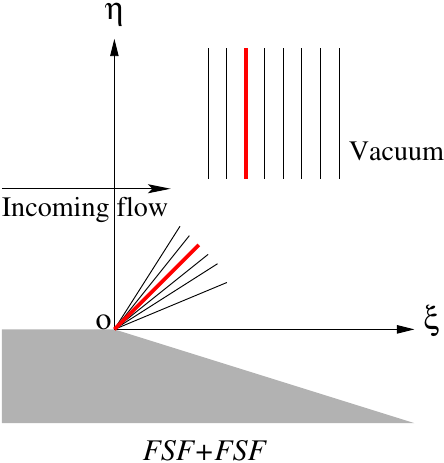}\\[2pt]
\vskip 6pt\includegraphics[scale=0.48]{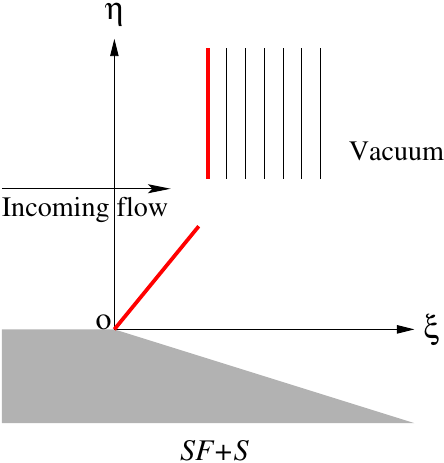}\qquad\qquad
 \includegraphics[scale=0.48]{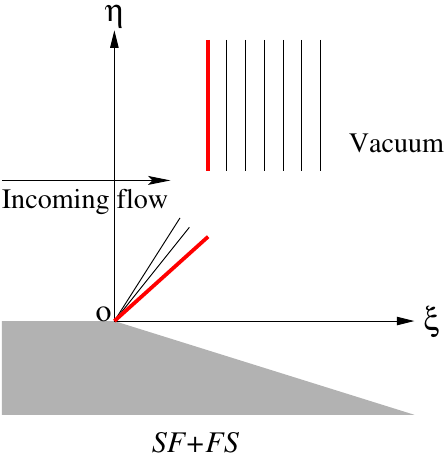}\qquad\qquad \includegraphics[scale=0.48]{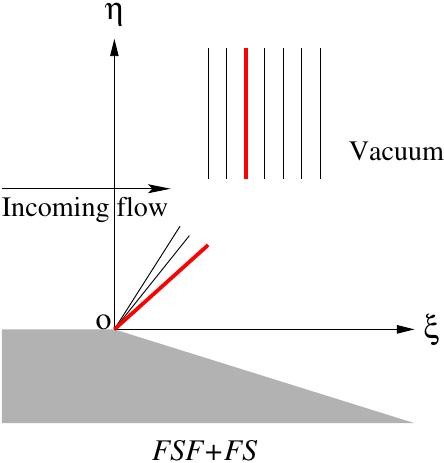}\\[2pt]
\vskip 6pt \includegraphics[scale=0.48]{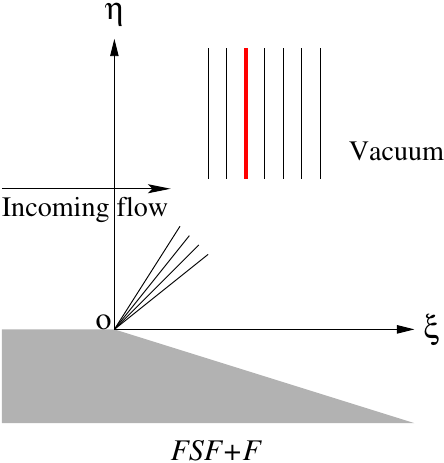}\qquad\qquad \includegraphics[scale=0.48]{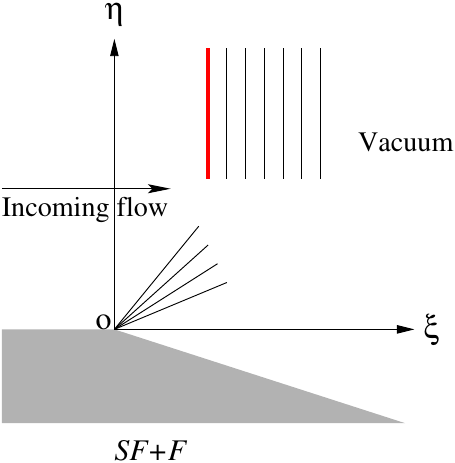}
\caption{\footnotesize Wave combinations in the 2D Riemann problem (\ref{Potential}), (\ref{initial})--(\ref{boundary}).}
\label{Combination}
\end{center}
\end{figure}

%In a recent paper \cite{VKG1},, Kluwick $\&$ Guardone constructed self-similar solutions to the 2D steady supersonic flow of BZT fluids past compressible and rarefactive ramps.
%For the BZT fluids, many additional wave configurations may possibly occur in the 2D steady supersonic ramp flows.
In this paper, we study the 2D Riemann problem (\ref{Potential}), (\ref{initial})--(\ref{boundary}) with a nonconvex equation of state $p=p(\tau)$ having two inflection points.
Assume that there exist $\tau_1^{i}$ and $\tau_2^{i}$ ($0<\tau_1^{i}<\tau_2^{i}$) such that the equation of state
$p=p(\tau)$ satisfies the following conditions:
\begin{equation}\label{81601}
\begin{aligned}
& \quad p''(\tau_1^i)=p''(\tau_2^i)=0; \quad p'(\tau)<0~\text{for}~\tau>0;\quad p''(\tau)<0~\text{for}~\tau_1^i<\tau<\tau_2^i;
\\[2pt]&  p''(\tau)>0~\text{for}~\tau\in(0, \tau_1^i)\cup(\tau_2^i, +\infty);\quad \lim\limits_{\tau\rightarrow 0}p'(\tau)=-\infty; \quad \lim\limits_{\tau\rightarrow +\infty}p'(\tau)=0.
 \end{aligned}
\end{equation}
Then, in the self-similar $(\xi, \eta)$-plane, the oblique waves emanating from the corner $\mathrm{O}$ may be a centered simple wave ${\it F}$, a ${\it FS}$ composite wave, a ${\it FSF}$ composite wave, a shock wave, or a ${\it SF}$ composite wave.
Moreover, the oncoming flow $(u_0, 0, \rho_0)$ and the vacuum can be connected by a planar fan wave, a ${\it SF}$ composite wave, or a ${\it FSF}$ composite wave.
For gases satisfying the assumptions in (\ref{81601}),
 the following eight distinct elementary wave interactions arise in the 2D Riemann problem (\ref{Potential}), (\ref{initial})--(\ref{boundary}) (see Figure \ref{Combination}):
\begin{description}
\item[(1)] interaction between two fan waves;
\item[(2)] interaction between two ${\it SF}$ composite waves;
\item[(3)] interaction between two ${\it FSF}$ composite waves;
\item[(4)] interaction between a  ${\it SF}$ composite wave and a shock wave;
\item[(5)] interaction between  ${\it SF}$ and ${\it FS}$ composite waves;
\item[(6)] interaction between ${\it FSF}$ and ${\it FS}$ composite waves;
\item[(7)] interaction between a ${\it FSF}$ composite wave and a fan wave;
\item[(8)] interaction between a ${\it SF}$ composite wave and a fan wave.
\end{description}

%=============
%==========

We aim to systematically study these wave interactions and construct global-in-time piecewise smooth supersonic solutions to the 2D Riemann problem (\ref{Potential}), (\ref{initial})--(\ref{boundary}).
To this end, we solve several boundary value problems for the system (\ref{PEC}), including Goursat problems, discontinuous Goursat problems, sonic shock free-boundary problems,  transonic shock free-boundary problems, and sonic-transonic shock free-boundary problems.
In contrast to ideal gases, for a nonconvex gas, the flow downstream (or upstream) of a shock may be sonic relative to the shock front.
Moreover, in certain wave interactions, the shock type cannot be determined a priori, so that the corresponding boundary conditions cannot be formulated in advance.
We establish a priori estimates for the shock types by comparing the shock and acoustic speeds along the shocks and by using Liu's extended entropy condition.
We also find that sonic shocks, which are envelopes of one family of wave characteristics, appear in certain oblique wave interactions.
This results in a fact that the flow downstream (or upstream) of a sonic shock is not $C^1$ smooth up to the shock boundary.
We use the hodograph transformation method to overcome the difficulty caused by the singularity on sonic shocks.

The rest of the paper is organized as follows.

In Section 2, we introduce the characteristic decomposition and the hodograph transformation for system (\ref{42501}). The notions of characteristic directions and characteristic angles $\alpha$ and $\beta$ are introduced in Section 2.1. In Section 2.2, we derive a set of characteristic equations for the variables $\alpha$, $\beta$, and $\rho$; these equations will be used extensively in establishing a priori estimates for the hyperbolicity of system (\ref{42501}). In Section 2.3, we derive several sets of characteristic decompositions for the variables $\rho$ and $c$; these decompositions will play a central role in controlling the bounds of the derivatives of the solutions. In Section 2.4, we present the hodograph transformation for system (\ref{42501}) and derive a set of characteristic decompositions in the hodograph plane. These decompositions are then used to establish a monotonicity property of the hodograph transformation, which in turn implies that the transformation is one-to-one.

In Section 3, we study the planar elementary wave connecting the oncoming flow with vacuum, as well as the oblique wave at the corner. We also identify eight distinct types of wave combinations in the 2D Riemann problem (\ref{Potential}), (\ref{initial})--(\ref{boundary}).

In Section 4, we study the interaction between two fan waves, which arise in  the 2D Riemann problem (\ref{Potential}), (\ref{initial})--(\ref{boundary}) for  $\tau_0>\tau_2^i$. We reduce the interaction problem to a standard Goursat problem and thereby construct a global-in-time classical solution to the 2D Riemann problem by the method of characteristics, provided that the ramp angle and the equation of state satisfy certain conditions. 
The key point for the global existence is to control the hyperbolicity of  system (\ref{Potential}) and the $C^1$ norm of the solution. To this end, 
we introduce a density-dependent invariant region for the characteristic angles to control the hyperbolicity of system (\ref{Potential}) in the interaction region, and thereby apply the characteristic decomposition method to control the bounds of the first-order derivatives of the solution. 

In Section 5, we study the interaction between two $SF$ composite waves, which arise  when $\tau_1^e<\tau_0<\tau_2^i$ and $\theta_w$ satisfies a certain condition. 
We reduce this interaction problem to a discontinuous Goursat-type problem (DGP), with the discontinuity located at the corner point of the DGP. 
Using phase plane analysis, we show that there are two centered waves emanating from the corner. We then apply the hodograph transformation to construct a local classical solution to the DGP. The solution in the remaining part of the interaction region is constructed as in Section 4.

In Section 6, we study the interaction between two $FSF$ composite waves, which arise when  $\tau_0<\tau_1^i$  and $\theta_w$ satisfies a certain condition. 
The main difficulty in analyzing this wave interaction is that the interaction region contains shocks. Moreover, the type of the shocks is not known a priori. By calculating the curvatures of the shocks and using Liu's extended entropy condition, we prove that the shocks in the interaction region must be post-sonic. We also prove that the shocks are envelopes of one of the two families of wave characteristics of the flow downstream of them, rather than characteristics themselves. This implies that the flow behind the shocks is not $C^1$-smooth up to the shock boundaries.
Hence, to construct the flow downstream of the sonic shocks, we need to solve two singular Cauchy problems. To overcome the difficulty caused by the singularity on the sonic shocks, we apply the hodograph transformation to transform these singular problems into regular Cauchy problems in the hodograph plane, and then show that this transformation can be inverted back to the physical plane. The solution in the remaining part of the interaction region is constructed as in Section 5.

In Section 7, we study the interaction between a $SF$ composite waves and an oblique shock wave, which arise  when  $\tau_1^e<\tau_0<\tau_2^i$  and $\theta_w$ satisfies a certain condition. This wave interaction can be reduced to a DGP with a straight $C_{-}$ boundary and a curved $C_{+}$ boundary.  We show that the interaction generates either a non-planar simple wave or a non-planar $SF$ composite wave, depending on the specific volume of the flow behind the oblique shock.
This non-planar simple wave or $SF$ composite wave will reflect at the ramp. We find that two types of reflection may occur: regular-like reflection and Mach-like reflection. In the case of regular-like reflection, a global solution to the reflection problem can be constructed by solving a Goursat or a discontinuous Goursat problem as in Sections 4 and 5, respectively. However, for Mach-like reflection, the method presented in this paper is not applicable; moreover, even the reflection configuration cannot be easily conjectured; we leave it open. 

In Section 8, we study the interaction between a $SF$ composite waves and a $FS$ composite wave, which arise  when  $\tau_1^e<\tau_0<\tau_1^i$  and $\theta_w$ satisfies a certain condition. The interaction also generates either a non-planar simple wave or a non-planar $SF$ composite wave.
The main contribution of this section is to characterize the propagation of a post-sonic shock in a rarefaction simple wave of a different characteristic family. 
We prove that when the shock and the simple wave belong to different characteristic families (e.g., a 1-shock and a 2-rarefaction wave), the post-sonic shock crosses the entire rarefaction simple wave and remains post-sonic. 

In Section 9, we study the interaction between a $FSF$ composite wave and a $FS$ composite wave, which arise  when  $\tau_0<\tau_1^e$  and $\theta_w$ satisfies a certain condition. The interaction also generates either a non-planar simple wave or a non-planar $SF$ composite wave.
The main contribution of this section is to characteristic the propagation of a pre-sonic shock in an expansion flow region. 
We prove that when a pre-sonic shock enters an expansion flow region, it becomes transonic. In some cases, this transonic shock subsequently becomes post-sonic. We present a transition case from a pre-sonic shock to a transonic shock and then to a post-sonic shock.

In Section 10, we study the interaction between a $FSF$ composite wave and a fan wave, which arise  when  $\tau_0<\tau_1^e$  and $\theta_w$ satisfies a certain condition.
The main contribution of this section is to characterize the transition from a post-sonic shock to a transonic shock, which arises in the interaction of a post-sonic shock with a rarefaction simple wave of the same characteristic family.
%We prove that when the shock and the simple wave belong to the same characteristic family (e.g., a 1-shock and a 1-rarefaction wave), the post-sonic shock becomes transonic as it enters the simple wave region. 
To this end, we consider (\ref{Potential}) in an angular domain between a shock free-boundary and a given characteristic boundary of another characteristic family. One main difficulty of this free-boundary problem is that the shock direction and the characteristic direction coincide at the corner of the angular domain, since the shock is post-sonic at the corner. Another difficulty is that the data on the characteristic boundary does not have a $C^1$ bound near the corner. We establish a delicate estimate of the solution near the corner and construct a local Lipschitz-continuous solution to the free boundary problem, and prove that the shock becomes transonic after leaving the corner.

In Section 11, we study the interaction between a $SF$ composite wave and a fan wave, which arise  when $\tau_1^e<\tau_0<\tau_1^i$ and $\theta_w$ satisfies a certain condition. Using the discussion in Section 8, we show that the interaction generates a non-planar shock-fan composite wave. 

In Section 12, we present, without rigorous proof, a brief discussion of the case in which the right initial state is a vacuum-free constant state for polytropic gases. 
%Some additional remarks, including the conclusions of this paper, vacuum-free Riemann initial data,  and several  open problems related to the 2D Riemann problem (\ref{Potential}), (\ref{initial})--(\ref{boundary}), are presented in Section 12. 

\section{Preliminaries}

\subsection{Characteristic angles and directions}
\begin{figure}[htbp]
\begin{center}
\includegraphics[scale=0.6]{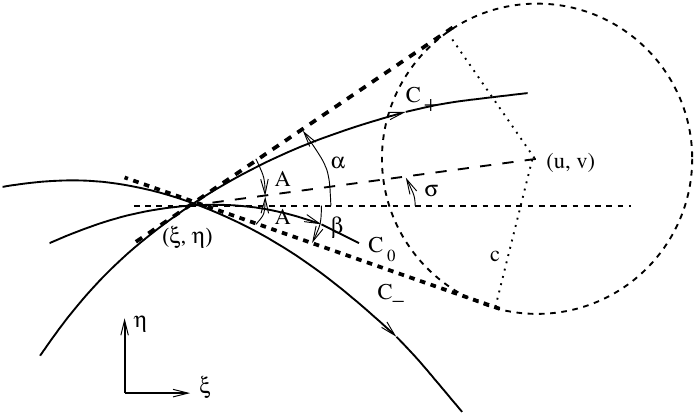}\qquad\qquad\includegraphics[scale=0.42]{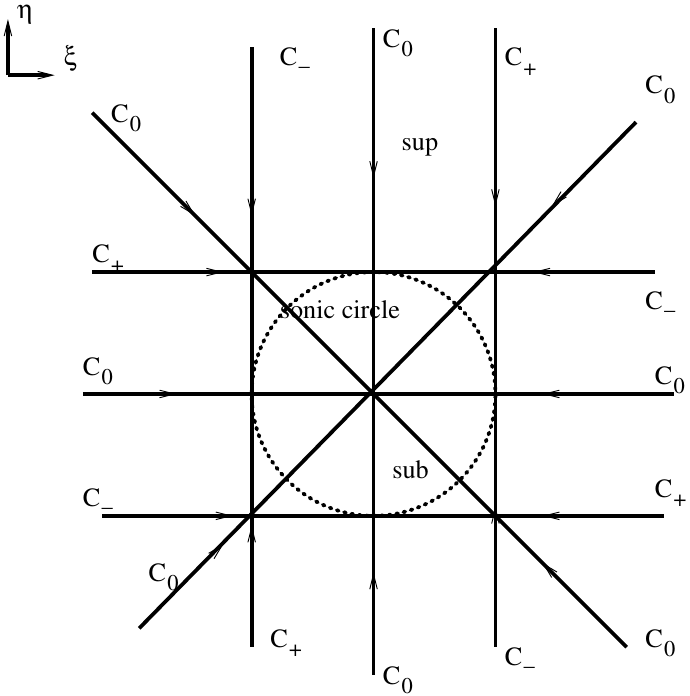}
\caption{\footnotesize Left: characteristic directions and characteristic angles; right: $C_{\pm}$ and $C_{0}$ characteristic lines for  a constant state.}
\label{Fig11}
\end{center}
\end{figure}

In this paper we shall confine ourselves to supersonic flow. The direction of a wave
characteristics is defined as the tangent direction that makes an
acute angle $A$ with the pseudoflow velocity
$(U, V)$. A direct computation gives
\begin{equation}
c^{2}=q^{2}\sin^{2}A.\label{210cqo}
\end{equation}
The angle $A$ is called the
 pseudo-Mach angle. For pseudosonic flow,
for which $q/c=1$, the direction of the
pseudoflow velocity
 is perpendicular to the directions of the wave characteristics.
A direct computation shows  that
the $C_+$
characteristic direction forms with the pseudoflow velocity
$(U, V)$ the angle $A$ from $C_{+}$ to  $(U, V)$ in
the clockwise direction, and the $C_-$ characteristic direction forms with the
pseudoflow direction the angle $A$ from $C_{-}$ to $(U, V)$
in the counterclockwise direction; see Figure \ref{Fig11} (left).

From (\ref{21}) we have
\begin{equation}\label{c}
c=\frac{|(U,V)\cdot(\lambda,-1)|}{|(\lambda,-1)|},
\end{equation}
which implies that the component of the pseudoflow velocity normal to the
direction of a $C_{+}$ ($C_{-}$) characteristic is equal to the sound speed.
 Equivalently, it can be stated that the tangent line of a $C_{+}$ ($C_{-}$)
characteristic at a point is tangent to the sonic circle of the
state at that point. Furthermore,  a $C_{+}$ ($C_{-}$) characteristic must be straight
if $(u, v, c)$ is constant along it; see Figure \ref{Fig11} (right).

 The $C_{+}$ ($C_{-}$) characteristic angle is defined as the angle measured counterclockwise from the positive $\xi$-axis to  the $C_{+}$ ($C_{-}$) characteristic direction.
We denote by $\alpha$ and
 $\beta$ the counterclockwise angles from the positive $\xi$-axis to the $C_{+}$ and $C_{-}$ characteristic directions, respectively, i.e.,
\begin{equation}\label{6803}
\lambda_{+}=\tan\alpha\quad \text{and}\quad \lambda_{-}=\tan\beta.
\end{equation}
We remark that $0<\alpha-\beta<\pi$ for supersonic flow.

%=============
%==========

 Let $\sigma$ be the angle measured counterclockwise from the positive $\xi$-axis to the pseudoflow velocity, i.e.,
 \begin{equation}\label{121407}
 \sigma=\arctan\left(\frac{V}{U}\right).
 \end{equation}
Hence,
\begin{equation}
\alpha=\sigma+A,\quad \beta=\sigma-A, \quad\sigma=\frac{\alpha+\beta}{2},\quad A=\frac{\alpha-\beta}{2}. \label{tau}
\end{equation}
Therefore, the relations between $(U, V, c)$ and $(\alpha, \beta, c)$ are given by
\begin{equation} \label{U}
U=\frac{c\cos\sigma}{\sin A}=q\cos\sigma,\quad V=\frac{c\sin\sigma}{\sin A}=q\sin\sigma.
\end{equation}
The mapping $(U, V, c)\rightarrow (\alpha, \beta, c)$ is bijective as long as system (\ref{PEC}) is hyperbolic.

\subsection{Characteristic equations}
We first introduce the following variables for later use:
\begin{equation}
\begin{aligned}
&\kappa(\tau):=\frac{-2p'(\tau)}{\tau p''(\tau)+2p'(\tau)}, \quad \omega(\tau):=-\frac{4p'(\tau)+\tau p''(\tau)}{\tau p''(\tau)},\quad \varpi(\tau):=\frac{2p'(\tau)+\tau p''(\tau)}{\tau p''(\tau)},\\[4pt]
&\qquad\qquad  \Omega(\tau, A):=~\omega(\tau)-\tan^2A,\quad \text{and}\quad  \Theta(\tau)=\arctan\sqrt{\omega(\tau)}.
\end{aligned}
\end{equation}

Multiplying (\ref{matrix1}) on the left by
$(1,\mp c\sqrt{U^{2}+V^{2}-c^{2}})$, we get the characteristic equations:
\begin{equation}
\left\{
  \begin{array}{ll}
  \displaystyle \bar{\partial}_{+}u+\lambda_{-}\bar{\partial}_{+}v
  =0,  \\[4pt]
    \displaystyle  \bar{\partial}_{-}u+\lambda_{+}\bar{\partial}_{-}v=0,
  \end{array}
\right.\label{form}
\end{equation}
where
\begin{equation}\label{3303}
\bar{\partial}_{+}=\cos\alpha\partial_{\xi}+\sin\alpha\partial_{\eta}\quad \text{and}\quad \bar{\partial}_{-}=\cos\beta\partial_{\xi}+\sin\beta\partial_{\eta}
\end{equation}
are the directional derivatives along the $C_{+}$ and $C_{-}$ characteristic directions, respectively.

From (\ref{U}) we have
\begin{equation}
\bar{\partial}_{\pm}u=\cos(\sigma\pm A)+\frac{\cos\sigma}{\sin A}\bar{\partial}_{\pm}c+\frac{c\cos\alpha\bar{\partial}_{\pm}\beta
-c\cos\beta\bar{\partial}_{\pm}\alpha}{2\sin^{2}A},\label{1}
\end{equation}
\begin{equation}
\bar{\partial}_{\pm}v=\sin(\sigma\pm A)+\frac{\sin\sigma}{\sin A}\bar{\partial}_{\pm}c
+\frac{c\sin\alpha\bar{\partial}_{\pm}\beta
-c\sin\beta\bar{\partial}_{\pm}\alpha}{2\sin^{2}A}.\label{2}
\end{equation}
Inserting (\ref{1}) and (\ref{2}) into (\ref{form}), we obtain
\begin{equation}\label{3}
\bar{\partial}_{+}c=-\frac{\cos2A}{\cot A}+\frac{c}{\sin2A}
(\bar{\partial}_{+}\alpha-\cos2A\bar{\partial}_{+}\beta),
\end{equation}
\begin{equation}
\bar{\partial}_{-}c=-\frac{\cos2A}{\cot A}+\frac{c}{\sin2A}
(\cos2A\bar{\partial}_{-}\alpha-\bar{\partial}_{-}\beta).\label{4}
\end{equation}

Differentiating Bernoulli's law (\ref{5802}) and
using (\ref{1}) and (\ref{2}), we get
\begin{equation}
\left(\frac{1}{\sin^{2}A}+\kappa\right)\bar{\partial}_{+}c=
\frac{c\cos A}{2\sin^{3}A}(\bar{\partial}_{+}\alpha-\bar{\partial}_{+}\beta)-\cot A.\label{bB1}
\end{equation}
Similarly, we have
\begin{equation}
\left(\frac{1}{\sin^{2}A}+\kappa\right)\bar{\partial}_{-}c=
\frac{c\cos A}{2\sin^{3}A}(\bar{\partial}_{-}\alpha-\bar{\partial}_{-}\beta)-\cot A.\label{bB2}
\end{equation}

Inserting (\ref{bB1}) into (\ref{3}), we obtain
\begin{equation}\label{6}
c\bar{\partial}_{+}\alpha=\Omega\cos^{2}A(c\bar{\partial}_{+}\beta+2\sin^{2}A).
\end{equation}
Similarly, inserting (\ref{bB2}) into (\ref{4}), we obtain
\begin{equation}\label{5}
c\bar{\partial}_{-}\beta=\Omega\cos^{2}A(c\bar{\partial}_{-}\alpha-2\sin^{2}A).
\end{equation}

Combining (\ref{3}) with (\ref{6}), we have
\begin{equation}
c\bar{\partial}_{+}\beta=-\frac{\tan A}{\varpi}\bar{\partial}_{+}c-2\sin^{2} A=-\frac{p''(\tau)}{2c\rho^4}\tan A\bar{\partial}_{+}\rho-2\sin^{2} A,\label{7}
\end{equation}
\begin{equation}
c\bar{\partial}_{+}\alpha=-\left(\frac{1+\kappa}{2}\right)\Omega\sin2 A\bar{\partial}_{+}c=-\frac{p''(\tau)}{4c\rho^4}\Omega\sin2 A\bar{\partial}_{+}\rho.\label{10}
\end{equation}
Combining (\ref{4}) with (\ref{5}), we have
\begin{equation}\label{82301}
c\bar{\partial}_{-}\alpha=\frac{\tan A}{\varpi}\bar{\partial}_{-}c+2\sin^{2} A=\frac{p''(\tau)}{2c\rho^4}\tan A\bar{\partial}_{-}\rho+2\sin^{2} A,
\end{equation}
\begin{equation}
c\bar{\partial}_{-}\beta=\left(\frac{1+\kappa}{2}\right) \Omega\sin2 A\bar{\partial}_{-}c=\frac{p''(\tau)}{4c\rho^4}\Omega\sin2 A\bar{\partial}_{-}\rho.\label{8}
\end{equation}

 By (\ref{7}) and (\ref{10}) we have
\begin{equation}
c\bar{\partial}_{\pm}A=\sin^2  A \left(1+\frac{p''(\tau)}{4c\rho^{4}}\big(\tan A-\Omega\sin A\cos A\big)\frac{\bar{\partial}_{\pm}\rho}{\sin^2  A}\right).
\end{equation}

From (\ref{1}), (\ref{2}), and (\ref{7})--(\ref{8}) we have
\begin{equation}\label{11}
\bar{\partial}_{+}u=\kappa\sin\beta\bar{\partial}_{+}c=c\tau\sin\beta\bar{\partial}_{+}\rho,\quad\bar{\partial}_{-}u=-\kappa\sin\alpha\bar{\partial}_{-}c=-c\tau\sin\alpha\bar{\partial}_{-}\rho,
\end{equation}
\begin{equation}\label{72804}
\bar{\partial}_{+}v=-\kappa\cos\beta\bar{\partial}_{+}c=-c\tau\cos\beta\bar{\partial}_{+}\rho,
\quad\bar{\partial}_{-}v=\kappa\cos\alpha\bar{\partial}_{-}c=c\tau\cos\alpha\bar{\partial}_{-}\rho.
\end{equation}

\subsection{Characteristic decompositions}
In this paper, we use the method of characteristic decomposition to construct global classical solutions to the 2D Riemann problem (\ref{Potential}), (\ref{initial})--(\ref{boundary}). This method was first introduced by Li, Zhang and Zheng \cite{LJ3} to study simple waves of the 2D pseudo-steady Euler equations, and has been widely used to study the global existence of smooth solutions to the compressible Euler equations in two independent variables.

\begin{prop}
(Commutator relation)  We have
\begin{equation}
\bar{\partial}_{-} \bar{\partial}_{+}- \bar{\partial}_{+} \bar{\partial}_{-}=
\frac{1}{\sin2 A}\Big[(\cos2 A \bar{\partial}_{+}\beta- \bar{\partial}_{-}\alpha) \bar{\partial}_{-}-
(\bar{\partial}_{+}\beta-\cos2 A \bar{\partial}_{-}\alpha) \bar{\partial}_{+}\Big].
\label{comm}
\end{equation}
\end{prop}
\begin{proof}
This commutator relation was first derived by Li, Zhang and Zheng \cite{LJ3}.
From (\ref{3303}) we have
\begin{equation}\label{61513}
\partial_{\xi}=-\frac{\sin\beta\bar{\partial}_{+}-\sin\alpha\bar{\partial}_{-}}{\sin2A}\quad \text{and}\quad
\partial_{\eta}=\frac{\cos\beta\bar{\partial}_{+}-\cos\alpha\bar{\partial}_{-}}{\sin2A}.
\end{equation}
Then the commutator relation can be obtained by a direct computation.
\end{proof}

\begin{prop}
For the variable $c$, we have the following characteristic decompositions:
\begin{equation}\label{cd}
\left\{
  \begin{array}{ll}
    \displaystyle c\bar{\partial}_{+}\bar{\partial}_{-}c=\left(\sin2 A+\frac{\bar{\partial}_{-}c}{2\varpi\cos^2 A}
+
\Big(\frac{\kappa\cos^2 2A+1}{2\cos^2 A}-\frac{c}{\kappa}\cdot\frac{{\rm d}\kappa}{{\rm d}c}\Big)\bar{\partial}_{+}c\right)\bar{\partial}_{-}c, \\[18pt]
   \displaystyle c\bar{\partial}_{-}\bar{\partial}_{+}c=\left(\sin2 A+\frac{\bar{\partial}_{+}c}{2\varpi\cos^2 A}+
\Big(\frac{\kappa\cos^2 2A+1}{2\cos^2 A}-\frac{c}{\kappa}\cdot\frac{{\rm d}\kappa}{{\rm d}c}\Big)\bar{\partial}_{-}c\right)\bar{\partial}_{+}c.
  \end{array}
\right.
\end{equation}
\end{prop}

\begin{proof}
%These decompositions are due to Li and Yang and Zheng \cite{Li-Yang-Zheng}. For the sake of completeness, we sketch the proof.
It follows from (\ref{11}) and (\ref{comm}) that
\begin{equation}
\begin{array}{rcl}
&&\bar{\partial}_{-} \bar{\partial}_{+}u- \bar{\partial}_{+} \bar{\partial}_{-}u\\[6pt]&=&%\\[4pt]&=&
\bar{\partial}_{+}\big(\kappa\sin\alpha\bar{\partial}_{-}c\big)+
\bar{\partial}_{-}\big(\kappa\sin\beta\bar{\partial}_{+}c\big)\\[4pt]
&=&-\displaystyle\frac{1}{\sin2 A}\Big((\bar{\partial}_{+}\beta-\cos2 A \bar{\partial}_{-}\alpha)\kappa\sin\beta \bar{\partial}_{+}c-( \bar{\partial}_{-}\alpha- \cos2 A\bar{\partial}_{+}\beta)\kappa\sin\alpha \bar{\partial}_{-}c\Big).
\end{array}
\end{equation}
Hence we have
\begin{equation}
\begin{array}{lll}
\displaystyle&&\sin\alpha\bar{\partial}_{+}\bar{\partial}_{-}c+
\sin\beta\bar{\partial}_{-}\bar{\partial}_{+}c+\displaystyle\frac{1}{\kappa}\cdot\frac{{\rm d}\kappa}{{\rm d}c}(\sin\alpha+\sin\beta)\bar{\partial}_{+}c\bar{\partial}_{-}c
\\[4pt]&=&-\displaystyle\frac{1}{\sin2 A}
\Big((\sin\beta \bar{\partial}_{+}\beta-\sin\beta \cos2 A \bar{\partial}_{-}\alpha+\cos\beta\sin2 A\bar{\partial}_{-}\beta) \bar{\partial}_{+}c\\[8pt]&&\qquad\displaystyle-( \sin\alpha\bar{\partial}_{-}\alpha- \sin\alpha\cos2 A\bar{\partial}_{+}\beta-\cos\alpha\sin2 A\bar{\partial}_{+}\alpha) \bar{\partial}_{-}c\Big).
\end{array}\label{12}
\end{equation}

Applying the commutator relation (\ref{comm}) for $c$, we get
\begin{equation}
\bar{\partial}_{-} \bar{\partial}_{+}c- \bar{\partial}_{+} \bar{\partial}_{-}c=
\frac{1}{\sin2 A}\Big((\cos2 A \bar{\partial}_{+}\beta- \bar{\partial}_{-}\alpha) \bar{\partial}_{-}c-
(\bar{\partial}_{+}\beta-\cos2 A \bar{\partial}_{-}\alpha) \bar{\partial}_{+}c\Big).\label{13}
\end{equation}
Inserting (\ref{13}) into (\ref{12}) we get
\begin{equation}
\begin{array}{rcl}
&&\displaystyle(\sin\alpha+\sin\beta)\bar{\partial}_{+}\bar{\partial}_{-}c
\\[4pt]&=&-\displaystyle\frac{1}{\sin2 A}
\Big(\cos\beta\sin2 A\bar{\partial}_{-}\beta\bar{\partial}_{+}c+( \sin\alpha+\sin\beta)\cos2 A\bar{\partial}_{+}\beta\bar{\partial}_{-}c\\[6pt]&&
\qquad-(\sin\alpha+\sin\beta)\bar{\partial}_{-}\alpha\bar{\partial}_{-}c+\cos\alpha\sin2 A\bar{\partial}_{+}
\alpha\bar{\partial}_{-}c\Big)-\displaystyle\frac{1}{\kappa}\cdot\frac{{\rm d}\kappa}{{\rm d}c}(\sin\alpha+\sin\beta)\bar{\partial}_{+}c\bar{\partial}_{-}c.
\end{array}\label{14}
\end{equation}
Thus, by (\ref{7})--(\ref{8}) we obtain
 the first equation of (\ref{cd}).

The proof for the other decomposition is similar. This completes the proof.
\end{proof}

\begin{figure}[htbp]
\begin{center}
\includegraphics[scale=0.45]{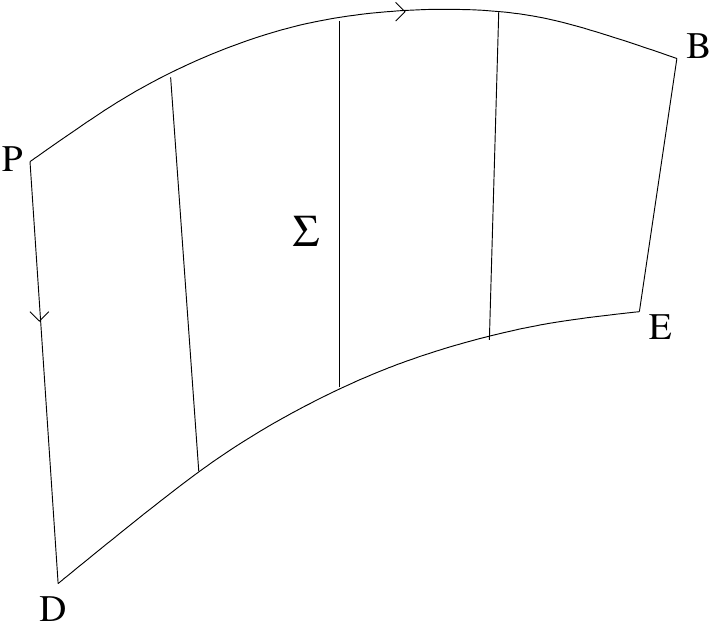}
\caption{\footnotesize Simple wave problem.}
\label{Simple}
\end{center}
\end{figure}

As an application of the characteristic decomposition, we consider the following simple wave problem, which will serve as a building block for solutions of the 2D Riemann problem.
Let $\widetilde{\mathrm{PB}}$ be a smooth curve and let $\overline{\mathrm{PD}}$ be a straight line; see Figure \ref{Simple}. We prescribe the boundary condition
\begin{equation}\label{Simplebd}
(u, v, \tau)=\left\{
               \begin{array}{ll}
                 (u_{+}, v_{+}, \tau_{+})(\xi, \eta), & \hbox{$(\xi, \eta)\in \widetilde{\mathrm{PB}}$;} \\[2pt]
                 (u_{_\mathrm{P}}, v_{_\mathrm{P}}, \tau_{_\mathrm{P}}), & \hbox{$(\xi, \eta)\in \overline{\mathrm{PD}}$,}
               \end{array}
             \right.
\end{equation}
where $(u_{_\mathrm{P}}, v_{_\mathrm{P}}, \tau_{_\mathrm{P}})=(u_{+}, v_{+}, \tau_{+})(\xi_{_\mathrm{P}}, \eta_{_\mathrm{P}})$,
such that $\widetilde{\mathrm{PB}}$ is a $C_{+}$ characteristic curve oriented from $\mathrm{P}$ to $\mathrm{B}$, and
$\overline{\mathrm{PD}}$ is a straight $C_{-}$ characteristic line oriented from $\mathrm{P}$ to $\mathrm{D}$.
Moreover, the first  characteristic equation of (\ref{form}) holds on $\widetilde{\mathrm{PB}}$.
Then problem (\ref{PEC}), (\ref{Simplebd}) is a standard Goursat problem.

%=============
%==========

\begin{lem}\label{simple}
Assume that the boundary data (\ref{Simplebd}) satisfy $\tau_{+}(\mathrm{P})>\tau_2^i$ and $\bar{\partial}_{+}\tau>0$ along $\widetilde{\mathrm{PB}}$. Then the Goursat problem (\ref{PEC}), (\ref{Simplebd}) admits a simple wave solution in a curvilinear quadrilateral domain $\Sigma$ bounded by $\widetilde{\mathrm{PB}}$, $\overline{\mathrm{PD}}$, $\widetilde{\mathrm{DE}}$, and $\overline{\mathrm{BE}}$, where $\overline{\mathrm{BE}}$ is a straight $C_{-}$ characteristic line emanating from $\mathrm{B}$ and $\widetilde{\mathrm{DE}}$ is a forward $C_{+}$ characteristic curve emanating from $\mathrm{D}$. Moreover, the solution satisfies $\bar{\partial}_{+}\rho<0$ on $\Sigma$.
\end{lem}

\begin{proof}
From (\ref{11}), (\ref{72804}), and (\ref{cd}) we have that
for the 2D pseudo-steady potential flow equations (\ref{PEC}), the flow in a region adjacent to a straight characteristic line must be either a simple wave or a constant state (see \cite{LJ3}). Thus, it remains only to prove that the Goursat problem admits a global classical solution. From the second equation of (\ref{cd}) we see that the forward straight $C_{-}$ characteristic lines emanating from $\widetilde{\mathrm{PB}}$ do not intersect before reaching their respective sonic points. From the point $\mathrm{P}$ we draw a cross $C_{+}$ characteristic curve of the simple wave, and this curve intersects with the straight $C_{-}$ characteristic line emanating from $\mathrm{B}$ at a point $\mathrm{E}$. This completes the proof.
\end{proof}

\begin{prop} We have the following characteristic decompositions:
\begin{equation}
\left\{
  \begin{array}{ll}
\begin{aligned}
c\bar{\partial}_{+}\left(\frac{\bar{\partial}_{-}c}{\sin^{2} A}\right)
   &=\bar{\partial}_{-}c\Bigg(\frac{1}{2\varpi\cos^{2} A}\Big(\frac{\bar{\partial}_{-}c}{\sin^2A}-\frac{\bar{\partial}_{+}c}{\sin^2A}\Big)
\\&\qquad\qquad\qquad+\left(\frac{1}{\varpi\cos^{2} A}-4\kappa\sin^{2} A-2-\frac{c}{\kappa}\cdot\frac{{\rm d}\kappa}{{\rm d}c}\right)
\frac{\bar{\partial}_{+}c}{\sin^{2} A}\Bigg),
\end{aligned}
\\[32pt]
\begin{aligned}
c\bar{\partial}_{-}\left(\frac{\bar{\partial}_{+}c}{\sin^{2} A}\right)
   &=\bar{\partial}_{+}c\Bigg(\frac{1}{2\varpi\cos^{2} A}\Big(\frac{\bar{\partial}_{+}c}{\sin^2A}-\frac{\bar{\partial}_{-}c}{\sin^2A}\Big)
\\&\qquad\qquad\qquad+\left(\frac{1}{\varpi\cos^{2} A}-4\kappa\sin^{2} A-2-\frac{c}{\kappa}\cdot\frac{{\rm d}\kappa}{{\rm d}c}\right)
\frac{\bar{\partial}_{-}c}{\sin^{2} A}\Bigg).
\end{aligned}
\end{array}
\right.\label{cd1}
\end{equation}
\end{prop}
\begin{proof}
These decompositions can be derived by (\ref{cd}). We omit the details.
%Similar characteristic decompositions were first obtained by Li and Zheng \cite{LMJ} in investigating semi-hyperbolic patches of solutions of the 2D self-similar system (\ref{PE}).
\end{proof}

%\begin{prop} We have the following decompositions:
%\begin{equation}
%\left\{
 % \begin{array}{ll}
 % \begin{aligned}
%\displaystyle &c\bar{\partial}_{+}\left(\frac{\bar{\partial}_{-}c}{\cos A}\right)
%=\frac{\bar{\partial}_{-}c}{\cos^{2} A}\Bigg(\sin2A\cos A+\sin^3A+\frac{1}{2\varpi}\frac{\bar{\partial}_{-}c}{\cos A}\\&\qquad\qquad\qquad\quad\qquad+\left(\frac{1}{2\varpi\cos^2 A}+\frac{(1+\sin^2A-\omega\cos^2 A)\sin^2 A}{2\varpi}\right)\frac{\bar{\partial}_{+}c}{\cos A}\Bigg),
%\end{aligned} \\[36pt]
%\begin{aligned}
%\displaystyle &c\bar{\partial}_{-}\left(\frac{\bar{\partial}_{+}c}{\cos A}\right)
%=\frac{\bar{\partial}_{+}c}{\cos^{2} A}\Bigg(\sin2A\cos A+\sin^3A+\frac{1}{2\varpi}\frac{\bar{\partial}_{+}c}{\cos A}\\&\qquad\qquad\qquad\quad\qquad+\left(\frac{1}{2\varpi\cos^2 A}+\frac{(1+\sin^2A-\omega\cos^2 A)\sin^2 A}{2\varpi}\right)\frac{\bar{\partial}_{-}c}{\cos A}\Bigg).
%\end{aligned}
%\end{array}
%\right.\label{cd2}
%\end{equation}
%\end{prop}
%\begin{proof}
%These decompositions can be proved directly by (\ref{cd}).
%\end{proof}

\begin{prop}
For the variable $\rho$, we have the following characteristic decomposition
\begin{equation}\label{81104}
\left\{
  \begin{array}{ll}
   \displaystyle c\bar{\partial}_{+}\bar{\partial}_{-}\rho=\sin2A\bar{\partial}_{-}\rho+
    \frac{\tau^4p''(\tau)}{4c\cos^{2}A }\Big[(\bar{\partial}_{-}\rho)^2+(\mathcal{F}-1)\bar{\partial}_{-}\rho\bar{\partial}_{+}\rho\Big],\\[12pt]
  \displaystyle c\bar{\partial}_{-}\bar{\partial}_{+}\rho=\sin2A\bar{\partial}_{+}\rho+
    \frac{\tau^4p''(\tau)}{4c\cos^{2}A }\Big[(\bar{\partial}_{+}\rho)^2+(\mathcal{F}-1)\bar{\partial}_{-}\rho\bar{\partial}_{+}\rho\Big],
  \end{array}
\right.
\end{equation}
where
$$
\mathcal{F}=2\sin^2A-\frac{8p'(\tau)\cos^4A}{\tau p''(\tau)}>0\quad \text{for}\quad~p''>0~~\text{and}~~p'<0.
$$
\end{prop}

\begin{proof}
This decomposition was first obtained in \cite{Lai5}. For the sake of completeness, we give the proof.

It follows from (\ref{11}) and (\ref{comm}) that
\begin{equation}
\begin{array}{rcl}
&&\bar{\partial}_{+}\big[c\tau\sin\alpha\bar{\partial}_{-}\rho\big]+
\bar{\partial}_{-}\big[c\tau\sin\beta\bar{\partial}_{+}\rho\big]
\\[6pt]&=&-\displaystyle\frac{1}{\sin2 A}\big[(\bar{\partial}_{+}\beta-\cos2 A \bar{\partial}_{-}\alpha)c\tau\sin\beta \bar{\partial}_{+}\rho-( \bar{\partial}_{-}\alpha- \cos2 A\bar{\partial}_{+}\beta)c\tau\sin\alpha \bar{\partial}_{-}\rho\big].
\end{array}
\end{equation}
Hence, we have
\begin{equation}
\begin{array}{rcl}
&&\displaystyle(\sin\alpha+\sin\beta)\frac{1}{c\tau }\frac{d(c\tau)}{d\rho }\bar{\partial}_{+}\rho\bar{\partial}_{-}\rho+\sin\alpha\bar{\partial}_{+}\bar{\partial}_{-}\rho+
\sin\beta\bar{\partial}_{-}\bar{\partial}_{+}\rho\\[8pt]&=&-\displaystyle\frac{1}{\sin2 A}
\Big[(\sin\beta \bar{\partial}_{+}\beta-\sin\beta \cos2 A \bar{\partial}_{-}\alpha+\cos\beta\sin2 A\bar{\partial}_{-}\beta) \bar{\partial}_{+}\rho\\[4pt]&&\qquad\qquad\qquad\qquad\displaystyle-( \sin\alpha\bar{\partial}_{-}\alpha- \sin\alpha\cos2 A\bar{\partial}_{+}\beta-\cos\alpha\sin2 A\bar{\partial}_{+}\alpha) \bar{\partial}_{-}\rho\Big].
\end{array}\label{12a}
\end{equation}

Applying the commutator relation (\ref{comm}) for $\rho$, we get
\begin{equation}
\bar{\partial}_{-} \bar{\partial}_{+}\rho- \bar{\partial}_{+} \bar{\partial}_{-}\rho=
\frac{1}{\sin2 A}\big[(\cos2 A \bar{\partial}_{+}\beta- \bar{\partial}_{-}\alpha) \bar{\partial}_{-}\rho-
(\bar{\partial}_{+}\beta-\cos2 A \bar{\partial}_{-}\alpha) \bar{\partial}_{+}\rho\big].\label{13a}
\end{equation}
Inserting this into (\ref{12a}), we get
\begin{equation}
\begin{array}{rcl}
&&\displaystyle(\sin\alpha+\sin\beta)\frac{1}{c\tau }\frac{d(c\tau)}{d\rho }\bar{\partial}_{+}\rho\bar{\partial}_{-}\rho+(\sin\alpha+\sin\beta)\bar{\partial}_{+}\bar{\partial}_{-}\rho\\[6pt]
&=&-\displaystyle\frac{1}{\sin2 A}
\Big[\cos\beta\sin2 A\bar{\partial}_{-}\beta\bar{\partial}_{+}\rho+( \sin\alpha+\sin\beta)\cos2 A\bar{\partial}_{+}\beta\bar{\partial}_{-}\rho\\[6pt]&&
\qquad\qquad\qquad\qquad-( \sin\alpha+\sin\beta)\bar{\partial}_{-}\alpha\bar{\partial}_{-}\rho+\cos\alpha\sin2 A\bar{\partial}_{+}
\alpha\bar{\partial}_{-}\rho\Big].
\end{array}\label{14a}
\end{equation}
Thus, by (\ref{7})--(\ref{8}) we get
\begin{equation}
 c\bar{\partial}_{+}\bar{\partial}_{-}\rho=\sin2A\bar{\partial}_{-}\rho+
    \frac{\tau^4p''(\tau)}{4c\cos^{2}A }\Big[(\bar{\partial}_{-}\rho)^2+(\mathcal{F}-1)\bar{\partial}_{-}\rho\bar{\partial}_{+}\rho\Big].
\end{equation}

The proof for the other decomposition is similar; we omit the details.
\end{proof}

%From (\ref{81104}),
%we have the following characteristic decomposition
%\begin{equation}\label{new}
%\left\{
 % \begin{array}{ll}
 % \begin{aligned}
  % \displaystyle c\bar{\partial}_{+}\Big(\frac{\bar{\partial}_{-}\rho}{\cos A}\Big)&=\left(\sin 2A+\frac{\sin^3A}{\cos A}\right)\frac{\bar{\partial}_{-}\rho}{\cos A}+
   % \frac{\tau^4p''(\tau)\sin A}{4c\cos A }\left(\frac{\bar{\partial}_{-}\rho}{\cos A}\right)^2\\&\qquad+  \frac{\tau^4p''(\tau)\sin A}{4c\cos A }\left(\mathcal{F}-1+\sin^2 A-\Omega \sin^2A\cos^2A\right)\frac{\bar{\partial}_{-}\rho\bar{\partial}_{+}\rho}{\cos^2 A},
   % \end{aligned}\\[32pt]
% \begin{aligned}
 %  \displaystyle c\bar{\partial}_{-}\Big(\frac{\bar{\partial}_{+}\rho}{\cos A}\Big)&=\left(\sin 2A+\frac{\sin^3A}{\cos A}\right)\frac{\bar{\partial}_{+}\rho}{\cos A}+
  %  \frac{\tau^4p''(\tau)\sin A}{4c\cos A }\left(\frac{\bar{\partial}_{+}\rho}{\cos A}\right)^2\\&\qquad+  \frac{\tau^4p''(\tau)\sin A}{4c\cos A }\left(\mathcal{F}-1+\sin^2 A-\Omega \sin^2A\cos^2A\right)\frac{\bar{\partial}_{-}\rho\bar{\partial}_{+}\rho}{\cos^2 A}.
   % \end{aligned}
  %\end{array}
%\right.
%\end{equation}

%From $\omega(\tau)>0$, one has
%\begin{equation}
%\begin{aligned}
%&\mathcal{F}-1+\sin^2 A-\Omega \sin^2A\cos^2A\\~>~&3\sin^2 A-1-\sin^2A\cos^2A-(1-\tan^2A)\sin^2A\cos^2A
%>0.
%\end{aligned}
%\end{equation}

\begin{prop}\label{pnew}
For any positive integer $n$, we have the following characteristic decomposition
\begin{equation}\label{new}
\left\{
  \begin{array}{ll}
   \displaystyle c\bar{\partial}_{+}\Big(\frac{\rho^n\bar{\partial}_{-}\rho}{\sin^2 A}\Big)=
    \frac{\tau^4p''(\tau)\bar{\partial}_{-}\rho}{4c\cos^{2}A }\Bigg[\frac{\rho^n\bar{\partial}_{-}\rho}{\sin^2 A}+\mathcal{H}\frac{\rho^n\bar{\partial}_{+}\rho}{\sin^2 A}\Bigg],\\[16pt]
  \displaystyle c\bar{\partial}_{-}\Big(\frac{\rho^n\bar{\partial}_{+}\rho}{\sin^2 A}\Big)=
    \frac{\tau^4p''(\tau)\bar{\partial}_{+}\rho}{4c\cos^{2}A }\Bigg[\frac{\rho^n\bar{\partial}_{+}\rho}{\sin^2 A}+\mathcal{H}\frac{\rho^n\bar{\partial}_{-}\rho}{\sin^2 A}\Bigg],
  \end{array}
\right.
\end{equation}
where
\begin{equation}\label{82803}
\mathcal{H}=2\sin^2A-\frac{8p'(\tau)\cos^4A}{\tau p''(\tau)}+\frac{4nc^2\cos^{2}A }{\tau^3p''(\tau)}-2\cos^2 A+2\Omega \cos^4 A-1.
\end{equation}
\end{prop}
\begin{proof}
This proposition can be obtained by using (\ref{7})--(\ref{8}) and (\ref{81104}).
\end{proof}

%===============================================================================================

\subsection{Hodograph transformation}

\subsubsection{Concept of hodograph transformation}
Consider the hodograph transformation for (\ref{PEC}):
$$
T:  (\xi, \eta)\rightarrow (u, v).
$$
If the Jacobian
$$
j(u, v;\xi, \eta)=\frac{\partial(u, v)}{\partial(\xi, \eta)}=u_{\xi}v_{\eta}-u_{\eta}v_{\xi}\neq 0
$$
for a solution $(u, v)(\xi, \eta)$ of (\ref{PEC}), we may consider $\xi$ and $\eta$ as functions of $u$ and $v$. From
\begin{equation}
u_{\xi}=j\eta_{v}, \quad u_{\eta}=-j\xi_{v}, \quad v_{\xi}=-j\eta_{u}, \quad v_{\eta}=j\xi_{u},
\end{equation}
and regarding $h$ as a function of $u$ and $v$,
we then see that $\xi(u, v)$ and $\eta(u, v)$ satisfy the equations
\begin{equation}\label{HT}
\left\{
  \begin{array}{ll}
  (c^{2}-U^{2})\eta_{v}+UV(\xi_{v}+\eta_{u})+(c^{2}-V^{2})\xi_{u}=0, \\[2pt]
  \xi_{v}-\eta_{u}=0.
  \end{array}
\right.
\end{equation}

%Differentiating (\ref{52201}) with respect to $u$ and $v$, one gets
%\begin{equation}\label{52203}
%\xi_{u}=1+h_{uu}, \quad \xi_v=h_{uv}, \quad \eta_u=h_{uv}, \quad \eta_{v}=h_{vv}+1.
%\end{equation}
%Inserting these into the first equation of (\ref{HT}), we obtain
%\begin{equation}\label{52202}
%(c^2-h_u^2)h_{vv}+2h_uh_vh_{uv}+(c^2-h_v^2)h_{uu}=h_u^2+h_v^2-c^2.
%\end{equation}
%This is a second order partial differential equation for $h$ alone.

The eigenvalues of (\ref{HT}) are
\begin{equation}
\Lambda=\Lambda_{\pm}=\frac{UV\pm
c\sqrt{U^{2}+V^{2}-c^{2}}}{c^{2}-V^{2}}.
\end{equation}
The characteristic curves of (\ref{HT}) in the $(u, v)$-plane are defined as the integral curves of
\begin{equation}\label{823301}
\Gamma_{\pm}:~\frac{{\rm d}v}{{\rm d}u}=\Lambda_\pm.
\end{equation}
It is easy to see that
\begin{equation}\label{52204}
\Lambda_{\pm}=-\frac{1}{\lambda_{\mp}},
\end{equation}
which implies that the directions of the $C$-characteristics of one kind are perpendicular to the $\Gamma$-characteristics of the other kind. Or, more precisely, the directions of $C_{+}$ and $\Gamma_{-}$ and of $C_{-}$ and $\Gamma_{+}$ through corresponding points $(\xi,\eta)$ and $(u, v)$ are perpendicular.

We introduce the normalized directional derivatives along the $\Gamma$-characteristics,
\begin{equation}\label{52503}
\hat{\partial}_{+}=-\sin\beta\partial_{u}+\cos\beta\partial_{v}, \quad
\hat{\partial}_{-}=\sin\alpha\partial_{u}-\cos\alpha\partial_{v}.
\end{equation}
The characteristic equations of (\ref{HT}) are
\begin{equation}\label{38}
\left\{
  \begin{array}{ll}
    \hat{\partial}_{+}\eta=\tan \alpha \hat{\partial}_{+}\xi, \\[2pt]
\hat{\partial}_{-}\eta=\tan \beta \hat{\partial}_{-}\xi.
  \end{array}
\right.
\end{equation}
This implies a correspondence between $C_{\pm}$ and $\Gamma_{\pm}$ characteristic curves, i.e., the image of a $C_{+}$ (resp. $C_{-}$) characteristic curve in the $(\xi, \eta)$-plane is a $\Gamma_{+}$ (resp. $\Gamma_{-}$) characteristic curve in  the hodograph plane.

\subsubsection{Phase space system of equations}
From Bernoulli's law (\ref{5802}), we have
\begin{equation}\label{52201}
h_u=\xi-u=-U, \quad h_{v}=\eta-v=-V.
\end{equation}
We obtain the following first order system for $(U, V, h)$:
\begin{equation}\label{52301}
\left(
  \begin{array}{ccc}
    c^2-V^2 & UV & 0 \\
    0 & 1 & 0 \\
    0 & 0 & 1 \\
  \end{array}
\right)\left(
         \begin{array}{c}
           U \\
           V \\
           h \\
         \end{array}
       \right)_u+ \left(
  \begin{array}{ccc}
    UV & c^2-U^2 & 0 \\
    -1 & 0 & 0 \\
    0 & 0 & 0 \\
  \end{array}
\right)\left(
         \begin{array}{c}
           U \\
           V \\
           h \\
         \end{array}
       \right)_v=\left(
                   \begin{array}{c}
                     2c^2-U^2-V^2 \\
                     0 \\
                     -U \\
                   \end{array}
                 \right).
\end{equation}
The system has three eigenvalues
\begin{equation}
\Lambda_0=0\quad \text{and} \quad \Lambda_{\pm}=\frac{UV\pm
c\sqrt{U^{2}+V^{2}-c^{2}}}{c^{2}-V^{2}}.
\end{equation}

Multiplying (\ref{52301}) on the left by the vector $(l_{+}, l_{-}, l_{0})^{T}$, where
$$
l_0=(0, 0, 1)\quad \text{and} \quad l_{\pm}=(1, \pm c\sqrt{U^{2}+V^{2}-c^{2}}, 0),
$$
 we obtain
\begin{equation}\label{52303}
\left\{
  \begin{array}{ll}
    \displaystyle\hat{\partial}_{+}U+\Lambda_{-}\hat{\partial}_{+}V=-\frac{2\cos 2A}{\sin\alpha}, \\[8pt]
   \displaystyle \hat{\partial}_{-}U+\Lambda_{+}\hat{\partial}_{-}V=\frac{2\cos 2A}{\sin\beta},\\[8pt]
    h_u=-U.
  \end{array}
\right.
\end{equation}
%where
%$$
%\partial_{+}=\partial_{u}+\Lambda_{+}\partial_{v}\quad \text{and}\quad \partial_{-}=\partial_{u}+\Lambda_{-}\partial_{v}.
%$$

By (\ref{U}) and (\ref{52201}) we have
\begin{equation}\label{52501}
\hat{\partial}_{\pm}c=-\frac{1}{\kappa}.
\end{equation}
From (\ref{U}) we have
\begin{equation}\label{6601}
\hat{\partial}_{\pm}U=\frac{\cos\sigma \hat{\partial}_{\pm}c}{\sin A}-\frac{c\sin\sigma}{2\sin A}(\hat{\partial}_{\pm}\alpha+\hat{\partial}_{\pm}\beta)-\frac{c\cos\sigma\cos A}{2\sin^2 A}(\hat{\partial}_{\pm}\alpha-\hat{\partial}_{\pm}\beta)
\end{equation}
and
\begin{equation}\label{6602}
\hat{\partial}_{\pm}V=\frac{\sin\sigma \hat{\partial}_{\pm}c}{\sin A}+\frac{c\cos\sigma}{2\sin A}(\hat{\partial}_{\pm}\alpha+\hat{\partial}_{\pm}\beta)-\frac{c\sin\sigma\cos A}{2\sin^2 A}(\hat{\partial}_{\pm}\alpha-\hat{\partial}_{\pm}\beta).
\end{equation}

Inserting (\ref{52501})--(\ref{6602}) into (\ref{52303}), one gets
\begin{equation}\label{52305}
\left\{
  \begin{array}{ll}
    \displaystyle\hat{\partial}_{+}\alpha=-\frac{\tau p''(\tau)}{4cp'(\tau)}\Omega\sin 2A,  \\[10pt]
    \displaystyle\hat{\partial}_{-}\beta=\frac{\tau p''(\tau)}{4cp'(\tau)}\Omega\sin 2A,  \\[12pt]
   \displaystyle \tau_u=\frac{\tau\cos\sigma}{c \sin A}.
  \end{array}
\right.
\end{equation}
Actually, system (\ref{52305}) is linearly degenerate in the sense of Lax \cite{Lax}.

\subsubsection{Characteristic decompositions in the hodograph plane}

\begin{prop}
(Commutator relation of $\hat{\partial}_{\pm}$) For system (\ref{52301}), the following commutator relation holds:
\begin{equation}\label{52401}
\hat{\partial}_{+}\hat{\partial}_{-}-\hat{\partial}_{-}\hat{\partial}_{+}=\tan A\hat{\partial}_{+}\alpha(\hat{\partial}_{+}-\hat{\partial}_{-}).
\end{equation}
\end{prop}
\begin{proof}
From (\ref{52503}) we have
\begin{equation}\label{317}
\partial_u=\frac{\cos\alpha\hat{\partial}_{+}+\cos\beta\hat{\partial}_{-}}{\sin 2A}, \quad \partial_v=\frac{\sin\alpha\hat{\partial}_{+}+\sin\beta\hat{\partial}_{-}}{\sin 2A}.
\end{equation}
By a direct computation one can obtain (\ref{52401}).
\end{proof}

Using the commutator relation, we easily derive:
\begin{prop}
For the variables $\alpha$ and $\beta$, we have the characteristic decompositions:
\begin{equation}\label{52502}
\left\{
  \begin{array}{ll}
    \hat{\partial}_{+}\hat{\partial}_{-}\alpha+\mathcal{W}(A, c)\hat{\partial}_{-}\alpha=\mathcal{Q}(A, c),  \\[4pt]
 \hat{\partial}_{-}\hat{\partial}_{+}\beta+\mathcal{W}(A, c)\hat{\partial}_{+}\beta=-\mathcal{Q}(A, c),
  \end{array}
\right.
\end{equation}
where
$$
\mathcal{W}(A, c)=-\frac{\tau p''(\tau)}{4cp'(\tau)}\Big[\Omega(4\sin^2 A-1)+2\tan^2A\Big],
$$
$$
\mathcal{Q}(A, c)=\left(\frac{\tau p''(\tau)}{4cp'(\tau)}\right)^2\Omega(1-3\tan^2A)\sin 2A-\frac{\tan A}{2\tau p'(\tau)}\cdot\frac{{\rm d}}{{\rm d}\tau}\left(\frac{\tau p''(\tau)}{p'(\tau)}\right).
$$
\end{prop}
\begin{proof}
We apply the commutator relation (\ref{52401}) to obtain
$$
\hat{\partial}_{+}\hat{\partial}_{-}\alpha-\hat{\partial}_{-}\hat{\partial}_{+}\alpha=\tan A\hat{\partial}_{+}\alpha(\hat{\partial}_{+}\alpha-\hat{\partial}_{-}\alpha).
$$
Using the expressions of $\hat{\partial}_{+}\alpha$ and $\hat{\partial}_{-}\beta$ in (\ref{52305}) and recalling (\ref{52501}), we directly obtain the first equation of (\ref{52502}). The second equation of (\ref{52502}) can be proved similarly.
\end{proof}

We need to show that the Jacobian $j^{-1}(u, v; \xi, \eta)$ does not vanish:
$$
j^{-1}(u, v; \xi, \eta)=\xi_u\eta_v-\xi_v\eta_u\neq0.
$$
By computation we obtain
$$
j^{-1}(u, v; \xi, \eta)=\frac{\sin^2 2A}{\cos\alpha\cos\beta}\hat{\partial}_{-}\xi\cdot\hat{\partial}_{+}\xi.
$$

A direct computation yields
\begin{equation}\label{61312}
\hat{\partial}_{+}\xi=-\hat{\partial}_{+}U-\sin\beta\quad \text{and}\quad \hat{\partial}_{-}\xi=-\hat{\partial}_{-}U+\sin\alpha.
\end{equation}
We compute
\begin{equation}\label{61311}
\begin{aligned}
\hat{\partial}_{+}U+\sin\beta&=\frac{\cos\sigma \hat{\partial}_{+}c}{\sin A}-\frac{c\sin\sigma}{2\sin A}(\hat{\partial}_{+}\alpha+\hat{\partial}_{+}\beta)-\frac{c\cos\sigma\cos A}{2\sin^2 A}(\hat{\partial}_{+}\alpha-\hat{\partial}_{+}\beta)+\sin\beta\\&=
\frac{\tau p''(\tau)}{2p'(\tau)}\cdot\frac{\cos\alpha}{\sin 2A}+\frac{c\cos\alpha}{2\sin^2 A}\hat{\partial}_{+}\beta
\\&=\frac{c\cos\alpha}{2\sin^2 A}\left(\hat{\partial}_{+}\beta+\frac{\tau p''(\tau)\tan A}{2cp'(\tau)}\right)
\end{aligned}
\end{equation}
and
\begin{equation}\label{61319}
\begin{aligned}
\hat{\partial}_{-}U-\sin\alpha&=\frac{\cos\sigma \hat{\partial}_{-}c}{\sin A}-\frac{c\sin\sigma}{2\sin A}(\hat{\partial}_{-}\alpha+\hat{\partial}_{-}\beta)-\frac{c\cos\sigma\cos A}{2\sin^2 A}(\hat{\partial}_{-}\alpha-\hat{\partial}_{-}\beta)-\sin\alpha\\&=
\frac{\tau p''(\tau)}{2p'(\tau)}\cdot\frac{\cos\beta}{\sin 2A}-\frac{c\cos\beta}{2\sin^2 A}\hat{\partial}_{-}\alpha
\\&=-\frac{c\cos\beta}{2\sin^2 A}\left(\hat{\partial}_{-}\alpha-\frac{\tau p''(\tau)\tan A}{2cp'(\tau)}\right).
\end{aligned}
\end{equation}
Therefore, we have
\begin{equation}\label{6805}
j^{-1}(u, v; \xi, \eta)=-(c\cot A)^2 \mathcal{Z}_{+}\mathcal{Z}_{-},
\end{equation}
where
\begin{equation}\label{62701}
\mathcal{Z}_{+}=\hat{\partial}_{+}\beta+\frac{\tau p''(\tau)\tan A}{2cp'(\tau)}\quad \text{and} \quad
\mathcal{Z}_{-}=\hat{\partial}_{-}\alpha-\frac{\tau p''(\tau)\tan A}{2cp'(\tau)}.
\end{equation}
Thus, to ensure that $j^{-1}$ does not vanish, we require $\mathcal{Z}_{\pm}\neq0$, $c\neq0$, and $A\neq\frac{\pi}{2}$.

From (\ref{61312}) and (\ref{62701}), we also have
\begin{equation}\label{26102402}
\hat{\partial}_{+}\xi=-\frac{c\cos\alpha}{2\sin^2 A}\mathcal{Z}_{+}\quad \text{and}\quad \hat{\partial}_{-}\xi=\frac{c\cos\alpha}{2\sin^2 A}\mathcal{Z}_{-}.
\end{equation}

To estimate $\mathcal{Z}_{\pm}$, we use (\ref{52502}) to derive the following characteristic equations:
\begin{prop}
For  $\mathcal{Z}_{\pm}$, we have
\begin{equation}\label{322}
\left\{
  \begin{array}{ll}
   \hat{\partial}_{+}\mathcal{Z}_{-}+\mathcal{W}\mathcal{Z}_{-}=\displaystyle\frac{\tau p''(\tau)}{4cp'(\tau)}(\tan^2 A+1)\mathcal{Z}_{+}, \\[10pt]
\hat{\partial}_{-}\mathcal{Z}_{+}+\mathcal{W}\mathcal{Z}_{+}=\displaystyle\frac{\tau p''(\tau)}{4cp'(\tau)}(\tan^2 A+1)\mathcal{Z}_{-}.
  \end{array}
\right.
\end{equation}
\end{prop}

 We compute
$$
\hat{\partial}_{+}c=c_{\xi}\hat{\partial}_{+}\xi+c_{\eta}\hat{\partial}_{+}\eta
=(c_{\xi}+c_{\eta}\tan\alpha)\hat{\partial}_{+}\xi=\sec\alpha\bar{\partial}_{+}c \hat{\partial}_{+}\xi.
$$
Thus, we have
$$
\bar{\partial}_{+}c=\frac{2\sin^2 A}{c\kappa\mathcal{Z}_{+}}.
$$
Combining this with $\frac{{\rm d}c}{{\rm d}\tau}=-\frac{c}{\tau \kappa}$, we get
\begin{equation}\label{6807}
\bar{\partial}_{+}\tau=-\frac{2\tau\sin^2 A}{c^2\mathcal{Z}_{+}}.
\end{equation}
Similarly, one has
\begin{equation}\label{6808}
\bar{\partial}_{-}\tau=\frac{2\tau\sin^2 A}{c^2\mathcal{Z}_{-}}.
\end{equation}

\subsubsection{One-to-one property of the hodograph transformation}
In this part we shall show that under some conditions the hodograph transformation can be inverted back to the physical plane.

\begin{lem}\label{onetoone1}
(One-to-one Property 1)
Assume that the system (\ref{52305}) admits a classical solution $(\alpha, \beta, \tau)=(\alpha, \beta, \tau)(u, v)$ in a domain  $\Pi$. Assume also that for any two points in $\Pi$, there exist two characteristic cures of different types connecting these two points.
Assume moreover that the solution satisfies
\begin{equation}\label{92502aa}
0<\alpha<\frac{\pi}{2}, \quad -\frac{\pi}{2}<\beta<0,
\end{equation}
\begin{equation}\label{92801e}
\rho>0,\quad \mathcal{Z}_{+}<0,\quad\text{and} \quad \mathcal{Z}_{-}>0
\end{equation}
in $\Pi$.
Then the mapping $(u, v)\rightarrow (\xi, \eta)$ defined by
\begin{equation}\label{92504a}
\xi=u-\frac{c(\tau)\cos(\sigma(u, v))}{\sin (A(u,v))}, \quad \eta=v-\frac{c(\tau)\sin(\sigma(u, v))}{\sin (A(u, v))}
\end{equation}
is invertible, where $\sigma(u, v)=\frac{\alpha(u, v)+\beta(u, v)}{2}$ and $A(u, v)=\frac{\alpha(u, v)-\beta(u, v)}{2}$.
\end{lem}

\begin{proof}
By (\ref{92502a}),
in the domain $\Pi$ we have
\begin{equation}\label{92503aa}
\hat{\partial}_{+}\xi=-\frac{c\cos\alpha}{2\sin^2A}\mathcal{Z}_{+}>0,\quad
\hat{\partial}_{-}\xi=\frac{c\cos\beta}{2\sin^2A}\mathcal{Z}_{-}>0.
\end{equation}
\begin{equation}\label{92503bv}
\hat{\partial}_{+}\eta=-\frac{c\sin\alpha}{2\sin^2A}\mathcal{Z}_{+}>0,\quad
\hat{\partial}_{-}\eta=\frac{c\sin\beta}{2\sin^2A}\mathcal{Z}_{-}<0.
\end{equation}

 By (\ref{92503aa}) and (\ref{92503bv}) it follows that at least one of $\xi$ and $\eta$ is monotone along the connecting $\Gamma_{\pm}$ characteristic path. Thus, no two distinct points in the $(u, v)$-plane are mapped to the same point in the $(\xi, \eta)$-plane.
This completes the proof.
\end{proof}

\begin{rem}\label{92505aa}
The condition (\ref{92502aa}) can be relaxed to $(\alpha, \beta)\in (\alpha_{_L}, \alpha_{_R})\times(\beta_{_L}, \beta_{_R})$,
where
\begin{equation}\label{92802f}
\alpha_{_R}-\alpha_{_L}<\frac{\pi}{2}, \quad \beta_{_R}-\beta_{_L}<\frac{\pi}{2},\quad \alpha_{_L}-\beta_{_R}>0\quad \text{and} \quad \alpha_{_R}-\beta_{_L}<\pi.
\end{equation}

When $\alpha_{_R}-\beta_{_R}>\frac{\pi}{2}$, we
let
$\alpha'=\alpha+\delta$ and $\beta'=\beta+\delta$,
where $\delta=\frac{\pi}{2}-\alpha_{_R}$.
Then  we have
$0<\alpha'<\frac{\pi}{2}$ and $-\frac{\pi}{2}<\beta'<0$.
When $\alpha_{_R}-\beta_{_R}<\frac{\pi}{2}$, we
let
$\alpha'=\alpha+\delta$ and $\beta'=\beta+\delta$,
where $\delta=-\beta_{_R}$.
Then  we have
$0<\alpha'<\frac{\pi}{2}$ and $-\frac{\pi}{2}<\beta'<0$.
Thus, by rotating the coordinates of the physical plane and using the fact that
 $\mathcal{Z}_{\pm}$ are invariant under the coordinate rotation,  we still have that the mapping (\ref{92504a}) is invertible.
\end{rem}

%=============
\begin{figure}[htbp]
\begin{center}
\includegraphics[scale=0.43]{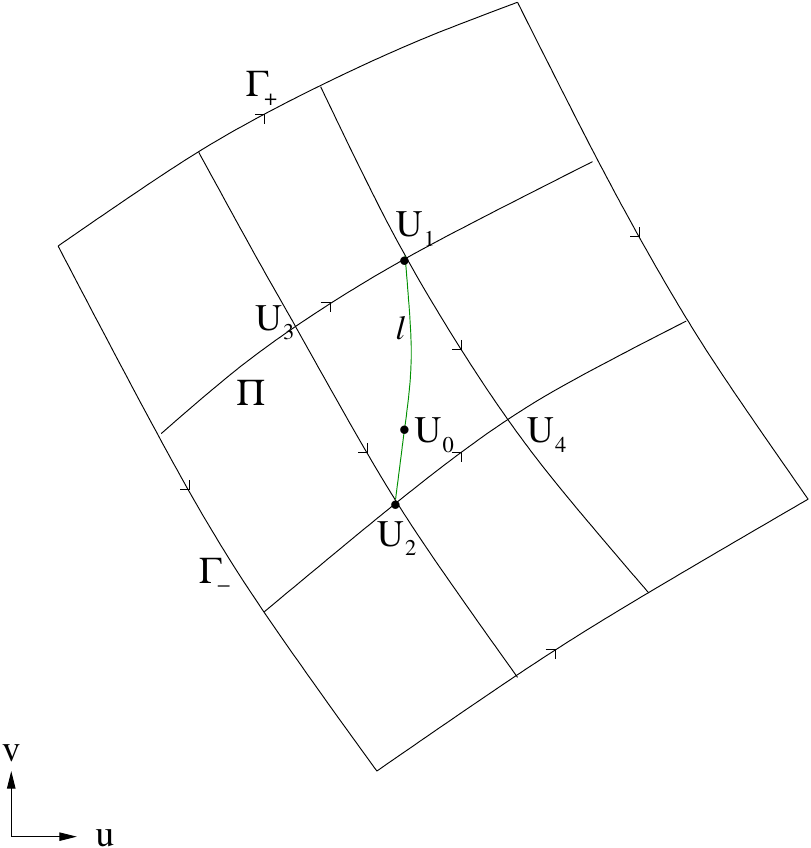}
\caption{\footnotesize A curvilinear quadrilateral region bounded by $\Gamma_{\pm}$ characteristic curves.}
\label{one-to-one}
\end{center}
\end{figure}
%==========

\begin{lem}\label{onetoone}
(One-to-one Property 2)
Assume that the system (\ref{52305}) admits a classical solution $(\alpha, \beta, \tau)=(\alpha, \beta, \tau)(u, v)$ in a curvilinear quadrilateral region $\Pi$ bounded by four characteristic curves in the $(u, v)$-plane; see Figure \ref{one-to-one}.
Assume furthermore that in $\Pi$, the solution satisfies (\ref{92801e}) and
\begin{equation}\label{92502a}
-\frac{\pi}{2}<\alpha<\frac{\pi}{2}, \quad -\frac{\pi}{2}<\beta<\frac{\pi}{2}, \quad  0<\alpha-\beta<\pi, \quad \rho>0.
\end{equation}
Then the mapping (\ref{92504a}) from $(u, v)$ to $(\xi, \eta)$ is invertible.
\end{lem}
\begin{proof}
From (\ref{92502a})
it follows that, in the domain $\Pi$,
\begin{equation}\label{92503a}
\hat{\partial}_{+}\xi=-\frac{c\cos\alpha}{2\sin^2A}\mathcal{Z}_{+}>0,\quad
\hat{\partial}_{-}\xi=\frac{c\cos\beta}{2\sin^2A}\mathcal{Z}_{-}>0.
\end{equation}

Let $\mathrm{U}_1$ and $\mathrm{U}_2$ be any two points in $\Pi$.
There exist two characteristic curves of different types connecting these two points; see Figure \ref{one-to-one}.

Suppose that $(\xi, \eta)(\mathrm{U}_1)=(\xi, \eta)(\mathrm{U}_2)$. Then by (\ref{92503a}) we see that
the forward $\Gamma_{+}$ characteristic curve emanating from $\mathrm{U}_2$ intersects the
forward $\Gamma_{-}$ characteristic curve emanating from $\mathrm{U}_1$ at a point $\mathrm{U}_4$, and the
backward $\Gamma_{+}$ characteristic curve emanating from $\mathrm{U}_1$ intersects the
backward $\Gamma_{-}$ characteristic curve emanating from $\mathrm{U}_2$ at a point $\mathrm{U}_3$.
Let $\Pi_{_\mathrm{U}}$ be the curvilinear quadrilateral region bounded by the $\Gamma$ characteristic curves $\widetilde{\mathrm{U_3U_1}}$,
 $\widetilde{\mathrm{U_1U_4}}$, $\widetilde{\mathrm{U_3U_2}}$, and $\widetilde{\mathrm{U_2U_4}}$.
 By (\ref{92503a}) we also see that
 there is a smooth non-characteristic curve ${\it l}$ in $\Pi_{_\mathrm{U}}$ such that
 $$
 \xi(u, v)=\xi(\mathrm{U}_1)\quad \text{on} \quad \ell.
$$

Let  $(t_1(u, v), t_2(u, v))$ be the unit tangent vector of $\ell$. We then define the derivative along $\ell$ by
$$
\partial_{\ell}=t_1(u, v)\partial_u+t_2(u, v)\partial_v \quad \text{along} \quad \ell.
$$
 Since $\ell$ is a non-characteristic curve, we have the decomposition
$$
\partial_{\ell}=a(u, v)\hat{\partial}_{+}+b(u, v)\hat{\partial}_{-} \quad \text{along} \quad \ell,
$$
where $a(u, v)\neq0$ and $b(u, v)\neq0$ for $(u, v)\in \ell$.

Since $\xi=\xi(U_1)$ on $\ell$, by (\ref{92503a}) we have
$$
\partial_{\ell}\xi=-\frac{ac\cos\alpha}{2\sin^2A}\mathcal{Z}_{+}+\frac{bc\cos\beta}{2\sin^2A}\mathcal{Z}_{-}=0\quad \text{along}\quad \ell.
$$
Since $\eta(\mathrm{U}_1)=\eta(\mathrm{U}_2)$, there must exist a $U_0\in\ell$ such that
$$
\partial_{\ell}\eta= -\frac{ac\sin\alpha}{2\sin^2A}\mathcal{Z}_{+}+\frac{bc\sin\beta}{2\sin^2A}\mathcal{Z}_{-}=0\quad\text{at}\quad \mathrm{U}_0.
$$
This leads to a contradiction, since $A=(\alpha-\beta)/2>0$, $c>0$, and $\mathcal{Z}_{\pm}\neq0$ on $\ell$.

So, there are no two distinct points in $\Pi$ that map to a single point in the $(\xi, \eta)$-plane.
\end{proof}

\begin{rem}\label{92505a}
The condition (\ref{92502a}) can be relaxed to
$$
\alpha_{_L}<\alpha<\alpha_{_R}, \quad \beta_{_L}<\beta<\beta_{_R}, \quad 0<\alpha-\beta<\pi,
$$
where
\begin{equation}\label{92806}
\alpha_{_R}-\alpha_{_L}<\pi, \quad \beta_{_R}-\beta_{_L}<\pi,\quad \alpha_{_L}-\beta_{_R}>-\pi, \quad
\text{and} \quad \alpha_{_R}-\beta_{_L}<\pi.
\end{equation}

When $\alpha_{_R}-\beta_{_R}>0$, we
let
$\alpha'=\alpha+\delta$ and $\beta'=\beta+\delta$,
where $\delta=\frac{\pi}{2}-\alpha_{_R}$.
Then by (\ref{92806}) we have
\begin{equation}\label{9281001}
-\frac{\pi}{2}<\alpha'<\frac{\pi}{2}, \quad -\frac{\pi}{2}<\beta'<\frac{\pi}{2}, \quad \text{and}\quad 0<\alpha'-\beta'<\pi.
\end{equation}
When $\alpha_{_R}-\beta_{_R}<0$, we
let
$\alpha'=\alpha+\delta$ and $\beta'=\beta+\delta$,
where $\delta=\frac{\pi}{2}-\beta_{_R}$.
Then by (\ref{92806}) we still have (\ref{9281001}).
Thus, by rotating the coordinates of the physical plane and using the fact that
 $\mathcal{Z}_{\pm}$ are invariant under the coordinate rotation,  we still have that the mapping (\ref{92504a}) is invertible.
\end{rem}

\subsection{Notations and some assumptions}
\subsubsection{Notations}
\begin{figure}[htbp]
\begin{center}
\includegraphics[scale=0.43]{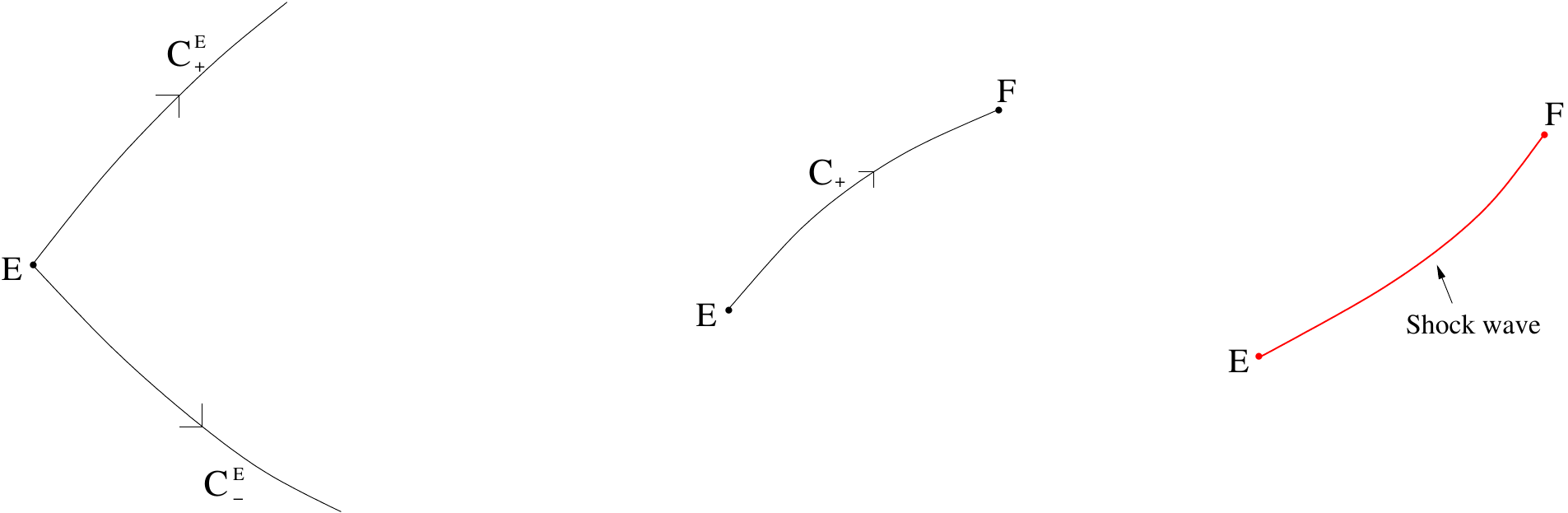}
\caption{\footnotesize Characteristic curves and shock wave curve.}
\label{Fig2.3}
\end{center}
\end{figure}

In this paper,
we use various curves to describe the flow structure. In order to distinguish between the different meanings of these curves, we shall introduce the following notation; see Figure \ref{Fig2.3}.
Let $\mathrm{E}$ and $\mathrm{F}$ be two points in the $(\xi,\eta)$-plane. We introduce the following notations:
\begin{itemize}
  \item $C_{+}^{\mathrm{E}}$ (resp. $C_{-}^{\mathrm{E}}$): a forward $C_{+}$ (resp. $C_{-}$) characteristic line emanating from point $\mathrm{E}$;
  \item $\widetilde{\mathrm{EF}}$: a $C_{+}$ (or $C_{-}$) characteristic curve segment from $\mathrm{E}$ to $\mathrm{F}$ along the corresponding characteristic direction;
  \item $\overline{\mathrm{EF}}$: a straight $C_{+}$ (or $C_{-}$) characteristic line segment from $\mathrm{E}$ to $\mathrm{F}$;
\item $\widehat{\mathrm{EF}}$: a shock curve connecting points $\mathrm{E}$ and $\mathrm{F}$;
\item $\wideparen{\mathrm{EF}}$: a vacuum boundary or a level curve of $\rho$ connecting points $\mathrm{E}$ and $\mathrm{F}$.
\end{itemize}

Let ${\it S}$ be an oblique shock curve in the $(\xi, \eta)$-plane. We introduce the following two notations:
\begin{itemize}
  \item $\bar{\partial}_{s}$: the derivative along ${\it S}$;
  \item $\phi$: the inclination angle of ${\it S}$.
\end{itemize}
Assume that ${\it S}$ can be represented by $\eta=\eta_s(\xi)$, $\xi\in (\xi', \xi'')$. Then, along ${\it S}$,
\begin{equation}\label{20101}
\phi:=\arctan\big(\eta_s'(\xi)\big), \quad \bar{\partial}_{s}u:=\cos\phi\frac{{\rm d}u(\xi, \eta_s(\xi))}{{\rm d}\xi},\quad \bar{\partial}_{s}v:=\cos\phi\frac{{\rm d}v(\xi, \eta_s(\xi))}{{\rm d}\xi}
\end{equation}
If the flow upstream or downstream of ${\it S}$ is $C^1$ smooth up to ${\it S}$, then we also have
\begin{equation}\label{110401}
\bar{\partial}_{s}(u , v)=\cos\phi\partial_{x}(u , v)+\sin\phi\partial_{y}(u , v)\quad\text{along}~{\it S}.
\end{equation}
However, if the upstream or downstream  flow is not $C^1$ smooth up to ${\it S}$ then (\ref{110401}) does not hold along ${\it S}$. This case occurs if the shock is a sonic shock.

%=============
%==========

\subsubsection{Some basic assumptions about the nonconvex equation of state}
To establish the existence of global-in-time solutions to the 2D Riemann problem (\ref{Potential}), (\ref{initial})--(\ref{boundary}), 
we make the following assumptions on the nonconvex equation of state with two inflection points:
\begin{description}
                 \item[(E1)] There exists $\gamma\in (1, 3)$ such that $\lim\limits_{\tau\rightarrow+\infty}\frac{\tau p''(\tau)}{p'(\tau)}=-(\gamma+1)$;
                 \item[(E2)] $\omega'(\tau)\leq 0$ for  $\tau>\tau_2^e$;
\vskip 2pt
                 \item[(E3)] $\frac{c}{\kappa}\cdot\frac{{\rm d}\kappa}{{\rm d}c}\rightarrow 0$ as $\tau\rightarrow +\infty$.
               \end{description}
Under the above assumptions, we have, as $\tau\to+\infty$,
$$
\kappa(\tau)\rightarrow \frac{2}{\gamma-1}, \quad \omega(\tau)\rightarrow\frac{3-\gamma}{\gamma+1}, \quad
\Theta(\tau)\rightarrow A_{\infty}:=\arctan\sqrt{\frac{3-\gamma}{\gamma+1}}\quad \text{and} \quad  \varpi(\tau)\rightarrow\frac{\gamma-1}{\gamma+1}.
$$
Moreover, there exist a sufficiently large $\tau_*>0$ and a sufficiently small $\varepsilon_*>0$ such that
\begin{equation}\label{51703}
\kappa>\frac{1}{\gamma-1}\quad\text{and}\quad
0<\varpi<\frac{2(\gamma-1)}{\gamma+1}\quad \text{for}\quad \tau\geq \tau_*;
\end{equation}
\begin{equation}\label{51703a}
\frac{\kappa\cos^2 2A+1}{2\cos^2 A}-\frac{c}{\kappa}\cdot\frac{{\rm d}\kappa}{{\rm d}c}>0\quad \text{for}\quad \tau\geq \tau_*\quad \text{and}\quad 0<A<\frac{\pi}{2};
\end{equation}
\begin{equation}\label{51704}
\frac{1}{\varpi\cos^{2} A}-4\kappa\sin^{2} A-2-\frac{c}{\kappa}\cdot\frac{{\rm d}\kappa}{{\rm d}c}>0\quad \text{for}\quad \tau\geq \tau_*\quad \text{and}\quad 0<A\leq \varepsilon_*.
\end{equation}

\begin{rem}
Consider the van der Waals gas equation of state
$$
p=\frac{S}{(\tau-1)^{\gamma}}-\frac{1}{\tau^2}.
$$
For suitable  $S$ and $\gamma>1$, the above assumptions are satisfied.
\end{rem}

\section{Elementary waves}

\subsection{Oblique shock wave}
We now discuss oblique shocks of the 2D pseudo-steady potential flow equations (\ref{PEC}).
Assume that there is a straight shock in the $(\xi, \eta)$-plane  and that the states on both sides of it are constant.
Let $\phi$ be the inclination angle of the shock front.
To formulate the shock relations, we denote the normal and tangential components of the flow vector by $N$ and $L$, respectively, i.e.,
\begin{equation}\label{102201}
N=\pm(U\sin\phi-V\cos\phi)\quad \text{and}\quad L=U\cos\phi+V\sin\phi.
\end{equation}
where the "$\pm$" are related to the classification of the shock wave, such as the 1-shock and 2-shock to be introduced next.

%=============
\begin{figure}[htbp]
\begin{center}
\includegraphics[scale=0.7]{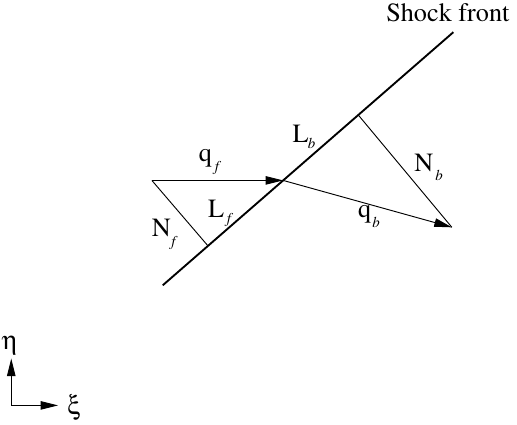}
\caption{\footnotesize A 2D self-similar oblique rarefaction shock front.}
\label{obliqueshock}
\end{center}
\end{figure}
%==========

We use subscripts `$f$' and `$b$' to denote the front-side and back-side states of the
shock front, respectively; see Figure \ref{obliqueshock}. Then
the shock relations for  (\ref{PEC}) are:
\begin{equation}\label{RH}
\left\{
  \begin{array}{ll}
    \rho_f N_f=\rho_b N_b =m>0, \\[4pt]
    L_{f}=L_b,\\[4pt]
   N_f^2+2h(\tau_f)= N_b^2+2h(\tau_b).
  \end{array}
\right.
\end{equation}
%In addition  the shock is also required to satisfy the following
%Liu's extended entropy condition (see \cite{Liu1}):
%\begin{equation}\label{43001}
%-\frac{2h(\tau_f)-2h(\tau_b)}{\tau_f^2-\tau_b^2}\geq -\frac{2h(\tau_f)-2h(\tau)}{\tau_f^2-\tau^2}\quad \text{for~ all}\quad\tau\in(\tau_f, \tau_b).
%\end{equation}
%for all $$.

From the first and third relations of (\ref{RH}), we have
\begin{equation}\label{RHm}
m^2=-\frac{2h(\tau_f)-2h(\tau_b)}{\tau_f^2-\tau_b^2}.
\end{equation}
Since the gas considered here is nonconvex, the shock is also required to satisfy the
Liu's extended entropy condition (Liu \cite{Liu1}):
\begin{equation}\label{43001}
-p'(\tau_b)\geq-\frac{2h(\tau_f)-2h(\tau_b)}{\tau_f^2-\tau_b^2}\geq -\frac{2h(\tau_f)-2h(\tau)}{\tau_f^2-\tau^2}\quad \text{for~ all}~ \tau\in\big(\min\{\tau_f, \tau_b\}, \max\{\tau_f, \tau_b\}\big).
\end{equation}
In this paper, we confine ourselves to rarefaction shocks, i.e., the shocks that satisfy $\tau_f<\tau_b$.

\begin{figure}[htbp]
\begin{center}
\includegraphics[scale=0.63]{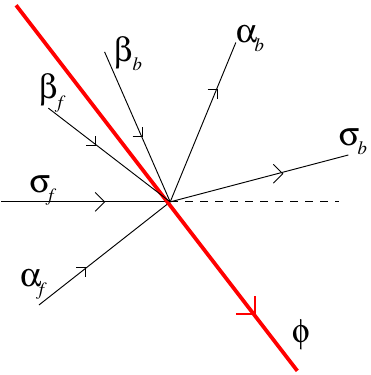}\qquad\qquad\qquad \qquad\qquad \includegraphics[scale=0.63]{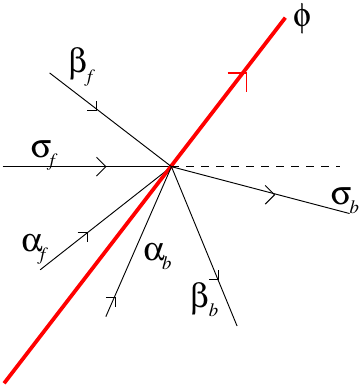}
\caption{\footnotesize Left: 1-shock; right: 2-shock.}
\label{Fig3.2}
\end{center}
\end{figure}

From (\ref{43001}) we have that along the shock,
\begin{equation}\label{20102}
N_f\geq c_f\quad \text{and}\quad N_b\leq c_b.
\end{equation}
We also assume that the flows upstream and downstream of the shock are supersonic, i.e., $q>c$.
Then by (\ref{20102}) we have that along the shock,
$\beta_b\leq \phi\leq \beta_f$ or $\alpha_f\leq \phi\leq \alpha_b$;
see Figure \ref{Fig3.2}.
 As in Lax \cite{Lax},  we call the oblique shock 1-shock (resp. 2-shock) if
\begin{equation}\label{83002a}
\beta_b\leq \phi\leq \beta_f\quad
 (\text{resp.} \alpha_f\leq \phi\leq \alpha_b).
\end{equation}

Thus, We  have
\begin{equation}\label{111306}
\left\{
  \begin{array}{ll}
    N=V\cos\phi-U\sin\phi\quad \text{and}\quad L=U\cos\phi+V\sin\phi, & \hbox{1-shock;} \\[4pt]
    N=U\sin\phi-V\cos\phi\quad \text{and}\quad L=U\cos\phi+V\sin\phi, & \hbox{2-shock.}
  \end{array}
\right.
\end{equation}
From (\ref{111306}) we also have
%\begin{equation}\label{102301}
%U=N\sin\phi+L\cos\phi\quad \text{and}\quad V=-N\cos\phi+L\sin\phi.
%\end{equation}
\begin{equation}\label{111306a}
\left\{
  \begin{array}{ll}
    U=L\cos\phi-N\sin\phi\quad \text{and}\quad V=N\cos\phi+L\sin\phi,  & \hbox{1-shock;} \\[4pt]
    U=N\sin\phi+L\cos\phi\quad \text{and}\quad V=L\sin\phi-N\cos\phi, & \hbox{2-shock.}
  \end{array}
\right.
\end{equation}

%=============
%==========

\begin{prop}\label{pro3.1.1}
For rarefaction 1-shocks,
we have $A_f\leq \sigma_f-\phi<\sigma_b-\phi\leq A_b$;
for rarefaction 2-shocks,
we have $A_f\leq\phi-\sigma_f<\phi-\sigma_b\leq A_b$.
\end{prop}
\begin{proof}
For rarefaction 1-shocks,
we have $L=q\cos(\sigma-\phi)$.
So, by the second relation of (\ref{RH}) we have
\begin{equation}\label{111311}
q_f\cos(\sigma_f-\phi)=q_b\cos(\sigma_b-\phi).
\end{equation}
By the definition of 1-shocks, we have
\begin{equation}\label{111312}
\sigma_b-\phi\leq \sigma_b-\beta_b=A_b<\frac{\pi}{2}\quad \text{and} \quad \sigma_{f}-\phi\geq \sigma_f-\beta_f=A_f.
\end{equation}
Since $\tau_f<\tau_b$, by Bernoulli's law (\ref{5802}) we have $q_b>q_f$. Combining this with (\ref{111311}) and (\ref{111312}), we have $A_f\leq \sigma_f-\phi<\sigma_b-\phi\leq A_b$.

The proof for the other case is similar; we omit the details.
\end{proof}

By (\ref{81601}),
there exist $\tau_{1}^a$ and $\tau_{2}^a$, with $\tau_1^a<\tau_{1}^{i}<\tau_{2}^{i}<\tau_{2}^a$, such that
\begin{equation}\label{5901}
p'(\tau_{1}^a)=p'(\tau_{2}^{i})\quad \text{and}\quad p'(\tau_{1}^{i})=p'(\tau_{2}^a).
 \end{equation}
We cite some propositions from Lai \cite{Laip1}.
\begin{prop}\label{5602}
({\bf Double-sonic})
There exists  $\tau_1^e$ and $\tau_2^e$, with $1<\tau_1^e<\tau_2^e$, such that
\begin{equation}\label{42301}
p'(\tau_1^e)=p'(\tau_2^e)=\frac{2h(\tau_1^e)-2h(\tau_2^e)}{(\tau_1^e)^2-(\tau_2^e)^2},
\end{equation}
and
\begin{equation}\label{42301aa}
-\frac{2h(\tau_1^e)-2h(\tau_2^e)}
{(\tau_1^e)^2-(\tau_2^e)^2}>-\frac{2h(\tau_1^e)-2h(\tau)}{(\tau_1^e)^2-\tau^2}\quad \text{for all}~ \tau\in(\tau_1^e, \tau_2^e).
\end{equation}
Moreover, $\tau_1^e\in (\tau_1^a, \tau_1^i)$ and $\tau_2^e\in (\tau_2^i, \tau_2^a)$.
\end{prop}

%=============
\begin{figure}[htbp]
\begin{center}
\includegraphics[scale=0.43]{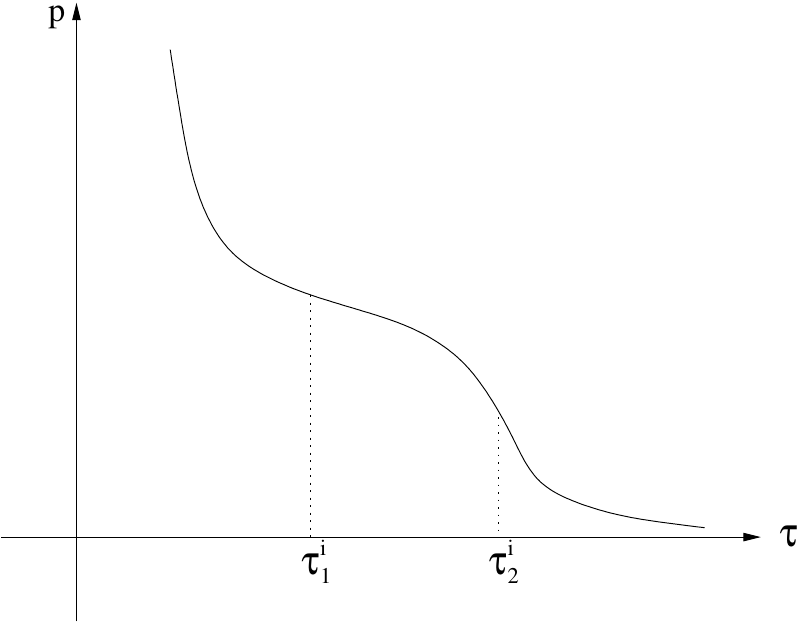}\qquad\qquad\includegraphics[scale=0.43]{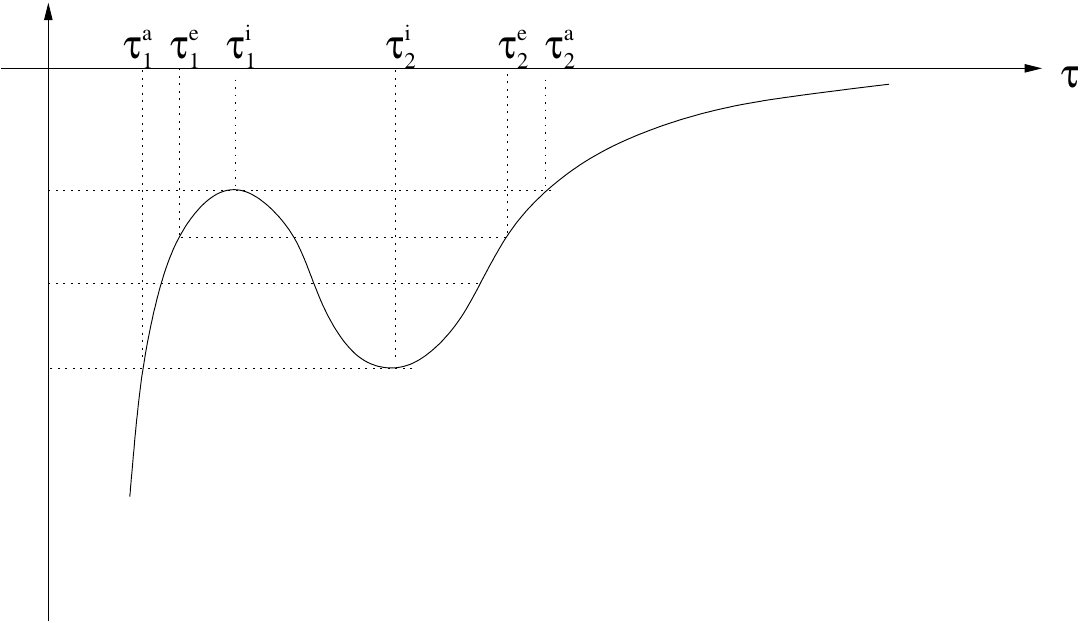}
\caption{\footnotesize Left: the isentrope $p=p(\tau)$; right: the graph of $p'(\tau)$.}
\label{EOS}
\end{center}
\end{figure}
%==========

%{\color{red} Why called double-sonic shock:}

\begin{prop}\label{51003}
({\bf Post-sonic}) For any $\tau_{f}\in(\tau_1^e, \tau_2^i)$, there exists one and only one $\psi_{po}(\tau_{f})\in(\tau_f,+\infty)$ such that
\begin{equation}\label{43003}
p'(\psi_{po}(\tau_f))=\frac{2h(\psi_{po}(\tau_{f}))-2h(\tau_{f})}{\psi_{po}^2(\tau_{f})-\tau_{f}^2}<p'(\tau_{f}),
\end{equation}
and
\begin{equation}\label{43004}
-\frac{2h(\psi_{po}(\tau_{f}))-2h(\tau_{f})}{\psi_{po}^2(\tau_{f})-\tau_{f}^2}>-\frac{2h(\tau)-2h(\tau_{f})}{\tau^2-\tau_{f}^2}\quad \text{for all}~\tau\in (\tau_{f}, \psi_{po}(\tau_{f})).
\end{equation}
Moreover, the function $\psi_{po}\in C^1[\tau_1^e, \tau_2^i)$ and satisfies
\begin{equation}\label{5906}
 \psi_{po}'(\tau_{f})<0\quad \text{and}\quad \tau_2^i<\psi_{po}(\tau_{f})<\tau_2^e\quad \text{for}~  \tau_{f}\in (\tau_1^e, \tau_2^i); %\lim\limits_{\tau\rightarrow \tau_2^i}s_{po}(\tau)=\tau_2^i.
\end{equation}
$$
\psi_{po}(\tau_{f})\rightarrow \tau_2^e\quad \text{as}~  \tau_{f}\rightarrow\tau_1^e; \qquad
\psi_{po}(\tau_{f})\rightarrow \tau_2^i\quad \text{as}~ \tau_{f}\rightarrow\tau_2^i.
$$
\end{prop}

\begin{prop}\label{pre}
({\bf Pre-sonic})
For any $\tau_{f}\in (\tau_1^e, \tau_1^i)$, there exists one and only one $\psi_{pr}(\tau_{f})\in (\tau_{f}, +\infty)$ such that
\begin{equation}\label{5701}
p'(\psi_{pr}(\tau_{f}))<\frac{2h(\psi_{pr}(\tau_{f}))-2h(\tau_{f})}{\psi_{pr}^2(\tau_{f})-\tau_{f}^2}=p'(\tau_{f})
\end{equation}
and
\begin{equation}\label{43007}
-\frac{2h(\psi_{pr}(\tau_{f}))-2h(\tau_{f})}{\psi_{pr}^2(\tau_{f})-\tau_{f}^2}>-
\frac{2h(\tau_{f})-2h(\tau)}{\tau_{f}^2-\tau^2}\quad \text{for all}~ \tau\in (\tau_{f}, \psi_{pr}(\tau_{f})).
\end{equation}
Moreover, the function $\psi_{pr}$ satisfies
\begin{equation}
\psi_{pr}'(\tau_{f})<0\quad \text{and}\quad
\psi_{pr}(\tau_{f})\in (\tau_1^i, \tau_2^e)\quad \text{for} ~\tau_{f}\in (\tau_1^e, \tau_1^i).
\end{equation}
\end{prop}

\begin{prop}\label{tran}
({\bf Transonic}) Assume $\tau_f\in (\tau_1^e, \tau_2^i)$. Then,
\begin{equation}\label{43012}
p'(\tau_b)<\frac{2h(\tau_b)-2h(\tau_f)}{\tau_b^2-\tau_f^2}<p'(\tau_f)\quad \text{for}~ \tau_b\in (\tau_*, \psi_{po}(\tau_f)),
\end{equation}
where $\tau_*=\left\{
                \begin{array}{ll}
                  \psi_{pr}(\tau_f), & \hbox{$\tau_f\leq \tau_1^i$;} \\
                  \tau_f, & \hbox{$\tau_f>\tau_1^i$.}
                \end{array}
              \right.
$ Moreover,
\begin{equation}\label{43013}
-\frac{2h(\tau_b)-2h(\tau_f)}{\tau_b^2-\tau_f^2}>-\frac{2h(\tau)-2h(\tau_f)}{\tau^2-\tau_f^2}\quad \text{for all}~ \tau\in (\tau_f, \tau_b).
\end{equation}
\end{prop}

We classify the oblique shocks of nonconvex gases as follows.
\begin{itemize}
  \item When $\tau_f=\tau_1^e$ and $\tau_b=\tau_2^e$, they are referred to as double-sonic shocks. For double-sonic shocks, one has
$$
N_f=c_f:=\sqrt{-\tau_f^2p'(\tau_f)} \quad \text{and}\quad N_b=c_b:=\sqrt{-\tau_b^2p'(\tau_b)}.
$$
  \item When $\tau_f\in (\tau_1^e, \tau_2^i)$ and $\tau_b=\psi_{po}(\tau_f)$, they are referred to as post-sonic shocks. For oblique post-sonic shocks, one has
$$
N_f>c_f\quad \text{and}\quad N_b=c_b.
$$
\item When $\tau_f\in (\tau_1^e, \tau_1^i)$ and $\tau_b=\psi_{pr}(\tau_f)$, they are referred to as pre-sonic shocks. For oblique pre-sonic shocks, one has
$$
N_f=c_f\quad \text{and}\quad N_b<c_b.
$$
\item For oblique transonic shocks\footnote{In this paper, the term ``transonic" refers to the flow velocity relative to the shock front.},  one has
$$
N_f>c_f\quad \text{and}\quad N_b<c_b.
$$
\end{itemize}

\subsection{Centered simple wave}
We now discuss centered simple waves emanating from the sharp corner $\mathrm{O}$.
In the following discussion, we use $\theta\in(-\pi, \pi)$ to denote the polar angle, i.e., the counterclockwise angle from the positive $\xi$-axis to the vector $(\xi, \eta)$.
\begin{lem}
Assume the functions $\bar{q}(\theta)$, $\bar{\tau}(\theta)$, and $\bar{\sigma}(\theta)$ satisfy the equations
\begin{equation}\label{102302}
\left\{
  \begin{array}{ll}
    \bar{q}'(\theta)\cos(\bar{A}(\theta))+\bar{q}(\theta)\bar{\sigma}'(\theta)\sin(\bar{A}(\theta))=0, \\[4pt]
    \bar{q}(\theta)\bar{q}'(\theta)+\bar{\tau}(\theta)p'(\bar{\tau}(\theta))\bar{\tau}'(\theta)=0,  \\[4pt]
    \bar{\sigma}(\theta)+\bar{A}(\theta)=\theta
  \end{array}
\right.
\end{equation}
and  $\bar{q}(\theta)>\bar{c}(\theta)>0$ for $\theta\in (\theta_1, \theta_2)$,
where
\begin{equation}\label{102303}
\bar{A}(\theta)=\arcsin\left(\frac{\bar{c}(\theta)}{\bar{q}(\theta)}\right)\quad \text{and}\quad \bar{c}(\theta)=\bar{\tau}(\theta)\sqrt{-p'(\bar{\tau}(\theta))}.
\end{equation}
Let
\begin{equation}\label{102306}
\bar{u}(\theta)=\bar{q}(\theta)\cos(\bar{\sigma}(\theta))\quad \text{and}\quad \bar{v}(\theta)=\bar{q}(\theta)\sin(\bar{\sigma}(\theta)).
\end{equation}
Then
\begin{equation}
(u, v, \tau)=(\bar{u}, \bar{v}, \bar{\tau})(\theta), \quad   \theta_1<\theta<\theta_2
\end{equation}
is a centered simple wave with straight $C_{+}$ characteristic lines of (\ref{PEC}) on the fan-shaped region
$$
\begin{aligned}
\Delta=\big\{(\xi, \eta)~\big|~\xi=r\cos\theta, ~\eta=r\sin\theta,~ r\in (0, \bar{r}(\theta)),
~\theta\in(\theta_1, \theta_2)
\big\},
\end{aligned}
$$
where $\bar{r}(\theta)=\sqrt{\bar{q}^2(\theta)-\bar{c}^2(\theta)}$.
The curve $\xi=\bar{r}(\theta)\cos\theta$, $\eta=\bar{r}(\theta)\sin\theta$, $\theta \in(\theta_1, \theta_2)$ is called the sonic edge of this centered simple wave.
\end{lem}
\begin{proof}
From the third equation of (\ref{102302}) we have
\begin{equation}\label{82104}
\big(\bar{u}(\theta), \bar{v}(\theta)\big)\cdot(\cos\theta, \sin\theta)=\bar{r}(\theta)\quad \text{and}\quad
\big(\bar{u}(\theta), \bar{v}(\theta)\big)\cdot(\sin\theta, -\cos\theta)=\bar{c}(\theta).
\end{equation}
Thus, for any fixed $\theta\in (\theta_1, \theta_2)$, the ray $C_{+}^{\theta}:$ $\xi=r\cos\theta$, $\eta=r\sin\theta$, where $0\leq r\leq  \bar{r}(\theta)$,
is tangent to the sonic circle $(\xi-\bar{u}(\theta))^2+(\eta-\bar{v}(\theta))^2=\bar{c}^2(\theta)$.
Therefore, by (\ref{c}) we see that for any $\theta\in (\theta_1, \theta_2)$, $C_{+}^{\theta}$
 is
a $C_{+}$ characteristic line. Hence, we have $$\alpha=\theta, \quad \bar{\partial}_{+}u=0,\quad \bar{\partial}_{+}v=0,\quad\text{and} \quad\bar{\partial}_{+}\tau=0\quad \text{in}\quad \Lambda.$$ Consequently,
\begin{equation}\label{82101}
(u-\xi)\bar{\partial}_{+}u+(v-\eta)\bar{\partial}_{+}v+\bar{\partial}_{+}h(\tau)=0\quad\text{and}\quad
\bar{\partial}_{+}u+\lambda_{-}\bar{\partial}_{+}v=0\quad\text{in}\quad \Delta.
\end{equation}
%Moreover, using (\ref{82104}) we have
%\begin{equation}\label{7103}
% A=\arctan\Big(\frac{\bar{c}(\theta)}{\bar{r}(\theta)-r}\Big)\quad \text{and}\quad\beta=\theta+\frac{\pi}{2}-2A\quad\text{in}\quad \Theta,
%\end{equation}
%where
%$$
%\theta=\arctan\frac{\eta}{\xi}-\frac{\pi}{2}, \quad r=\sqrt{\xi^2+\eta^2}.
%$$

By the first equation of (\ref{102302}) and $\alpha=\theta$, we have
\begin{equation}\label{82102}
\begin{aligned}
\bar{\partial}_{-}u+\lambda_{+}\bar{\partial}_{-}v
&=\big( \bar{u}'(\theta)+\tan\theta  \bar{v}'(\theta)\big)\bar{\partial}_{-}\theta\\&=\big(\bar{q}(\theta)'\cos(\bar{A}(\theta))
+\bar{q}(\theta)\bar{\sigma}'(\theta)\sin(\bar{A}(\theta))\big)\frac{\bar{\partial}_{-}\theta}{\cos\theta}=0 \quad\text{in}\quad \Delta.
\end{aligned}
\end{equation}

By the second equation of (\ref{102302}) we have
$$
\begin{aligned}
&\big(\bar{u}(\theta)-\xi\big)\bar{u}'(\theta)+\big(\bar{v}(\theta)-\eta\big)\bar{v}'(\theta)
+h'(\bar{\tau}(\theta))\bar{\tau}'(\theta)\\
~=~&\big(\bar{u}(\theta)+r\cos\theta\big)\bar{u}'(\theta)+\big(\bar{v}(\theta)+r\sin\theta\big)\bar{v}'(\theta)
+h'(\bar{\tau}(\theta))\bar{\tau}'(\theta)\\~=~&
\bar{u}(\theta)\bar{u}'(\theta)+\bar{v}(\theta)\bar{v}'(\theta)
+h'(\bar{\tau}(\theta))\bar{\tau}'(\theta)\\  ~=~& \bar{q}(\theta)\bar{q}'(\theta)+\bar{\tau}(\theta)p'(\bar{\tau}(\theta))\bar{\tau}'(\theta)~=~0.
\end{aligned}
$$
Thus,
\begin{equation}\label{82103}
(u-\xi)\bar{\partial}_{-}u+(v-\eta)\bar{\partial}_{-}v+\bar{\partial}_{-}h(\tau)=0\quad\text{in}\quad \Delta.
\end{equation}
This completes the proof of the lemma.
\end{proof}

By (\ref{102303}), system (\ref{102302}) can be written as
\begin{equation}\label{102304}
\left\{
  \begin{array}{ll}
    \displaystyle\bar{q}'=\frac{2\bar{c}p'(\bar{\tau})\cos \bar{A}}{\bar{\tau}p''(\bar{\tau})},\\[10pt]
     \displaystyle\bar{\tau}'=-\frac{2\bar{q}\bar{c}\cos \bar{A}}{\bar{\tau}^2p''(\bar{\tau})},\\[10pt]
 \displaystyle\bar{\sigma}'=-\frac{2 p'(\bar{\tau})\cos^2 \bar{A}}{\bar{\tau}p''(\bar{\tau})}.
  \end{array}
\right.
\end{equation}
System (\ref{102304}) is a nonlinear system of ordinary differential equations. From the second equation of (\ref{102304}) we also have
\begin{equation}\label{81501}
\bar{c}'=\varpi(\tau)\bar{q}\cos \bar{A}.
\end{equation}

\subsection{Elementary waves for $\tau_0>\tau_2^{i}$}

\begin{figure}[htbp]
\begin{center}
  \includegraphics[scale=0.6]{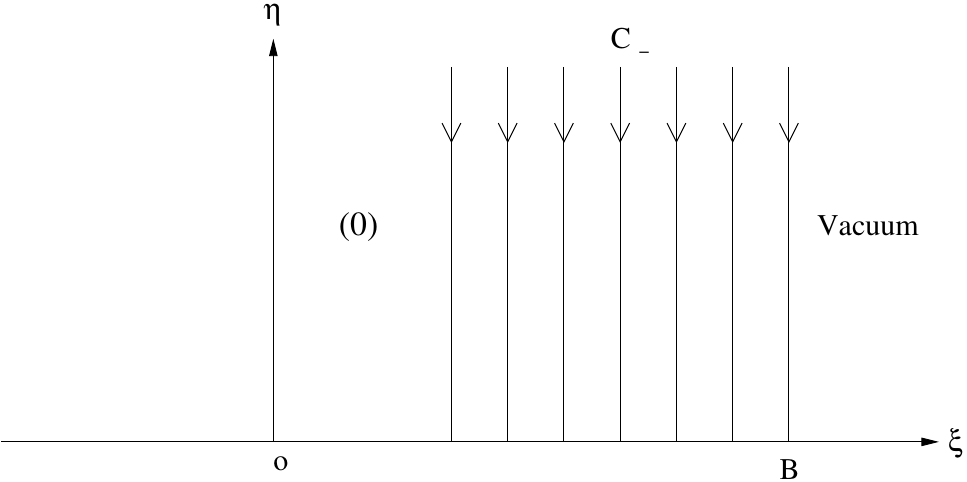}
\caption{\footnotesize Straight $C_{-}$ characteristic lines in the planar rarefaction wave.}
\label{Fig1.4.1.1}
\end{center}
\end{figure}

From $\tau_0>\tau_2^i$ we have $p''(\tau)>0$ for $\tau>\tau_0$. Then by the classical theory for hyperbolic systems of conservation laws, we know that away from the corner
the supersonic oncoming flow $(0)$ is connected to the vacuum by a plane rarefaction wave (see Figure \ref{Fig1.4.1.1}). The rarefaction wave can be
represented by
\begin{equation}\label{71501}
(u, v, \tau)=(u_{-}, v_{-}, \tau_{-})(\xi), \quad\xi_{_P}<\xi<\xi_{_B}, \quad \eta>0.
\end{equation}
Here, $\xi_{_\mathrm{P}}=u_0-c_0$, $\xi_{_\mathrm{B}}=u_0+\displaystyle\int_{\tau_0}^{+\infty}\sqrt{-p'(\tau)}{\rm d}\tau$, $v_{-}(\xi)\equiv0$, and the functions $u_{-}(\xi)$ and $\tau_{-}(\xi)$ are determined by
\begin{equation}\label{71601A}
u_{-}(\xi)-\int_{\tau_0}^{\tau_{-}(\xi)}\sqrt{-p'(\tau)}{\rm d}\tau=u_0\quad \text{and} \quad \xi=u_{-}(\xi)-c_{-}(\xi),
\end{equation}
where $c_{-}=\sqrt{-\tau_{-}^2p'(\tau_{-})}$;
see \cite{S2}.

%=============
%==========

Due to the sharp corner, there is a centered simple wave emanating from the origin.
In order to construct the centered simple wave, we consider (\ref{102304}) with data
\begin{equation}\label{102305}
(\bar{q}, \bar{\tau}, \bar{\sigma})(\alpha_0)=(u_0, \tau_0, 0),
\end{equation}
where $\alpha_0=\arcsin(\frac{c_0}{u_0})$.

Let
$$
\hat{q}(\tau)=\sqrt{u_0^2+2h(\tau_0)-2h(\tau)}, \quad \tau\geq\tau_0.
$$
Then the function $q=\hat{q}(\tau)$ is strictly increasing for $\tau>\tau_0$ and has an inverse function denoted by  $\tau=\hat{\tau}(q)$, where $q>u_0$.

Let
$$
\alpha_v=-\int_{u_0}^{q_{_{lim}}}\frac{\sqrt{q^2-\hat{c}^2(q)}}{q\hat{c}(q)}{\rm d}q,
$$
where
\begin{equation}\label{3506}
q_{_{lim}}=\sqrt{u_0^2+2h(\tau_0)-h_{\infty}}, \quad h_{\infty}=\lim\limits_{\tau\rightarrow +\infty}h(\tau), \quad \text{and}\quad \hat{c}(q)=\sqrt{-\hat{\tau}^2(q)p'(\hat{\tau}(q))}.
\end{equation}
%Let
%$$
%R_{+}(\sigma, q)=\sigma+\int_{u_0}^{q}\frac{\sqrt{q^2-\hat{c}^2(q)}}{q\hat{c}(q)}{\rm d}q.
%$$
%From  the Bernoulli's law (\ref{5802}) we know that the solution of the initial value problem (\ref{102304}), (\ref{102305}) satisfies
%$\hat{c}(\bar{q}(\theta))=\bar{c}(\theta)$ and
%$\frac{{\rm d}}{{\rm d}\theta}R_{+}(\bar{\sigma}(\theta), \bar{q}(\theta))=0$ for $\theta<\alpha_0$.
 Then  the initial value problem (\ref{102304}), (\ref{102305}) admits a solution for $\theta\in (\alpha_v, \alpha_0)$. By (\ref{102304}) and (\ref{81501}) the solution satisfies
\begin{equation}\label{81801}
\bar{q}'(\theta)<0, \quad  \bar{\sigma}'(\theta)>0, \quad \bar{\tau}'(\theta)<0, \quad \text{and}\quad 0<\bar{A}(\theta)<\frac{\pi}{2}
\end{equation}
for $\alpha_v< \theta<\alpha_0$.
Moreover,
$$
\lim\limits_{\theta\rightarrow \alpha_v}\bar{\tau}(\theta)=+\infty \quad \text{and}\quad \lim\limits_{\theta\rightarrow \alpha_v}\bar{q}(\theta)=q_{_{lim}}.
$$
%Moreover, by the first equation of  we have
%$$

%$$

%Let $\bar{u}(\theta)=\bar{q}(\theta)\cos(\bar{\sigma}(\theta))$ and $\bar{v}(\theta)=\bar{q}(\theta)\sin(\bar{\sigma}(\theta))$.
%Then the function
%$(u, v, \tau)=(\bar{u}, \bar{v}, \bar{\tau})(\theta)$
%is a complete expansion fan wave of the system (\ref{PEC}) in the region
%$
%\big\{(\xi, \eta)~\big|~\xi=r\cos\theta, ~\eta=r\sin\theta,~ r\in (0, \bar{r}(\theta)),
%~\theta\in(\alpha_v, \alpha_0)
%\big\}
%$.

%=============
\begin{figure}[htbp]
\begin{center}
\includegraphics[scale=0.58]{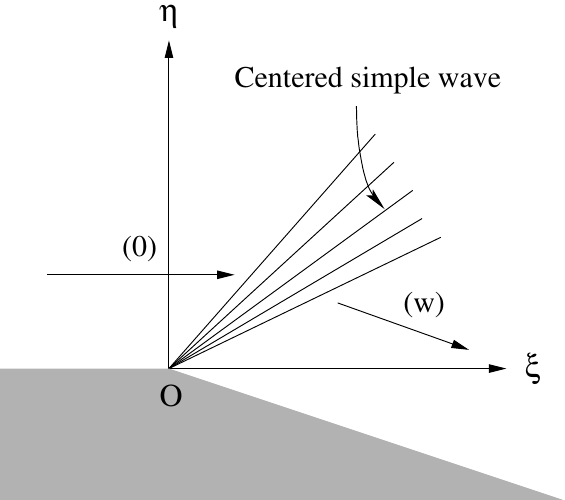}\qquad\qquad \qquad\includegraphics[scale=0.58]{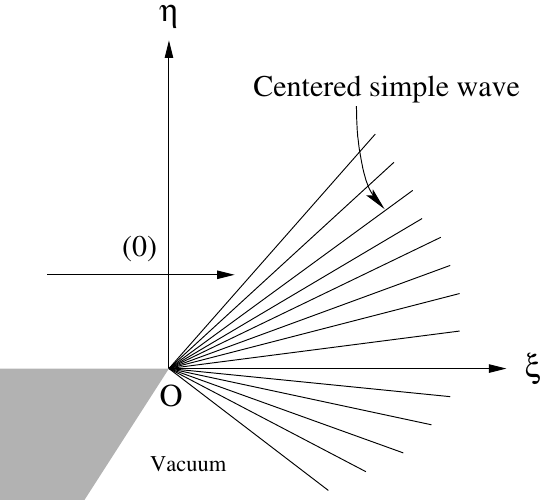}
\caption{\footnotesize Centered simple wave near the corner for $\tau_0>\tau_2^{i}$. }
\label{Fig1.4.1.2}
\end{center}
\end{figure}
%==========
%Therefore, we obtain the following Prandtl-Meyer fan wave solution near the sharp corner.

  When
 $\pi+{\alpha}_v<\theta_w<\pi$, there exists an $\alpha_w\in (\alpha_v, \alpha_0)$ such that
$\bar{v}(\alpha_w)=\bar{u}(\alpha_w)\tan\theta_w$. In this case, at the corner
the oncoming flow $(0)$ is deflected through a centered simple wave to a constant state $(u_w, v_w, \tau_w)$, where  $(u_w, v_w, \tau_w)=(\bar{u}, \bar{v}, \bar{\tau})(\alpha_w)$. The flow near the corner can be represented by
\begin{equation}\label{820a}
(u, v, \tau)=\left\{
               \begin{array}{ll}
                 (u_0, 0, \tau_0), & \hbox{$\theta>\alpha_0$;} \\[2pt]
                 (\bar{u}, \bar{v}, \bar{\tau})(\theta), & \hbox{${\alpha}_w<\theta<\alpha_0$;} \\[2pt]
                 (u_w, v_w, \tau_w), & \hbox{$\theta_w-\pi<\theta<{\alpha}_w$;}
               \end{array}
             \right.
\end{equation}
see Figure \ref{Fig1.4.1.2} (left).
When $\theta_w<\pi+{\alpha}_v$, the oncoming flow $(0)$ is deflected through a complete centered simple wave to a vacuum state. The flow near the corner can be represented by
\begin{equation}\label{820b}
(u, v, \tau)=\left\{
               \begin{array}{ll}
                 (u_0, 0, \tau_0), & \hbox{$\theta>\alpha_{0}$;} \\[2pt]
                 (\bar{u}, \bar{v}, \bar{\tau})(\theta), & \hbox{${\alpha}_v<\theta<\alpha_{0}$;} \\[2pt]
                 \text{vacuum}, & \hbox{$\theta_w-\pi<\theta<{\alpha}_v$;}
               \end{array}
             \right.
\end{equation}
see Figure \ref{Fig1.4.1.2} (right).

\subsection{Elementary waves for $\tau_1^{i}<\tau_0<\tau_2^{i}$}
Since $\tau_{1}^{i}<\tau_{0}<\tau_{2}^{i}$, by
Proposition \ref{51003} one has $\psi_{po}(\tau_0)\in (\tau_2^i, \tau_2^e)$ such that
\begin{equation}\label{5601}
p'(\psi_{po}(\tau_0))=\frac{2h(\psi_{po}(\tau_0))-2h(\tau_0)}{\psi_{po}^2(\tau_0)-\tau_0^2}.
\end{equation}
For convenience, we let $\tau_{po}=\psi_{po}(\tau_0)$, $\rho_{po}=\frac{1}{\tau_{po}}$, and $c_{po}=\sqrt{-\tau_{po}^2p'(\tau_{po})}$.

%=============
\begin{figure}[htbp]
\begin{center}
  \includegraphics[scale=0.5]{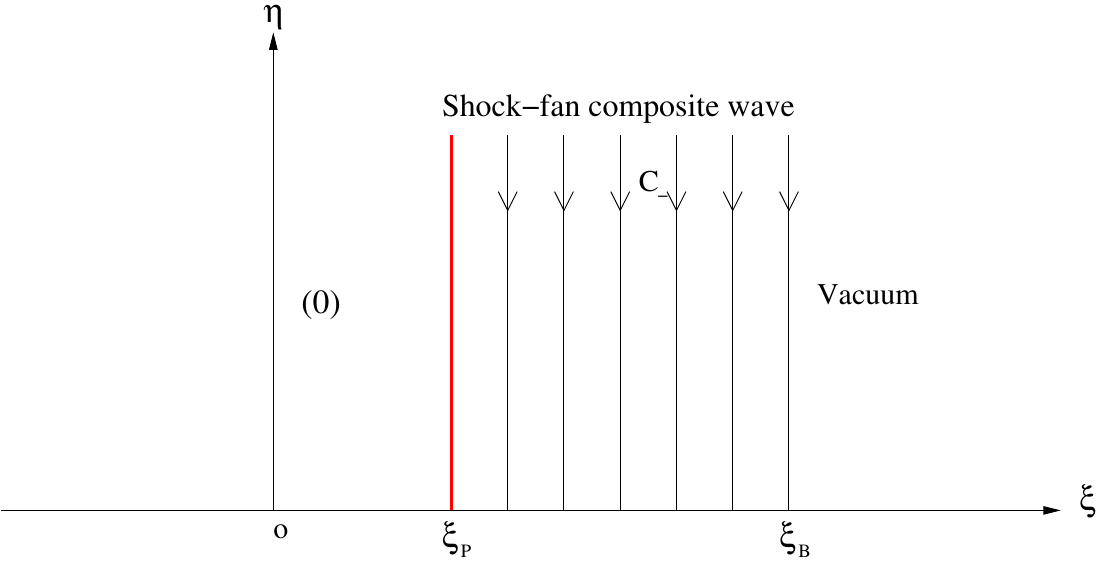}
\caption{\footnotesize A planar shock–fan composite wave with straight $C_{-}$ characteristic curves.}
\label{Fig1.4.2.1}
\end{center}
\end{figure}
%==========

By the result on the 1D Riemann problem for nonconvex gases (see \cite{Wen1}), we know that away from the corner
the supersonic oncoming flow $(0)$ is connected to the vacuum by a planar shock–fan  composite wave; see Figure \ref{Fig1.4.2.1}.
The shock wave of the composite wave  is located at $\xi=\xi_{_\mathrm{P}}$, where $\xi_{_\mathrm{P}}=u_0-\frac{\rho_{po}c_{po}}{\rho_0}$.
The fan wave of the composite wave can be
represented by
\begin{equation}\label{71501c}
(u, v, \tau)=(u_{-}, v_{-}, \tau_{-})(\xi), \quad\xi_{_\mathrm{P}}\leq\xi<\xi_{_\mathrm{B}}.
\end{equation}
Here,
$\xi_{_\mathrm{B}}=u_{-}(\xi_{_\mathrm{P}})+\displaystyle\int_{\tau_{po}}^{+\infty}\sqrt{-p'(\tau)}{\rm d}\tau$, $v_{-}(\xi)\equiv0$, and the functions $u_{-}(\xi)$ and $\tau_{-}(\xi)$ are determined by
\begin{equation}\label{71601}
u_{-}(\xi)-\int_{\tau_{po}}^{\tau_{-}(\xi)}\sqrt{-p'(\tau)}{\rm d}\tau=u_{-}(\xi_{_\mathrm{P}}), \quad \xi=u_{-}(\xi)-\sqrt{-\tau_{-}^2(\xi)p'(\tau_{-}(\xi))},
\end{equation}
and $(u_{-}, \tau_{-})(\xi_{_\mathrm{P}})=(\xi_{_\mathrm{P}}+c_{po}, \tau_{po})$.

Due to the existence of the sharp corner, there is an oblique wave emanating from the origin.
We shall show that at the corner an oblique shock or shock–fan composite wave consisting of a post-sonic shock and a centered simple wave may occur.

Assume that there is an oblique shock emanating from the origin and the state on its front side is $(0)$.
We use  $(u_b, v_b, \tau_b)$ to denote the state on the back side of the oblique shock, where $\tau_b\in (\tau_0, \tau_{po}]$. Then by (\ref{111306a}) we have
\begin{equation}\label{102401}
u_b=N_b\sin\phi_b+L_b\cos\phi_b, \quad v_b=-N_b\cos\phi_b+L_b\sin\phi_b,
\end{equation}
where
$$
N_b=\tau_b\sqrt{-\frac{2h(\tau_b)-2h(\tau_0)}{\tau_b^2-\tau_0^2}},\quad L_b=u_0\cos\phi_b,
$$
and $\phi_b=\arcsin \left(\frac{\tau_0}{u_0}\sqrt{-\frac{2h(\tau_b)-2h(\tau_0)}{\tau_b^2-\tau_0^2}}\right)$ is the inclination angle of the oblique shock.
This implies that $u_b$, $v_b$, and $\phi_b$ can be seen as functions of $\tau_b\in (\tau_0, \tau_{po}]$.

Let
\begin{equation}\label{4901}
\left\{
  \begin{array}{ll}
  \displaystyle \phi_{po}=\arcsin \Big(\frac{\rho_{po}c_{po}}{\rho_0u_0}\Big), \\[4pt]
  u_{po}=c_{po}\sin\phi_{po}+u_0\cos^2\phi_{po}, \\[4pt]
v_{po}=u_0\cos\phi_{po}\sin\phi_{po}-c_{po}\cos\phi_{po}.
  \end{array}
\right.
\end{equation}
Then we have $(\phi_b, u_b, v_b)=(\phi_{po}, u_{po}, v_{po})$ for $\tau_{b}=\tau_{po}$.

\begin{figure}[htbp]
\begin{center}
\quad \includegraphics[scale=0.5]{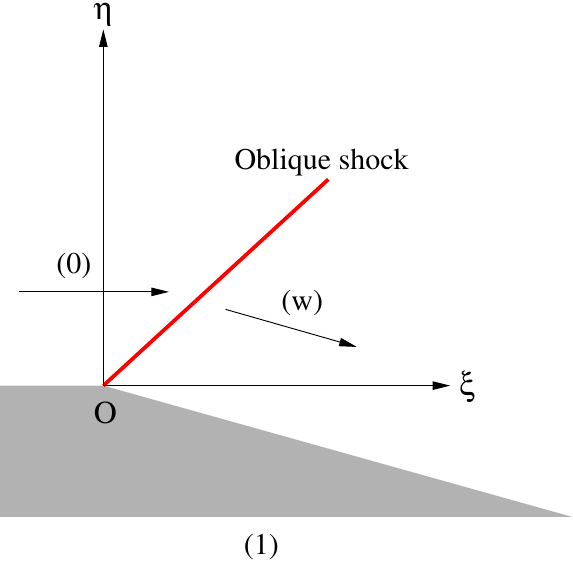}\qquad \includegraphics[scale=0.5]{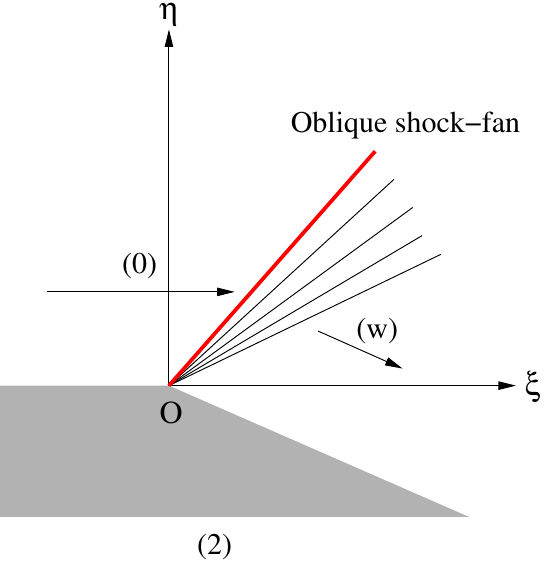}\qquad\includegraphics[scale=0.5]{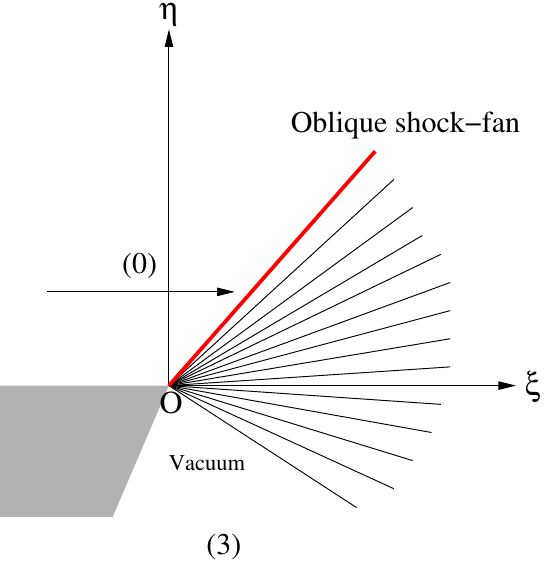}
\caption{\footnotesize Oblique wave at the corner for  $\tau_1^{i}<\tau_0\leq \tau_2^{i}$. }
\label{Fig1.4.2.2}
\end{center}
\end{figure}

We define
\begin{equation}\label{102408}
\sigma_{po}=\arctan\Big(\frac{v_{po}}{u_{po}}\Big).
\end{equation}
Then when
$\pi+\sigma_{po}<\theta_w<\pi$  there exists a $\tau_w\in (\tau_0, \psi_{po}(\tau_0))$ such that
$\arctan(\frac{v_w}{u_w})=\theta_w-\pi$, where  $(u_w, v_w)=(u_b(\tau_w), v_b(\tau_w))$,  and the oncoming flow $(0)$ can be deflected through an oblique shock wave to the constant state $(u_w, v_w, \tau_w)$. In this case, the flow near the corner can be represented by
\begin{equation}\label{41801}
(u, v, \tau)=\left\{
               \begin{array}{ll}
                 (u_0, 0, \tau_0), & \hbox{$\theta>\phi_w$;} \\[2pt]
                 (u_w, v_w, \tau_w), & \hbox{$\theta_w-\pi<\theta<\phi_w$,}
               \end{array}
             \right.
\end{equation}
where $\phi_w=\phi_b(\tau_w)$ and $(u_w, v_w)=(u_b, v_b)(\tau_w)$;
see Figure \ref{Fig1.4.2.2} (1).

%=============
%==========

When $\theta_w<\pi+\sigma_{po}$, at the corner the oncoming flow $(0)$ can be deflected through an oblique shock–fan composite wave consisting of an oblique post-sonic shock and a centered simple wave.
The inclination angle of the oblique shock is $\phi_{po}$, and the states on the front-side and back-side of the shock are $(0)$ and $(u_{po}, v_{po}, \tau_{po})$, respectively. In order to obtain the centered simple wave adjacent to the oblique post-sonic shock, we solve (\ref{102304}) with data
\begin{equation}\label{102402}
(\bar{q}, \bar{\tau}, \bar{\sigma})(\phi_{po})=(q_{po}, \tau_{po}, \sigma_{po}),
\end{equation}
where $q_{po}=\sqrt{u_{po}^2+v_{po}^2}$.

Let
$$
\alpha_v=\sigma_{po}-\int_{q_{po}}^{{q}_{_{lim}}}\frac{\sqrt{q^2-\hat{c}^2(q)}}{q\hat{c}(q)}{\rm d}q,
$$
where $\hat{c}(q)$ and $q_{_{lim}}$ are the same as those defined in (\ref{3506}).
Then the initial value problem (\ref{102304}), (\ref{102402}) admits a solution for $\theta\in (\alpha_v, \phi_{po})$. Moreover, the solution satisfies
(\ref{81801}) for $\alpha_v< \theta<\phi_{po}$.

 When
 $\pi+\alpha_v<\theta_w<\pi+\sigma_{po}$, there exists an $\alpha_w\in (\alpha_v, \phi_{po})$ such that
$\bar{v}(\alpha_w)=\bar{u}(\alpha_w)\tan\theta_w$. In this case,
the oncoming flow $(0)$ is deflected through an oblique shock–fan composite wave to a constant state $(u_w, v_w, \tau_w)$, where $(u_w, v_w, \tau_w)=(\bar{u}, \bar{v}, \bar{\tau})(\alpha_w)$. The flow near the corner can be represented by
\begin{equation}\label{8230a}
(u, v, \tau)=\left\{
               \begin{array}{ll}
                 (u_0, 0, \tau_0), & \hbox{$\theta>\phi_{po}$;} \\[2pt]
                 (\bar{u}, \bar{v}, \bar{\tau})(\theta), & \hbox{${\alpha}_w<\theta<\phi_{po}$;} \\[2pt]
                 (u_w, v_w, \tau_w), & \hbox{$\theta_w-\pi<\theta<{\alpha}_w$;}
               \end{array}
             \right.
\end{equation}
see Figure \ref{Fig1.4.2.2} (2).
 When $\theta_w<\pi+\alpha_v$, the oncoming flow $(0)$ is deflected through an oblique shock–fan composite wave to a vacuum. The flow near the corner can be represented by
\begin{equation}\label{8230}
(u, v, \tau)=\left\{
               \begin{array}{ll}
                 (u_0, 0, \tau_0), & \hbox{$\theta>\phi_{po}$;} \\[2pt]
                 (\bar{u}, \bar{v}, \bar{\tau})(\theta), & \hbox{$\alpha_v<\theta<\phi_{po}$;} \\[2pt]
                 \text{vacuum}, & \hbox{$\theta_w-\pi<\theta<\alpha_v$;}
               \end{array}
             \right.
\end{equation}
see Figure \ref{Fig1.4.2.2} (3).

\subsection{Elementary waves for $\tau_1^{e}\leq\tau_0\leq \tau_1^{i}$}
As in Section 3.4, when $\tau_1^{e}\leq\tau_0\leq \tau_1^{i}$,
away from the corner
the constant state  is connected to the vacuum by a planar shock–fan  composite wave.
The shock wave of the planar composite wave is located at $\xi=\xi_{_\mathrm{P}}$, and the fan wave of the planar composite wave can be
represented by (\ref{71501c}) and (\ref{71601}), where  $\tau_{po}$, $c_{po}$, and $\rho_{po}$ are the same as those defined in Section 3.4.

%=============
\begin{figure}[htbp]
\begin{center}
\includegraphics[scale=0.43]{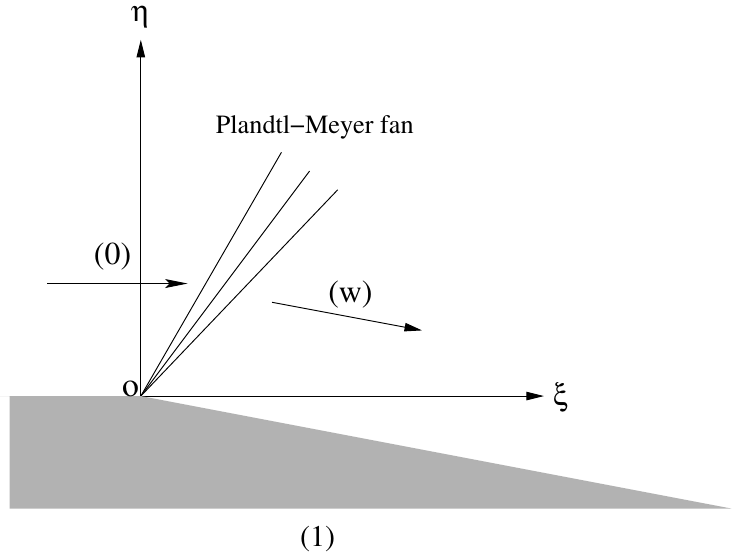} \qquad\quad\includegraphics[scale=0.43]{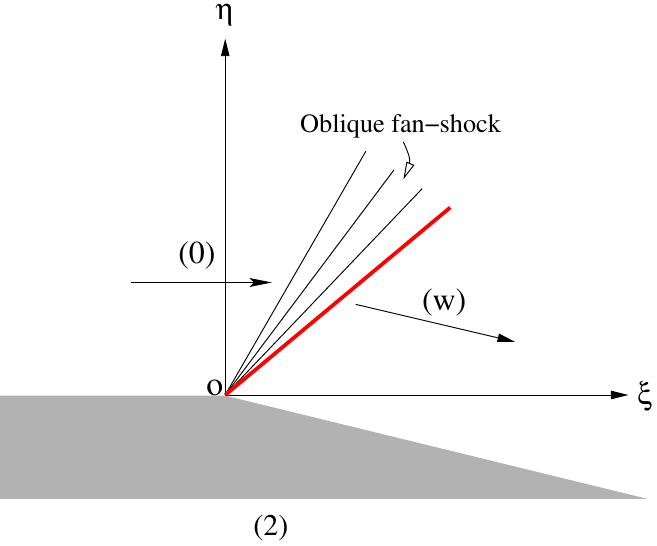}
\\
\includegraphics[scale=0.43]{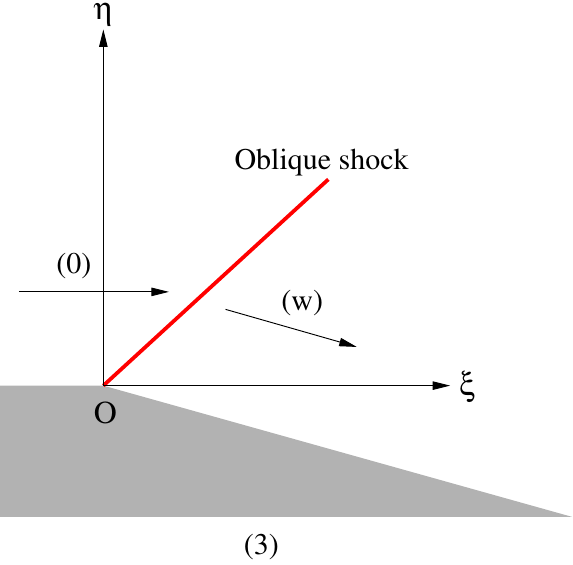}\qquad \includegraphics[scale=0.43]{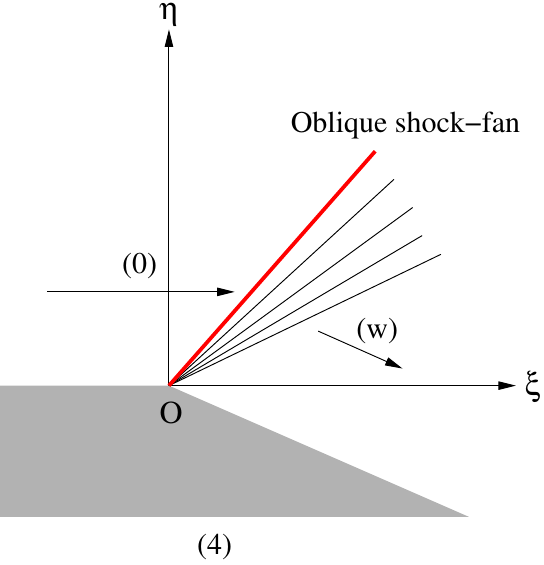}\qquad\includegraphics[scale=0.43]{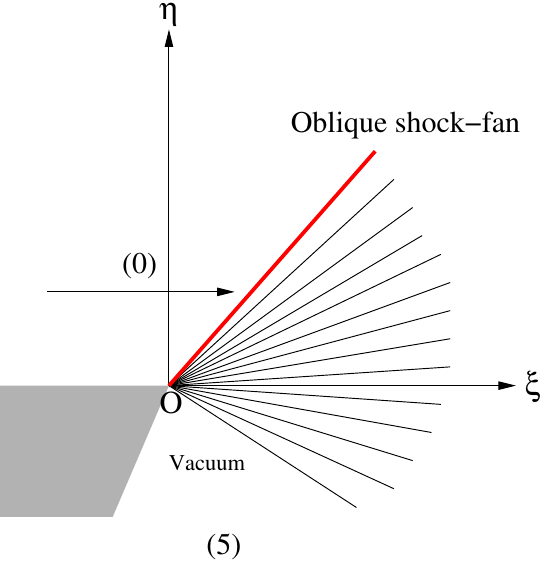}
\caption{\footnotesize Oblique wave at the corner for  $\tau_1^{e}<\tau_0\leq \tau_1^{i}$.}
\label{Fig4.1.3.1}
\end{center}
\end{figure}
%==========

We next discuss the oblique wave at the corner.
In contrast to the case of $\tau_1^i<\tau_0<\tau_2^i$, two additional oblique wave configurations may occur at the corner.

 We first solve (\ref{102304}) with data
\begin{equation}\label{102407}
(\bar{q}, \bar{\tau}, \bar{\sigma})(\alpha_0)=(u_0, \tau_0, 0).
\end{equation}
 Since $\tau_1^e<\tau_0<\tau_1^i$, there exists an $\alpha_i<\alpha_0$ such that the
 initial value problem (\ref{102304}), (\ref{102407}) admits a solution for $\theta\in [\alpha_i, \alpha_{0})$.
The solution satisfies (\ref{81801}) for $\theta\in [\alpha_i, \phi_{0})$.
Moreover, $\bar{\tau}(\alpha_i)=\tau_1^i$.

Let
 \begin{equation}\label{81802}
\sigma_i=\arctan \left(\frac{\bar{v}(\alpha_i)}{\bar{u}(\alpha_i)}\right).
\end{equation}
When $\theta_w>\pi+\sigma_i$, there exists an $\alpha_w \in (\alpha_i, \alpha_0)$ such that $\theta_w-\pi=\arctan\left(\frac{\bar{v}(\alpha_w)}{\bar{u}(\alpha_w)}\right).$
Hence, at the corner
the oncoming flow can be deflected through a centered simple wave to a constant state $(u_w, v_w, \tau_w)$, where $(u_w, v_w, \tau_w)= (\bar{u}, \bar{v}, \bar{\tau})(\alpha_w)$. In this case,
the flow near the corner can be represented by
\begin{equation}\label{100901}
(u, v, \tau)=\left\{
               \begin{array}{ll}
                 (u_0, 0, \tau_0), & \hbox{$\theta\geq\alpha_0$;} \\[2pt]
                 (\bar{u}, \bar{v}, \bar{\tau})(\theta), & \hbox{$\alpha_w<\theta<\alpha_0$;} \\[2pt]
                (u_w, v_w, \tau_w), & \hbox{$\theta_w-\pi<\theta<{\alpha}_w$;}
               \end{array}
             \right.
\end{equation}
 see Figure \ref{Fig4.1.3.1} (1).

 We next consider the possibility of the occurrence of an oblique fan–shock composite wave consisting of a centered simple wave and a pre-sonic shock at the corner.
We assume that the inclination angle of the pre-sonic shock of the composite wave is
$\phi\in (\alpha_i, \alpha_0]$.
Then the state on the front side of the pre-sonic shock is
$(\bar{u}, \bar{v}, \bar{\tau})(\phi)$.
We denote by $(u_s, v_s, \tau_s)$ the  state on the back side of the pre-sonic shock. Then, we have
\begin{equation}\label{81803}
\tau_s=\psi_{pr}(\bar{\tau}(\phi)), \quad u_{s}=N_{{s}}\sin\phi+L_{s}\cos\phi, \quad v_{s}=-N_{s}\cos\phi+L_{s}\sin\phi,
\end{equation}
where
$$
N_{s}=\frac{\psi_{pr}(\bar{\tau}(\phi))}{\bar{\tau}(\phi)}\big(\bar{u}(\phi)
\sin\phi-\bar{v}(\phi)\cos\phi\big)\quad \text{and}\quad L_{s}=\bar{u}(\phi)\cos\phi+\bar{v}(\phi)\sin\phi.
$$

Let
$$
\sigma_{_M}=\sup\limits_{\phi\in (\alpha_i, \alpha_0)}\arctan\bigg(\frac{v_{s}(\phi)}{u_{s}(\phi)}\bigg)\quad \text{and}\quad \sigma_{m}=\inf\limits_{\phi\in (\alpha_i, \alpha_0)}\arctan\bigg(\frac{v_{s}(\phi)}{u_{s}(\phi)}\bigg).
$$
When $\theta_w-\pi\in(\sigma_{m}, \sigma_{_M})$, there exists a $\phi_{pr}\in(\alpha_i, \alpha_0)$ such that $\theta_w-\pi=\arctan\big(\frac{v_{s}(\phi_{pr})}{u_{s}(\phi_{pr})}\big)$.
Hence, at the corner the oncoming flow $(0)$ can be deflected by a fan–shock composite wave.
In this case the flow near the corner can be represented by
\begin{equation}\label{2662601}
(u, v, \tau)=\left\{
               \begin{array}{ll}
                 (u_0, 0, \tau_0), & \hbox{$\theta\geq\alpha_0$;} \\[2pt]
                 (\bar{u}, \bar{v}, \bar{\tau})(\theta), & \hbox{$\phi_{pr}<\theta<\alpha_0$;} \\[2pt]
                (u_w, v_w, \tau_w), & \hbox{$\theta_w-\pi<\theta<{\phi}_{pr}$,}
               \end{array}
             \right.
\end{equation}
where  $(u_w, v_w, \tau_w)=(u_s, v_s, \tau_s)(\phi_{pr})$.

Let
$
\sigma_{*}=\arctan\Big(\frac{v_{s}(\alpha_0)}{u_{s}(\alpha_0)}\Big)
$,
 and let $\sigma_{po}$ be as defined in (\ref{102408}). Then when $\pi+\sigma_{po}<\theta_w<\pi+\sigma_{*}$, the elementary wave near the corner can be an oblique shock; see Figure \ref{Fig4.1.3.1} (3); when $\theta_w<\pi+\sigma_{po}$ the oblique wave at the corner can be an oblique shock–fan composite wave; see Figure \ref{Fig4.1.3.1} (4) and (5). The discussion is similar to that in Section 3.4.

\subsection{Elementary waves for $\tau_0<\tau_1^{e}$}
\begin{figure}[htbp]
\begin{center}
  \includegraphics[scale=0.5]{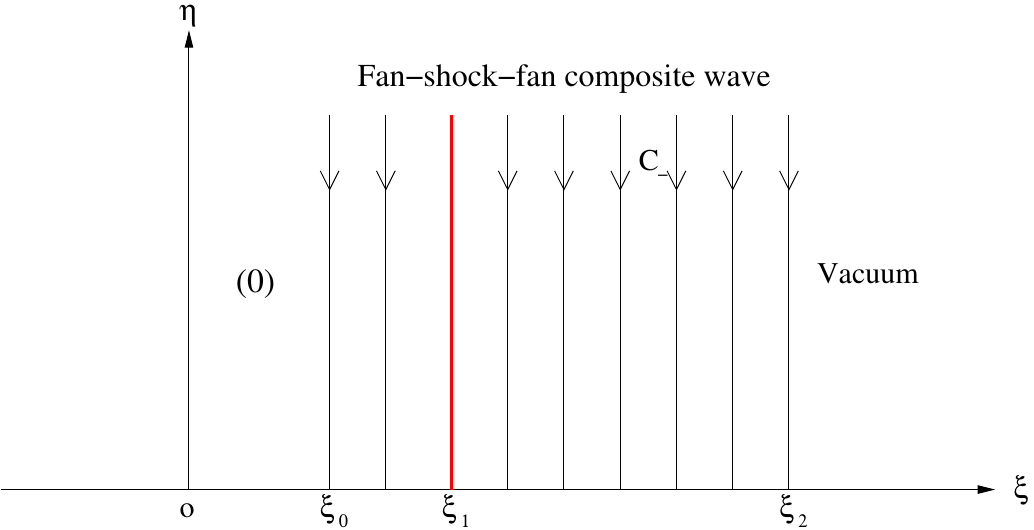}
\caption{\footnotesize A planar fan–shock–fan composite wave with straight $C_{-}$ characteristic curves.}
\label{Fig1.4.4.1}
\end{center}
\end{figure}

By the result of the 1D Riemann problem for nonconvex gases (see \cite{Wen1}), we know that 
the oncoming flow $(0)$ is connected to the vacuum by  a planar fan–shock–fan composite wave away from the corner; see Figure \ref{Fig1.4.4.1}. This planar fan–shock–fan composite wave can be represented by
\begin{equation}\label{pfsf}
(u, v, \tau)=\left\{
               \begin{array}{ll}
                 (u_0, 0, \tau_0), & \hbox{$\xi<\xi_0$;} \\[2pt]
                 (\hat{u}_l, \hat{v}_l, \hat{\tau}_l)(\xi), & \hbox{$\xi_0\leq \xi\leq \xi_1$;} \\[2pt]
                 (\hat{u}_r, \hat{v}_r, \hat{\tau}_r)(\xi), & \hbox{$\xi_1\leq \xi\leq \xi_2$;} \\[2pt]
                 \text{vacuum}, & \hbox{$\xi>\xi_2$.}
               \end{array}
             \right.
\end{equation}
Here, $\xi_{0}=u_0-c_0$; $\hat{v}_l(\xi)\equiv0$; $\hat{v}_r(\xi)\equiv0$; $\hat{u}_l(\xi)$ and $\hat{\tau}_l(\xi)$ are determined by
$$
\hat{u}_l(\xi)-\int_{\tau_0}^{\hat{\tau}_l(\xi)}\sqrt{-p'(\tau)}{\rm d}\tau=u_0\quad \text{and} \quad \xi=\hat{u}_l(\xi)-c(\hat{\tau}_l(\xi));
$$
$\xi_{1}$ is determined by $\hat{\tau}_l(\xi)=\tau_1^e$;  $\hat{u}_r(\xi)$ and $\hat{\tau}_r(\xi)$ are determined by
$$
\hat{u}_r(\xi)-\int_{\tau_2^e}^{\hat{\tau}_l(\xi)}\sqrt{-p'(\tau)}{\rm d}\tau=\hat{u}_2\quad \text{and} \quad \xi=\hat{u}_r(\xi)-c(\hat{\tau}_r(\xi));
$$
$\hat{u}_1=\hat{u}_{l}(\xi_1)$; $\hat{u}_2=\frac{\tau_2^e}{\tau_1^e}(\hat{u}_1-{\xi_1})+\xi_1$;  $\xi_{2}=\hat{u}_2+\displaystyle\int_{\tau_2^e}^{+\infty}\sqrt{-p'(\tau)}{\rm d}\tau$.

%=============
%==========

We next discuss the oblique wave at the corner.
We shall show that at the corner a centered simple wave, a fan–shock composite wave, or a fan–shock–fan composite wave may occur.

\begin{figure}[htbp]
\begin{center}
 \includegraphics[scale=0.48]{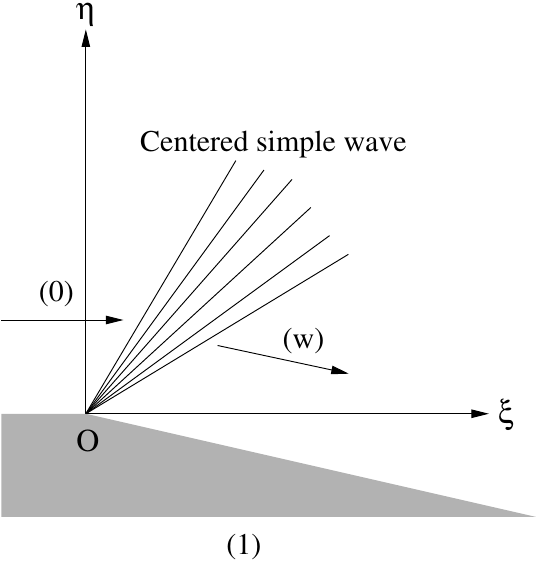}\qquad\qquad \includegraphics[scale=0.48]{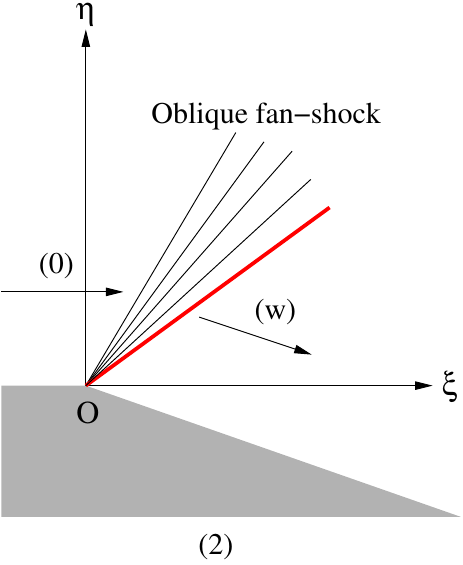}\\[2pt] \includegraphics[scale=0.48]{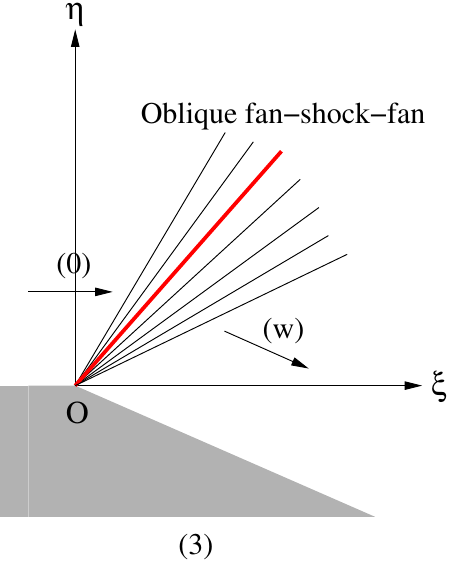}\qquad\qquad \includegraphics[scale=0.48]{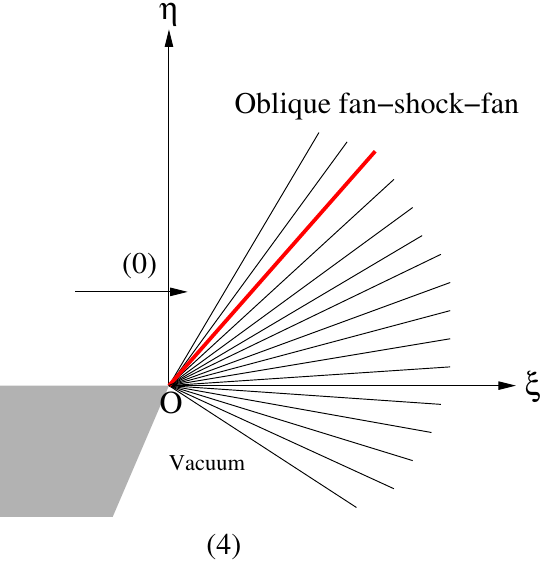}
\caption{\footnotesize Oblique wave at the sharp corner for  $\tau_1^{e}<\tau_0\leq \tau_1^{i}$.}
\label{Fig1.4.4.2}
\end{center}
\end{figure}

We first solve the initial value problem (\ref{102304}), (\ref{102407}).
 Since $\tau_0<\tau_1^e$, there exist $\alpha_i$ and $\alpha_e$ with $\alpha_i<\alpha_e<\alpha_0$ such that $\bar{\tau}(\alpha_i)=\tau_1^i$ and $\bar{\tau}(\alpha_e)=\tau_1^e$.
Let $\sigma_i$ be as defined in (\ref{81802}), and let
$$
\bar{\sigma}_{_M}=\sup\limits_{\phi\in (\alpha_i, \alpha_e)}\arctan\bigg(\frac{v_{s}(\phi)}{u_{s}(\phi)}\bigg)\quad \text{and}\quad \bar{\sigma}_{m}=\inf\limits_{\phi\in (\alpha_i, \alpha_e)}\arctan\bigg(\frac{v_{s}(\phi)}{u_{s}(\phi)}\bigg),
$$
where the functions $u_s(\phi)$ and $v_s(\phi)$ are the same as those defined in (\ref{81803}).
As in Section 3.5, when $\theta_w>\pi+\sigma_i$ the oncoming flow can be deflected by a centered simple wave;  when $\theta_w-\pi\in(\bar{\sigma}_{m}, \bar{\sigma}_{_M})$ the oncoming flow can be deflected by a fan–shock composite wave; see Figure \ref{Fig1.4.4.2} (1) and (2).

%=============
%==========

 We now consider the possibility of the occurrence of an oblique fan–shock–fan composite wave consisting of a centered simple wave, an oblique double-sonic shock, and another centered simple wave at the corner.
 The centered simple wave adjacent to the constant state $(0)$ can be represented by
 $$
 (u, v, \tau)=(\bar{u}, \bar{v}, \bar{\tau})(\theta), \quad \alpha_e\leq \theta\leq \alpha_0.
 $$
 The inclination angle of the oblique double-sonic shock is $\alpha_e$, and the states on the front side and backside of the shock
 are $(\bar{u}, \bar{v}, \bar{\tau})(\alpha_e)$ and $(u_s, v_s, \tau_s)(\alpha_e)$, respectively.
Since $\bar{\tau}(\alpha_e)=\tau_1^e$, we have $\tau_s(\alpha_e)=\tau_2^e$.

 In order to construct the other centered simple wave adjacent to the shock, we solve (\ref{102304}) with data
 \begin{equation}\label{81805}
 (\bar{q}, \bar{\tau}, \bar{\sigma})(\alpha_e)=(q_2^e, \tau_2^e, \sigma_2^e),
 \end{equation}
 where $q_2^e=\sqrt{u_s^2(\alpha_e)+v_s^2(\alpha_e)}$ and $\sigma_2^e=\arctan\left(\frac{v_s(\alpha_e)}{u_s(\alpha_e)}\right)$.
Let
$$
\alpha_v=\sigma_{2}^e-\int_{q_{2}^e}^{{q}_{_{lim}}}\frac{\sqrt{q^2-\hat{c}^2(q)}}{q\hat{c}(q)}{\rm d}q,
$$
where $\hat{c}(q)$ and $q_{_{lim}}$ are the same as those defined in (\ref{3506}).
Then the initial value problem (\ref{102304}), (\ref{81805}) admits a solution for $\theta\in (\alpha_v, \alpha_e)$. Moreover, the solution satisfies
(\ref{81801}).

For convenience we use the subscripts `$l$' and `$r$' to denote the centered simple waves upstream and downstream of the oblique double-sonic shock, respectively.
When
 $\alpha_v<\theta_w-\pi<\sigma_2^e$, there exists an $\alpha_w\in (\alpha_v, \alpha_e)$ such that
$\bar{v}_r(\alpha_w)=\bar{u}_r(\alpha_w)\tan\theta_w$. In this case, at the sharp corner
the oncoming flow $(0)$ can be deflected through an oblique fan–shock–fan composite wave to a constant state $(u_w, v_w, \tau_w)$, where $(u_w, v_w, \tau_w)=(\bar{u}_r, \bar{v}_r, \bar{\tau}_r)(\alpha_w)$. The flow near the corner can be represented by
\begin{equation}\label{ofsf1}
(u, v, \tau)=\left\{
               \begin{array}{ll}
                 (u_0, 0, \tau_0), & \hbox{$\theta>\alpha_0$;} \\[2pt]
                 (\bar{u}_l, \bar{v}_l, \bar{\tau}_l)(\theta), & \hbox{${\alpha}_e<\theta\leq \alpha_0$;} \\[2pt]
                  (\bar{u}_r, \bar{v}_r, \bar{\tau}_r)(\theta), & \hbox{${\alpha}_w\leq \theta<\alpha_e$;} \\[2pt]
                 (u_w, v_w, \tau_w), & \hbox{$\theta_w-\pi<\theta<{\alpha}_w$;}
               \end{array}
             \right.
\end{equation}
see Figure \ref{Fig1.4.4.2} (3).
 When $\theta_w<\pi+\alpha_v$, the oncoming flow $(0)$ can be deflected through an oblique fan–shock–fan composite wave to the vacuum. The flow near the corner can be represented by
\begin{equation}\label{ofsf2}
(u, v, \tau)=\left\{
               \begin{array}{ll}
                 (u_0, 0, \tau_0), & \hbox{$\theta>\alpha_0$;} \\[2pt]
                (\bar{u}_l, \bar{v}_l, \bar{\tau}_l)(\theta), & \hbox{${\alpha}_e<\theta\leq \alpha_0$;} \\[2pt]
                  (\bar{u}_r, \bar{v}_r, \bar{\tau}_r)(\theta), & \hbox{${\alpha}_v\leq \theta<\alpha_e$;} \\[2pt]
                 \text{vacuum}, & \hbox{$\theta_w-\pi<\theta<\alpha_v$;}
               \end{array}
             \right.
\end{equation}
see Figure \ref{Fig1.4.4.2} (4).

\subsection{Uniqueness of the oblique wave at the corner}
It has been shown that when $u_0$ is sufficiently large, the oblique wave at the corner can be uniquely determined by the wedge angle $\theta_w$; see Lai \cite{Laip2}.

\section{Interaction between a planar fan wave and a centered simple wave}

In this section, we study the 2D Riemann problem (\ref{Potential}), (\ref{initial})--(\ref{boundary}) for $\tau_0>\tau_2^e$.

\subsection{Characteristic boundaries of the wave interaction region} 
\begin{figure}[htbp]
\begin{center}
\includegraphics[scale=0.48]{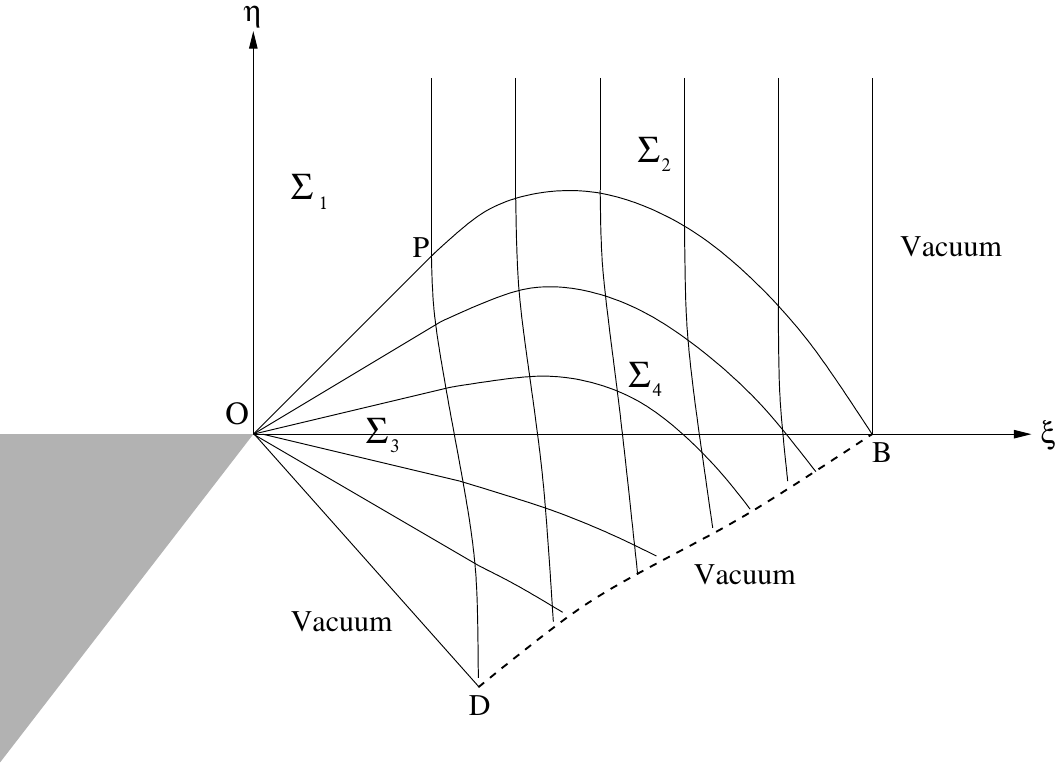}
\caption{\footnotesize Interaction between two fan waves. }
\label{Fig4.1}
\end{center}
\end{figure}

By the result of Section 3.3, in the $(\xi, \eta)$-plane there is a planar rarefaction wave with straight $C_{-}$ characteristic lines away from the origin and a centered simple wave  with straight $C_{+}$ characteristic lines emanating from the origin. The two elementary waves start to interact with each other at the point $\mathrm{P}(\xi_{_\mathrm{P}}, \xi_{_\mathrm{P}}\tan\alpha_0)$, where $\xi_{_\mathrm{P}}=u_0-c_0$ and $\alpha_0=\arcsin(\frac{c_0}{u_0})$. Through $\mathrm{P}$ we draw a forward $C_{+}$ cross characteristic curve of the planar rarefaction wave and a forward $C_{-}$ cross characteristic curve of the centered simple wave; see Figure \ref{Fig4.1}.

In order to construct a global-in-time classical solution to the 2D Riemann problem (\ref{Potential}), (\ref{initial})--(\ref{boundary}) for $\tau_0>\tau_2^e$,
we make the following assumption:
\begin{description}
  \item[($h_{11}$)] $\alpha_0+\frac{\pi}{2}>2\Theta(\tau_0)$.
       % \item[(H$_{11}$)]  $\theta>\alpha_0-\pi+2\Theta(\bar{\tau}(\theta))$ for $\alpha_v\leq \theta\leq\alpha_0$.
\end{description}

%=============
%==========

\begin{lem}\label{lem3.1}
Under assumption ($h_{11}$), the forward $C_{+}$ characteristic curve emanating from the point $\mathrm{P}$ crosses the entire planar rarefaction wave and terminates at the point $\mathrm{B}(\xi_{_\mathrm{B}}, 0)$, where $\xi_{_B}=u_0+\int_{\tau_0}^{+\infty}\sqrt{-p'(\tau)}{\rm d}\tau$. Moreover, along this $C_{+}$ cross characteristic curve there holds
\begin{equation}\label{7201e}
2\Theta(\tau)-\frac{\pi}{2}<\alpha\leq \alpha_0.
\end{equation}
\end{lem}
\begin{proof}
For convenience, we use $C_{+}^{\mathrm{P}}$ to denote the $C_{+}$ characteristic curve emanating from $\mathrm{P}$.

By the definition of $\beta$ we have
\begin{equation}\label{71602}
\beta(\xi, \eta)=-\frac{\pi}{2}\quad \text{on}\quad C_{+}^{\mathrm{P}}.
\end{equation}

From (\ref{210cqo}) and (\ref{71601}), we have
\begin{equation}\label{7102}
 A~=~\arcsin\Big(\frac{c}{q}\Big)~=~\arctan\Bigg(\frac{\tau_{-}(\xi)\sqrt{-p'(\tau_{-}(\xi))}}{\eta}\Bigg)\quad \text{on}\quad C_{+}^{\mathrm{P}},
\end{equation}
where the function $\tau_{-}(\xi)$ is defined in (\ref{71601A}).
Moreover, we have
$$
A(\mathrm{P})=\frac{\alpha(\mathrm{P})-\beta(\mathrm{P})}{2}=\frac{\alpha_0}{2}+\frac{\pi}{4}.
$$

Then the $C_{+}$ characteristic curve $C_{+}^{\mathrm{P}}$
can be represented by a function $\eta=\eta_{+}(\xi)$ which satisfies
\begin{equation}\label{7101}
\left\{
  \begin{array}{ll}
    \displaystyle\frac{{\rm d}\eta_{+}}{{\rm d}\xi}=\tan\bigg(\underbrace{-\frac{\pi}{2}+2\arctan\Big(\frac{\tau_{-}(\xi)\sqrt{-p'(\tau_{-}(\xi))}}{\eta_{+}}\Big)}_{\alpha}\bigg),\quad  \xi_{_P}<\xi<\xi_{_B},  \\[12pt]
  \eta_{+}(\xi_{_\mathrm{P}})=\xi_{_\mathrm{P}}\tan\alpha_0.
  \end{array}
\right.
\end{equation}

From (\ref{71602}), one has $\bar{\partial}_{+}\beta=0$ along $C_{+}^{\mathrm{P}}$.
Hence, by (\ref{7}), (\ref{71602}), and $\tau_0>\tau_2^i$, one has
\begin{equation}\label{62201}
\bar{\partial}_{+}\rho=-\frac{2c\rho^4\sin2 A}{p''(\tau)}<0\quad \text{along}\quad C_{+}^{\mathrm{P}}.
\end{equation}
Inserting this into (\ref{10}) we get
$$
c\bar{\partial}_{+}A=(\tan^2 \Theta-\tan^2 A)\sin^2A\cos^2A\quad \text{along}\quad C_{+}^{\mathrm{P}}.
$$
Consequently, by assumption {\bf (E2)} we have
\begin{equation}\label{3501}
c\bar{\partial}_{+}(A-\Theta)=(\tan^2 \Theta-\tan^2 A)\sin^2A\cos^2A-\frac{c\omega'(\tau)\bar{\partial}_{+}\tau}{2\sqrt{\omega(\tau)}(1+\omega(\tau))}\quad  \text{along}\quad C_{+}^{\mathrm{P}}.
\end{equation}
Combining this with (\ref{62201}) and $A(\mathrm{P})>\Theta(\tau_0)$, we get
\begin{equation}\label{121301}
A-\Theta>0\quad  \text{on}\quad C_{+}^{\mathrm{P}}.
\end{equation}
In addition, by (\ref{10}) we also have
$$
c\bar{\partial}_{+}\alpha=-\frac{p''(\tau)}{4c\rho^4}\Omega\sin2 A\bar{\partial}_{+}\rho<0\quad  \text{along}\quad C_{+}^{\mathrm{P}}.
$$
Thus we get (\ref{7201e}).

By (\ref{7201e}) and (\ref{7102}) one has $\eta_{+}(\xi_{_B})=0$.
This completes the proof of Lemma \ref{lem3.1}.
\end{proof}

By the result of Section 3.3 we have that the complete centered simple wave adjacent to the constant state $(0)$ can be represented by
$$
(u, v, \tau)=(\bar{u}, \bar{v}, \bar{\tau})(\theta), \quad \alpha_v\leq \theta\leq\alpha_0,
$$
where $\bar{u}(\theta)=\bar{q}(\theta)\cos ( \bar{\sigma}(\theta)))$, $\bar{v}(\theta)=\bar{q}(\theta)\sin( \bar{\sigma}(\theta)))$, and $(\bar{q}(\theta), \bar{\tau}(\theta), \bar{\sigma}(\theta))$ is the solution of the problem (\ref{102304}), (\ref{102305}).  Moreover, $\alpha=\theta$ in the centered simple wave region.

\begin{lem}\label{lem3.3}
Under assumption ($h_{11}$),
the forward $C_-$ characteristic emanating from the point $\mathrm{P}$ can cross the whole centered simple wave and end up at the point
$
\mathrm{D}\big(\bar{r}(\alpha_v)\cos \alpha_{v}, ~~\bar{r}(\alpha_v)\sin \alpha_{v}\big),
$
where the function $\bar{r}(\theta)=\sqrt{\bar{q}^2(\theta)+\bar{\tau}^2(\theta)p'(\bar{\tau}(\theta))}$.
Moreover, there exists a small $\delta_{0}>0$ such that along this $C_{-}$ cross characteristic curve,
\begin{equation}\label{7106}
\delta_{0}~\leq~ A~\leq~ \frac{\alpha_0}{2}+\frac{\pi}{4}.
\end{equation}
\end{lem}
\begin{proof}
For convenience,
we use $C_{-}^{\mathrm{P}}$ to denote the forward $C_{-}$ cross characteristic curve  emanating from $\mathrm{P}$. Assume that $C_{-}^{\mathrm{P}}$ can be described by a one-parameter form: $\xi=\xi_{-}(\theta)=r_{-}(\theta)\cos\theta$, $\eta=\eta_{-}(\theta)=r_{-}(\theta)\sin\theta$.

Since $C_{-}^{\mathrm{P}}$  is a $C_{-}$ characteristic curve, we have $$\cos\big(\alpha(\xi_{-}, \eta_{-})-2A(\xi_{-}, \eta_{-})\big) \frac{{\rm d}\eta_{-}}{{\rm d}\theta}-\sin\big(\alpha(\xi_{-}, \eta_{-})-2A(\xi_{-}, \eta_{-})\big) \frac{{\rm d}\xi_{-}}{{\rm d}\theta}~=~0,$$
where
\begin{equation}\label{32401}
A(\xi_{-}(\theta),\eta_{-}(\theta))=\arctan \left(\frac{\bar{\tau}(\theta)\sqrt{-p'(\bar{\tau}(\theta))}}{\bar{r}(\theta)-\sqrt{\xi_{-}^2(\theta)+\eta_{-}^2(\theta)}}\right).
\end{equation}
Combining this with $\alpha=\theta$ on $C_{-}^{\mathrm{P}}$, we have
that $r_{-}(\theta)$ satisfies
\begin{equation}\label{62302}
\left\{
  \begin{array}{ll}
   r_{-}'(\theta)=-\cot \big(2A(\xi_{-}(\theta), \eta_{-}(\theta))\big) r_{-}(\theta), \quad \alpha_v<\theta<\alpha_0,\\[8pt]
   \displaystyle r_{-}(\alpha_0)=\frac{u_0-c_0}{\cos\alpha_0}.
  \end{array}
\right.
\end{equation}

 Since
\begin{equation}\label{80801}
\Theta(\tau_0)<A\big(\xi_{-}(\alpha_0), \eta_{-}(\alpha_0)\big)=\frac{\alpha_0}{2}+\frac{\pi}{4}<\frac{\pi}{2},
\end{equation}
the initial value problem (\ref{62302}) admits a local solution.

Let
\begin{equation}\label{new81901}
\begin{aligned}
&A_{-}(\theta)=A\big(\xi_{-}(\theta), \eta_{-}(\theta)\big),\quad \Omega_{-}(\theta)=\tan^2(\Theta_{-}(\theta))-\tan^2(A_{-}(\theta)),\quad \\& \Theta_{-}(\theta)=\arctan\sqrt{\omega(\bar{\tau}(\theta))},  \quad \text{and}\quad
\beta_{-}(\theta)=\theta-2A_{-}(\theta).
\end{aligned}
\end{equation}
Then, by (\ref{8}) we have
\begin{equation}\label{62401}
\beta_{-}'(\theta)=\frac{p''(\bar{\tau}(\theta))\Omega_{-}(\theta)}{4\bar{c}^2(\theta)\bar{\rho}^4(\theta)}\sin (2A_{-}(\theta))\bar{\rho}'(\theta),
\end{equation}
where $\bar{\rho}(\theta)=\bar{\tau}^{-1}(\theta)$ and $\bar{c}(\theta)=\sqrt{-\bar{\tau}^2(\theta)p'(\bar{\tau}(\theta))}$.

From the second equation of (\ref{102304}) and $\bar{\tau}(\alpha_0)=\tau_0>\tau_2^i$, we get
\begin{equation}\label{80802}
\bar{\rho}'(\theta)>0\quad \text{and}\quad \bar{\Theta}_{-}'(\theta)>0 \quad \text{for}\quad \alpha_v<\theta<\alpha_0.
\end{equation}

From (\ref{80801}), (\ref{62401}),(\ref{80802}), and ($h_{11}$) one knows
$$
A_{-}'(\theta)=\frac{1}{2}\big(1-\beta_{-}'(\theta)\big)>0
\quad
\text{suppose}\quad A_{-}(\theta)=\frac{\alpha_0}{2}+\frac{\pi}{4}.$$
Then by $A_{-}(\alpha_0)=\frac{\alpha_0}{2}+\frac{\pi}{4}$  and  (\ref{32401}) one knows that the problem (\ref{62302}) has a solution for $\alpha_v<\alpha<\alpha_0$.
Moreover,
\begin{equation}\label{3601}
0<A_{-}(\theta)<\frac{\alpha_0}{2}+\frac{\pi}{4}\quad \text{for}\quad  \alpha<\alpha_0.
\end{equation}
 This implies that the $C_{-}$ characteristic emanating from $\mathrm{P}$ can cross the whole centered simple wave.

In what follows,
we are going to prove $$\inf\limits_{\theta\in (\alpha_v, \alpha_0)}A_{-}(\theta)>0.$$
We prove this by contradiction.
Suppose $\inf\limits_{\theta\in (\alpha_v, \alpha_0)}A_{-}(\theta)=0$. Then in view of assumptions {\bf (E1)} and {\bf (E2)},
one of the following two cases must occur.
\begin{itemize}
  \item[1.] There exists an $\alpha_m\in (\alpha_v, \alpha_0)$ such that $A_{-}(\theta)<\Theta_{-}(\theta)<\frac{\pi}{4}$ for $\alpha_v<\theta<\alpha_m$.
  \item[2.] For any $\alpha_m\in (\alpha_v, \alpha_0)$, there exists a $\alpha'\in (\alpha_v, \alpha_m)$ and $\alpha''\in (\alpha_v, \alpha')$  such that $A_{-}(\alpha')=\Theta_{-}(\alpha')$, $A_{-}(\alpha'')<\frac{\Theta_\infty}{2}$, and
$A_{-}(\theta)<\Theta_{-}(\theta)$ for $\theta\in (\alpha'', \alpha')$.
\end{itemize}

\noindent%
{1.} By (\ref{62302}) and (\ref{62401}),  we have
\begin{equation}\label{62403}
\beta_{-}'(\theta)>0,\quad r_{-}'(\theta)<0,\quad\text{and}\quad 0<r_{-}(\theta)<\bar{r}(\theta)
\quad\text{for}\quad \alpha_{v}<\theta\leq\alpha_m.
\end{equation}
Hence
\begin{equation}\label{62404}
0<r_{-}(\alpha_v)=\lim\limits_{\alpha\rightarrow \alpha_{v}^{+}}r_{-}(\theta)\leq \bar{r}(\alpha_v)\quad\text{and}\quad \lim\limits_{\alpha\rightarrow \alpha_{v}^{+}}A_{-}(\theta)\quad \text{exists}.
\end{equation}

From (\ref{62404}) and the differential mean value theorem we know that
there exists an $\alpha^{m}\in (\alpha_v, \alpha_m)$ such that $\beta_{-}(\alpha^m)-\alpha_v=-A_{\varepsilon}$, where
$A_{\varepsilon}>0$ is determined by
$$
\frac{r_{-}(\alpha_m)}{\sin A_{\varepsilon}}=\frac{r_{-}(\alpha_v)}{\sin (\pi-A_{\varepsilon}-(\alpha_m-\alpha_v))};
$$
see Figure \ref{Fig4.2}.
Then according to $\beta_{-}'(\theta)>0$ for $\theta\in (\alpha_v, \alpha^m)$, we have
$$
2\lim\limits_{\theta\rightarrow \alpha_{v}^{+}}A_{-}(\theta)=\lim\limits_{\theta\rightarrow \alpha_{v}^{+}}\big(\theta-\beta_{-}(\theta)\big)>\alpha_v-\beta_{-}(\alpha^m)>A_{\varepsilon}.
$$
This leads to a contradiction.

%=============
\begin{figure}[htbp]
\begin{center}
\includegraphics[scale=0.43]{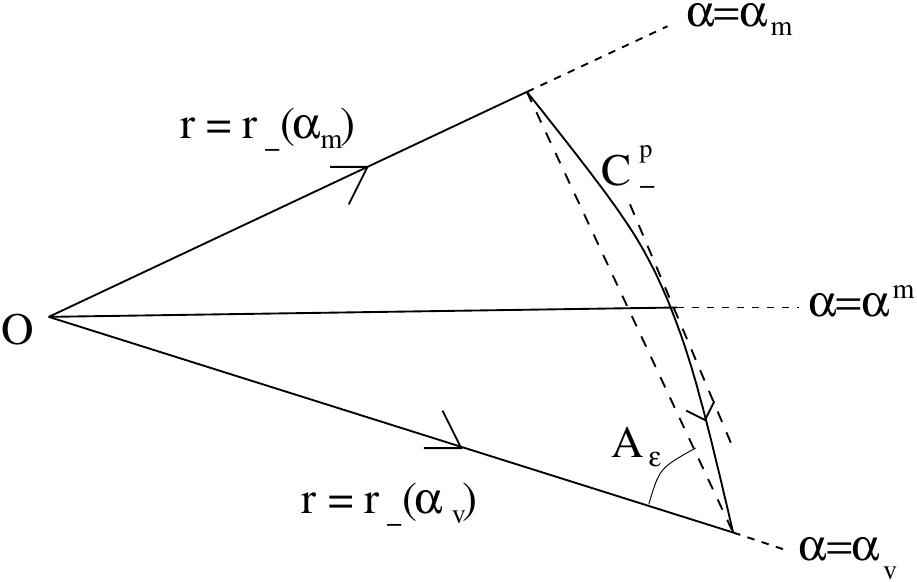}
\caption{\footnotesize Convexity of $C_{-}^{\mathrm{P}}$ near the vacuum.}
\label{Fig4.2}
\end{center}
\end{figure}
%==========

\noindent%
{2.}
By the mean value theorem we know that there exists an $\alpha'''\in (\alpha'', \alpha')$ such that
$A_{-}'(\alpha''')=\frac{A_{-}(\alpha')-A_{-}(\alpha'')}{\alpha'-\alpha''}>1$ as $\alpha_m-\alpha_v$ is sufficiently small.
While, by (\ref{62401}), (\ref{80802}), and $A_{-}(\alpha''')<\Theta_{-}(\alpha''')$ one has
$$
A_{-}'(\alpha''')=\frac{1-\beta_{-}'(\alpha''')}{2}<1.
$$
This leads to a contradiction.
Then we have $\inf\limits_{\theta\in (\alpha_v, \alpha_0)}A_{-}(\theta)>0$.

From (\ref{210cqo}) and $\inf\limits_{\theta\in (\alpha_v, \alpha_0)}A_{-}(\theta)>0$ we have $$\lim\limits_{\theta\rightarrow a_v}\xi_{-}(\theta)=\bar{r}(\alpha_v)\cos \alpha_{v}\quad \text{and}\quad
\lim\limits_{\theta\rightarrow a_v}\eta_{-}(\theta)=\bar{r}(\alpha_v)\sin \alpha_{v}.$$
Hence, the $C_{-}$ characteristic curve emanating from $\mathrm{P}$ ends up at the point $\mathrm{D}$.
This completes the proof of Lemma \ref{lem3.3}.
\end{proof}
We now use $\widetilde{\mathrm{PD}}$ and $\widetilde{\mathrm{PB}}$ to denote the $C_{-}$  and $C_{+}$ cross characteristic curves  emanating from $\mathrm{P}$, respectively.
By (\ref{102304}), Lemma \ref{lem3.3}, $\bar{\rho}'(\theta)>0$, and $\alpha=\theta$ for the centered simple wave,  we also have
\begin{equation}\label{62303}
\bar{\partial}_-\rho~=~\bar{\rho}'(\theta)\big(\cos\beta\theta_{\xi}+\sin\beta\theta_{\eta}\big)
~=~\frac{4p'(\tau)\cos^2A}{r\tau^2 p''(\tau)}<0
\quad \text{on}\quad \widetilde{\mathrm{PD}}
\end{equation}
and
\begin{equation}\label{51701}
\bar{\partial}_-c
~=~-\frac{\sin 2A}{r}\cdot \bar{c}'(\theta)~=~-\frac{4\varpi c \cos^2A}{r}
\quad \text{on}\quad \widetilde{\mathrm{PD}}.
\end{equation}
From (\ref{51701}) we see that $\frac{\bar{\partial}_{-}c}{\sin ^2 A}$ is finite on $\widetilde{\mathrm{PD}}$.

Let
$
\theta=\bar{\alpha}(\tau)
$, $\tau\geq \tau_0$ be the inverse function of $\tau=\bar{\tau}(\theta)$, $\alpha_v\leq\theta\leq\alpha_0$.
Clearly,
\begin{equation}
\bar{\alpha}(\tau_0)=\alpha_0, \quad \lim\limits_{\tau\rightarrow+\infty}\bar{\alpha}(\tau)=\alpha_v, \quad \text{and}\quad \bar{\alpha}'(\tau)<0\quad  \text{for}\quad \tau\geq \tau_0.
\end{equation}

%We shall divide the discussion into  the following two cases: $0<\theta_w\leq\alpha_v+\pi$; $\alpha_v+\pi<\theta_w<\pi$.

\subsection{Wave interaction for the case $0<\theta_w\leq\alpha_v+\pi$}
%We first study the case of $0<\theta_w\leq\alpha_v+\pi$.
\subsubsection{Main theorem}
When $0<\theta_w\leq\alpha_v+\pi$, the centered simple wave emanating from the corner $\mathrm{O}$ is a complete simple wave.
In order to construct a global-in-time classical solution to the 2D Riemann problem (\ref{Potential}), (\ref{initial})--(\ref{boundary}),
we make the following assumption:
\begin{description}
        \item[($h_{12}$)]  $\theta>\alpha_0-\pi+2\Theta(\bar{\tau}(\theta))$ for all $\alpha_v\leq \theta\leq\alpha_0$.
\end{description}

\begin{thm}\label{thmf+f1}
Assume $\tau_0>\tau_2^i$ and $0<\theta_w\leq \alpha_v+\pi$. Assume moreover that assumptions ($h_{11}$) and ($h_{12}$) hold.
Then the 2D Riemann problem (\ref{Potential}), (\ref{initial})--(\ref{boundary})  admits a global-in-time self-similar solution.
\end{thm}

To construct a self-similar solution to (\ref{Potential}), (\ref{initial})--(\ref{boundary}) for $\tau_0>\tau_2^i$ and $0<\theta_w\leq \alpha_v+\pi$,
it only remains to construct the solution in a region $\Sigma_{4}$ bounded by $\widetilde{\mathrm{PB}}$, $\widetilde{\mathrm{PD}}$, and a gas–vacuum interface.
To this end,
we consider (\ref{PEC}) with the boundary conditions
\begin{equation}\label{I1}
(u, v, \tau)=\left\{
               \begin{array}{ll}
                 (u_{-}, v_{-}, \tau_{-})(\xi), & \hbox{$(\xi, \eta)\in \widetilde{\mathrm{PB}}$;} \\[2pt]
                 (\bar{u}, \bar{v}, \bar{\tau})(\theta), & \hbox{$(\xi, \eta)\in\widetilde{\mathrm{PD}}$.}
               \end{array}
             \right.
\end{equation}
The problem (\ref{PEC}), (\ref{I1}) is a standard Goursat problem.

\subsubsection{Local solution}
%We have
%\begin{equation}\label{52703}
% \beta\mid_{\widetilde{PB}}~=~-\frac{\pi}{2}.
%\end{equation}
%Combining this with (\ref{7}) we get
%\begin{equation}
%\bar{\partial}_{+}\rho\mid_{\widetilde{PB}}~=~-\frac{4c\sin^{2} A}{\tau^{4}p''(\tau)\tan A}~<~0.\label{52701}
%\end{equation}

For each $s>0$, we let
$\mathrm{B}_s$ and $\mathrm{D}_s$ be the points on $\widetilde{\mathrm{PB}}$ and $\widetilde{\mathrm{PD}}$, respectively, satisfying  $\tau(\mathrm{B}_s)=\tau(\mathrm{D}_s)=\tau_0+s$.

%=============
\begin{figure}[htbp]
\begin{center}
\includegraphics[scale=0.52]{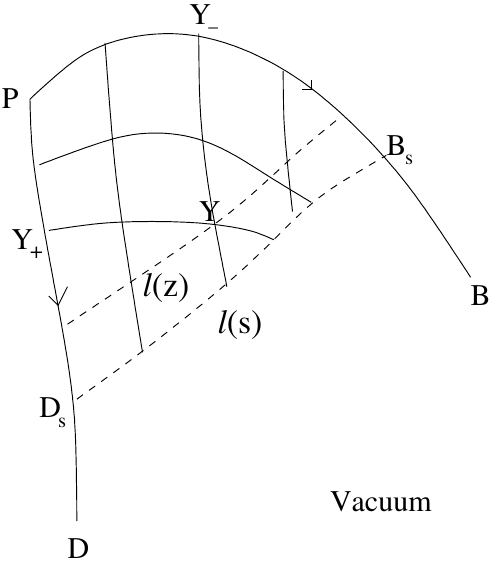}
\caption{\footnotesize Domain $\Sigma_4(s)$ and the level curve $\it{l}(s)$.}
\label{Fig4.3}
\end{center}
\end{figure}
%==========

\begin{lem} \label{loc54}
(Local solution)
For sufficiently small $s>0$, the Goursat problem (\ref{PEC}), (\ref{I1}) admits a unique $C^{1}$ solution on a curvilinear triangular region $\overline{\Sigma_{4}(s)}$ bounded by $\widetilde{\mathrm{PB_s}}$, $\widetilde{\mathrm{PD_s}}$, and the level curve ${\it l}(s)$ defined by $\tau(\xi,\eta)=\tau_0+s$; see Figure \ref{Fig4.3}. Moreover, the solution satisfies
\begin{equation}\label{022501}
 \bar{\partial}_{\pm}\rho<0\quad \text{in}\quad \overline{\Sigma_4(s)}.
\end{equation}
\end{lem}
\begin{proof}
The local existence of a classical solution follows by a standard argument  from Li and Yu \cite{Li-Yu}. Using (\ref{62201}), (\ref{62303}) and the characteristic decompositions (\ref{81104}), we deduce that the solution satisfies
$\bar{\partial}_{\pm}\rho<0$. Hence, for $s-s_0$ sufficiently small, the level curve ${\it l}(s)$ exists.
This completes the proof.
\end{proof}

\subsubsection{Hyperbolicity and $C^1$ norm estimate}
In what follows, we shall extend the local solution to the whole determinate region $\Sigma_4$. To this end, it remains to establish some uniform a priori $C^{1}$ estimates for the classical solution.

It is easy to verify that
there exists a sufficiently large integer $n>0$ such that
\begin{equation}\label{113002}
\mathcal{H}>0\quad \text{for}\quad \tau>\tau_0\quad \text{and}\quad 0<A<\frac{\pi}{2},
\end{equation}
where the variable $\mathcal{H}$ is defined in (\ref{82803}).
We define
$$
\mathcal{M}_{1}=\min \left\{\inf \limits_{\widetilde{\mathrm{PB}}}\frac{\rho^n\bar{\partial}_{+}\rho}{\sin^2 A}, ~
\inf \limits_{\widetilde{\mathrm{PD}}}\frac{\rho^n\bar{\partial}_{-}\rho}{\sin^2 A}\right\}\quad
\text{and}\quad \mathcal{N}_{1}(\tau)=\max\limits_{\tau_0\leq \varsigma\leq\tau}\frac{\varsigma^{4+n}p''(\varsigma)}{4c(\varsigma)}.
$$

\begin{figure}[htbp]
\begin{center}
\includegraphics[scale=0.6]{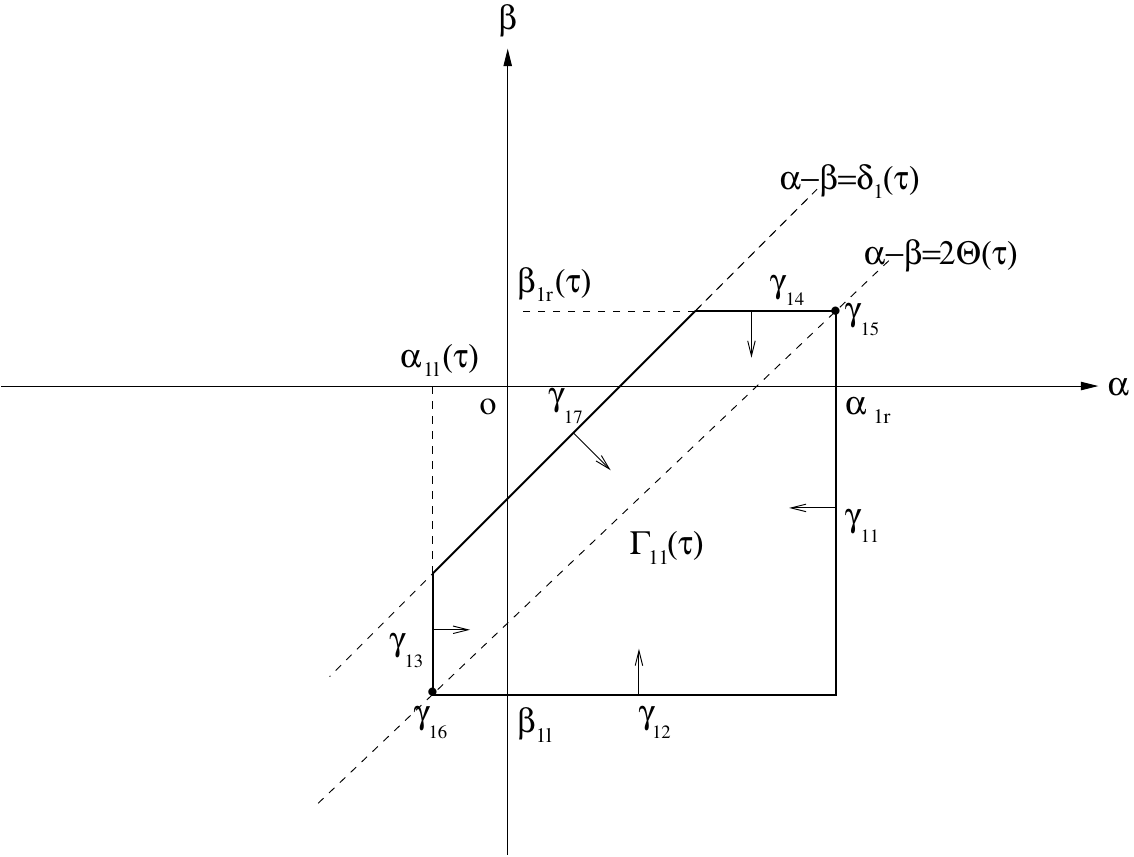}
\caption{\footnotesize Invariant region of $(\alpha, \beta)$.}
\label{Fig4.4}
\end{center}
\end{figure}

Let
$$
 \varepsilon_{11}=\inf\limits_{\theta\in(\alpha_v, \alpha_0)}\Big\{\theta-\big(\alpha_0-\pi+2\Theta(\bar{\tau}(\theta))\big)\Big\}, \quad \varepsilon_{12}=\frac{\pi}{2}-\alpha_0, \quad\text{and}\quad\varepsilon_1=\frac{1}{4}\min\{\varepsilon_{11}, \varepsilon_{12}\}.
$$
For $\tau\geq\tau_0$,
we introduce
$$
\Gamma_{11}(\tau)=\big\{(\alpha, \beta)\mid \alpha_{1l}(\tau)<\alpha<\alpha_{1r},~\beta_{1l}<\beta<\beta_{1r}(\tau),~\alpha-\beta>\delta_1(\tau)\big\},
$$
where
$$
\begin{aligned}
&\qquad \alpha_{1l}(\tau)=\alpha_0-\pi+2\Theta(\tau)+2\varepsilon_1, \quad \alpha_{1r}=\alpha_0+\varepsilon_1, \quad \beta_{1l}=\alpha_0-\pi+2\varepsilon_1, \\&\beta_{1r}(\tau)=\alpha_0+\varepsilon_1-2\Theta(\tau),\quad \text{and} \quad \delta_1(\tau)=\min\left\{\arctan\bigg(\frac{-1}{8\mathcal{M}_{1}\mathcal{N}_{1}(\tau)}\bigg),~ A_{\infty},~ \delta_0\right\};
\end{aligned}
$$
see Figure \ref{Fig4.4}.

%=============
%==========

By  assumption {\bf (E2)} we have
\begin{equation}\label{32402aa}
\Gamma_{11}(\tau')\subset\Gamma_{11}(\tau'')\quad \text{for~any}\quad  \tau_0\leq \tau'<\tau''.
\end{equation}

From Lemma \ref{lem3.1} we immediately obtain the following lemma.
\begin{lem}\label{112601}
Under assumption (${h}_{12}$),
we have
 $$(\alpha, \beta)(\xi, \eta_{+}(\xi))\in \Gamma_{11} \big(\tau_{-}(\xi)\big)\quad \text{for}\quad \xi\in(\xi_{_\mathrm{P}}, \xi_{_\mathrm{B}}).$$
\end{lem}

\begin{lem}\label{112602}
Under assumption (${h}_{12}$),
we have
\begin{equation}\label{81902}
(\theta, \beta_{-}(\theta))=(\alpha_0, -\frac{\pi}{2})\in \Gamma_{11} (\bar{\tau}(\theta))\quad \text{for}\quad \theta\in(\alpha_v, \alpha_0),
\end{equation}
where the function $\beta_{-}(\theta)$ is defined in (\ref{new81901}).
\end{lem}
\begin{proof}
It is easy to verify that $(\alpha_0, \beta_{-}(\alpha_0))\in \Gamma_{11}(0)$.
We next
prove the following assertion:
\begin{itemize}
  \item  For any fixed $\theta'\in(\alpha_v, \alpha_0)$, if $(\theta, \beta_{-}(\theta))\in \Gamma_{11}(\bar{\tau}(\theta))$ for all $\theta\in (\theta', \alpha_0)$ then
$(\theta', \beta_{-}(\theta'))\in \Gamma_{11}(\bar{\tau}(\theta')) $.
\end{itemize}

From Lemma \ref{lem3.3} we immediately have
$\theta'-\beta_{-}(\theta')\geq 2\delta_0>\delta_1(\bar{\tau}(\theta'))$.
By assumption (${h}_{12}$)  and the definition of the constant $\varepsilon_1$, one has
$\theta'>\alpha_{1l}(\bar{\tau}(\theta'))$. So, it only remains  to prove $\beta_{1l}<\beta_{-}(\theta')<\beta_{1r}(\bar{\tau}(\theta'))$.

Suppose $\beta_{-}(\theta')=\beta_{1l}$. Then
by the assumption of the assertion we have $\beta_{-}'(\theta')\geq 0$.
We also have
$$
\begin{aligned}
\Omega_{-}(\theta')&=\tan^2\Big(\Theta(\bar{\tau}(\theta'))\Big)-\tan^2\Big(\frac{\theta'-\beta_{1l}}{2}\Big)\\&<
\tan^2\Big(\Theta(\bar{\tau}(\theta'))\Big)-\tan^2\Big(\frac{\alpha_{1l}(\bar{\tau}(\theta'))-\beta_{1l}}{2}\Big)=0.
\end{aligned}
$$
However,
by (\ref{62401}) and (\ref{80802}) one has $\beta_{-}'(\theta')<0$.
This leads to a contradiction. Hence, we get $\beta_{-}(\theta')>\beta_{1l}$.

Suppose $\beta_{-}(\theta')=\beta_{1r}(\bar{\tau}(\theta'))$. Then
by the assumption of the assertion we have $\beta_{-}'(\theta')\leq 0$.
We also have
$$
\begin{aligned}
\Omega_{-}(\theta')&=\tan^2\Big(\Theta(\bar{\tau}(\theta'))\Big)-\tan^2\Big(\frac{\theta'-\beta_{1r}(\bar{\tau}(\theta'))}{2}\Big)\\&>
\tan^2\Big(\Theta(\bar{\tau}(\theta'))\Big)-\tan^2\Big(\frac{\alpha_{1r}-\beta_{1r}(\bar{\tau}(\theta'))}{2}\Big)=0.
\end{aligned}
$$
While,
by (\ref{62401}) one has $\beta_{-}'(\theta')>0$.
This leads to a contradiction. So, we get $\beta_{-}(\theta')<\beta_{1r}(\bar{\tau}(\theta'))$.
This completes the proof of the assertion.

Therefore, by a continuity argument we get (\ref{81902}). This completes the proof.
\end{proof}

%\begin{rem}
%Lemmas \ref{112601} and  \ref{112602} imply the following hyperbolicity estimate
%$$
%(\alpha, \beta)\in \Gamma_{11}(\tau)\quad \text{ono}\quad \widetilde{\mathrm{PB}}\cup \widetilde{\mathrm{PD}}.
%$$
%\end{rem}

\begin{lem}\label{51708}
Suppose that the Goursat problem (\ref{PEC}), (\ref{I1}) admits a classical solution in $\overline{\Sigma_4(s)}$ for some $s>0$. Then, under assumptions (${h}_{11}$) and (${h}_{12}$), the solution satisfies
\begin{equation}\label{51602}
(\alpha, \beta)\in \Gamma_{11}(\tau) \quad \text{and}\quad 2\mathcal{M}_{1}<\frac{\rho^n\bar{\partial}_{\pm}\rho}{\sin ^2A}<0 \quad \text{in}\quad \overline{\Sigma_4(s)}.
\end{equation}
\end{lem}
\begin{proof}
First, by Lemmas \ref{112601} and \ref{112602} and the definition of $\mathcal{M}_{1}$ we know that there exists a sufficiently small $\epsilon>0$ such that
$(\alpha, \beta)\in\Gamma_{11}(\tau)$
and $2\mathcal{M}_{1}<\frac{\rho^n\bar{\partial}_{\pm}\rho}{\sin^2 A}<0$ on $\overline{\Sigma_4(\epsilon)}$.
In order to use a continuity argument to prove this lemma,
we shall prove the following assertion.
\begin{itemize}
  \item For any fixed $z\in(0, s]$, suppose
$(\alpha, \beta)\in \Gamma_{11}(\tau)$ and $2\mathcal{M}_{1}<\frac{\rho^n\bar{\partial}_{+}\rho}{\sin^2 A}<0$ hold for every point in $\Sigma_4(z)\setminus{\it l}(z)$, then
they also hold on ${\it l}(z)$.
\end{itemize}

For any point $\mathrm{Y}\in {\it l}(z)$ the backward $C_{+}$ and $C_{-}$ characteristic curves emanating from $\mathrm{Y}$ remain in $\Sigma_4(z)$ until they intersect $\widetilde{\mathrm{PD}}$ and $\widetilde{\mathrm{PB}}$ at some points $\mathrm{Y}_{+}$ and $\mathrm{Y}_{-}$, respectively; see Figure \ref{Fig4.3}. (Remark: if $\mathrm{Y}\in\widetilde{\mathrm{PD}}$ then $\mathrm{Y}_{+}=\mathrm{Y}$ and $\mathrm{Y}_{-}=\mathrm{P}$; if $\mathrm{Y}\in\widetilde{\mathrm{PB}}$ then $\mathrm{Y}_{-}=\mathrm{Y}$  and $\mathrm{Y}_{+}=\mathrm{P}$.)

 Integrating (\ref{new}) along the $C_{\pm}$ characteristic curves from $\mathrm{Y}_{\pm}$ to $\mathrm{Y}$ and recalling the definition of $\mathcal{M}_{1}$, we have
\begin{equation}\label{82602}
2\mathcal{M}_{1}<\left(\frac{\rho^n\bar{\partial}_{\pm}\rho}{\sin^2 A}\right)(\mathrm{Y})<0.
\end{equation}

 If $\mathrm{Y}\in\widetilde{\mathrm{PD}}\cup\widetilde{\mathrm{PB}}$, then by Lemmas \ref{112601} and \ref{112602} we have $(\alpha, \beta)(\mathrm{Y})\in\Gamma_{11}(z)$. It only remains to consider the case of $\mathrm{Y}\notin\widetilde{\mathrm{PD}}\cup\widetilde{\mathrm{PB}}$. If $\mathrm{Y}\notin\widetilde{\mathrm{PD}}\cup\widetilde{\mathrm{PB}}$, then $\mathrm{Y}_{+}\neq \mathrm{Y}$ and $\mathrm{Y}_{-}\neq \mathrm{Y}$, and consequently $\widetilde{\mathrm{Y}_{+}\mathrm{Y}}$ and $\widetilde{\mathrm{Y_{-}Y}}$ exist. We next prove  $(\alpha, \beta)(\mathrm{Y})\notin \partial\Gamma_{11}(\tau)$, where  $ \partial\Gamma_{11}(\tau)$ is boundary of the domain  $\Gamma_{11}(\tau)$.
For convenience, we divide $\partial\Gamma_{11}(\tau)$ into the following seven parts:
\begin{itemize}
  \item $\gamma_{11}(\tau)=\big\{(\alpha, \beta)|\alpha=\alpha_{1r}, \beta_{1l}\leq\beta<\beta_{1r}(\tau)\big\}$;
  \item $\gamma_{12}(\tau)=\big\{(\alpha, \beta)|\beta=\beta_{1l}, \alpha_{1l}(\tau)<\alpha<\alpha_{1r}\big\}$;
  \item $\gamma_{13}(\tau)=\big\{(\alpha, \beta)|\alpha=\alpha_{1l}(\tau), \beta_{1l}<\beta\leq\beta_{1r}(\tau), \alpha-\beta>\delta_1(\tau)\big\}$;
  \item $\gamma_{14}(\tau)=\big\{(\alpha, \beta)|\beta=\beta_{1r}(\tau), \alpha_{1l}(\tau)<\alpha<\alpha_{1r}, \alpha-\beta>\delta_1(\tau)\big\}$;
  \item $\gamma_{15}(\tau)=\big\{(\alpha, \beta)|\alpha=\alpha_{1r}, \beta=\beta_{1r}(\tau)\big\}$;
  \item $\gamma_{16}(\tau)=\big\{(\alpha, \beta)|\alpha=\alpha_{1l}(\tau), \beta=\beta_{1l}\big\}$;
  \item $\gamma_{17}(\tau)=\big\{(\alpha, \beta)|\alpha_{1l}(\tau)\leq\alpha\leq \alpha_{1r}, \beta_{1l}\leq\beta\leq\beta_{1r}(\tau), \alpha-\beta=\delta_1(\tau)\big\}$;
\end{itemize}
see Figure \ref{Fig4.4}.
If $\alpha_{1l}(\tau)-\beta_{1r}(\tau)>\delta_1(\tau)$ then $\gamma_{17}(\tau)$ is an empty set.
We shall use the method of contradiction to prove $(\alpha, \beta)(\mathrm{Y})\notin \partial\Gamma_{11}(\tau)$.

Suppose $(\alpha, \beta)\in\gamma_{11}(\tau)$ at the point $\mathrm{Y}$. Then
by the assumption that $(\alpha, \beta)\in \Gamma_{11}(\tau)$ in $\Sigma_4(z)\setminus{\it l}(z)$ one has
$\bar{\partial}_{+}\alpha(\mathrm{Y})\geq0$.
However, by (\ref{10}) and (\ref{82602}) we have
\begin{equation}\label{113003}
c\bar{\partial}_{+}\alpha=-\frac{p''(\tau)}{4c\rho^4}\big(\underbrace{\tan^2(\Theta(\tau))-\tan^2A}_{<0}\big)\sin2 A\bar{\partial}_{+}\rho<0\quad\text{at}\quad \mathrm{Y}.
\end{equation}
This leads to a contradiction. Hence, we have $(\alpha, \beta)\notin\gamma_{11}(\tau)$  at the point $\mathrm{Y}$. Similarly, using (\ref{8}) we have $(\alpha, \beta)\notin\gamma_{12}(\tau)$  at the point $\mathrm{Y}$.

Suppose $(\alpha, \beta)\in\gamma_{13}(\tau)$  at the point $\mathrm{Y}$.  Then by the assumption that $(\alpha, \beta)\in \Gamma_{11}(\tau)$ in $\Sigma_4(z)\setminus{\it l}(z)$ one has
$\bar{\partial}_{+}\alpha(\mathrm{Y})\leq 0$.
However, by (\ref{10}) and (\ref{82602}) we have
$$
c\bar{\partial}_{+}\alpha=-\frac{p''(\tau)}{4c\rho^4}\big(\underbrace{\tan^2(\Theta(\tau))-\tan^2A}_{>0}\big)\sin2 A\bar{\partial}_{+}\rho>0\quad\text{at}\quad\mathrm{Y}.
$$
This leads to a contradiction. Hence, we get $(\alpha, \beta)\notin\gamma_{13}(\tau)$  at $\mathrm{Y}$. Similarly, we have $(\alpha, \beta)\notin\gamma_{14}(\tau)$  at $\mathrm{Y}$.

Suppose $(\alpha, \beta)\in \gamma_{16}(\tau)$  at $\mathrm{Y}$.
We let $\hat{\alpha}$ be a function defined on $\widetilde{\mathrm{Y_{+}Y}}$ such that
\begin{equation}\label{112905}
\left\{
  \begin{array}{ll}
   c\bar{\partial}_{+}\hat{\alpha}=\displaystyle-\frac{p''(\tau)}{4c\rho^4}\left[\tan^2 \big(\Theta(\tau_0+z)\big)-\tan^2\Big(\frac{\hat{\alpha}
   -\beta_{1l}}{2}\Big)\right]\sin2 A\bar{\partial}_{+}\rho\quad\text{along}\quad \widetilde{\mathrm{Y_{+}Y}}, \\[8pt]
    \hat{\alpha}(\mathrm{Y}_{+})=\alpha(\mathrm{Y}_{+}).
  \end{array}
\right.
\end{equation}
Then, by $\alpha(\mathrm{Y}_{+})>\alpha_{1l}(\tau_0+z)$ we have
\begin{equation}\label{82604}
\hat{\alpha}(\mathrm{Y})
>\alpha_{1l}(\tau(\mathrm{Y}_{+}))>\alpha_{1l}(\tau_0+z).
\end{equation}

From (\ref{10}) and (\ref{112905}) we have
\begin{equation}\label{82603}
\begin{array}{rcl}
  && c\bar{\partial}_{+}(\alpha-\hat{\alpha})=\underbrace{-\frac{p''(\tau)}{4c\rho^4}\Big[\tan^2(\Theta(\tau))-\tan^2 \big(\Theta(\tau_0+z)\big)\Big]\sin2 A\bar{\partial}_{+}\rho}_{>0, ~\text{since}~ \Theta>\Theta(\tau_0+z)~\text{along}~\widetilde{\mathrm{Y_{+}Y}}. }\\[32pt]&&\displaystyle+\frac{p''(\tau)}
{4c\rho^4}\Big[\tan^2\Big(\frac{\alpha-\beta}{2}\Big)-\tan^2\Big(\frac{\hat{\alpha}
-\beta_{1l}}{2}\Big)\Big]\sin2 A\bar{\partial}_{+}\rho
\quad\text{along}\quad \widetilde{\mathrm{Y_{+}Y}}.
  \end{array}
\end{equation}
By the assumption of the assertion we have $(\alpha,\beta)(\xi,\eta)\in\Gamma_{11}(\tau)$ for $(\xi,\eta)\in \widetilde{\mathrm{Y_{+}Y}}\setminus \{\mathrm{Y}\}$. Thus,
 $$\tan^2\Big(\frac{\alpha-\beta}{2}\Big)-\tan^2\Big(\frac{\hat{\alpha}
-\beta_{1l}}{2}\Big)<0\quad \text{if}\quad \alpha=\hat{\alpha}.$$
Therefore, by integrating (\ref{82603}) along $\widetilde{\mathrm{Y_{+}Y}}$ from $\mathrm{Y}_{+}$ to $\mathrm{Y}$ and using (\ref{82604}), we obtain $$\alpha(\mathrm{Y})>\hat{\alpha}(\mathrm{Y})>\alpha_{1l}(\tau_0+z).$$ This leads to a contradiction. Hence, we have $(\alpha, \beta)\notin \gamma_{16}(\tau)$ at $\mathrm{Y}$.
 Similarly, we have $(\alpha, \beta)\notin \gamma_{15}(\tau)$  at $\mathrm{Y}$.

Suppose $(\alpha, \beta)\in\gamma_{17}(\tau)$  at the point $\mathrm{Y}$. Then
by the assumption that $(\alpha, \beta)(\xi, \eta)\in \Gamma_{11}(\tau(\xi, \eta))$ for all $(\xi, \eta)\in \Sigma_4(z)\setminus{\it l}(z)$ we have $\bar{\partial}_{+}A\leq 0$ at $\mathrm{Y}$.
However, by (\ref{82301}) and (\ref{8}) we have
\begin{equation}\label{113005}
\begin{aligned}
c\bar{\partial}_{+} A=c\bar{\partial}_{-}\Big(\frac{\alpha-\beta}{2}\Big)&=\sin^2  A \left(1+\frac{p''(\tau)\tan A}{4c\rho^{4+n}}\frac{\rho^n\bar{\partial}_{+}\rho}{\sin^2 A}-\frac{p''(\tau)\Omega\sin^2 A}{8c\rho^{4+n}}\frac{\rho^n\bar{\partial}_{+}\rho}{\sin^2 A}\right)\\&>
\sin^2  A \left(1+\frac{p''(\tau)\tan A}{4c\rho^{4+n}}\frac{\rho^n\bar{\partial}_{-}\rho}{\sin^2 A}+\frac{p''(\tau)\tan A\sin^2 A}{4c\rho^{4+n}}\frac{\rho^n\bar{\partial}_{-}\rho}{\sin^2 A}\right)\\&>
\sin^2  A \Big(1+2\mathcal{M}_{1}\mathcal{N}_{1}(\tau)\tan(\delta_1(\tau))\Big)> 0\quad\text{at}\quad \mathrm{Y}.
\end{aligned}
\end{equation}
This leads to a contradiction.
Then, we get $(\alpha, \beta)\notin\gamma_{17}(\tau)$  at $\mathrm{Y}$.

This completes the proof of the assertion.

Therefore, by (\ref{32402aa}) and a continuity argument we know that the solution satisfies (\ref{51602}). This completes the proof of the lemma.
\end{proof}

The estimates in (\ref{51602}) are not enough to guarantee that the solution has a uniform  $C^1$ norm.
To this end, we shall seek a uniform bound for $\frac{\bar{\partial}_{\pm}c}{\sin^2 A}$ and
 a non-zero lower bound of $A$.
 The existence of a  non-zero lower bound of $A$ is also crucial for the regularity of the vacuum boundary of the wave interaction region $\Sigma_4$.

We define the following constants:
$$
 \mathcal{C}_1=\min\left\{\inf\limits_{\widetilde{\mathrm{PD}}}\frac{2\bar{\partial}_{-}c}{\sin^2 A},~ -\frac{8(\gamma-1)}{(\gamma+1)\tan A_{\infty}},~ \frac{\tau_{*}^{2+n}\big(2p'(\tau_*)+\tau_*p''(\tau_*)\big)\mathcal{M}_1}{\sqrt{-p'(\tau_*)}}\right\},\quad  \mathcal{D}_1=\frac{2\mathcal{C}_1}{\sin^2 \varepsilon_*},
$$
where the constants $\varepsilon_*$ and $\tau_*$ are the same as those defined in Section 2.5.

Let $\delta_*<\min\{\delta_1(\tau_*), \varepsilon_*\}$ be a sufficiently small positive constant such that
\begin{equation}\label{51601}
2+\frac{\sin A}{\varpi}\Big(\frac{\sin^2 A}{\cos A}+\frac{1}{\cos A}-\omega(\tau)\cos A\Big)\mathcal{D}_1>0\quad \text{for}\quad 0<A\leq \frac{\delta_*}{2}\quad\text{and} \quad\tau\geq\tau_*,
\end{equation}
where the constant $\tau_*$ is the same as that defined in Section 2.5.

\begin{lem}\label{61510}
Suppose that the Goursat problem (\ref{PEC}), (\ref{I1}) admits a classical solution on $\overline{\Sigma_4(s)}$ for some $s>\tau_*-\tau_0$. Then, the solution satisfies
\begin{equation}\label{61509}
\alpha-\beta>\delta_*, \quad \mathcal{C}_1<\bar{\partial}_{\pm}c<0,\quad\text{and}\quad \mathcal{D}_1<\frac{\bar{\partial}_{\pm}c}{\sin^2 A}<0 \quad \text{in}\quad \overline{\Sigma_4(s)}\setminus\overline{\Sigma_4(\tau_*-\tau_0)}.
\end{equation}
\end{lem}

\begin{proof}
First, from Lemma \ref{51708} we have
\begin{equation}\label{72906}
\alpha-\beta>\delta_*\quad \text{on}\quad {\it l}(s_*)
\end{equation}
By a direct computation we have
$$
\frac{\bar{\partial}_{\pm}c}{\sin ^2A}=\frac{\rho^n\bar{\partial}_{\pm}\rho}{\sin ^2A}\cdot\frac{2\tau^3p'(\tau)+\tau^4p''(\tau)}{2c\rho^n}.
$$
Hence, by (\ref{51602}) we have
\begin{equation}\label{72902}
 \mathcal{C}_1<\frac{\bar{\partial}_{\pm}c}{\sin ^2A}<0\quad \text{on}\quad {\it l}(s_*).
\end{equation}

Let $\mathrm{B}_{*}$ ($\mathrm{D}_{*}$, resp.) and $\mathrm{B}_{s}$ ($\mathrm{D}_{s}$, resp.) be the points on $\widetilde{\mathrm{PB}}$ ($\widetilde{\mathrm{PD}}$, resp.) such that $\tau(\mathrm{B}_{*})=\tau_*$ ($\tau(\mathrm{D}_{*})=\tau_*$, resp.) and $\tau(\mathrm{B}_{s})=\tau_s$ ($\tau(\mathrm{D}_{*})=\tau_s$, resp.); see Figure \ref{Fig4.5}.
%Let $\mathrm{B}_{s}$ and $\mathrm{D}_{s}$ be the point on $\widetilde{\mathrm{PB}}$ and $\widetilde{\mathrm{PD}}$ such that $\tau(\mathrm{B}_{s})=\tau(\mathrm{D}_{s})=\tau_0+s$.

%=============
\begin{figure}[htbp]
\begin{center}
\includegraphics[scale=0.52]{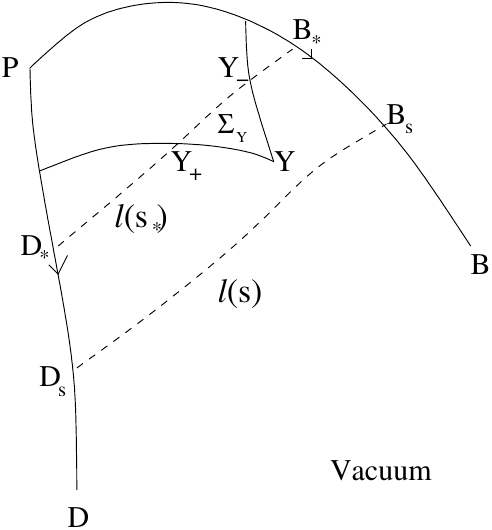}
\caption{\footnotesize Domain $\Sigma_Y$.}
\label{Fig4.5}
\end{center}
\end{figure}
%==========

From (\ref{7}) and (\ref{71602}) we have
$$
\frac{\bar{\partial}_{+}c}{\sin^2 A}=\frac{-2\varpi}{\tan A}\quad \text{on}\quad \widetilde{\mathrm{B_{*}B}}.
$$
Thus, by (\ref{51703}) and (\ref{121301}) one has
\begin{equation}\label{72903}
-\frac{8(\gamma-1)}{(\gamma+1)\tan A_{\infty}}<\frac{\bar{\partial}_{+}c}{\sin^2 A}<0\quad \text{on}\quad \widetilde{\mathrm{B_{*}B}}.
\end{equation}
By the the definition of $\mathcal{C}_1$ we have
\begin{equation}\label{72904}
\mathcal{C}_1<\frac{\bar{\partial}_{-}c}{\sin^2 A}<0\quad \text{on}\quad \widetilde{\mathrm{D_{*}D}}.
\end{equation}

See Figure \ref{Fig4.5}. Let $\mathrm{Y}$ be an arbitrary point in $\Sigma_4(s)\setminus\Sigma_4(s_*)$. The backward $C_{+}$ and $C_{-}$ characteristic curves intersect $\widetilde{\mathrm{B_*B_s}}\cup{\it l}(s_*)$ and $\widetilde{\mathrm{D}_*\mathrm{D}_s}\cup{\it l}(s_*)$ at some points $\mathrm{Y}_{+}$ and $\mathrm{Y}_{-}$, respectively. We denote by $\Sigma_{_\mathrm{Y}}$ a closed domain bounded by $\widetilde{\mathrm{Y_{-}Y}}$, $\widetilde{\mathrm{Y_{+}Y}}$, and a portion of $\widetilde{\mathrm{B_*B}_s}\cup\widetilde{\mathrm{D_*D}_s}\cup{\it l}(s_*)$.
In order to use the continuity argument to prove (\ref{61509}), we establish the following assertion:
if the inequalities in (\ref{61509}) hold at every point of $\Sigma_{_\mathrm{Y}}\backslash\{\mathrm{Y}\}$, then they also hold at $\mathrm{Y}$.

Suppose $(\bar{\partial}_{-}c)(\mathrm{Y})=\mathcal{C}_1$. Then by (\ref{72903}) we know $\mathrm{Y}\notin\widetilde{\mathrm{D_{*}D}}$. This implies that the characteristic curve $\widetilde{\mathrm{Y_{+}Y}}$ exists. Hence, by the assumption of the assertion  one has $\bar{\partial}_{+}\bar{\partial}_{-}c\leq 0$ at $\mathrm{Y}$.
However, by the first equation of (\ref{cd}) and (\ref{51703}) we obtain
$$
\begin{aligned}
c\bar{\partial}_{+}\bar{\partial}_{-}c~=~&\bar{\partial}_{-}c\left(\sin2 A+\frac{\bar{\partial}_{-}c}{2\varpi\cos^2 A}
+
\Big(\frac{1-\frac{\kappa\sin^2 2A}{\kappa+1}}{2\varpi\cos^2 A}-\frac{c}{\kappa}\cdot\frac{{\rm d}\kappa}{{\rm d}c}\Big)\bar{\partial}_{+}c\right)\\~>~&\bar{\partial}_{-}c\left(\sin2 A+\frac{\bar{\partial}_{-}c}{2\varpi\cos^2 A}\right)>0\quad \text{at}\quad \mathrm{Y}.
\end{aligned}
$$
This leads to a contradiction. So, we get $(\bar{\partial}_{-}c)(\mathrm{Y})>\mathcal{C}_1$. Similarly, we have
$(\bar{\partial}_{+}c)(\mathrm{Y})>\mathcal{C}_1$.

Suppose $\big(\frac{\bar{\partial}_{-}c}{\sin^{2} A}\big)(\mathrm{Y})=\mathcal{D}_1$.
Then by $(\bar{\partial}_{-}c)(\mathrm{Y})>\mathcal{C}_1$ one has $A(\mathrm{Y})<\varepsilon_*$.
Moreover, by the assumption of the assertion one has $\bar{\partial}_{+}\left(\frac{\bar{\partial}_{-}c}{\sin^{2} A}\right)\leq 0$ at $\mathrm{Y}$.
However, by the first equation of (\ref{cd1}) and (\ref{51704}) we have
\begin{equation}\label{7205}
\begin{aligned}
 &c\bar{\partial}_{+}\left(\frac{\bar{\partial}_{-}c}{\sin^{2} A}\right)\\=&~\bar{\partial}_{-}c\Bigg(\!\frac{1}{2\varpi\cos^{2} A}\bigg(\underbrace{\frac{\bar{\partial}_{-}c}{\sin^{2} A}-\frac{\bar{\partial}_{+}c}{\sin^{2} A}}_{\leq 0}\bigg)
+\bigg(\underbrace{\frac{1}{\varpi\cos^{2} A}-4\kappa\sin^{2} A-2-\frac{c}{\kappa}\cdot\frac{{\rm d}\kappa}{{\rm d}c}}_{>0}\bigg)
\frac{\bar{\partial}_{+}c}{\sin^{2} A}\!\Bigg)\\
>&~0\quad \text{at}\quad \mathrm{Y}.
\end{aligned}
\end{equation}
This leads to a contradiction. So we get $\big(\frac{\bar{\partial}_{-}c}{\sin^{2} A}\big)(\mathrm{Y})>\mathcal{D}_1$. Similarly, we have $\big(\frac{\bar{\partial}_{+}c}{\sin^{2} A}\big)(\mathrm{Y})>\mathcal{D}_1$.

Suppose $(\alpha-\beta)(\mathrm{Y})=\delta_*$.
Then $\mathrm{Y}\notin \widetilde{\mathrm{B_{*}B}_s}\cup \widetilde{\mathrm{D_{*}D}_s}$.
Thus, the characteristic curves $\widetilde{\mathrm{Y_{+}Y}}$ and $\widetilde{\mathrm{Y_{-}Y}}$ exist.
Hence, by the assumption of the assertion and $(\alpha-\beta)(\mathrm{Y})=\delta_*$ we have $\partial_{-}A\leq 0$ at $\mathrm{Y}$. However,
 by (\ref{82301}), (\ref{8}), and (\ref{51601}) we have
$$
\begin{aligned}
c\bar{\partial}_{-}A&=\frac{\sin^2 A}{2}\left(2+\frac{\sin A}{\varpi}\left(\frac{\sin^2 A}{\cos A}+\frac{1}{\cos A}-\omega(\tau)\cos A\right)\frac{\bar{\partial}_{-}c}{\sin^2 A}\right)\\[2pt]&\geq
\frac{\sin^2 A}{2}\left(2+\frac{\sin A}{\varpi}\left(\frac{\sin^2 A}{\cos A}+\frac{1}{\cos A}-\omega(\tau)\cos A\right)\mathcal{D}_1\right)~>~0\quad\text{at}\quad \mathrm{Y}.
\end{aligned}
$$
This leads to a contradiction. Hence, we get $(\alpha-\beta)(\mathrm{Y})>\delta_*$.

%By (\ref{cd}) and
%(\ref{51703}) one has
%\begin{equation}\label{7204}
%\begin{aligned}
%c\bar{\partial}_{+}\bar{\partial}_{-}c&=\bar{\partial}_{-}c\left(\sin2 A+\frac{\bar{\partial}_{-}c}{2\varpi\cos^2 A}
%+
%\Big(\frac{1-\frac{\kappa\sin^2 2A}{\kappa+1}}{2\varpi\cos^2 A}-\frac{c}{\kappa}\cdot\frac{{\rm d}\kappa}{{\rm d}c}\Big)\bar{\partial}_{+}c\right)\\&>
%\bar{\partial}_{-}c\left(\sin2 A+\frac{\bar{\partial}_{-}c}{2\varpi\cos^2 A}\right)>
%\bar{\partial}_{-}c\left(\sin2 A-\frac{1}{\cos^2 A}\right)>0\quad \text{along}\quad \widetilde{F_{+}F}.
%\end{aligned}
%\end{equation}
%Integrating this along $\widetilde{F_{+}F}$ we get
%$
%\mathcal{C}_1<\bar{\partial}_{-}c(S)<0
%$. Similarly, one has $\mathcal{C}_1<\bar{\partial}_{+}c(S)<0$.

%We now estimate $\frac{\bar{\partial}_{\pm}c}{\sin^2 A}$.

We then complete the proof of the assertion.
Therefore, by a continuity argument and (\ref{72906})--(\ref{72902}) we obtain (\ref{61509}). This completes the proof.
\end{proof}

\begin{lem}\label{lem3}
(Uniform a priori $C^1$ estimate) Suppose that the Goursat problem (\ref{PEC}), (\ref{I1})  admits a classical solution on $\overline{\Sigma_{4}(s)}$ for some $s>0$. Then, under assumption ($h_{12}$), there exists a positive constant $\mathcal{Q}$,   independent of $s$, such that the solution satisfies
$$
|\nabla u|+|\nabla v|+|\nabla \rho|<\mathcal{Q}\quad \text{and}\quad |u|+|v|<u_0+2\int_{\tau_0}^{+\infty}\frac{c(\tau)}{\tau}{\rm d}\tau.
$$
\end{lem}
\begin{proof}
By Lemma \ref{51708} we know that there is a constant $\mathcal{Q}_1$ such that
$$
| \nabla \rho|<\mathcal{Q}_1\quad \text{for}\quad \tau\leq\tau_*.
$$

From assumption {\bf(E1)} we know that $\max\limits_{\tau\geq\tau_*}\left|\frac{2c}{2\tau^3p'(\tau)+\tau^4p''(\tau)}\right|$ is finite.
So,
by Lemma \ref{51708}  we have
$$
\begin{aligned}
&|\nabla \rho|<\left|\frac{\bar{\partial}_{\pm}\rho}{\sin A\cos A}\right|\leq\left|\frac{\bar{\partial}_{\pm}c}{\sin ^2A}\cdot\frac{2c}{2\tau^3p'(\tau)+\tau^4p''(\tau)}\cdot\frac{\sin A}{\cos A}\right|\\&\qquad\leq \frac{\mathcal{D}_1}{\sin (\frac{\varepsilon_1}{2})}\cdot\max\limits_{\tau\geq\tau_*}\left|\frac{2c}{2\tau^3p'(\tau)+\tau^4p''(\tau)}\right|\quad \text{for}\quad \tau>\tau_*.
\end{aligned}
$$
Therefore, by (\ref{11}), (\ref{72804}), and (\ref{61513}), there is a constant $\mathcal{Q}$ independent of $s$ such that the solution satisfies
$$
|\nabla u|+|\nabla v|+|\nabla \rho|<\mathcal{Q} \quad\text{in}\quad \overline{\Sigma_4(s)}.
$$

Using  (\ref{11}) and (\ref{72804}) and recalling $\bar{\partial}_{\pm}\rho<0$ and $(u, v)(\mathrm{P})=(u_0, 0)$, we have
$$
|u|+|v|<u_0+2\int_{\tau_0}^{+\infty}\frac{c(\tau)}{\tau}{\rm d}\tau \quad\text{in}\quad \overline{\Sigma_4(s)}.
$$
This completes the proof.
\end{proof}

\subsubsection{Global solvability and regularity of the vacuum boundary}
\begin{thm}\label{thm1}
Under assumptions ($h_{11}$) and ($h_{12}$), the Goursat problem (\ref{PEC}), (\ref{I1}) admits a global classical solution on a bounded domain bounded by $\widetilde{\mathrm{PD}}$, $\widetilde{\mathrm{PB}}$, and a vacuum boundary $\widetilde{\mathrm{BD}}$.
Moreover, if $\alpha_0+A_\infty>\frac{\pi}{2}$, then the vacuum boundary $\wideparen{\mathrm{BD}}$ can be represented by a Lipschitz-continuous function $\xi=\mathcal{V}(\eta)$ defined on $[\bar{r}(\alpha_v)\sin\alpha_v, 0]$ which satisfies
\begin{equation}\label{198903}
\sup \limits _{\eta_1\neq\eta_2}\frac{|\mathcal{V}(\eta_1)-\mathcal{V}(\eta_2)|}{|\eta_1-\eta_2|}~\leq~\tan(\pi-\alpha_0-A_{\infty}).
\end{equation}
\end{thm}
\begin{proof}
Since $\bar{\partial}_{\pm}\tau>0$, each level curve ${\it l}(s)$ has nowhere a characteristic direction.
As in \cite{LJ2}, we take level curves of $\tau$ as the Cauchy support and solve a
mixed initial and characteristic boundary value problem in each extension step.
Hence, by Lemma \ref{lem3} we know that for any $s>0$, the local classical solution of  (\ref{PEC}), (\ref{I1}) can be extended to $\overline{\Sigma_4(s)}$. Then the problem admits a solution on $\overline{\Sigma_4}=\bigcup\limits_{s>0}\overline{\Sigma_4(s)}$.

By Lemmas \ref{51708} and \ref{61510} 
we know that the solution satisfies
 $\frac{\delta_*}{2}<A<\frac{\pi-\varepsilon_1}{2}$.
Combining this with  $q=\frac{c}{\sin A}$ and Lemma \ref{lem3} we know that the domains $\Sigma_4(s)$ are uniformly bounded  with respect to $s>0$.

Using Bernoulli's law (\ref{5802}), we have
\begin{equation}
h_{\xi}=-(u-\xi)u_{\xi}-(v-\eta)v_{\xi}\quad \text{and}\quad   h_{\eta}=-(u-\xi)u_{\eta}-(v-\eta)v_{\eta}.
\end{equation}
Thus, along each level curve ${\it l}(s)$,
\begin{equation}\label{62409}
-h_{\eta}u_{\xi}+h_{\xi}u_{\eta}=-(v-\eta)(v_{\xi}u_{\eta}-u_{\xi}v_{\eta})
\end{equation}
and
\begin{equation}\label{62410}
-h_{\eta}v_{\xi}+h_{\xi}v_{\eta}=(u-\xi)(v_{\xi}u_{\eta}-u_{\xi}v_{\eta}).
\end{equation}

Using (\ref{11}), (\ref{72804}), and (\ref{61513}), we have
$$
u_{\xi}=
-\frac{\sin^2\beta\bar{\partial}_{+}\rho+\sin^2\alpha\bar{\partial}_{-}\rho}{\sin2A}\cdot c\tau,\quad
v_{\xi}
=\frac{\sin\beta\cos\beta\bar{\partial}_{+}\rho+\sin\alpha\cos\alpha\bar{\partial}_{-}\rho}{\sin2A}\cdot c\tau,
$$
$$
v_{\eta}=-
\frac{\cos^2\beta\bar{\partial}_{+}\rho+\cos^2\alpha\bar{\partial}_{-}\rho}{\kappa\sin2A}\cdot c\tau,
\quad
u_{\eta}=
\frac{\sin\beta\cos\beta\bar{\partial}_{+}\rho+\sin\alpha\cos\alpha\bar{\partial}_{-}\rho}{\kappa\sin2A}\cdot c\tau.
$$
Thus by Lemma \ref{51708} we have
$$
v_{\xi}u_{\eta}-u_{\xi}v_{\eta}=-c\tau \sin 2A \bar{\partial}_{+}\rho\bar{\partial}_{-}\rho\neq0.
$$
This implies that the image of each level curve ${\it l}(s)$ in the $(u, v)$-plane is a smooth curve.
Moreover, by (\ref{62409}) and (\ref{62410}) we have that along each level curve ${\it l}(s)$,
\begin{equation}\label{2401}
\frac{{\rm d}u}{{\rm d}v}=-\frac{v-\eta}{u-\xi}=-\tan\Big(\frac{\alpha+\beta}{2}\Big).
\end{equation}
We denote by $u=\mathcal{V}_{s}(v)$ the image of  ${\it l}(s)$ in the $(u, v)$-plane.
By Lemmas \ref{51708} and \ref{61510} and recalling that $\varepsilon_1$ can be arbitrary small, we know that the solution satisfies
$$
\alpha_0-\pi+2A_{\infty}\leq \alpha\leq \alpha_0,\quad
\alpha_0-\pi\leq\beta<\alpha_0-2A_{\infty}.
$$
Thus, we have
$$
|\mathcal{V}_{s}'(v)|\leq \tan(\pi-\alpha_0-A_{\infty}).
$$
Therefore, by  the Arzela-Ascoli theorem  we know that there exists a Lipschitz-continuous function $u=\mathcal{V}(v)$, $v\in[\bar{r}(\alpha_v)\sin\alpha_v, 0]$ such that the image of the vacuum boundary in the $(u, v)$-plane can be represented by this function.

From $q=\frac{c}{\sin A}$ and $\frac{\delta_*}{2}<A<\frac{\pi-\varepsilon_1}{2}$ we know that on the vacuum boundary $q=0$. That is to say $u=\xi$ and $v=\eta$ on the vacuum boundary.
Hence, the vacuum boundary in the physical plane can be represented by $\xi=\mathcal{V}(\eta)$, $\eta\in[\bar{r}(\alpha_v)\sin\alpha_v, 0]$.
This completes the proof.
\end{proof}

For convenience, we denote by $(u, v, \tau)=(u_4, v_4, \tau_4)(\xi, \eta)$ the solution of the problem (\ref{PEC}, \ref{I1}).
Let
\begin{equation}\label{21603}
\overline{\Sigma_1}~=~\big\{(\xi, \eta)~\big|~0< \xi\leq \xi_{_\mathrm{P}}, ~~\eta\geq\xi\tan\alpha_0\big\}\cup
\big\{(\xi, \eta)~\big|~\xi\leq 0, ~~\eta\geq0\big\},
\end{equation}
\begin{equation}\label{21604}
\Sigma_2~=~\Big\{(\xi, \eta)~\big|~ \xi_{_\mathrm{P}}<\xi<\xi_{_\mathrm{B}}, ~\eta>\eta_{+}(\xi)\Big\},
\end{equation}
and
$$
\Sigma_3~=~\big\{(\xi, \eta)~\big|~ \xi=r\cos s, ~\eta=r\sin s,~r\in (0, r_{-}(s)),
~s\in(\alpha_v, \alpha_0)\big\};
$$
see Figure \ref{Fig4.1}.
We define
\begin{equation}\label{71201a}
(u, v, \tau)=\left\{
                              \begin{array}{ll}
                                (u_0, 0, \tau_0), & \hbox{$(\xi, \eta)\in\Sigma_1$,} \\[2pt]
                                (u_-, v_{-}, \tau_{-})(\xi, \eta), & \hbox{$(\xi, \eta)\in\Sigma_2$,} \\[2pt]
(\bar{u}, \bar{v}, \bar{\tau})(\theta), & \hbox{$(\xi, \eta)\in\Sigma_3$,} \\[2pt]
                               (u_4, v_4, \tau_4)(\xi, \eta), & \hbox{$(\xi, \eta)\in\overline{\Sigma_4}$,} \\[2pt]
                                \text{vacuum}, & \hbox{otherwise}.
                              \end{array}
                            \right.
\end{equation}
Then the triple $(u, v, \tau)$ defined in (\ref{71201a}) is a self-similar solution to the 2D Riemann problem (\ref{Potential}), (\ref{initial})--(\ref{boundary}) for $\tau_0>\tau_2^e$ and $0<\theta_w<\alpha_v+\pi$.
This completes the proof of Theorem \ref{thmf+f1}.

\subsection{Wave interaction for the case $\alpha_v+\pi<\theta_w<\pi$}
This part treats the case $\alpha_v+\pi<\theta_w<\pi$.
As shown in Section 3.3, for $\theta_w>\alpha_v+\pi$
the oncoming flow $(0)$ is deflected through an incomplete centered simple wave into a constant state $(u_w, v_w, \tau_w)$, and the flow near the corner is represented by (\ref{820a}).
\subsubsection{Main theorem}
\begin{figure}[htbp]
\begin{center}
 \includegraphics[scale=0.5]{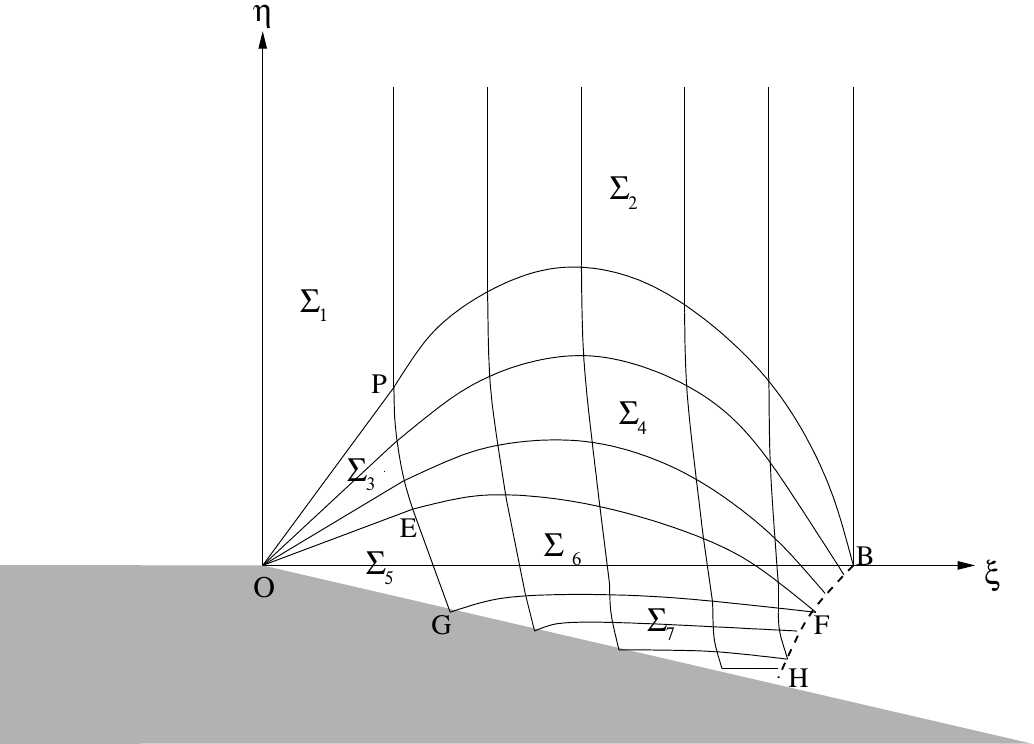}
\caption{\footnotesize Interaction of fan waves and reflection of a non-planar simple wave at the ramp.}
\label{Fig4.6}
\end{center}
\end{figure}

 For this case, to construct a global-in-time solution to the 2D Riemann problem (\ref{Potential}), (\ref{initial})--(\ref{boundary})
we also make the following assumptions:
\begin{description}
  \item[($h_{13}$)] $\theta>-\frac{\pi}{2}+2\Theta(\bar{\tau}(\theta))$ for $\alpha_w\leq \theta\leq\alpha_0$.
\item[($h_{14}$)] $\theta_w>\frac{\pi}{2}+2\Theta(\tau_0)$.
\end{description}
\begin{thm}\label{thmf+f2}
Assume $\tau_0>\tau_2^e$ and $\theta_w>\alpha_v+\pi$. Assume further that assumptions ($h_{11}$), ($h_{13}$), and ($h_{14}$) hold.
Then the 2D Riemann problem (\ref{Potential}), (\ref{initial})--(\ref{boundary})  admits a global-in-time self-similar solution; see Figure \ref{Fig4.6}.
\end{thm}

%=============
%==========

\subsubsection{Solution in domain $\overline{\Sigma_4}$}
Through $\mathrm{P}$ we draw a forward $C_{-}$ cross characteristic curve in the incomplete centered simple wave.
 This characteristic curve intersects the straight $C_{+}$ characteristic line $\eta=\xi\tan\alpha_w$ at a point $E$, where the constant $\alpha_w$ is defined in Section 3.3.
 The curve $\widetilde{\mathrm{PE}}$ is a portion of the $C_{-}$ cross characteristic curve $\widetilde{\mathrm{PD}}$ constructed in Section 4.2, and can be represented by the one-parameter form $\xi=r_{-}(\theta)\cos\theta$, $\eta=r_{-}(\theta)\sin\theta$, $\alpha_w\leq \theta\leq \alpha_0$, where the function $r_{-}(\theta)$ is defined in the proof of Lemma \ref{lem3.3}.

We first consider (\ref{PEC}) with the boundary conditions
\begin{equation}\label{I2}
(u, v, \tau)=\left\{
               \begin{array}{ll}
                 (u_{-}, v_{-}, \tau_{-})(\xi), & \hbox{$(\xi, \eta)\in\widetilde{\mathrm{PB}}$;} \\[3pt]
                 (\bar{u}, \bar{v}, \bar{\tau})(\theta), & \hbox{$(\xi, \eta)\in\widetilde{\mathrm{PE}}$.}
               \end{array}
             \right.
\end{equation}
Problem (\ref{PEC}), (\ref{I2}) is a standard Goursat problem and illuminate the interaction between the planar rarefaction wave and the incomplete centered simple wave.

Since assumption (${h}_{13}$) is stronger than assumption (${h}_{12}$), Theorem \ref{thm1} implies that the Goursat problem (\ref{PEC}), (\ref{I2}) admits a global $C^1$ solution on the domain ${\Sigma_4}$ bounded by $\widetilde{\mathrm{PB}}$, $\widetilde{\mathrm{PE}}$, $\widetilde{\mathrm{EF}}$, and $\wideparen{\mathrm{BF}}$, where $\widetilde{\mathrm{EF}}$ is a forward $C_{+}$ characteristic curve emanating from $E$ and  $\wideparen{\mathrm{BF}}$ is a vacuum boundary.
Moreover, the solution satisfies $\bar{\partial}_{\pm}\rho<0$ in $\overline{\Sigma_4}$.

Let $\varepsilon_2$ be a sufficiently small positive constant satisfying $\varepsilon_2<\frac{1}{8}\min\big\{\frac{\pi}{2}-\alpha_0, \pi-\theta_w\big\}$.
For any $\tau\geq\tau_0$, we define
$$
\Gamma_{12}(\tau)=\Big\{(\alpha, \beta)\mid \bar{\alpha}_{1l}(\tau)<\alpha<\bar{\alpha}_{1r},~\bar{\beta}_{1l}<\beta<\bar{\beta}_{1r}(\tau),~\alpha-\beta>\delta_2\Big\},
$$
where
$$
\begin{aligned}
&\qquad\qquad\bar{\alpha}_{1l}(\tau)=2\Theta(\tau)-\frac{\pi}{2}-\varepsilon_2, \quad \bar{\beta}_{1l}=-\frac{\pi}{2}-\varepsilon_2, \\& \bar{\alpha}_{1r}=\max\Big\{\frac{\pi}{2}+2(\theta_w-\pi), \alpha_0\Big\}+2\varepsilon_2, \quad \bar{\beta}_{1r}(\tau)=\bar{\alpha}_{1r}-2\Theta(\tau),
\end{aligned}
$$
and $\delta_2$ is a small positive constant, yet to be determined, satisfying  $\delta_2<\delta_*$. Here, the constant $\delta_*$ is as defined in (\ref{51601}).

By assumption {\bf(E2)} we have
\begin{equation}\label{32402a}
\Gamma_{12}(\tau')\subset\Gamma_{12}(\tau'')\quad \text{for}\quad  \tau_0\leq \tau'<\tau''.
\end{equation}
We divide $\partial\Gamma_{12}(\tau)$ into the following seven parts:
\begin{itemize}
  \item $\bar{\gamma}_{11}(\tau)=\big\{(\alpha, \beta)|\alpha=\bar{\alpha}_{1r}, \bar{\beta}_{1l}\leq\beta<\bar{\beta}_{1r}(\tau)\big\}$;
  \item $\bar{\gamma}_{12}(\tau)=\big\{(\alpha, \beta)|\beta=\bar{\beta}_{1l}, \bar{\alpha}_{1l}(\tau)<\alpha<\bar{\alpha}_{1r}\big\}$;
  \item $\bar{\gamma}_{13}(\tau)=\big\{(\alpha, \beta)|\alpha=\bar{\alpha}_{1l}(\tau), \bar{\beta}_{1l}<\beta\leq \bar{\beta}_{1r}(\tau), \alpha-\beta>\delta_2\big\}$;
  \item $\bar{\gamma}_{14}(\tau)=\big\{(\alpha, \beta)|\beta=\bar{\beta}_{1r}(\tau), \bar{\alpha}_{1l}(\tau)<\alpha<\bar{\alpha}_{1r}, \alpha-\beta>\delta_2\big\}$;
  \item $\bar{\gamma}_{15}(\tau)=\big\{(\alpha, \beta)|\alpha=\bar{\alpha}_{1r}, \beta=\bar{\beta}_{1r}(\tau)\big\}$;
  \item $\bar{\gamma}_{16}(\tau)=\big\{(\alpha, \beta)|\alpha=\bar{\alpha}_{1l}(\tau), \beta=\bar{\beta}_{1l}\big\}$;
  \item $\bar{\gamma}_{17}(\tau)=\big\{(\alpha, \beta)|\bar{\alpha}_{1l}(\tau)\leq\alpha\leq \bar{\alpha}_{1l}, \bar{\beta}_{1l}\leq\beta\leq\bar{\beta}_{1r}(\tau), \alpha-\beta=\delta_2\big\}$.
\end{itemize}

\begin{figure}[htbp]
\begin{center}
\includegraphics[scale=0.6]{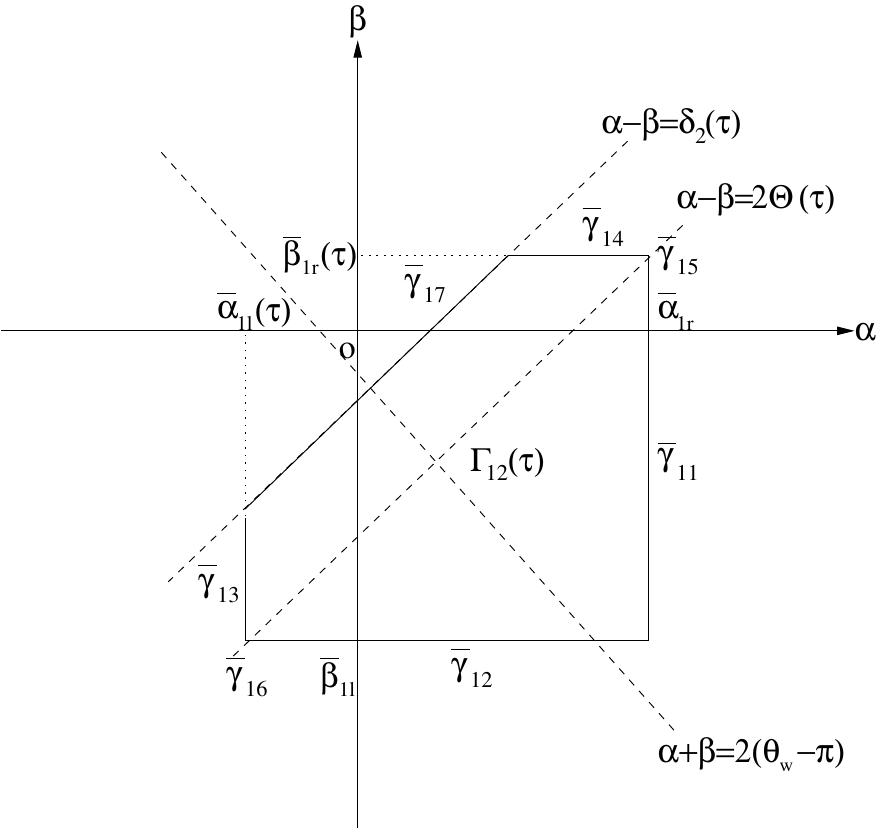}
\caption{\footnotesize Invariant region $\Gamma_{12}(\tau)$.}
\label{Fig4.7}
\end{center}
\end{figure}

From assumptions (${h}_{13}$) and (${h}_{14}$), we have
\begin{equation}\label{113006a}
\bar{\alpha}_{1l}(\tau)<2(\theta_w-\pi)-\bar{\beta}_{1l}<\bar{\alpha}_{1r}
\end{equation}
and
\begin{equation}\label{113007a}
\bar{\alpha}_{1l}(\tau)<\theta_w-\pi<\bar{\beta}_{1r}(\tau)
\end{equation}
for any $\tau\geq\tau_0$.
\begin{rem}\label{BZT61901}
From (\ref{113006a}) we have that in the $(\alpha, \beta)$-plane, the straight line $\alpha+\beta=2(\theta_w-\pi)$ and the line segment $\bar{\gamma}_{12}$ intersect.
From (\ref{113007a}) we have that when $\delta_2>0$ is sufficiently small, the straight line $\alpha+\beta=2(\theta_w-\pi)$ and the line segment $\bar{\gamma}_{17}$ intersect; see Figure \ref{Fig4.7}.
\end{rem}

%=============
%==========
From Lemma \ref{lem3.1} we immediately obtain the following result.
\begin{lem}
Under assumption (${h}_{11}$), 
 $$(\alpha, \beta)(\xi, \eta_{+}(\xi))\in \Gamma_{12}\big(\tau_{-}(\xi)\big)\quad \text{for}\quad \xi\in(\xi_{_\mathrm{P}}, \xi_{_\mathrm{B}}).$$
\end{lem}

\begin{lem}\label{51502}
Under assumptions (${h}_{11}$) and (${h}_{13}$), 
\begin{equation}\label{51705}
(\theta, \beta_{-}(\theta))\in \Gamma_{12}(\bar{\tau}(\theta))\quad \text{for}\quad \theta\in[\alpha_w, \alpha_0].
\end{equation}
\end{lem}
\begin{proof}
It is easy to check $(\alpha_0, -\frac{\pi}{2})\in\Gamma_{12}(0)$.
From assumption (${h}_{12}$), we have
$$
\bar{\alpha}_{1l}(\bar{\tau}(\theta))<\theta<\bar{\alpha}_{1r}\quad \text{for}\quad \theta\in[\alpha_w, \alpha_0].
$$
Then the rest of the proof is similar to that of Lemma \ref{112601}. We omit the details.
\end{proof}

As in Lemmas \ref{51708} and \ref{61510}, we have the following estimate:
\begin{lem}\label{113001}
Under assumptions (${h}_{11}$) and (${h}_{13}$), the solution of the Goursat problem (\ref{PEC}), (\ref{I2}) satisfies
\begin{equation}\label{51101}
(\alpha, \beta)\in \Gamma_{12}(\tau)\quad \text{and}\quad 2\mathcal{M}_{1}<\frac{\rho^n\bar{\partial}_{\pm}\rho}{\sin ^2A}<0
\text{in}\quad \overline{\Sigma_4}.
\end{equation}
Moreover,
\begin{equation}
\mathcal{D}_1<\frac{\bar{\partial}_{\pm}c}{\sin^2 A}<0\quad \text{for}\quad \tau>\tau_*.
\end{equation}
Here, the constants $\mathcal{M}_{1}$ and $\mathcal{D}_1$ are the same as those defined in Section 4.2.
\end{lem}

For convenience,
we denote the solution of the Goursat problem (\ref{PEC}), (\ref{I2}) by $(u, v, \tau)=(u_4, v_4, \tau_4)(\xi, \eta)$, $(\xi, \eta)\in \overline{\Sigma_4}$.

\subsubsection{Constant state solution on $\overline{\Sigma_5}$}
Let $\overline{\mathrm{OE}}$ be the line segment connecting $\mathrm{O}$ and $\mathrm{E}$. Then
 $\overline{\mathrm{OE}}$ is a straight $C_{+}$ characteristic line, and the state on it is the constant state $(u_w, v_w, \tau_w)$.
 We consider (\ref{PEC}) with the boundary conditions
\begin{equation}\label{BD3a}
\left\{
                       \begin{array}{ll}
                         (u, v, \tau)= (u_{w}, v_{w}, \tau_{w}) & \hbox{on $\overline{\mathrm{OE}}$,} \\[4pt]
                              v= u\tan\theta_w & \hbox{on $\partial\mathcal{W}_{+}$.}
                       \end{array}
                     \right.
\end{equation}

 %=============
\begin{figure}[htbp]
\begin{center}
\includegraphics[scale=0.46]{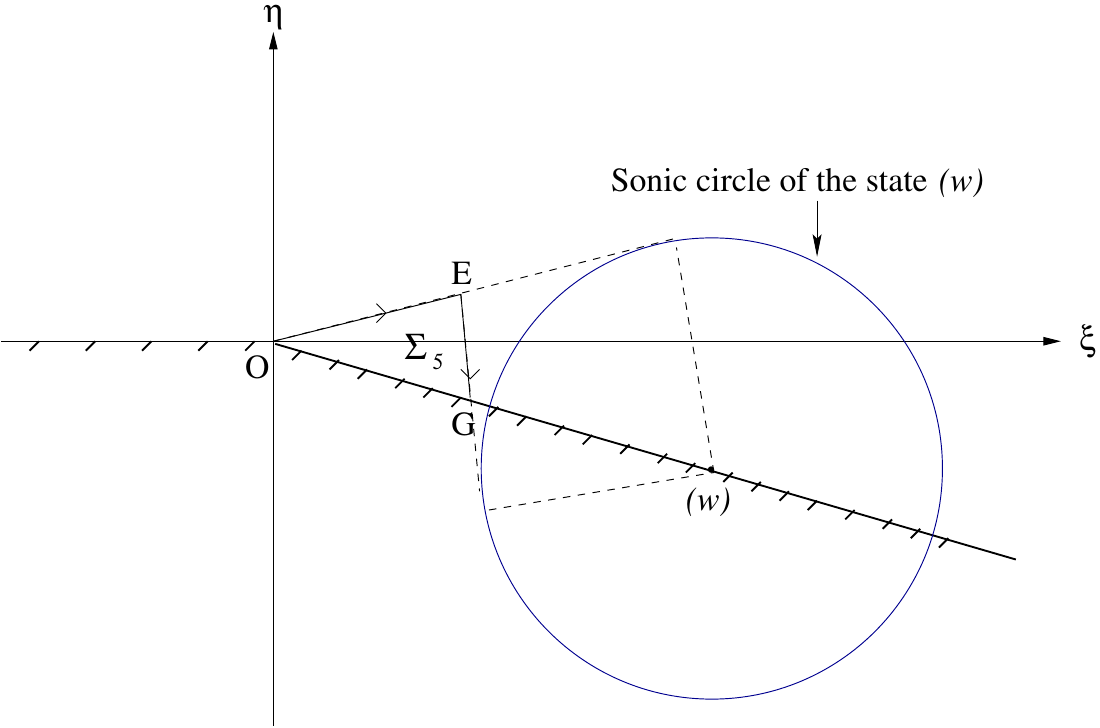}
\caption{\footnotesize Sonic circle of the state $(u_w, v_w, \tau_w)$ and the straight characteristic lines $\overline{\mathrm{OE}}$ and $\overline{\mathrm{EG}}$.}
\label{Fig4.8}
\end{center}
\end{figure}
%==========

From (\ref{7}), we have $c\bar{\partial}_{+}(\alpha+\beta)=c\bar{\partial}_{+}\beta=-2\sin^2A<0$ along $\overline{\mathrm{OE}}$.
Hence, by the limit $\lim\limits_{r\rightarrow 0+}(\alpha+\beta)\big(r\cos\alpha_w, r\sin\alpha_w\big)=2(\theta_w-\pi)$  we have $(\alpha+\beta)(\mathrm{E})<2(\theta_w-\pi)$. Consequently, since $\alpha\geq\beta$, we have
\begin{equation}\label{1981002}
\beta(\mathrm{E})<\theta_w-\pi.
\end{equation}
By (\ref{51705}) we have
\begin{equation}\label{821a}
\beta(\mathrm{E})=\beta_{-}(\alpha_w)>\bar{\beta}_{1l}>-\frac{\pi}{2}-\frac{1}{8}(\pi-\theta_w).
\end{equation}
Thus, from (\ref{1981002}) and (\ref{821a}) we have
\begin{equation}\label{821b}
0<(\theta_w-\pi)-\beta(\mathrm{E})<\frac{\pi}{2}.
\end{equation}

We draw a forward $C_{-}$ characteristic line through $\mathrm{E}$.
Using (\ref{82301}) and (\ref{8}), we have that along this characteristic line,
\begin{equation}\label{1981006}
\beta\equiv\beta(\mathrm{E})\quad \text{and} \quad c\bar{\partial}_{-}(\alpha+\beta)=c\bar{\partial}_{-}\alpha=2\sin^2A>0.
\end{equation}
Therefore, by (\ref{821b})  this $C_{-}$ characteristic line intersects $\partial\mathcal{W}_{+}$ at a point $\mathrm{G}$ before reaching the sonic circle of the state $(u_w, v_w, \tau_w)$.
Let $\overline{\Sigma_5}$ be the triangular domain bounded by $\overline{\mathrm{OE}}$, $\overline{\mathrm{EG}}$, and $\overline{\mathrm{OG}}$.
Then, the problem (\ref{PEC}), (\ref{BD3a}) has a constant state solution $(u, v, \tau)=(u_w, v_w, \tau_w)$ on  $\overline{\Sigma_5}$.
Moreover, from (\ref{1981006}) we have $(\alpha, \beta)\in \Gamma_{12}(\tau_w)$ on $\overline{\mathrm{EG}}$.

\subsubsection{Simple wave solution on $\overline{\Sigma_6}$}
We next consider (\ref{PEC}) with data
\begin{equation}\label{BD3}
(u, v, \tau)=\left\{
                       \begin{array}{ll}
                          (u_{4}, v_{4}, \tau_{4})(\xi, \eta) & \hbox{on $\widetilde{\mathrm{EF}}$,} \\[2pt]
                      (u_w, v_w, \tau_w) & \hbox{on $\overline{\mathrm{EG}}$.}
                       \end{array}
                     \right.
\end{equation}
%The problem (\ref{PEC}), (\ref{BD3}) is a standard Goursat problem.

\begin{figure}[htbp]
\begin{center}
\includegraphics[scale=0.55]{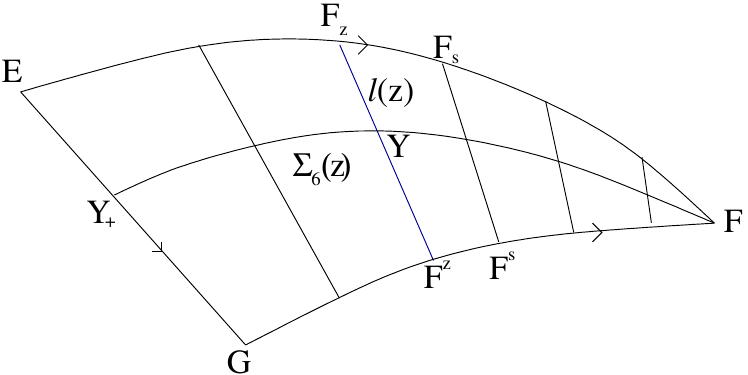}
\caption{\footnotesize Simple wave solution in the domain $\Sigma_6$.}
\label{Fig4.9}
\end{center}
\end{figure}

By Lemma \ref{simple}, the Goursat problem (\ref{PEC}), (\ref{BD3}) admits a simple wave solution with straight $C_{-}$ characteristic lines emanating from $\widetilde{\mathrm{EF}}$.
For any $s>s_{w}$ where $s_w=\tau_w-\tau_0$, let $\mathrm{F}_{s}$ be a point on $\widetilde{\mathrm{EF}}$ such that $\tau(\mathrm{F}_{s})=\tau_0+s$. The straight $C_{-}$ characteristic line emanating from $\mathrm{F}_{s}$ intersects  the forward $C_{+}$ characteristic curve emanating from $\mathrm{G}$ at a point $\mathrm{F}^{s}$. We denote by ${\Sigma_6(s)}$ the domain bounded by the characteristic lines $\widetilde{\mathrm{EF}_s}$, $\overline{\mathrm{EG}}$, $\widetilde{\mathrm{GF}^s}$, and $\overline{\mathrm{F}_s\mathrm{F}^s}$, which is a straight $C_{-}$ line from $\mathrm{F}_s$ to $\mathrm{F}^s$; see Figure \ref{Fig4.9}.

 %=============
%==========

\begin{lem}\label{lem5.2}
Suppose that assumptions (${h}_{11}$), (${h}_{13}$), and (${h}_{14}$) hold.  Then for any $s>s_w$, the simple wave solution satisfies
\begin{equation}
(\alpha, \beta)\in \Gamma_{12}(\tau), \quad \bar{\partial}_{-}\rho=0, \quad \text{and}\quad 2\mathcal{M}_1<\frac{\rho^n\bar{\partial}_{+}\rho}{\sin ^2A}<0 \quad \text{on}\quad \overline{\Sigma_6(s)}.
\end{equation}
Moreover, when $s>\tau_*-\tau_0$, there holds
\begin{equation}
\mathcal{D}_1<\frac{\bar{\partial}_{+}c}{\sin^2 A}<0 \quad \text{in}\quad \overline{\Sigma_6(s)}\setminus\overline{\Sigma_6(\tau_*-\tau_0)}.
\end{equation}
\end{lem}
\begin{proof}
The proof proceeds in two steps.

\noindent
{Step 1.} Using $\bar{\partial}_{-}\rho\mid_{\overline{EG}}~=~0$ and the first equation of (\ref{cd}),  the solution satisfies
$
\bar{\partial}_{-}\rho\equiv0.
$

From Lemma \ref{113001} we have that for any $s>s_w$,
\begin{equation}\label{51706}
(\alpha, \beta)\in \Gamma_{12}(\tau)\quad \text{and}\quad 2\mathcal{M}_1<\frac{\rho^n\bar{\partial}_{+}\rho}{\sin ^2A}<0
\quad \text{on}\quad \widetilde{\mathrm{EF}_s}.
\end{equation}
Moreover,
\begin{equation}\label{51707}
\mathcal{D}_1<\frac{\bar{\partial}_{+}c}{\sin^{2} A}<0\quad \text{on}\quad \widetilde{\mathrm{F_{*}F}},
\end{equation}
where $\mathrm{F}_{*}\in\widetilde{\mathrm{EF}}$ is determined by $\tau(\mathrm{F}_{*})=\tau_*$.

Integrating the second equation of (\ref{new}) along $C_{-}$ characteristic curves emanating from $\widetilde{\mathrm{EF}}$ and recalling (\ref{51706}), we have $2\mathcal{M}_1<\frac{\rho^n\bar{\partial}_{+}\rho}{\sin ^2A}<0$ in $\overline{\Sigma_6(s)}$. Integrating the second equation of (\ref{cd1}) along $C_{-}$ characteristic curves emanating from $\widetilde{\mathrm{F_{*}F}}$ and recalling (\ref{51707}), we have $\mathcal{D}_1<\frac{\bar{\partial}_{+}c}{\sin^2 A}<0$ in
$\overline{\Sigma_6(s)}\setminus\overline{\Sigma_6(\tau_*-\tau_0)}$ for $s>\tau_*-\tau_0$.

\vskip 2pt

\noindent
{Step 2.} We assert that
for any fixed $z\in(s_w, s]$, if
$(\alpha, \beta)\in \Gamma_{12}(\tau)$  holds for every point in $\Sigma_6(z)\setminus{\it l}(z)$ then
it also holds on ${\it l}(z)$. Indeed, ${\it l}(z)=\overline{\mathrm{F}_{z}\mathrm{F}^{z}}$ for this Goursat problem.

First, since $(\alpha, \beta)(\mathrm{F}_{z})\in \Gamma_{12}(\tau)$ and $\bar{\partial}_{-}\beta=0$ on ${\it l}(z)$, we have
$\bar{\beta}_{1l}<\beta<\bar{\beta}_{1r}(\tau_0+z)$  on ${\it l}(z)$.

Let $\mathrm{Y}$ be an arbitrary point on ${\it l}(z)$ and $\mathrm{Y}\neq \mathrm{F}_{z}$. The backward $C_{+}$ characteristic curve emanating from $\mathrm{Y}$ intersects $\widetilde{\mathrm{EG}}$ at a point $\mathrm{Y}_{+}$; see Figure \ref{Fig4.9}.
Suppose $\alpha(\mathrm{Y})=\bar{\alpha}_{1r}$. Then by the assumption in the assertion we have $\bar{\partial}_{+}\alpha\geq 0$ at $\mathrm{Y}$. However, as in (\ref{113003}), we have $\bar{\partial}_{+}\alpha<0$ at $\mathrm{Y}$. This leads to a contradiction. Hence, we have $\alpha(\mathrm{Y})<\bar{\alpha}_{1r}$. Similarly, we have  $\alpha(\mathrm{Y})>\bar{\alpha}_{1l}(\tau(\mathrm{Y}))$.
Suppose $A(\mathrm{Y})=\delta_2$. Then by the assumption in the assertion, we have $\bar{\partial}_{+}A\leq 0$ at $\mathrm{Y}$. However, as in (\ref{113005}), we have
$\bar{\partial}_{+}A>0$ at $\mathrm{Y}$. This leads to a contradiction. Hence, we have $A(\mathrm{Y})>\delta_2$.
This completes the proof of the assertion.

Therefore, by a continuity argument, we get $(\alpha, \beta)\in \Gamma_{12}(\tau)$ on $\overline{\Sigma_6(s)}$.
This completes the proof of the Lemma.
\end{proof}

By Lemma \ref{lem5.2}, we know that the Goursat problem (\ref{PEC}), (\ref{BD3}) admits a global classical solution in a region $\Sigma_6=\bigcup\limits_{s>\tau_w-\tau_0}\Sigma_6(s)$.
We next show that the forward $C_{+}$ characteristic issuing from $\mathrm{G}$  ends up at $\mathrm{F}$.

 From $\bar{\partial}_{-}\rho\equiv0$, we have $\bar{\partial}_{-}u\equiv0$ and $\bar{\partial}_{-}v\equiv0$. Hence, $(u, v)(\mathrm{F}^{s})=(u, v)(\mathrm{F}_{s})$. So, when $s\rightarrow+\infty$, we have $$(u, v)(\mathrm{F}^{s})=(u, v)(\mathrm{F}_{s})\rightarrow(u_4(\mathrm{F}), v_4(\mathrm{F})).$$ From $A=\arcsin(\frac{c}{q})$ and $A>\delta_2>0$, we have
$$
\xi(\mathrm{F}^{s})-u(\mathrm{F}^{s})\rightarrow 0\quad\text{and}\quad \eta(\mathrm{F}^{s})-v(\mathrm{F}^{s})\rightarrow 0\quad\text{as}\quad s\rightarrow+\infty.
$$
Hence, $\mathrm{F}^{s}\rightarrow \mathrm{F}$ as $s\rightarrow+\infty$.
So, the domain $\Sigma_6$ is bounded by $\widetilde{\mathrm{EF}}$, $\overline{\mathrm{EG}}$, and  $\widetilde{\mathrm{GF}}$.
For convenience, we denote the simple wave solution in $\overline{\Sigma_6}$ by $(u, v, \tau)=(u_{6}, v_{6}, \tau_{6})(\xi,\eta)$.

\subsubsection{Solution on $\overline{\Sigma_7}$}
\begin{figure}[htbp]
\begin{center}
\includegraphics[scale=0.56]{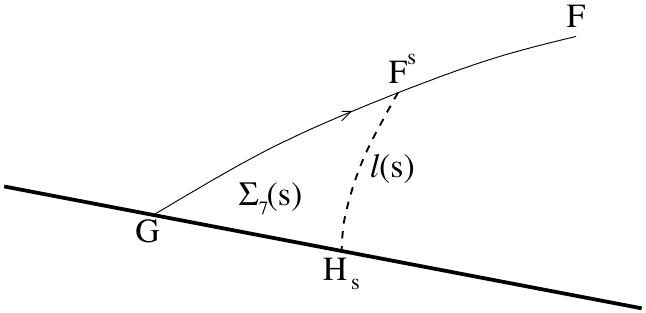}\qquad\qquad\qquad \includegraphics[scale=0.56]{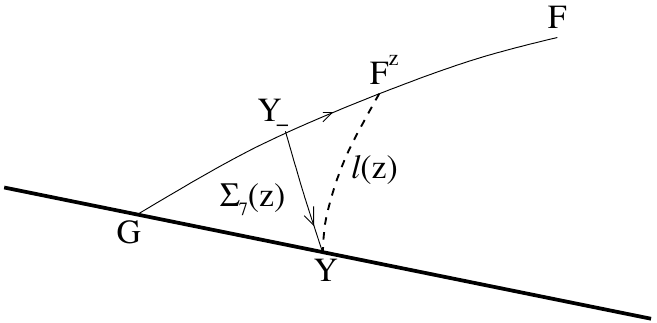}
\caption{\footnotesize Domain $\Sigma_7(s)$.}
\label{Fig4.10}
\end{center}
\end{figure}

We now consider (\ref{PEC}) with the boundary condition:
\begin{equation}\label{I3}
\left\{
               \begin{array}{ll}
                 (u, v, \tau)=(u_{6}, v_{6}, \tau_{6})(\xi,\eta) & \hbox{on $\widetilde{\mathrm{GF}}$;} \\[3pt]
                v= u\tan\theta_w & \hbox{on $\partial\mathcal{W}_{+}$,}
               \end{array}
             \right.
\end{equation}
where 
$$
\partial\mathcal{W}_{+}=\big\{(\xi, \eta)\mid \xi=\eta\cot\theta_w, \eta<0\big\}.
$$
\begin{lem}
When $s-s_w>0$ is sufficiently small, the boundary value problem (\ref{PEC}), (\ref{I3}) admits a classical solution on the domain $\overline{\Sigma_7(s)}$ bounded by $\widetilde{\mathrm{GF}^s}$, $\partial \mathcal{W}_{+}$, and the level curve ${\it l}(s)$ defined by $\tau(\xi, \eta)=\tau_0+s$; see Figure \ref{Fig4.10} (left).
\end{lem}
\begin{proof}
Problem (\ref{PEC}), (\ref{I3}) is a mixed characteristic boundary and slip boundary value problem. Existence of a local $C^1$ solution is known via the method of characteristics; see \cite{Li-Yu} (Chapter 2).
Since $\bar{\partial}_{+}\rho<0$ on $\widetilde{\mathrm{GF}}$, by the second equation of (\ref{cd}) we have $\bar{\partial}_{+}\rho<0$.

From the boundary condition on $\partial\mathcal{W}_{+}$, we have
 $$\bar{\partial}_{0}\sigma=0\quad \text{on} \quad\partial\mathcal{W}_{+}.$$
Hence,
 $$\left(\frac{\bar{\partial}_{+}+\bar{\partial}_{-}}{2}\right)\left(\frac{\alpha+\beta}{2}\right)=0\quad \text{on} \quad\partial\mathcal{W}_{+}.$$
Inserting (\ref{7})--(\ref{8}) into this equation, we get
\begin{equation}\label{62902}
\bar{\partial}_{+}\rho=\bar{\partial}_{-}\rho\quad \text{on}\quad \partial\mathcal{W}_{+}.
\end{equation}
Thus, we have $\bar{\partial}_{-}\rho<0$ on $\partial\mathcal{W}_{+}$. Combining this with the first equation of (\ref{cd}), we obtain that the solution satisfies
$\bar{\partial}_{-}\rho<0$.
Therefore, when $s-s_w>0$ is sufficiently small, the level curve of $\tau(\xi,\eta)=\tau_0+s$ exists. This completes the proof of the lemma.
\end{proof}

%=============
%==========

\begin{lem}\label{1001a}
Suppose that assumptions (${h}_{11}$),  (${h}_{13}$), and (${h}_{14}$) hold.
Suppose further that the problem (\ref{PEC}), (\ref{I3}) admits a classical solution on $\overline{\Sigma_7(s)}$ for some $s>s_w$. Then the solution satisfies
\begin{equation}\label{51801}
(\alpha, \beta)\in \Gamma_{12}(\tau)\quad \text{and}\quad  2\mathcal{M}_1<\frac{\rho^n\bar{\partial}_{+}\rho}{\sin ^2A}<0 \quad \text{in}\quad \overline{\Sigma_7(s)}.
\end{equation}
Moreover, when $s>\tau_*-\tau_0$ there also holds
\begin{equation}
\mathcal{D}_1<\frac{\bar{\partial}_{+}c}{\sin^2 A}<0 \quad \text{in}\quad \overline{\Sigma_7(s)}\setminus\overline{\Sigma_7(\tau_*-\tau_0)}.
\end{equation}
\end{lem}

\begin{proof}
There exists a sufficiently small $\varepsilon>0$ such that
$(\alpha, \beta)\in\Gamma_{12}(\tau)$
and $2\mathcal{M}<\frac{\rho^n\bar{\partial}_{\pm}\rho}{\sin^2 A}<0$ in $\overline{\Sigma_7(s_w+\varepsilon)}$.
In order to use a continuity argument to prove this lemma, we shall prove the following assertion.
\begin{itemize}
  \item For any fixed $z\in(s_w, s]$, if the desired estimates
 hold for every point in $\overline{\Sigma_7(z)}\setminus{\it l}(z)$ then
they also hold on ${\it l}(z)$.
\end{itemize}

For any point $\mathrm{Y}\in {\it l}(z)$, the backward $C_{+}$ and $C_{-}$ characteristic curves through $\mathrm{Y}$ remain in $\Sigma_7(z)$ until they intersect $\partial\mathcal{W}_{+}$ and $\widetilde{\mathrm{FG}}$ at some points $\mathrm{Y}_{+}$ and $\mathrm{Y}_{-}$, respectively; see Figure \ref{Fig4.10} (right).

If $\mathrm{Y}\notin\partial\mathcal{W}_{+}$ and $\mathrm{Y}\neq \mathrm{F}^{z}$, then the backward $C_{\pm}$ characteristic curves emanating from $\mathrm{Y}$ exist.
Hence, as in the proofs of Lemmas \ref{51708} and \ref{61510}, we know that the desired estimates hold at $\mathrm{Y}$.
Hence, it only remains to prove that when $\{\mathrm{Y}\}=\partial\mathcal{W}_{+}\cap {\it l}(z)$, the desired estimates hold at $\mathrm{Y}$.
When $\{\mathrm{Y}\}=\partial\mathcal{W}_{+}\cap {\it l}(z)$ the backward $C_{-}$ characteristic curves emanating from $\mathrm{Y}$ exists.

Integrating the second equation of (\ref{new}) along $\widetilde{\mathrm{Y}_{-}\mathrm{Y}}$ from $\mathrm{Y}_{-}$ to $\mathrm{Y}$ and recalling that $2\mathcal{M}_1<(\frac{\rho^n\bar{\partial}_{+}\rho}{\sin ^2A})(\mathrm{Y}_{-})<0$, we have
$2\mathcal{M}_1<(\frac{\rho^n\bar{\partial}_{+}\rho}{\sin ^2A})(\mathrm{Y})<0$. Combining this with (\ref{62902}), we get
$2\mathcal{M}_1<(\frac{\rho^n\bar{\partial}_{-}\rho}{\sin ^2A})(\mathrm{Y})<0$.

Suppose that $s>\tau_*-\tau_0$ and $(\frac{\bar{\partial}_{+}c}{\sin^2 A})(\mathrm{Y})=\mathcal{D}_1$. Then, by the assumption in the assertion we have $\bar{\partial}_{-}(\frac{\bar{\partial}_{+}c}{\sin^2 A})\leq 0$ at $\mathrm{Y}$. However, as in (\ref{7205}), we have  $\bar{\partial}_{-}(\frac{\bar{\partial}_{+}c}{\sin^2 A})> 0$ at $\mathrm{Y}$. This leads to a contradiction. Hence, we have $(\frac{\bar{\partial}_{+}c}{\sin^2 A})(\mathrm{Y})>\mathcal{D}_1$.

We next prove that $(\alpha, \beta)\notin \Gamma_{12}(\tau)$ at the point $\mathrm{Y}$.
Suppose that $(\alpha, \beta)(\mathrm{Y})\in \partial\Gamma_{12}(\tau(\mathrm{Y}))$.
Then, by $\sigma(\mathrm{Y})=(\frac{\alpha+\beta}{2})(\mathrm{Y})=\theta_w-\pi$ and Remark \ref{BZT61901}, we have $(\alpha, \beta)(\mathrm{Y})\in\bar{\gamma}_{12}(\tau(\mathrm{Y}))$ or $(\alpha, \beta)(\mathrm{Y})\in\bar{\gamma}_{17}(\tau(\mathrm{Y}))$.
\begin{itemize}
  \item Suppose that $(\alpha, \beta)(\mathrm{Y})\in\bar{\gamma}_{12}(\tau(\mathrm{Y}))$. Then, by the assumption in the assertion, we have $\bar{\partial}_{-}\beta\leq 0$ at $\mathrm{Y}$. However,
by (\ref{8}) and (\ref{82602}), we have
$$
c\bar{\partial}_{-}\beta=\frac{p''(\tau)}{4c\rho^4}\big(\underbrace{\tan^2(\Theta(\tau))-\tan^2A}_{<0}\big)\sin2 A\bar{\partial}_{+}\rho>0\quad\text{at}\quad\mathrm{Y}.
$$
This leads to a contradiction.
  \item
Suppose that $(\alpha, \beta)(\mathrm{Y})\in\bar{\gamma}_{17}(\tau(\mathrm{Y}))$.
Then, by the assumption in the assertion we have $\bar{\partial}_{-}A\leq 0$ at $\mathrm{Y}$. However, as in (\ref{113005}), we have
$\bar{\partial}_{-}A>0$ at $\mathrm{Y}$. This leads to a contradiction.
\end{itemize}
Thus, we get $(\alpha, \beta)(\mathrm{Y})\notin \partial\Gamma_{12}(\tau(\mathrm{Y}))$. This implies that $(\alpha, \beta)(\mathrm{Y})\in \Gamma_{12}(\tau(\mathrm{Y}))$.
This completes the proof of the assertion. Therefore, by a continuity argument we obtain the desired estimates. This completes the proof.
\end{proof}

\begin{lem}\label{92301a}
Suppose that assumptions (${h}_{11}$), (${h}_{13}$), and (${h}_{14}$) hold.
Suppose further that the problem (\ref{PEC}), (\ref{I3}) admits a classical solution on $\overline{\Sigma_7(s)}$ for some $s>s_w$. Then there exists a positive constant $\mathcal{Q}'$   independent of $s$ such that the solution satisfies
$$
|\nabla u|+|\nabla v|+|\nabla \rho|<\mathcal{Q}'\quad \text{and}\quad |u|+|v|<(2+|\cot\theta_w|)\int_{\tau_0}^{+\infty}\frac{c(\tau)}{\tau}{\rm d}\tau\quad \mbox{in}\quad \overline{\Sigma_7(s)}.
$$
\end{lem}
\begin{proof}
As in the proof of  Lemma \ref{lem3}, we can prove
the estimate for $|\nabla u|+|\nabla v|+|\nabla \rho|$. Integrating (\ref{72804}) along the forward $C_{-}$ characteristic curves emanating from $\widetilde{\mathrm{PB}}$ and noting $v=0$ on $\widetilde{\mathrm{PB}}$ we have
$|v|<\int_{\tau_0}^{+\infty}\frac{c(\tau)}{\tau}{\rm d}\tau$ in $\overline{\Sigma_7(s)}$. By the boundary condition, we have
$|u|<|\cot\theta_w|\int_{\tau_0}^{+\infty}\frac{c(\tau)}{\tau}{\rm d}\tau$ on $\partial\mathcal{W}_{+}\cap\overline{\Sigma_7(s)}$.
Therefore, by integrating (\ref{11}) along the forward $C_{+}$ characteristic curves emanating from $\partial\mathcal{W}_{+}\cap\overline{\Sigma_7(s)}$, we have $|u|<(1+|\cot\theta_w|)\int_{\tau_0}^{+\infty}\frac{c(\tau)}{\tau}{\rm d}\tau$. This completes the proof of the lemma.
\end{proof}

By Lemma \ref{92301a}, we immediately have the following lemma.
\begin{lem}\label{lem419}
Suppose that assumptions  (${h}_{11}$), (${h}_{13}$), and (${h}_{14}$) hold.
Then problem (\ref{PEC}), (\ref{I3}) admits a classical solution in a curvilinear triangle domain $\Sigma_7$ bounded by $\widetilde{\mathrm{GF}}$, $\partial\mathcal{W}_{+}$, and a vacuum boundary $\wideparen{\mathrm{FH}}$.
\end{lem}

\subsubsection{Global solvability and regularity of the vacuum boundary}
%We denote the solution of the BVP (\ref{PEC}, \ref{BC1}) by $(u, v, \tau)=(\hat{u}, \hat{v}, \hat{\tau})(\xi, \eta)$.
Let the domains $\Sigma_1$ and $\Sigma_2$ be the same as those defined in (\ref{21603}) and (\ref{21604}), respectively.
Let
$$
\Sigma_3~=~\big\{(\xi, \eta)~\big|~ \xi=r\cos s, ~\eta=r\sin s,~r\in (0, r_{-}(s)),
~s\in(\alpha_w, \alpha_0)\big\}.
$$
We define
\begin{equation}\label{71201aa}
(u, v, \tau)=\left\{
                              \begin{array}{ll}
                                (u_0, 0, \tau_0), & \hbox{$(\xi, \eta)\in\Sigma_1$,} \\[2pt]
                                (u_-, v_{-}, \tau_{-})(\xi, \eta), & \hbox{$(\xi, \eta)\in\Sigma_2$,} \\[2pt]
(\bar{u}, \bar{v}, \bar{\tau})(\arctan(\eta/\xi)), & \hbox{$(\xi, \eta)\in\Sigma_3$,} \\[2pt]
                               (u_w, v_w, \tau_w), & \hbox{$(\xi, \eta)\in\Sigma_5$,} \\[2pt]
(u_i, v_i, \tau_i)(\xi, \eta), & \hbox{$(\xi, \eta)\in\overline{\Sigma_i}$, $i=4, 6, 7$,} \\[2pt]
                                \text{vacuum}, & \hbox{otherwise}.
                              \end{array}
                            \right.
\end{equation}
Then the $(u, v, \tau)$ defined in (\ref{71201aa}) is a self-similar solution to the 2D Riemann problem (\ref{Potential}), (\ref{initial})--(\ref{boundary}) for $\tau_0>\tau_2^i$ and $\theta_w>\alpha_v+\pi$.
This completes the proof of Theorem \ref{thmf+f2}.

\begin{thm}
When $\alpha_0+A_\infty>\frac{\pi}{2}$,
the vacuum boundary $\wideparen{\mathrm{BF}}\cup\wideparen{\mathrm{FH}}$ can be represented by a Lipschitz-continuous function $\xi=\widetilde{\mathcal{V}}(\eta)$ defined on $[v({\mathrm{H}}), 0]$.  Moreover, this function satisfies
\begin{equation}\label{71011}
\sup \limits _{\eta_1\neq\eta_2}\frac{|\widetilde{\mathcal{V}}(\eta_1)-\widetilde{\mathcal{V}}(\eta_2)|}{|\eta_1-\eta_2|}\leq\tan\Big(\frac{\pi}{2}-A_\infty\Big).
\end{equation}
\end{thm}
\begin{proof}

Since $\varepsilon_2$ can be arbitrary small,  the solution satisfies
\begin{equation}\label{121401}
(\alpha, \beta)\in \big\{(\alpha, \beta) \mid \tilde{\alpha}_{1l}\leq \alpha\leq \tilde{\alpha}_{1r},~\tilde{\beta}_{1l}\leq \beta\leq \tilde{\beta}_{1r}, ~\alpha-\beta\geq \delta_2\big\}\quad \text{in}\quad \Sigma_4\cup\Sigma_6\cup\Sigma_7,
\end{equation}
where
$$
\tilde{\alpha}_{1l}=2A_\infty-\frac{\pi}{2}, \quad \tilde{\beta}_{1l}=-\frac{\pi}{2},\quad\tilde{\alpha}_{1r}=\max\Big\{\frac{\pi}{2}+2(\theta_w-\pi), \alpha_0\Big\}, \quad \text{and}\quad\tilde{\beta}_{1r}=\tilde{\alpha}_{1r}-2A_\infty.
$$
Hence, we have
$$
A_\infty-\frac{\pi}{2}\leq\frac{\alpha+\beta}{2}\leq\max\Big\{\frac{\pi}{2}+2(\theta_w-\pi), \alpha_0\Big\}-A_{\infty}\quad \text{in}\quad \Sigma_4\cup\Sigma_6\cup\Sigma_7.
$$

In view of Theorem \ref{thm1},
it only remains to discuss the vacuum boundary $\wideparen{\mathrm{FH}}$.

Using the first equation in (\ref{72804}) and (\ref{121401}), we have
\begin{equation}\label{121402}
\bar{\partial}_{+}v\geq 0\quad \text{on}\quad \widetilde{\mathrm{GF}}.
\end{equation}
Using (\ref{72804}) and (\ref{62902}), we get
\begin{equation}\label{121403}
\bar{\partial}_{0}v=\frac{\bar{\partial}_{+}v+\bar{\partial}_{-}v}{2}=2c\tau\sin\theta_w\bar{\partial}_{+}\rho<0
\quad \text{on}\quad \widetilde{\mathrm{GH}}.
\end{equation}

Let ${\mathrm{H}}_s$ be the point on $\widetilde{\mathrm{GH}}$ such that $\tau({\mathrm{H}}_s)=\tau_0+s$. Then, by (\ref{121402}) and (\ref{121403}) we have
\begin{equation}\label{121404}
-\int_{\tau_0}^{+\infty}\frac{c(\tau)}{\tau}{\rm d}\tau<v({\mathrm{H}}_s)<v(\mathrm{G})<v(\mathrm{F}^s)=v(\mathrm{F}_s)<0.
\end{equation}

As in the proof of Theorem \ref{thm1}, there exists a function $u=\widehat{\mathcal{V}}_s(v)$ defined on $[v({\mathrm{H}}_s), v(\mathrm{F}^s)]$ such that the image of the level curve $\wideparen{\mathrm{F}^s{H}_s}$ can be represented by this function. Meanwhile, $|\widehat{\mathcal{V}}_s'(v)|<\tan(\frac{\pi}{2}-A_{\infty})$ for $v\in [v({\mathrm{H}}_s), v(\mathrm{F}^s)]$.

Using the Arzela-Ascoli theorem and recalling (\ref{121402})--(\ref{121404}), we know that $v(\mathrm{H})$ exists and there is a Lipschitz-continuous function $u=\widehat{\mathcal{V}}(v)$ defined on $[v({\mathrm{H}}), v(\mathrm{F})]$ such that  the image of the vacuum boundary $\wideparen{\mathrm{F}{\mathrm{H}}}$ can be represented by this function. Combining this with $q=\frac{c}{\sin A}$ and $A>\delta_2$, the vacuum boundary $\wideparen{\mathrm{F}{\mathrm{H}}}$ can be represented by $\xi=\hat{\mathcal{V}}(\eta)$, $\eta\in [v({\mathrm{H}}), v(\mathrm{F})]$.
This completes the proof.
\end{proof}

%============================================================================================================

\section{Interaction between two shock–fan composite waves}
In this section, we study the 2D Riemann problem (\ref{Potential}), (\ref{initial})--(\ref{boundary}) for $\tau_1^e\leq\tau_0<\tau_2^i$.
In the following discussion, the  constants $\tau_{po}$, $\rho_{po}$, $c_{po}$, $\xi_{_\mathrm{P}}$, $\xi_{_\mathrm{B}}$, $\phi_{po}$, $u_{po}$, $v_{po}$, $\sigma_{po}$, $\alpha_v$, and the functions $(u_{-}, v_{-}, \tau_{-})(\xi)$, $(\bar{u}, \bar{v}, \bar{\tau})(\theta)$ are as defined in Section 3.4.
Assume $u_0>\frac{\rho_{po}c_{po}}{\rho_0}$ and $\theta_w-\pi<\sigma_{po}$. Then, away from the corner, the oncoming flow $(0)$ and the vacuum
are connected by a planar shock–fan composite wave, which can be represented by (\ref{71501c}).
At the corner, the oncoming flow is deflected by an oblique shock–fan composite wave to a vacuum or a constant state $(u_w, v_w, \tau_w)$. The flow near the corner can be represented by (\ref{8230}) or (\ref{8230a}), respectively.

% where $\rho_{po}$ and $c_{po}$ are defined in Section 4.4.
%=============
\begin{figure}[htbp]
\begin{center}
\includegraphics[scale=0.55]{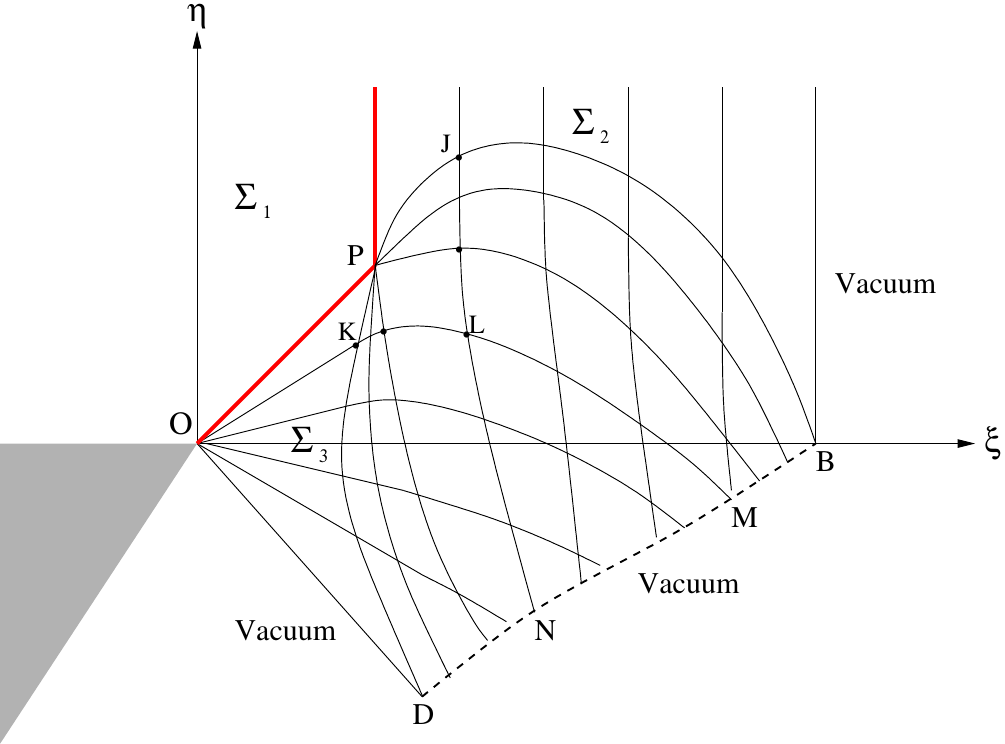}
\caption{\footnotesize Interaction of shock–fan composite waves for $\tau_1^e\leq\tau_0<\tau_2^i$ and $\theta_w<\alpha_v+\pi$.}
\label{Fig5.1}
\end{center}
\end{figure}
%==========

\subsection{Characteristic angles at the intersection point}
 These two composite waves begin to interact with each other at the point $\mathrm{P}(\xi_{_\mathrm{P}}, \eta_{_\mathrm{P}})$,
 where 
$$
\xi_{_\mathrm{P}}=u_0-\frac{\rho_{po}c_{po}}{\rho_0}, \quad \eta_{_\mathrm{P}}=\xi_{_\mathrm{P}}\tan\phi_{po}=u_0(1-\sin\phi_{po})\tan\phi_{po}=
\frac{\rho_{po}c_{po}(1-\sin\phi_{po})}{\rho_0\cos\phi_{po}},
$$
and $\phi_{po}$ is the inclination angle of the oblique shock (see Section 3.4).

For convenience,  we use the subscripts $2$ and $3$ to denote the flow states of the planar rarefaction wave and the centered simple wave, respectively; that is,
$$
(u_2, v_2, \tau_2)(\xi, \eta)=(u_{-}, v_{-}, \tau_{-})(\xi)\quad \text{and} \quad (u_3, v_3, \tau_3)(\xi, \eta)=(\bar{u}, \bar{v}, \bar{\tau})(\theta).
$$

To obtain a global solution of the 2D Riemann problem, we first make the following assumption:
\begin{description}
\item[($h_{21}$)] $\min\big\{\alpha_2(\mathrm{P})+\frac{\pi}{2},  \phi_{po}-\beta_3(\mathrm{P})\big\}>2\Theta(\tau_{po})$, $\alpha_2(\mathrm{P})-\beta_3(\mathrm{P})<\pi$, $\sigma_2(\mathrm{P})>\sigma_3(\mathrm{P})$.
\end{description}

\begin{figure}[htbp]
\begin{center}
\includegraphics[scale=0.5]{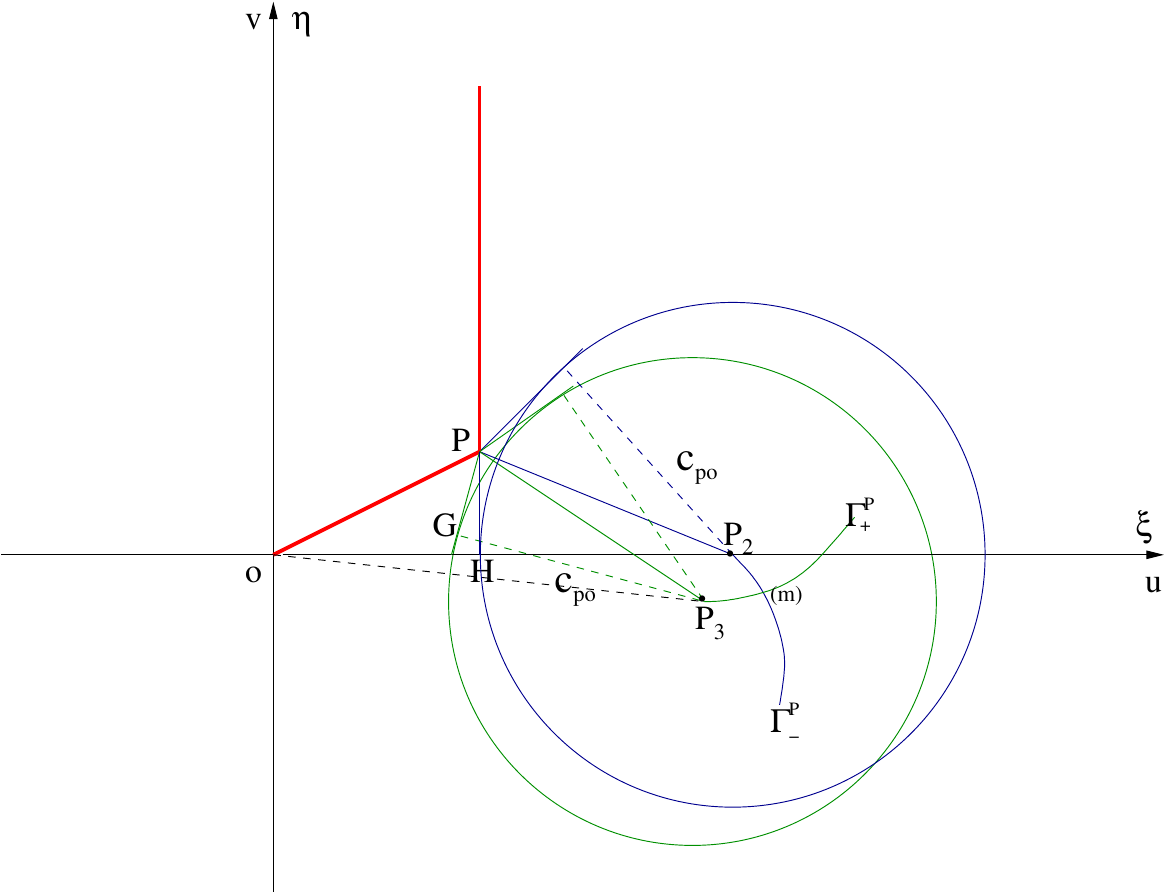}
\caption{\footnotesize Local structure at $\mathrm{P}$ and $\Gamma_{\pm}$ characteristic curves, where the coordinates of the points  $\mathrm{P}_{2}$ and $\mathrm{P}_3$ are $(u_{-}(\xi_{_\mathrm{P}}), 0)$ and $(u_{po}, v_{po})$, respectively.}
\label{Fig5.2}
\end{center}
\end{figure}

By computation, we have
$$
\alpha_{2}(\mathrm{P})=\beta_{2}(\mathrm{P})+2A_{2}(\mathrm{P})=-\frac{\pi}{2}+2\arctan \Big(\frac{c_{po}}{\eta_{_P}}\Big)=-\frac{\pi}{2}+2\arctan \left(\frac{\rho_0(1+\sin\phi_{po})}{\rho_{po}\cos\phi_{po}}\right),
$$
\begin{equation}\label{821c}
\sigma_2(\mathrm{P})=\frac{1}{2}(\alpha_2(\mathrm{P})+\beta_2(\mathrm{P}))=
-\frac{\pi}{2}+\arctan \left(\frac{\rho_0(1+\sin\phi_{po})}{\rho_{po}\cos\phi_{po}}\right),
\end{equation}
and
$$
\rho_{po}u_{-}(\xi_{_\mathrm{P}})-\rho_0u_0=\rho_{po}(u_{-}(\xi_{_\mathrm{P}})-\xi_{_\mathrm{P}})
-\rho_0(u_0-\xi_{_\mathrm{P}})+(\rho_{po}-\rho_0)\xi_{_\mathrm{P}}=(\rho_{po}-\rho_0)\xi_{_\mathrm{P}}<0.
$$
Hence,
$$
u_{-}(\xi_{_\mathrm{P}})\sin\phi_{po}=u_{-}(\xi_{_\mathrm{P}})\frac{\rho_{po} c_{po}}{\rho_0 u_0}<c_{po}.
$$
This implies
that the straight line $\eta=\xi\tan\phi_{po}$ intersects the sonic circle
$(\xi-u_{-}(\xi_{_\mathrm{P}}))^2+\eta^2=c_{po}^2$; see Figure \ref{Fig5.2}.
Furthermore, since the straight line $\eta-\eta_{_P}=\tan(\alpha_{2}(\mathrm{P})) (\xi-\xi_{_P})$ is tangent to the sonic circle
$(\xi-u_{-}(\xi_{_\mathrm{P}}))^2+\eta^2=c_{po}^2$,
we have
\begin{equation}\label{4906}
\alpha_{2}(\mathrm{P})>\phi_{po}.
\end{equation}

%Let $q_s=\sqrt{(u_b-\xi_{_P})^2+(v_b-\eta_{_P})^2}$ and $A_s=\arcsin\frac{c_{po}}{q_b}$.

%The straight line $\eta=\xi\tan\phi_{po}$ is tangent to the sonic circle $(\xi-u_{po})^2+(\eta-v_{po})^2=c_{po}^2$ at a point $\mathrm{G}=(\xi_{_\mathrm{G}}, \eta_{_\mathrm{G}})$. By computation, one has
Since $\tau_2(\mathrm{P})=\tau_3(\mathrm{P})=\tau_{po}$, from Bernoulli's law (\ref{5802}) we have
\begin{equation}\label{82701a}
A_{3}(\mathrm{P})=A_{2}(\mathrm{P})\quad \text{and}\quad (\xi_{_\mathrm{P}}-u_{-}(\xi_{_\mathrm{P}}))^2+\eta_{_\mathrm{P}}^2=
(\xi_{_\mathrm{P}}-u_{po})^2+(\eta_{_\mathrm{P}}-v_{po})^2.
\end{equation}
%$$
%\begin{aligned}
%A_{3}(\mathrm{P})&=\arctan\Big(\frac{c_{po}}{|\overline{\mathrm{PG}}|}\Big)=\arctan\Big(\frac{c_{po} \cos\phi_{po}}{\xi_{_\mathrm{G}}-\xi_{_\mathrm{P}}}\Big)\\&=\arctan\left(\frac{c_{po} \cos\phi_{po}}{u_0(1-\sin^2\phi_{po})-u_0(1-\sin\phi_{po})}\right)=\arctan\left(\frac{\rho_0(1+\sin\phi_{po})}{\rho_{po}\cos\phi_{po}}\right).
%\end{aligned}
%$$
Hence, 
\begin{equation}\label{4903}
\beta_{3}(\mathrm{P})=\alpha_{3}(\mathrm{P})-2A_{3}(\mathrm{P})=\phi_{po}-2\arctan \left(\frac{\rho_0(1+\sin\phi_{po})}{\rho_{po}\cos\phi_{po}}\right)
\end{equation}
and
\begin{equation}\label{821bb}
\sigma_3(\mathrm{P})=\phi_{po}-\arctan \left(\frac{\rho_0(1+\sin\phi_{po})}{\rho_{po}\cos\phi_{po}}\right).
\end{equation}

%=============
%==========

From (\ref{4901}) we have
$$
\begin{aligned}
u_{po}-c_{po}-\xi_{_\mathrm{P}}&= c_{po}(\sin\phi_{po}-1)+u_0(1-\sin^2\phi_{po})-u_0(1-\sin\phi_{po})
\\&=(\sin\phi_{po}-1)c_{po}\Big(1-\frac{\rho_{po}}{\rho_0}\Big)<0.
\end{aligned}
$$
Hence, we get
\begin{equation}\label{4905}
\beta_{3}(\mathrm{P})<-\frac{\pi}{2}.
\end{equation}
%since the straight line $\eta=\xi\tan\beta_{3}(P)$ is tangent to the sonic circle $(\xi-u_{po})^2+(\eta-v_{po})^2=c_{po}^2$.

It is easy to verify that
\begin{equation}
\alpha_{2}(\mathrm{P})\rightarrow \frac{\pi}{2}, \quad \beta_{3}(\mathrm{P})\rightarrow -\frac{\pi}{2},\quad \text{and}\quad \sigma_{po}\rightarrow 0\quad\text{as}\quad \phi_{po}\rightarrow \frac{\pi}{2}.
\end{equation}
By computation, we obtain that
if  $\frac{\rho_{po}}{\rho_0}>\frac{1}{2}$,
then
 \begin{equation}\label{4902}
 \frac{{\rm d}}{{\rm d}\phi_{po}}\big(\alpha_{2}(\mathrm{P})-\beta_{3}(\mathrm{P})\big)\bigg|_{\phi_{po}=\frac{\pi}{2}}~=~\frac{2\rho_{po}}{\rho_0}-1>0.
\end{equation}
Therefore, when $\frac{\rho_{po}}{\rho_0}>\frac{1}{2}$ and $\phi_{po}$ is close to $\frac{\pi}{2}$, there holds
\begin{equation}\label{4904}
\alpha_{2}(\mathrm{P})-\beta_{3}(\mathrm{P})<\pi.
\end{equation}
%This condition is crucial in the following discussion.

Form (\ref{821c}) and (\ref{821bb}) we also have that when $\phi_{po}$ is close to $\frac{\pi}{2}$, 
\begin{equation}\label{821aa}
\sigma_2(\mathrm{P})>\sigma_3(\mathrm{P}).
\end{equation}

%In this section
Since $\phi_{po}=\arcsin \big(\frac{\rho_{po}c_{po}}{\rho_0u_0}\big)$, when $u_0$ is close to $\frac{\rho_{po}c_{po}}{\rho_0}$, assumption ($h_{21}$) holds.

\begin{rem}\label{26102201}
From (\ref{82701a}) we also have
\begin{equation}\label{82702a}
\arctan\Big(\frac{v_{po}}{u_{po}}\Big)-\beta_3(\mathrm{P})>\frac{\pi}{2}.
\end{equation}
See Figure \ref{Fig5.2}. The straight line $\eta-\eta_{_\mathrm{P}}=(\xi-\xi_{_\mathrm{P}})\tan\beta_3(\mathrm{P})$ is tangent to the sonic circle $(\xi-u_{po})^2+(\eta-v_{po})^2=c_{po}^2$ at the point $\mathrm{G}(\xi_{_\mathrm{G}}, \eta_{_\mathrm{G}})$.
Let $\mathrm{H}=(\xi_{_\mathrm{P}}, 0)$. Then by (\ref{82701a}) we have $|\mathrm{PG}|=|\mathrm{PH}|$, and hence $\eta_{_\mathrm{G}}>0$. This leads to (\ref{82702a}).
The relation (\ref{82702a}) will not be used in this section, but will be used in Section 7.
\end{rem}

\subsection{Cross characteristic curves}
Through $\mathrm{P}$, we draw a forward $C_{+}$ cross characteristic curve in the planar fan wave and a forward $C_{-}$ cross characteristic curve in the centered fan wave.
As in Section 4.1,
the $C_{+}$ and $C_{-}$ cross characteristic curves can cross the entire fan waves.
Moreover, we have the following results. 

\begin{lem}
Assume $\alpha_2(\mathrm{P})+\frac{\pi}{2}>2\Theta(\tau_{po})$. Then the forward $C_{+}$ characteristic curve emanating from the point $\mathrm{P}$ can cross the entire planar rarefaction wave and terminate at the point $\mathrm{B}(\xi_{_\mathrm{B}}, 0)$. Moreover, along this $C_{+}$ cross characteristic curve there holds
\begin{equation}\label{7201}
2\Theta(\tau_{-}(\xi))-\frac{\pi}{2}<\alpha\leq \alpha_2(\mathrm{P}).
\end{equation}
%where the function $\tau_{-}(\xi)$ is defined in (\ref{71601}).
\end{lem}

%The $C_{+}$ cross characteristic curve $\widetilde{\mathrm{PB}}$ can be represented by a function $\eta=\eta_{+}(\xi)$, $\xi_{_\mathrm{P}}\leq \xi<\xi_{_\mathrm{B}}$ which satisfies
%$$
%\frac{{\rm d}\eta_{+}(\xi)}{{\rm d}\xi}=\alpha_2(\xi, \eta_{+}(\xi)), \quad \eta_{+}(\xi_{_\mathrm{P}})=0.
%$$

\begin{lem}
Assume $\phi_{po}-\beta_3(\mathrm{P})>2\Theta(\tau_{po})$. Then
the $C_-$ characteristic emanating from  $\mathrm{P}$ can cross the entire centered simple wave and terminate at the point
$
\mathrm{D}=\big(\bar{r}(\alpha_v)\cos \alpha_{v}, ~~\bar{r}(\alpha_v)\sin \alpha_{v}\big),
$
where the function $\bar{r}(\theta)=\sqrt{\bar{q}^2(\theta)+\bar{\tau}^2(\theta)p'(\bar{\tau}(\theta))}$ is as defined in Section 3.3.
Moreover, there exists a sufficiently small $\delta_{0}>0$ such that along this $C_{-}$ cross characteristic curve,
\begin{equation}
\delta_{0}~\leq~ A~\leq~ A_3(\mathrm{P}).
\end{equation}
\end{lem}

From now on we use $\widetilde{\mathrm{PB}}$ and $\widetilde{\mathrm{PD}}$ to denote the forward $C_{+}$ and $C_{-}$ cross characteristic curves emanating from point $\mathrm{P}$.

\subsection{Centered waves}
We shall show that the interaction of the shocks generates two centered waves emanating from $\mathrm{P}$.
We first present the definition of centered waves for the quasilinear hyperbolic system (\ref{PEC}); see Li and Yu \cite{Li-Yu} (pp. 188-190).

\subsubsection{$C_{+}$ centered wave}
\begin{defn}
  Let $\Delta_{+}(r)$ be an angular domain with curved boundaries:
\begin{equation}
\Delta_{+}(r)=\big\{(\xi,\eta)\mid\xi_{_\mathrm{P}}\leq \xi\leq\xi_{_\mathrm{P}}+r,~ \eta_{1}(\xi)\leq \eta\leq \eta_{2}(\xi)\big\},
\end{equation}
where
$\eta_{1}(\xi_{_\mathrm{P}})=\eta_{2}(\xi_{_\mathrm{P}})=\eta_{_\mathrm{P}}$
and
$\eta_{1}'(\xi_{_\mathrm{P}})<\eta_{2}'(\xi_{_\mathrm{P}})$.

A function $(u, v, \tau)(\xi,\eta)$ is called a $C_{+}$  centered wave solution of (\ref{PEC}) with center at $\mathrm{P}$ if the following properties hold:
\begin{enumerate}
  \item $(u, v, \tau)$ can be determined implicitly by the functions
$\eta=\hat{\eta}(\xi,\nu)$ and
$(u, v, \tau)=(\hat{u},\hat{v}, \hat{\tau})(\xi,\nu)$
defined on a rectangular domain
$T_{+}(r)=\{(\xi,\nu)\mid\xi_{_\mathrm{P}}\leq \xi\leq\xi_{_\mathrm{P}}+r, \nu_1\leq \nu\leq \nu_2\}$, where $\nu_1=\eta_{1}'(\xi_{_\mathrm{P}})$ and $\nu_2=\eta_{2}'(\xi_{_\mathrm{P}})$.
Moreover, $\hat{\eta}$ and $(\hat{u}, \hat{v}, \hat{\tau})$ belong to $C^{1}\big(T_{+}(r)\big)$, and for any $(\xi,\nu)\in T_{+}(r)\setminus\{\xi=\xi_{_\mathrm{P}}\}$ we have
$\hat{\eta}_{\nu}(\xi,\nu)>0$.
\item The function $(u, v, \tau)(\xi,\eta)$ defined above satisfies (\ref{PEC}) on $\Delta_{+}(r)\setminus \mathrm{P}$.
\item For any fixed $\nu\in[\nu_1,\nu_2]$, $\eta=\hat{\eta}(\xi,\nu)$ gives a $C_{+}$ characteristic curve emanating from $\mathrm{P}$ with the slope $\nu$ at $\mathrm{P}$; that is,
\begin{equation}\label{pgz1}
\partial_{\xi}\hat{\eta}(\xi, \nu)=\hat{\lambda}_{+}(\xi, \nu),\quad
    \hat{\eta}(\xi_{_\mathrm{P}},\nu)=\eta_{_\mathrm{P}},
\quad
\partial_{\xi}\hat{\eta}(\xi_{_\mathrm{P}},\nu)=\nu.
\end{equation}
Here, $\hat{\lambda}_{\pm}(\xi, \nu)=\lambda_{\pm}\big(\hat{u}(\xi, \nu)-\xi, \hat{v}(\xi, \nu)-\hat{\eta}(\xi, \nu), \hat{c}(\xi, \nu)\big)$ and $\hat{c}(\xi, \nu)=c(\hat{\tau}(\xi, \nu))$.
\item $\nu=\nu_1$ and $\nu=\nu_2$ correspond to $\eta=\eta_1(\xi)$ and $\eta=\eta_2(\xi)$, respectively.
\end{enumerate}
Let $(\widetilde{u}_{+},\widetilde{v}_{+}, \widetilde{\tau}_{+})(\nu)=(\hat{u},\hat{v}, \hat{\tau})(\xi_{_\mathrm{P}},\nu)$ for $\nu_1\leq \nu\leq \nu_2$.
 $(\widetilde{u}_{+},\widetilde{v}_{+}, \widetilde{\tau}_{+})(\nu)$ is called the principal part of this $C_{+}$ centered wave, and $\nu_2-\nu_1$ is called the amplitude of the centered wave.
\end{defn}

We
set
\begin{equation}\label{1802}
\tan(\hat{\alpha}(\xi, \nu))=\hat{\lambda}_{+}(\xi, \nu),\quad
\tan(\hat{\beta}(\xi, \nu))=\hat{\lambda}_{-}(\xi, \nu).
\end{equation}
Thus,
from (\ref{pgz1}) we  have that for the $C_{+}$ centered wave defined above,
\begin{equation}
\frac{\partial \hat{\eta}(\xi,\nu)}{\partial \nu}=\int_{\xi_{_\mathrm{P}}}^{\xi}\sec^{2}\hat{\alpha}
(\xi,\nu)\frac{\partial \hat{\alpha}(\xi, \nu)}{\partial \nu} d\xi.
\end{equation}

Introducing the transformation $\xi=\xi$, $\eta=\hat{\eta}(\xi, \nu)$,  we have
\begin{equation}\label{transd}
\frac{\partial}{\partial \xi}= \frac{\partial}{\partial \xi}-\frac{\partial \hat{\eta}}{\partial \xi} \Big(\frac{\partial \hat{\eta}}{\partial \nu}\Big)^{-1}\frac{\partial}{\partial \nu}\quad\text{and} \quad
\frac{\partial}{\partial \eta}= \Big(\frac{\partial \hat{\eta}}{\partial \nu}\Big)^{-1}\frac{\partial}{\partial \nu}.
\end{equation}
Inserting this into (\ref{PEC}) and (\ref{5802}) and letting $\xi\rightarrow\xi_{_\mathrm{P}}$, we know that the principal part of the $C_{+}$ centered wave satisfies
\begin{equation}\label{principal}
\left\{
  \begin{array}{ll}
   \displaystyle -\nu\widetilde{\rho}_{+}\frac{{\rm d} \widetilde{U}_{+}}{{\rm d} \nu}-
\nu\widetilde{U}_{+}\frac{{\rm d} \widetilde{\rho}_{+}}{{\rm d} \nu} +\widetilde{\rho}_{+}\frac{{\rm d} \widetilde{V}_{+}}{{\rm d} \nu}+
\widetilde{V}_{+}\frac{{\rm d} \widetilde{\rho}_{+}}{{\rm d} \nu}=0,\\[6pt]
 \displaystyle \frac{{\rm d} \widetilde{U}_{+}}{{\rm d}\nu}+\nu\frac{{\rm d} \widetilde{V}_{+}}{{\rm d}\nu}=0,
  \end{array}
\right.
\end{equation}
and
\begin{equation}\label{principal1}
 \frac{\widetilde{U}_{+}^2+\widetilde{V}_{+}^2}{2}+h(\widetilde{\tau}_{+})=\text{Const.}.
\end{equation}
where $\widetilde{U}_{+}=\widetilde{u}_{+}(\nu)-\xi_{_\mathrm{P}}$, $\widetilde{V}_{+}=\widetilde{v}_{+}(\nu)-\eta_{_\mathrm{P}}$, and $\widetilde{c}_{+}(\nu)=c(\widetilde{\tau}_{+}(\nu))$.

Inserting (\ref{transd}) into (\ref{82301})--(\ref{8}) and letting $\xi\rightarrow\xi_{_\mathrm{P}}$, we also have
\begin{equation}\label{8304}
\widetilde{c}_{+}\frac{{\rm d}  \widetilde{\alpha}_{+}}{{\rm d} \nu}=\frac{p''(\widetilde{\tau}_{+})}{2\widetilde{c}_{+}\widetilde{\rho}_{+}^4}\tan \widetilde{A}_{+}\frac{{\rm d}  \widetilde{\rho}_{+}}{{\rm d} \nu}
\end{equation}
and
\begin{equation}\label{8303}
\widetilde{c}_{+}\frac{{\rm d}  \widetilde{\beta}_{+}}{{\rm d} \nu}=\frac{p''(\widetilde{\tau}_{+})}{4\widetilde{c}_{+}\widetilde{\rho}_{+}^4}
\big(\omega(\widetilde{\tau}_{+})-\tan^2\widetilde{A}_{+}\big)\sin2 \widetilde{A}_{+}\frac{{\rm d}  \widetilde{\rho}_{+}}{{\rm d} \nu},
\end{equation}
where $(\widetilde{\alpha}_{+},\widetilde{\beta}_{+})(\nu)=(\hat{\alpha}, \hat{\beta})(\xi_{_\mathrm{P}},\nu)$ and $\widetilde{A}_{+}=\frac{\widetilde{\alpha}_{+}-\widetilde{\beta}_{+}}{2}$.

Since
$\tan(\widetilde{\alpha}_{+}(\nu))=\nu$,
 by (\ref{8304}) we have
\begin{equation}\label{120601}
 \frac{{\rm d}\widetilde{\rho}_{+}}{{\rm d} \nu}>0\quad \text{if}\quad \widetilde{\tau}_{+}>\tau_{po}.
\end{equation}
In addition, we have
\begin{equation}\label{82502}
\widetilde{c}_{+}\frac{{\rm d}  \widetilde{\sigma}_{+}}{{\rm d} \nu}=\frac{p''(\widetilde{\tau}_{+})}{8\widetilde{c}_{+}\widetilde{\rho}_{+}^4}
\Big[\big(\omega(\widetilde{\tau}_{+})-\tan^2\widetilde{A}_{+}\big)\sin2 \widetilde{A}_{+}+2\tan \widetilde{A}_{+}\Big]\frac{{\rm d}  \widetilde{\rho}_{+}}{{\rm d} \nu}>0 \quad \text{if}\quad \widetilde{\tau}_{+}>\tau_{po},
\end{equation}
where $\widetilde{\sigma}_{+}=\frac{\widetilde{\alpha}_{+}+\widetilde{\beta}_{+}}{2}$.

\subsubsection{$C_{-}$ centered wave}
\begin{defn}
  Let $\Delta_{-}(r)$ be an angular domain with curved boundaries:
\begin{equation}
\Delta_{-}(r)=\big\{(\xi,\eta)\mid\eta_{_\mathrm{P}}-r\leq \eta\leq\eta_{_\mathrm{P}},~ \xi_{1}(\eta)\leq \xi\leq \xi_{2}(\eta)\big\},
\end{equation}
where
$\xi_{1}(\eta_{_\mathrm{P}})=\xi_{2}(\eta_{_\mathrm{P}})=\xi_{_\mathrm{P}}$
and
$\xi_{1}'(\eta_{_\mathrm{P}})>\xi_{2}'(\eta_{_\mathrm{P}})$.

A function $(u, v, \tau)(\xi,\eta)$ is called a $C_{-}$  centered wave solution of (\ref{PEC}) with center at $\mathrm{P}$ if the following properties hold:
\begin{enumerate}
  \item $(u, v, \tau)$ can be determined implicitly by the functions
$\xi=\breve{\xi}(\eta,\mu)$ and
$(u, v, \tau)=(\breve{u},\breve{v}, \breve{\tau})(\eta,\mu)$
defined on a rectangular domain
$T_{-}(r)=\{(\eta,\mu)\mid\eta_{_\mathrm{P}}-r\leq \eta\leq\eta_{_\mathrm{P}}, \mu_2\leq \mu\leq \mu_1\}$ where $\mu_1=\xi_{1}'(\eta_{_\mathrm{P}})$ and $\mu_2=\xi_{2}'(\eta_{_\mathrm{P}})$.
Moreover, $\breve{\xi}$ and $(\breve{u}, \breve{v}, \breve{\tau})$ belong to $C^{1}\big(T_{-}(r)\big)$, and for any $(\eta,\mu)\in T_{-}(r)\setminus\{\eta=\eta_{_\mathrm{P}}\}$ there holds
$\breve{\xi}_{\mu}(\eta,\mu)<0$.
\item The function $(u, v, \tau)(\xi,\eta)$ defined above satisfies (\ref{PEC}) on $\Delta_{-}(r)\setminus \mathrm{P}$.
\item For any fixed $\mu\in[\mu_2,\mu_1]$, $\xi=\breve{\xi}(\eta,\mu)$ gives a $C_{-}$ characteristic curve emanating from $\mathrm{P}$ with the slope $\mu^{-1}$ at $\mathrm{P}$; that is,
\begin{equation}\label{pgz2}
\partial_{\eta}\breve{\xi}(\eta, \mu)=\breve{\lambda}_{-}^{-1}(\eta, \mu),\quad
    \breve{\xi}(\eta_{_\mathrm{P}},\mu)=\xi_{_\mathrm{P}},
\quad
\partial_{\eta}\breve{\xi}(\eta_{_\mathrm{P}},\mu)=\mu.
\end{equation}
Here, $\breve{\lambda}_{\pm}(\eta, \mu)=\lambda_{\pm}\big(\breve{u}(\eta, \mu)-\breve{\xi}(\eta, \mu), \breve{v}(\eta, \mu)-\eta, \breve{c}(\eta, \mu)\big)$ and $\breve{c}(\eta, \mu)=c(\breve{\tau}(\eta, \mu))$.
\item $\mu=\mu_1$ and $\mu=\mu_2$ correspond to $\xi=\xi_1(\xi)$ and $\xi=\xi_2(\xi)$, respectively.
\end{enumerate}
Let $(\widetilde{u}_{-},\widetilde{v}_{-}, \widetilde{\tau}_{-})(\mu)=(\breve{u},\breve{v}, \breve{\tau})(\eta_{_\mathrm{P}},\mu)$ for $\mu_2\leq \mu\leq \mu_1$.
$(\widetilde{u}_{-},\widetilde{v}_{-}, \widetilde{\tau}_{-})(\mu)$ is called the principal part of this $C_{-}$ centered wave and $\mu_1-\mu_2$ the amplitude of the centered wave.
\end{defn}

We
set
$$
\tan(\breve{\alpha}(\eta, \nu))=\breve{\lambda}_{+}(\eta, \nu),\quad
\tan(\breve{\beta}(\eta, \nu))=\breve{\lambda}_{-}(\eta, \nu).
$$
For the $C_{-}$ centered wave defined above, we have
\begin{equation}
\frac{\partial \breve{\xi}(\eta,\mu)}{\partial \mu}=\int_{\eta}^{\eta_{_\mathrm{P}}}\csc^{2}\breve{\beta}
(\eta,\mu)\frac{\partial \breve{\beta}(\eta, \mu)}{\partial \mu} d\eta.
\end{equation}

Introducing the transformation $\eta=\eta$, $\xi=\breve{\xi}(\eta, \mu)$,  we have
\begin{equation}\label{transdd}
\frac{\partial}{\partial \eta}= \frac{\partial}{\partial \eta}-\frac{\partial \breve{\xi}}{\partial \eta} \Big(\frac{\partial \breve{\xi}}{\partial \mu}\Big)^{-1}\frac{\partial}{\partial \mu}\quad\text{and} \quad
\frac{\partial}{\partial \xi}= \Big(\frac{\partial \breve{\xi}}{\partial \mu}\Big)^{-1}\frac{\partial}{\partial \mu}.
\end{equation}
Inserting this into (\ref{PEC}) and (\ref{5802}) and letting $\eta\rightarrow\eta_{_\mathrm{P}}$, we know that the principal part of the $C_{-}$ centered wave satisfies
\begin{equation}\label{82306}
\left\{
  \begin{array}{ll}
   \displaystyle -\mu\widetilde{\rho}_{-}\frac{{\rm d} \widetilde{U}_{-}}{{\rm d} \mu}-
\mu\widetilde{U}_{-}\frac{{\rm d} \widetilde{\rho}_{-}}{{\rm d} \mu} +\widetilde{\rho}_{-}\frac{{\rm d} \widetilde{V}_{-}}{{\rm d} \mu}+
\widetilde{V}_{-}\frac{{\rm d} \widetilde{\rho}_{-}}{{\rm d} \mu}=0,\\[6pt]
 \displaystyle \frac{{\rm d} \widetilde{U}_{-}}{{\rm d}\mu}+\mu\frac{{\rm d} \widetilde{V}_{-}}{{\rm d}\mu}=0,
  \end{array}
\right.
\end{equation}
and
\begin{equation}\label{82307}
 \frac{\widetilde{U}_{-}^2+\widetilde{V}_{-}^2}{2}+h(\widetilde{\tau}_{-})=\text{Const.},
\end{equation}
where $\widetilde{U}_{-}=\widetilde{u}_{-}(\mu)-\xi_{_\mathrm{P}}$, $\widetilde{V}_{-}=\widetilde{v}_{-}(\mu)-\eta_{_\mathrm{P}}$, and $\widetilde{c}_{-}(\mu)=c(\widetilde{\tau}_{-}(\mu))$.
Moreover, from (\ref{7}) and (\ref{10}) we have
\begin{equation}\label{8306}
\widetilde{c}_{-}\frac{{\rm d}  \widetilde{\beta}_{-}}{{\rm d} \mu}=-\frac{p''(\widetilde{\tau}_{-})}{2\widetilde{c}_{-}\widetilde{\rho}_{-}^4}\tan \widetilde{A}_{-}\frac{{\rm d}  \widetilde{\rho}_{-}}{{\rm d} \mu}
\end{equation}
and
\begin{equation}\label{8305}
\widetilde{c}_{-}\frac{{\rm d}  \widetilde{\alpha}_{-}}{{\rm d} \mu}=-\frac{p''(\widetilde{\tau}_{-})}{4\widetilde{c}_{-}\widetilde{\rho}_{-}^4}
\big(\omega(\widetilde{\tau}_{-})-\tan^2\widetilde{A}_{-}\big)\sin2 \widetilde{A}_{-}\frac{{\rm d}  \widetilde{\rho}_{-}}{{\rm d} \mu},
\end{equation}
where  $\widetilde{A}_{-}=\frac{\widetilde{\alpha}_{-}-\widetilde{\beta}_{-}}{2}$.

As in (\ref{1802}) we have $\cot(\widetilde{\beta}_{-}(\mu))=\mu$.
So, by (\ref{8306}) we have
\begin{equation}\label{010901}
 \frac{{\rm d}\widetilde{\rho}_{-}}{{\rm d} \mu}>0\quad \text{if}\quad \widetilde{\tau}_{-}>\tau_{po}.
\end{equation}
In addition, we have
\begin{equation}
\widetilde{c}_{-}\frac{{\rm d}  \widetilde{\sigma}_{-}}{{\rm d} \mu}=-\frac{p''(\widetilde{\tau}_{-})}{8\widetilde{c}_{-}\widetilde{\rho}_{-}^4}
\Big[\big(\omega(\widetilde{\tau}_{-})-\tan^2\widetilde{A}_{-}\big)\sin2 \widetilde{A}_{-}+2\tan \widetilde{A}_{-}\Big]\frac{{\rm d}  \widetilde{\rho}_{-}}{{\rm d} \mu}<0 \quad \text{if}\quad \widetilde{\tau}_{-}>\tau_{po},
\end{equation}
where $\widetilde{\sigma}_{-}=\frac{\widetilde{\alpha}_{-}+\widetilde{\beta}_{-}}{2}$.

\subsubsection{Principal parts of the $C_{\pm}$ centered waves}
From (\ref{823301}), (\ref{52204}), and the second equation of (\ref{principal}),
the image of the principal part of a $C_{+}$  centered wave in the $(u, v)$-plane is a $\Gamma_{-}$ characteristic curve.
Similarly, the image of the principal part of a $C_{-}$  centered wave in the $(u, v)$-plane is a $\Gamma_{+}$ characteristic curve.

Since the leading characteristic of  the $C_{+}$ centered wave is $\widetilde{\mathrm{PB}}$, we
set
 \begin{equation}\label{82904bzt}
\nu_2=\tan(\alpha_2(\mathrm{P}))
\end{equation}
and
\begin{equation}\label{82304}
(\widetilde{u}_{+},\widetilde{v}_{+}, \widetilde{\rho}_{+})(\nu_2)=(u_{-}(\xi_{_P}), 0, \rho_{po}).
\end{equation}
%where $\rho_g=1/\tau_g$.
To obtain the principal part of the $C_{+}$ centered wave, we  consider (\ref{principal})--(\ref{principal1}) with  initial data (\ref{82304}).

Since the leading characteristic of  the $C_{-}$ centered wave is $\widetilde{\mathrm{PD}}$, we
set
 \begin{equation}\label{82904}
\mu_1=\cot(\beta_3(\mathrm{P}))
\end{equation}
and
\begin{equation}\label{82304a}
(\widetilde{u}_{-},\widetilde{v}_{-}, \widetilde{\rho}_{-})(\mu_1)=(u_{po}, v_{po}, \rho_{po}).
\end{equation}
%where $\rho_g=1/\tau_g$.
To obtain the principal part of the $C_{-}$ centered wave, we consider (\ref{82306})--(\ref{82307}) with  initial data (\ref{82304a}).

Assume that
\begin{equation}\label{120302}
\int_{q_{_\mathrm{P}}}^{q_{_{lim}}}\frac{\sqrt{q^2-c^2(\grave{\tau}(q))}}{qc(\grave{\tau}(q))}{\rm d}q>\frac{\sigma_2(\mathrm{P})-\sigma_3(\mathrm{P})}{2},
\end{equation}
where $q_{_{\mathrm{P}}}
=\sqrt{(u_{po}-\xi_{_\mathrm{P}})^2+(v_{po}-\eta_{_\mathrm{P}})^2}$, the function $\grave{\tau}(q)$ is determined by
$$
\frac{1}{2}q^2+\int_{\tau_{po}}^{\grave{\tau}(q)}\tau p'(\tau)~{\rm d}\tau=\frac{1}{2}q_{_{\mathrm{P}}}^2,
$$
and $q_{_{lim}}$ is determined by
$$
\frac{1}{2}q_{_{lim}}^2+\int_{\tau_{po}}^{+\infty}\tau p'(\tau)~{\rm d}\tau=\frac{1}{2}q_{_{\mathrm{P}}}^2.
$$
Then, by assumption ($h_{21}$), there exists a $q_{m}\in (q_{_{\mathrm{P}}}, q_{_{lim}})$ such that
$$
\int_{q_{_\mathrm{P}}}^{q_{m}}\frac{\sqrt{q^2-c^2(\grave{\tau}(q))}}{qc(\grave{\tau}(q))}{\rm d}q=\frac{\sigma_2(\mathrm{P})-\sigma_3(\mathrm{P})}{2}
.
$$
%In what follows we only discuss the first case which can occur when $u_0$ is close to $\frac{\rho_{po}c_{po}}{\rho_0}$.

We introduce the following constants:
$$
\begin{aligned}
&\sigma_m=\frac{\sigma_2(\mathrm{P})+\sigma_3(\mathrm{P})}{2}, \quad \tau_m=\grave{\tau}(q_m), \quad A_m=\arcsin\Big(\frac{c(\tau_m)}{q_m}\Big), \\ &\qquad \qquad\qquad~ \alpha_m=\sigma_m+A_m, \quad \text{and}\quad \beta_m=\sigma_m-A_m.
\end{aligned}
$$
We also define the constants
\begin{equation}\label{121405}
\nu_1=\tan\alpha_m\quad \text{and} \quad \mu_2=\cot\beta_m.
\end{equation}

\begin{lem}\label{8250c}
Under assumption (\ref{120302}), the initial value problem (\ref{principal})--(\ref{principal1}), (\ref{82304}) admits a solution on $[\nu_1, \nu_2]$, the initial value problem (\ref{82306})--(\ref{82307}), (\ref{82304a}) admits a solution on $[\mu_2, \mu_1]$, and the solutions satisfy
$$
(\widetilde{u}_{-},\widetilde{v}_{-}, \widetilde{\rho}_{-})(\mu_2)=(\widetilde{u}_{+},\widetilde{v}_{+}, \widetilde{\rho}_{+})(\nu_1).
$$
\end{lem}

\begin{proof}
We introduce the Riemann invariants:
$$
R_{\pm}(\widetilde{U}, \widetilde{V})=\widetilde{\sigma}\pm \int_{q_{_\mathrm{P}}}^{\widetilde{q}}\frac{\sqrt{q^2-c^2(\grave{\tau}(q))}}{qc(\grave{\tau}(q))}{\rm d}q,
$$
where $\widetilde{q}=\sqrt{\widetilde{U}^2+\widetilde{V}^2}$ and $\widetilde{\sigma}=\arctan\Big(\frac{\widetilde{V}}{\widetilde{U}}\Big)$.

By the second equation of (\ref{principal}), we have
$$
\begin{array}{rcl}
&&\displaystyle\frac{{\rm d}}{{\rm d} \nu} R_{+}\big(\widetilde{U}_{+}(\nu), \widetilde{V}_{+}(\nu)\big)
\\[8pt]&=&\displaystyle\frac{\widetilde{U}_{+}\frac{{\rm d}\widetilde{V}_{+}}{{\rm d} \nu}-\widetilde{V}_{+}\frac{{\rm d}\widetilde{U}_{+}}{{\rm d} \nu}}{\widetilde{q}_{+}^2}
+\frac{\sqrt{\widetilde{q}_{+}^2-\widetilde{c}_{+}^2}}{\widetilde{q}_{+}\widetilde{c}_{+}}\cdot\frac{\widetilde{U}_{+}\frac{{\rm d}\widetilde{U}_{+}}{{\rm d} \nu}+\widetilde{V}_{+}\frac{{\rm d}\widetilde{V}_{+}}{{\rm d} \nu}}{\widetilde{q}_{+}}
\\[8pt]&=&\displaystyle\frac{1}{\widetilde{q}_{+}}\left(\cos\widetilde{\sigma}_{+}\frac{{\rm d}\widetilde{V}_{+}}{{\rm d} \nu}-\sin\widetilde{\sigma}_{+}\frac{{\rm d}\widetilde{U}_{+}}{{\rm d} \nu}+\cot\widetilde{A}_{+}\Big(
\cos\widetilde{\sigma}_{+}\frac{{\rm d}\widetilde{U}_{+}}{{\rm d} \nu}+\sin\widetilde{\sigma}_{+}\frac{{\rm d}\widetilde{V}_{+}}{{\rm d} \nu}\Big)\right)
\\[8pt]&=&\displaystyle\frac{1}{\widetilde{q}_{+}\sin\widetilde{A}_{+}}\bigg(\cos\widetilde{\alpha}_{+}\frac{{\rm d}\widetilde{U}_{+}}{{\rm d} \nu}+\sin\widetilde{\alpha}_{+}\frac{{\rm d}\widetilde{V}_{+}}{{\rm d} \nu}\bigg)~=~0,
\end{array}
$$
where $\widetilde{q}_{+}=\sqrt{\widetilde{U}_{+}^2+\widetilde{V}_{+}^2}$ and $\widetilde{\sigma}_{+}=\arctan\Big(\frac{\widetilde{V}_{+}}{\widetilde{U}_{+}}\Big)$.
Similarly, we have
$$
\frac{{\rm d}}{{\rm d} \mu} R_{-}\big(\widetilde{U}_{-}(\mu), \widetilde{V}_{-}(\mu)\big)=0.
$$

To determined the principal parts of the $C_{\pm}$ centered waves, we solve
\begin{equation}\label{120301}
\left\{
  \begin{array}{ll}
    \displaystyle\widetilde{\sigma}+ \int_{q_{_\mathrm{P}}}^{\widetilde{q}}\frac{\sqrt{q^2-c^2(\grave{\tau}(q))}}{qc(\grave{\tau}(q))}{\rm d}q=\sigma_2(\mathrm{P}),\\[12pt]
   \displaystyle \widetilde{\sigma}- \int_{q_{_\mathrm{P}}}^{\widetilde{q}}\frac{\sqrt{q^2-c^2(\grave{\tau}(q))}}{qc(\grave{\tau}(q))}{\rm d}q=\sigma_3(\mathrm{P}),
  \end{array}
\right.
\end{equation}
where $\widetilde{q}\in (q_{_{\mathrm{P}}}, q_{_{lim}})$.
Under assumption (\ref{120302}), the nonlinear equations (\ref{120301}) admits a unique solution $(\widetilde{q}, \widetilde{\sigma})=(q_m, \sigma_m)$, with $\widetilde{q}\in (q_{_{\mathrm{P}}}, q_{_{lim}})$.

From (\ref{principal1}), (\ref{120601}), and (\ref{82502}), we have
\begin{equation}\label{82903}
\widetilde{q}_{+}'(\nu)<0,\quad \widetilde{\sigma}_{+}'(\nu)>0,\quad  \text{and}\quad \widetilde{\rho}_{+}'(\nu)>0\quad  \text{for}\quad \nu<\nu_2.
\end{equation}

Hence, by $R_{+}(\widetilde{U}_{+}(\nu), \widetilde{V}_{+}(\nu))\equiv \sigma_2(\mathrm{P})$ we know that the initial value problem (\ref{principal})--(\ref{principal1}), (\ref{82304}) admits a solution on $[\nu_1, \nu_2]$, and
$(\widetilde{U}_{+},\widetilde{V}_{+}, \widetilde{\tau}_{+})(\nu_1)=(q_m\cos\sigma_m, q_m\sin\sigma_m, \tau_m)$.

Similarly, the initial value problem (\ref{82306})--(\ref{82307}), (\ref{82304a}) admits a solution on $[\mu_2, \mu_1]$, and
$(\widetilde{U}_{-},\widetilde{V}_{-}, \widetilde{\tau}_{-})(\mu_2)=(q_m\cos\sigma_m, q_m\sin\sigma_m, \tau_m)$.
This completes the proof.
\end{proof}

\begin{rem}
When $u_0$ is close to $\frac{\rho_{po}c_{po}}{\rho_0}$, assumption (\ref{120302}) holds.
The assumption (\ref{120302}) is not necessary. If this assumption is not satisfied, then there is a vacuum between the two centered wave.
\end{rem}

Let $\Gamma_{-}^{\mathrm{P}}$ be the curve $(u, v)=(\widetilde{u}_{+},\widetilde{v}_{+})(\nu)$, $\nu_1\leq \nu\leq\nu_2$, and let $\Gamma_{+}^{\mathrm{P}}$ be the curve $(u, v)=(\widetilde{u}_{-},\widetilde{v}_{-})(\mu)$, $\mu_2\leq\mu\leq\mu_1$ in the hodograph $(u, v)$-plane.
Let
\begin{equation}
\big(\widetilde{\alpha}_{h},\widetilde{\beta}_{h}, \widetilde{\tau}_{h}\big)(u, v)
=\left\{
                                                                    \begin{array}{ll}
      (\widetilde{\alpha}_{-}, \widetilde{\beta}_{-}, \widetilde{\tau}_{-})(\mu), & \hbox{$(u, v)=(\widetilde{u}_{-},\widetilde{v}_{-})(\mu)$, $\mu_2\leq\mu\leq\mu_1$;} \\[4pt]
     (\widetilde{\alpha}_{+}, \widetilde{\beta}_{+}, \widetilde{\tau}_{+})(\nu), & \hbox{$(u, v)=(\widetilde{u}_{+},\widetilde{v}_{+})(\nu)$, $\nu_1\leq \nu\leq\nu_2$.}
                                                                             \end{array}
                                                                           \right.
\end{equation}

%Since $\xi\sin\theta-\eta\cos\theta=\zeta_1$ ($S_1$) is tangent to the sonic circle $(\xi-w_2\sin\theta)^2+(\eta+w_2\cos\theta)^2=c_2^2$, we have
%\begin{equation}\label{122006}
%\theta+\sigma_1=\arcsin\frac{c_{b}}{q_{b}}.
%\end{equation}

%Actually, the principal part of the $C_{+}$ centered wave to the DGP  (\ref{42501}, \ref{61602}) is
%\begin{equation}\label{9104}
%(u, v, \tau)=\big(\widetilde{u}_{+}(\nu), \widetilde{v}_{+}(\nu), \widetilde{\tau}_{+}(\nu)\big),\quad\nu_1\leq\nu\leq\nu_{2},
%\end{equation}
%and the  principal part of the $C_{-}$ centered wave is
%\begin{equation}\label{9105}
%(u, v, \tau)=\big(\widetilde{u}_{-}(\nu), \widetilde{v}_{-}(\nu), \widetilde{\tau}_{-}(\nu)\big),\quad-\nu_{2}\leq\nu\leq -\nu_1.
%\end{equation}\breve{}

\vskip 4pt
From (\ref{120601}), we see that the function $\tau=\widetilde{\tau}_{+}(\nu)$, $\nu\in[\nu_1,\nu_2]$ has an inverse function $\nu=\widetilde{\nu}_{+}(\tau)$, $\tau\in[\tau_{po}, \tau_{m}]$.
Let $(\check{\alpha}_{+}, \check{\beta}_{+}, \check{\sigma}_{+}, \check{q}_{+})(\tau)=(\widetilde{\alpha}_{+}, \widetilde{\beta}_{+}, \widetilde{\sigma}_{+}, \widetilde{q}_{+})(\widetilde{\nu}_{+}(\tau))$. Then by (\ref{8304}) and (\ref{8303}) we have
\begin{equation}\label{61902}
\left\{
  \begin{array}{ll}
\displaystyle \frac{{\rm d}  \check{\alpha}_{+}}{{\rm d}  \tau}=-\frac{\tau^2p''(\tau)}{2c^2(\tau)}\tan \check{A}_{+},  \\[10pt]
  \displaystyle  \frac{{\rm d}  \check{\beta}_{+}}{{\rm d}  \tau }=-\frac{\tau^2p''(\tau)}{4 c^2(\tau)
  }\big(\omega(\tau)-\tan^2\check{A}_{+}\big)\sin2 \check{A}_{+},
  \end{array}
\right.
\end{equation}
where $\check{A}_{+}=\frac{\check{\alpha}_{+}-\check{\beta}_{+}}{2}$.
From (\ref{82904bzt}) and (\ref{121405}), we obtain
\begin{equation}\label{61803A}
(\check{\alpha}_{+},\check{\beta}_{+})(\tau_{po})=\left(\alpha_2(\mathrm{P}), -\frac{\pi}{2}\right)
\end{equation}
and
\begin{equation}\label{102501a}
   (\check{\alpha}_{+},\check{\beta}_{+})(\tau_m)=(\alpha_m, \beta_m).
\end{equation}

From (\ref{8306}), we obtain that the function $\tau=\widetilde{\tau}_{-}(\mu)$ has an inverse function $\mu=\widetilde{\mu}_{-}(\tau)$.
Let  $(\check{\alpha}_{-}(\tau), \check{\beta}_{-}(\tau))=\big(\widetilde{\alpha}_{-}(\widetilde{\mu}_{-}(\tau)), \widetilde{\beta}_{-}(\widetilde{\mu}_{-}(\tau))\big)$. Then by (\ref{8306})--(\ref{8305}) we have
\begin{equation}\label{120602}
\left\{
  \begin{array}{ll}
\displaystyle \frac{{\rm d}  \check{\beta}_{-}}{{\rm d} \tau}=\frac{\tau^2p''(\tau)}{2c^2(\tau)}\tan \check{A}_{-},  \\[8pt]
  \displaystyle \frac{{\rm d}  \check{\alpha}_{-}}{{\rm d}  \tau}=\frac{\tau^2p''(\tau)}{4c^2(\tau)
}\big(\omega(\tau)-\tan^2\check{A}_{-}\big)\sin2 \check{A}_{-},
  \end{array}
\right.
\end{equation}
where $\check{A}_{-}=\frac{\check{\alpha}_{-}-\check{\beta}_{-}}{2}$.
Similarly, we have
\begin{equation}\label{120603}
(\check{\alpha}_{-},\check{\beta}_{-})(\tau_{po})=(\phi_{po}, \beta_3(\mathrm{P}))
\end{equation}
and
\begin{equation}\label{61401}
   (\check{\alpha}_{-},\check{\beta}_{-})(\tau_m)=(\alpha_m, \beta_m).
\end{equation}

 In the following, we divide the discussion into  the following two cases: $0<\theta_w\leq\alpha_v+\pi$; $\alpha_v+\pi<\theta_w<\pi+\sigma_{po}$.

%From (\ref{72803}) one can see
%\begin{equation}\label{82401}
%A_m\rightarrow A_d\quad \text{as}\quad \tau_1^e-\tau_0\rightarrow 0,
%\end{equation}
%where the constant $A_d$ is defined in (\ref{101801}).

%\subsubsection{\bf Some important estimates for the principal parts}

\subsection{Wave interaction for the case $0<\theta_w\leq \alpha_v+\pi$}
%We first study the case of $0<\theta_w\leq \alpha_v+\pi$.
\subsubsection{Main theorem}
To obtain a global solution to the 2D Riemann problem, we also make the following assumptions:
\begin{description}
  \item[($h_{22}$)] $\theta>\beta_3(\mathrm{P})+2\Theta(\bar{\tau}(\theta))$ for  $\theta\in[\alpha_v, \alpha_0]$;
\item[($h_{23}$)] $\check{\alpha}_{+}(\tau)>\beta_3(\mathrm{P})+2\Theta(\tau)$ for  $\tau\in[\tau_{po}, \tau_m]$;
   \item[($h_{24}$)] $\check{\beta}_{-}(\tau)<\alpha_2(\mathrm{P})-2\Theta(\tau)$ for  $\tau\in[\tau_{po}, \tau_m]$.
\end{description}
The main result is stated as follows.
\begin{thm}\label{thmfs+fs1}
Assume $\tau_1^i<\tau_0<\tau_2^i$ and $0<\theta_w\leq \alpha_v+\pi$. Assume further that assumptions ($h_{21}$)--($h_{24}$) hold.
Then the 2D Riemann problem (\ref{Potential}), (\ref{initial})--(\ref{boundary})  admits a global-in-time self-similar solution.
\end{thm}

The problem can be reduced to the interaction between the two shock–fan composite waves in the $(\xi, \eta)$-plane.
Hence, we consider (\ref{PEC}) with the boundary conditions
\begin{equation}\label{I21}
(u, v, \tau)=\left\{
               \begin{array}{ll}
                 (u_{-}, v_{-}, \tau_{-})(\xi), & \hbox{$(\xi, \eta)\in \widetilde{\mathrm{PB}}$;} \\[2pt]
                 (\bar{u}, \bar{v}, \bar{\tau})(\theta), & \hbox{$(\xi, \eta)\in\widetilde{\mathrm{PD}}$.}
               \end{array}
             \right.
\end{equation}
 Problem (\ref{PEC}), (\ref{I21}) is a discontinuous Goursat problem (DGP).
In what follows, we are going to construct a continuous and piecewise smooth solution to this DGP.

\subsubsection{Invariant region}
Let $\varepsilon$ be a small positive constant such that
$$
\alpha_2(\mathrm{P})-\beta_3(\mathrm{P})+2\varepsilon<\pi.
$$
%\begin{equation}\label{8220a}
%\check{\alpha}_{+}(\tau)>\beta_3(\mathrm{P})+2\Theta(\tau)+\varepsilon\quad \text{and}\quad
%\check{\beta}_{-}(\tau)<\alpha_2(\mathrm{P})-2\Theta(\tau)-\varepsilon \quad \text{for}\quad \tau\in[\tau_{po}, \tau_m];
%\end{equation}

For $\tau\geq \tau_{po}$, we define
$$
\Gamma_{20}(\tau)=\Big\{(\alpha, \beta)\mid \alpha_{2l}(\tau)<\alpha<\alpha_{2r},~\beta_{2l}<\beta<\beta_{2r}(\tau),~\alpha-\beta>\delta_m\Big\},
$$
where
$$
\begin{aligned}
&\alpha_{2l}(\tau)=\beta_3(\mathrm{P})-\varepsilon+2\Theta(\tau), \quad \beta_{2l}=\beta_3(\mathrm{P})-\varepsilon,\quad \alpha_{2r}=\alpha_2(\mathrm{P})+\varepsilon, \\& \beta_{2r}(\tau)=\alpha_2(\mathrm{P})+\varepsilon-2\Theta(\tau),\quad \text{and} \quad \delta_m=\min\limits_{\tau\in [\tau_{po}, \tau_m]}\arcsin\Big(\frac{c(\tau)}{q_m}\Big).
\end{aligned}
$$

%From (\ref{113006}) we know that in the $(\alpha, \beta)$-plane, the straight line $\alpha+\beta=2(\theta_w-\pi)$ and the line segment $\bar{\gamma}_{12}$ interSection
%From (\ref{113007}) we know that when $\delta_2>0$ is sufficiently small, the straight line $\alpha+\beta=2(\theta_w-\pi)$ and the line segment $\bar{\gamma}_{17}$ interSection

\begin{lem}\label{62901}
Suppose that assumptions ($h_{23}$) and ($h_{24}$) hold.
Then,
$$
(\check{\alpha}_{-},\check{\beta}_{-})(\tau)\in \Gamma_{20}(\tau)\quad \text{and}\quad (\check{\alpha}_{+},\check{\beta}_{+})(\tau)\in \Gamma_{20}(\tau)\quad \text{for}\quad \tau\in [\tau_{po}, \tau_m].
$$
\end{lem}
\begin{proof}
From (\ref{61803A}) and (\ref{4905}), we know that there exists a sufficiently small $\epsilon>0$ such that
$(\check{\alpha}_{+},\check{\beta}_{+})(\tau)\in \Gamma_{20}(\tau)$ for $\tau\in [\tau_{po}, \tau_{po}+\epsilon)$.

We next prove the following assertion.
For any fixed $\tau_s\in(\tau_{po}, \tau_m]$, if $(\check{\alpha}_{+},\check{\beta}_{+})(\tau)\in \Gamma_{20}(\tau)$ for $\tau\in [\tau_{po}, \tau_s)$ then
$(\check{\alpha}_{+},\check{\beta}_{+})(\tau_s)\in \Gamma_{20}(\tau_s)$.

By the definition of $\delta_m$, we have
 $$(\check{\alpha}_{+}-\check{\beta}_{+})(\tau_{s})=2\arcsin\left( \frac{c(\tau_{s})}{\check{q}_{+}(\tau_{s})}\right)\geq 2\delta_m>\delta_m.$$
By assumption ($h_{23}$), we have
\begin{equation}\label{62501}
\alpha_{2l}(\tau_s)<\check{\alpha}_{+}(\tau_{s})<\alpha_{2r}.
\end{equation}
We next  prove $\beta_{2l}<\check{\beta}_{+}(\tau_{s})<\beta_{2r}(\tau_s)$.

Suppose $\check{\beta}_{+}(\tau_{s})=\beta_{2r}(\tau_s)$. Then by the assumption in the assertion we have
$
\check{\beta}_{+}'(\tau_{s})\geq 0
$.
Moreover, by (\ref{62501}) we have
 $$\check{A}_{+}(\tau_{s})
=\frac{\check{\alpha}_{+}(\tau_{s})-(\alpha_2(\mathrm{P})+\varepsilon-2\Theta(\tau_{s}))}{2}<\Theta(\tau_{s}).$$
However, by the second equation of (\ref{61902})  we have
 $$
\frac{{\rm d}  \check{\beta}_{+}}{{\rm d}\tau}=-\frac{\tau^2p''(\tau)}{4 c^2(\tau)
  }\big(\tan^2 \Theta-\tan^2\check{A}_{+}\big)\sin2 \check{A}_{+}<0\quad \text{for}\quad \tau=\tau_{s}.
$$
This leads to a contradiction. Hence, we have
$\check{\beta}_{+}(\tau_{s})<\beta_{2r}(\tau_s)$.

 Suppose $\check{\beta}_{+}(\tau_{s})=\beta_{2l}$. Then, by the assumption in the assertion,  we have
$
\check{\beta}_{+}'(\tau_{s})\leq 0
$.
Moreover, by (\ref{62501}) we have
 $$\check{A}_{+}(\tau_{s})
=\frac{\check{\alpha}_{+}(\tau_{s})-(\beta_3(\mathrm{P})-\varepsilon)}{2}>\Theta(\tau_{s}).$$
However, by the second equation of (\ref{61902}) we have
 $$
\frac{{\rm d}  \check{\beta}_{+}}{{\rm d}  \tau }=-\frac{\tau^2p''(\tau)}{4 c^2(\tau)
  }\big(\tan^2 \Theta-\tan^2\check{A}_{+}\big)\sin2 \check{A}_{+}>0\quad \text{at}\quad \tau=\tau_{s}.
$$
This leads to a contradiction.
Hence, we have $\check{\beta}_{+}(\tau_{s})>\beta_{2l}$.
This completes the prove of the assertion.

 Therefore, by a continuity argument we obtain $(\check{\alpha}_{+},\check{\beta}_{+})(\tau)\in \Gamma_{20}(\tau)$ for $\tau\in [\tau_{po}, \tau_m]$.
The proof for the other estimate is similar.
This completes the proof.
\end{proof}

\subsubsection{Local solution to the DGP  (\ref{PEC}), (\ref{I21})}
To overcome the difficulty caused by the multi-valued singularity at the point $\mathrm{P}$,
we apply the hodograph transformation method to find a local solution to the DGP (\ref{PEC}), (\ref{I21}).

%=============
\begin{figure}[htbp]
\begin{center}
\includegraphics[scale=0.65]{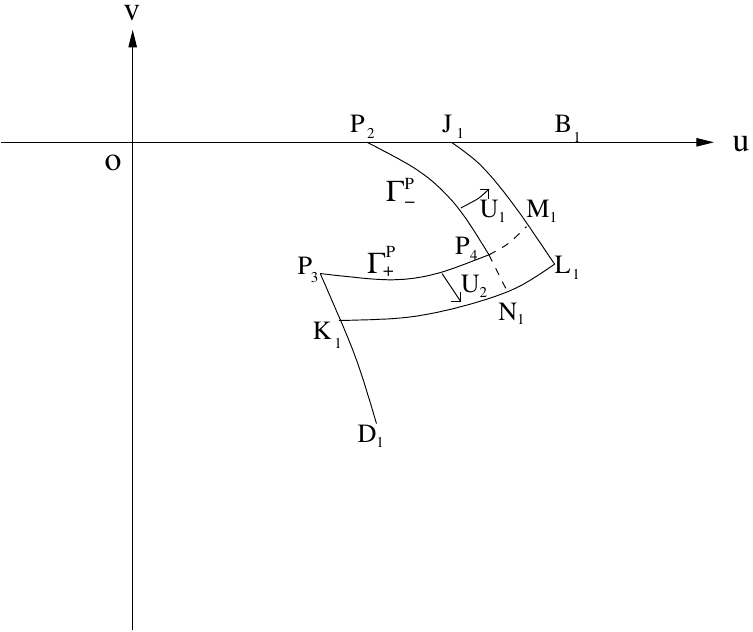}
\caption{\footnotesize Domain $\Sigma_4^h$ and $\Gamma_{\pm}$ characteristic curves in the hodograph plane.}
\label{Fig5.3}
\end{center}
\end{figure}
%==========

Let
$$
\widetilde{\mathrm{P_3D_1}}=\big\{(u, v)~|~u=\bar{u}(\theta), v=\bar{v}(\theta), {\alpha}_v\leq \theta\leq\phi_{po}\big\},
$$
with $\mathrm{P_3}(\bar{u}(\phi_{po}), \bar{v}(\phi_{po}))$ and $\mathrm{D_1}(\bar{u}({\alpha}_v), \bar{v}({\alpha}_v))$.
Let
$$
({\alpha}_{3}^{h}, {\beta}_{3}^{h}, {\tau}_{3}^{h})(u, v)=(\bar{\alpha}, \bar{\beta}, \bar{\tau})(\theta)\quad \text{for}\quad (u, v)=(\bar{u}(\theta), \bar{v}(\theta)),  \quad {\alpha}_v\leq \theta\leq\phi_{po}.
$$

Let
$$
\overline{\mathrm{P_2B_1}}=\big\{(u, v)~|~u=u_{-}(\xi), ~v=0,~ \xi_{_\mathrm{P}}\leq \xi\leq\xi_{_\mathrm{B}}\big\},
$$
with ${\mathrm{P_2}}(\xi_{_P}+c_{po},0)$ and ${\mathrm{B_1}}(u_{-}(\xi_{_B}),0)$.
Let
$$
({\alpha}_{2}^{h}, {\beta}_{2}^{h}, {\tau}_{2}^{h})(u, 0)=\Big(\alpha_2(\xi, \eta_{+}(\xi)), -\frac{\pi}{2}, \tau_{-}(\xi)\Big)\quad \text{for}\quad u=u_{-}(\xi),  \quad  \xi_{_\mathrm{P}}\leq \xi\leq\xi_{_\mathrm{B}}.
$$

We now consider (\ref{52305}) with the boundary conditions
\begin{equation}\label{61701}
 (\alpha, \beta, \tau)=\left\{
  \begin{array}{ll}
(\widetilde{\alpha}_{h},\widetilde{\beta}_{h}, \widetilde{\tau}_{h})(u, v), & \hbox{$(u, v)\in\Gamma_{\pm}^{P}$;} \\[2pt]
   ({\alpha}_{2}^{h}, {\beta}_{2}^{h}, {\tau}_{2}^{h})(u, v), & \hbox{$(u, v)\in\overline{\mathrm{P_2B_1}}$;} \\[2pt]
   ({\alpha}_{3}^{h}, {\beta}_{3}^{h}, {\tau}_{3}^{h})(u, v), & \hbox{$(u, v)\in\widetilde{\mathrm{P_3D_1}}$.}
  \end{array}
\right.
\end{equation}
The problem (\ref{52305}), (\ref{61701}) is a characteristic boundary value problem; see Figure \ref{Fig5.3}.

\begin{lem}
Assume $\epsilon>0$ is sufficiently small, and let the points $\mathrm{K_{1}}=(\bar{u}(\phi_{po}-\epsilon), \bar{v}(\phi_{po}-\epsilon))$ and $\mathrm{J_{1}}=(u_{-}(\xi_{_\mathrm{P}})+\varepsilon, 0)$. The characteristic boundary value problem  (\ref{52305}), (\ref{61701}) admits a continuous and piecewise smooth solution denoted by $(\alpha, \beta, \tau)=(\alpha_h, \beta_h, \tau_h)(u, v)$ on a closed domain $\overline{\Sigma_{4}^{h}}$ bounded by $\Gamma_{\pm}^{P}$, $\overline{\mathrm{P_2J_1}}$, $\widetilde{\mathrm{P_3K_{1}}}$, $\widetilde{\mathrm{J_{1}L_{1}}}$, and $\widetilde{\mathrm{K_{1}L_{1}}}$. Here, $\widetilde{\mathrm{K_{1}L_{1}}}$ is a forward $\Gamma_{+}$ characteristic curve emanating from $\mathrm{K_{1}}$, $\widetilde{\mathrm{J_{1}L_{1}}}$ is a forward $\Gamma_{-}$ characteristic curve emanating from $\mathrm{J_{1}}$, and the two characteristic curves meet at $\mathrm{L_{1}}$. Moreover, the solution satisfies
\begin{equation}\label{62502}
\big(\alpha_{h}(u, v), \beta_h(u, v)\big)\in  \Gamma_{20}\big(\tau_h(u, v)\big)\quad \mbox{for}\quad (u, v)\in \overline{\Sigma_4^h}.
\end{equation}
\end{lem}
\begin{proof}
The existence of a continuous and piecewise $C^1$ solution in $\Sigma_4^h$ follows from  Li and Yu \cite{Li-Yu}. By Lemma \ref{62901}, for sufficiently small $\epsilon>0$, the estimate (\ref{62502}) holds.
\end{proof}

For $(u, v)\in \overline{\Sigma_4^h}$,
we define
$$\xi_4(u, v)=u-\frac{c_h\cos\sigma_h}{\sin A_h}\quad \text{and} \quad \eta_4(u, v)=v-\frac{c_h\sin\sigma_h}{\sin A_h},$$
where $A_h=\frac{\alpha_h-\beta_h}{2}$, $\sigma_h=\frac{\alpha_h+\beta_h}{2}$, and $c_h=c(\tau_h)$.

\begin{lem}
 The mapping $(\xi, \eta)=(\xi_{4}(u, v), \eta_4(u, v))$ is one-to-one from $\overline{\Sigma_4^h}\setminus\Gamma_{\pm}^{P}$ to $\overline{\Sigma_4}\setminus\{\mathrm{P}\}$, where
$$
\overline{\Sigma_4}=\Big\{(\xi,\eta)~\big|~ (\xi, \eta)=(\xi_4(u, v), \eta_4(u,v)), (u, v)\in \overline{\Sigma_4^h}\Big\}.
$$
\end{lem}
\begin{proof}
Since $(\xi, \eta)\equiv (\xi_{_\mathrm{P}}, \eta_{_\mathrm{P}})$ on $\Gamma_{\pm}^{\mathrm{P}}$, we have
\begin{equation}\label{121801}
\mathcal{Z}_{\pm}=0\quad \text{on}\quad \Gamma_{+}^{\mathrm{P}}\cup \Gamma_{-}^{\mathrm{P}}.
\end{equation}

Since $\xi(\bar{u}(\theta), \bar{v}(\theta))=\xi_{-}(\theta)$ for $\bar{\alpha}_v\leq \theta\leq\phi_{po}$, we have
$$
\begin{aligned}
\hat{\partial}_{-}\xi&=\sin\alpha\xi_{u}(u, v)-\cos\alpha\xi_v(u, v)\\[4pt]&
=-\frac{\bar{u}'(\theta)\xi_{u}(\bar{u}(\theta), \bar{v}(\theta))}{\sqrt{(\bar{u}'(\theta))^2+(\bar{v}'(\theta))^2}}
-\frac{\bar{v}'(\theta)\xi_{v}(\bar{u}(\theta), \bar{v}(\theta))}{\sqrt{(\bar{u}'(\theta))^2+(\bar{v}'(\theta))^2}}\\&
=-\frac{\bar{\xi}_{-}'(\theta)}{\sqrt{(\bar{u}'(\theta))^2+(\bar{v}'(\theta))^2}}=
\frac{\cos\beta r_{-}(\theta)}{\sin(2A)\sqrt{(\bar{u}'(\theta))^2+(\bar{v}'(\theta))^2}}\quad \text{on}\quad \widetilde{\mathrm{P_3D_{1}}}.
\end{aligned}
$$
Hence, we get
\begin{equation}\label{121802}
\mathcal{Z}_{-}=\frac{2\sin^2A\hat{\partial}_{-}\xi}{c\cos\beta}=
\frac{\cot A r_{-}(\theta)}{c\sqrt{(\bar{u}'(\theta))^2+(\bar{v}'(\theta))^2}}>0\quad \text{on}\quad \widetilde{\mathrm{P_3D_{1}}}.
\end{equation}

By a direct computation, we have
\begin{equation}\label{121803}
\mathcal{Z}_{+}=-\frac{2\sin^2 A}{c\cos\alpha}\hat{\partial}_{+}\xi=
-\frac{2\sin^2 A}{c\cos\alpha u_{-}'(\xi)}<0\quad \text{on}\quad \overline{\mathrm{P_2J_{1}}}.
\end{equation}

Integrating (\ref{322}) along $\Gamma_{\pm}$ characteristic curves emanating from points on the characteristic boundaries of  the problem (\ref{52305}), (\ref{61701}) and noticing (\ref{121801})--(\ref{121803}), we know that the solution satisfies
\begin{equation}\label{102506}
\mathcal{Z}_{-}>0\quad \text{and}\quad \mathcal{Z}_{+}<0\quad \text{on}\quad \overline{\Sigma_4^h}\setminus\Gamma_{\pm}^{P}.
\end{equation}

Referring to Figure \ref{Fig5.3}, we see that the $\Gamma_{+}$ characteristic curve emanating from $\mathrm{P}_4$ intersects $\widetilde{\mathrm{J_{1}L_{1}}}$ at a point $\mathrm{M_{1}}$, and the $\Gamma_{-}$ characteristic curve emanating from $\mathrm{P}_4$ intersects $\widetilde{\mathrm{K_{1}L_{1}}}$ at a point $\mathrm{N_{1}}$.
Let $\overline{\Sigma_{41}^h}$ be the closed curvilinear quadrilateral domain bounded by $\Gamma_{\pm}$ characteristic curves
$\overline{\mathrm{P_2J_{1}}}$, $\widetilde{\mathrm{P_2P_4}}$, $\widetilde{\mathrm{P_4M_{1}}}$, and $\widetilde{\mathrm{J_{1}M_{1}}}$. Let $\overline{\Sigma_{42}^h}$ be the closed curvilinear quadrilateral domain bounded by $\Gamma_{\pm}$ characteristic curves
$\widetilde{\mathrm{P_3K_{1}}}$, $\widetilde{\mathrm{P_3P_4}}$, $\widetilde{\mathrm{P_4N_{1}}}$, and $\widetilde{\mathrm{K_{1}N_{1}}}$. let $\overline{\Sigma_{43}^h}$ be the closed curvilinear quadrilateral domain bounded by $\Gamma_{\pm}$ characteristic curves
$\widetilde{\mathrm{P_4N_{1}}}$, $\widetilde{\mathrm{M_{1}L_{1}}}$, $\widetilde{\mathrm{P_4M_{1}}}$, and $\widetilde{\mathrm{N_{1}L_{1}}}$. Obviously, $\overline{\Sigma_4^h}=\overline{\Sigma_{41}^h}\cup\overline{\Sigma_{42}^h}\cup\overline{\Sigma_{43}^h}$.
For any two points $\mathrm{U}_1$, $\mathrm{U}_2\in \overline{\Sigma_4^h}\setminus\Gamma_{\pm}^{P}$, one of the following three cases happens:
\begin{enumerate}
  \item $U_1,~U_2\in\big(\overline{\Sigma_{41}^h}\setminus\Gamma_{-}^{P}\big)\cup\big(\overline{\Sigma_{43}^h}\setminus\widetilde{\mathrm{P_4N_{1}}}\big)$;
  \item $U_1,~U_2\in\big(\overline{\Sigma_{42}^h}\setminus\Gamma_{+}^{P}\big)\cup\big(\overline{\Sigma_{43}^h}\setminus\widetilde{\mathrm{P_4M_{1}}}\big)$;
  \item $U_1\in \overline{\Sigma_{41}^h}\setminus\Gamma_{-}^{P}$ and $U_2\in \overline{\Sigma_{42}^h}\setminus\Gamma_{+}^{P}$.
\end{enumerate}
By Lemma \ref{onetoone}, Remark \ref{92505a}, and (\ref{62502}), we have $(\xi, \eta)(\mathrm{U}_1)\neq(\xi, \eta)(\mathrm{U}_2)$ for cases (i) and (ii).

We now consider the third case.
Since  $\mathrm{U}_1\in \overline{\Sigma_{41}^h}\setminus\Gamma_{-}^{\mathrm{P}}$, the backward $\Gamma_{+}$ characteristic emanating from $\mathrm{U}_1$ stays in $\Sigma_{41}^h$ until it intersects $\Gamma_{-}^{\mathrm{P}}$ at a point $\mathrm{U}_{+}$.
Since $\mathrm{U}_2\in \overline{\Sigma_{42}^h}\setminus\Gamma_{+}^{\mathrm{P}}$, the backward $\Gamma_{-}$ characteristic emanating from $\mathrm{U}_2$ stays in $\Sigma_{42}^h$ until it intersects $\Gamma_{+}^{\mathrm{P}}$ at a point $\mathrm{U}_{-}$.
We introduce the coordinate rotation
$$
\xi'=\xi\cos\sigma_m+\eta\sin\sigma_m\quad \text{and} \quad \eta'=\eta\cos\sigma_m-\xi\sin\sigma_m.
$$
Then by $\alpha_2(\mathrm{P})-\beta_3(\mathrm{P})<\pi$
we have that when  $\varepsilon$  is sufficiently small,
$$
0<\alpha<\pi\quad \text{in}\quad \Sigma_{41}^{h}\quad \text{and} \quad -\pi<\beta<0\quad \text{in}\quad \Sigma_{42}^{h}.
$$
Combining this with the fact that $\mathcal{Z}_{\pm}$ are invariant under coordinate rotation,  we have $$\eta'(\mathrm{U_1})>\eta'(\mathrm{U_{+}})=\eta'(\mathrm{U_{-}})>\eta'(\mathrm{U}_2).$$
This implies $(\xi, \eta)(\mathrm{U}_1)\neq(\xi, \eta)(\mathrm{U}_2)$ for the third case.

So, there are no two distinct points in $\overline{\Sigma_4^h}\setminus\Gamma_{\pm}^{P}$ that map to the same point in $\overline{\Sigma_4}\setminus \mathrm{P}$.
This completes the proof of the lemma.
\end{proof}

%Let and $\mathrm{N}_{1}=(\xi_4(\mathrm{N_h}), \eta_4(\mathrm{N_h}))$ be points on the $(\xi, \eta)$-plane.
Let
 $$
\widetilde{\mathrm{JL}}=\Big\{(\xi,\eta)~\big|~ (\xi, \eta)=(\xi_4(u, v), \eta_4(u,v)), (u, v)\in \widetilde{\mathrm{J_{1}L_{1}}}\Big\}
$$
and
 $$
\widetilde{\mathrm{KL}}=\Big\{(\xi,\eta)~\big|~ (\xi, \eta)=(\xi_4(u, v), \eta_4(u,v)), (u, v)\in \widetilde{\mathrm{K_{1}L_{1}}}\Big\}
$$
be two curves on the $(\xi, \eta)$-plane, where $\mathrm{J}=(\xi_4(\mathrm{J}_{1}), \eta_4(\mathrm{J}_{1}))$,  $\mathrm{L}=(\xi_4(\mathrm{L}_{1}), \eta_4(\mathrm{L}_{1}))$, and
$\mathrm{K}=(\xi_4(\mathrm{K_{1}}), \eta_4(\mathrm{K_{1}}))$.
Then by the one-to-one property of the hodograph transformation we obtain a local classical solution to the DGP (\ref{PEC}), (\ref{I21}).
\begin{lem}\label{dgploc}
 The DGP (\ref{PEC}), (\ref{I21}) admits a continuous and piecewise $C^1$ solution on the curvilinear quadrilateral domain $\Sigma_4$ bounded by characteristic curves $\widetilde{\mathrm{PJ}}$, $\widetilde{\mathrm{PK}}$, $\widetilde{\mathrm{JL}}$, and $\widetilde{\mathrm{KL}}$.
\end{lem}

For convenience, we denote the solution by $(u, v, \tau)=(u_4, v_4, \tau_4)(\xi, \eta)$, $(\xi,\eta)\in\Sigma_4$.
%\subsubsection{\bf Hyperbolicity and a priori $C^1$ norm estimate}

\subsubsection{Global solution to the DGP  (\ref{PEC}), (\ref{I21})}
To obtain a global solution to the DGP, it only remains to consider  (\ref{PEC}) with the boundary conditions
\begin{equation}\label{CBD1}
(u, v, \tau)=\left\{
               \begin{array}{ll}
                 (u_{-}, v_{-}, \tau_{-})(\xi), & \hbox{$(\xi, \eta)\in \widetilde{\mathrm{JB}}$;} \\[2pt]
                 (\bar{u}, \bar{v}, \bar{\tau})(\theta), & \hbox{$(\xi, \eta)\in\widetilde{\mathrm{KD}}$;}\\[2pt]
 (u_4, v_4, \tau_4)(\xi, \eta), & \hbox{$(\xi, \eta)\in\widetilde{\mathrm{KL}}\cup \widetilde{\mathrm{JL}}$.}\\
               \end{array}
             \right.
\end{equation}

There exists a sufficiently large integer $n>0$ such that
$$
\mathcal{H}>0\quad \text{for}\quad \tau>\tau_{po}\quad \text{and}\quad 0<A<\frac{\pi}{2},
$$
where the variable $\mathcal{H}$ is as defined in (\ref{82803}).
We define the following constants:
$$
\begin{aligned}
&\mathcal{C}_2=\min\left\{\inf\limits_{\widetilde{\mathrm{KD}}\cup \widetilde{\mathrm{JL}}}\frac{2\bar{\partial}_{-}c}{\sin^2 A},~ \inf\limits_{\widetilde{\mathrm{KL}}\cup \widetilde{\mathrm{JB}}}\frac{2\bar{\partial}_{+}c}{\sin^2 A},~ -\frac{8(\gamma-1)}{(\gamma+1)\tan A_{\infty}},~ \frac{\tau_{*}^{2+n}(2p'(\tau_*)+\tau_*p''(\tau_*))\mathcal{M}_1}{\sqrt{-p'(\tau_*)}}\right\},\\& \mathcal{D}_2=\frac{2\mathcal{C}_2}{\sin^2 \varepsilon_*}, \quad \mathcal{M}_{2}=\min \left\{\inf \limits_{\widetilde{\mathrm{JB}}\cup \widetilde{\mathrm{KL}}}\frac{\rho^n\bar{\partial}_{+}\rho}{\sin^2 A}, ~
\inf \limits_{\widetilde{\mathrm{KD}}\cup \widetilde{\mathrm{JL}}}\frac{\rho^n\bar{\partial}_{-}\rho}{\sin^2 A}\right\},\quad \mathcal{N}_{*}=\max\limits_{\tau_{po}\leq\tau\leq\tau_{*}}\frac{\tau^{4+n}p''(\tau)}{4c(\tau)}.
\end{aligned}
$$
Let $\delta_*<\min\big\{\arctan\big(\frac{-1}{8\mathcal{M}_{1}\mathcal{N}_{*}}\big),~ A_{\infty},~ \delta_0, ~ \delta_m\big\}$ be a small positive constant such that
$$
2+\frac{\sin A}{\varpi}\Big(\frac{\sin^2 A}{\cos A}+\frac{1}{\cos A}-\omega(\tau)\cos A\Big)\mathcal{D}_2>0\quad \text{for}\quad 0<A\leq \frac{\delta_*}{2}\quad\text{and} \quad\tau\geq\tau_*.
$$

For $\tau\geq\tau_{po}$,
we define
$$
\Gamma_{21}(\tau)=\Big\{(\alpha, \beta)\mid \alpha_{2l}(\tau)<\alpha<\alpha_{2r}(\tau),~\beta_{2l}(\tau)<\beta<\beta_{2r}(\tau),~\alpha-\beta>\delta_*\Big\}.
$$
As in Section 4.2, we have the following  estimates:
\begin{itemize}
  \item $(\alpha, \beta)(\xi, \eta_{+}(\xi))\in \Gamma_{21} \big(\tau_{-}(\xi)\big)$ for $\xi\in(\xi_{_\mathrm{P}}, \xi_{_\mathrm{B}})$;
  \item $(\theta, \beta_{-}(\theta))\in \Gamma_{21} (\bar{\tau}(\theta))$ for $\theta\in(\alpha_v, \alpha_0)$.
\end{itemize}
Combining these with (\ref{62502}), we obtain that
\begin{equation}\label{92401a}
(\alpha, \beta)\in \Gamma_{21}(\tau)\quad \text{on}\quad \widetilde{\mathrm{JB}}\cup \widetilde{\mathrm{KD}}\cup  \widetilde{\mathrm{KL}}\cup \widetilde{\mathrm{JL}}.
\end{equation}

Using (\ref{92401a}) and the monotonicity condition $\bar{\partial}_{\pm}\rho<0$ on the $C_{\pm}$ characteristic boundaries, we obtain the following result.
\begin{lem}\label{8230b}
Suppose that assumptions ($h_{21}$)--($h_{24}$) hold.
Then the characteristic boundary value (\ref{PEC}), (\ref{CBD1}) admits a classical solution on the domain $\Sigma_5$
bounded by $\widetilde{\mathrm{JB}}$, $\widetilde{\mathrm{KD}}$, $\widetilde{\mathrm{KL}}$, $\widetilde{\mathrm{JL}}$, and a vacuum boundary connecting $\mathrm{B}$ and $\mathrm{D}$. The solution satisfies
\begin{equation}\label{8220b}
(\alpha, \beta)\in \Gamma_{21}(\tau)\quad \text{and} \quad  2\mathcal{M}_2<\frac{\rho^n\bar{\partial}_{+}\rho}{\sin ^2A}<0\quad \mbox{in}\quad \Sigma_5.
\end{equation}
The solution also satisfies
\begin{equation}\label{8220c}
\mathcal{D}_2<\frac{\bar{\partial}_{+}c}{\sin^2 A}<0\quad  \text{for}\quad  \tau>\tau_*.
\end{equation}
\end{lem}

For convenience we denote by  $(u, v, \tau)=(u_5, v_5, \tau_5)(\xi, \eta)$, $(\xi,\eta)\in\Sigma_5$ the solution to the characteristic boundary value problem (\ref{PEC}), (\ref{CBD1}).
Therefore, the DGP (\ref{PEC}), (\ref{I21}) admits a classical solution in the domain $\Sigma_4\cup \Sigma_5$.

Let
$$
\Sigma_1~=~\big\{(\xi, \eta)~\big|~0\leq \xi\leq \xi_{_\mathrm{P}}, ~~\eta>\xi\tan\phi_{po}\big\}\cup
\big\{(\xi, \eta)~\big|~\xi<0, ~~\eta>0\big\},
$$
$$
\Sigma_2~=~\Big\{(\xi, \eta)~\big|~ \xi_{_\mathrm{P}}<\xi<\xi_{_\mathrm{B}}, ~\eta>\eta_{+}(\xi)\Big\},
$$
$$
\Sigma_3~=~\big\{(\xi, \eta)~\big|~ \xi=r\cos s, ~\eta=r\sin s,~r\in (0, r_{-}(s)),
~s\in(\alpha_v, \alpha_0)\big\}.
$$
We define
\begin{equation}\label{8220d}
(u, v, \tau)=\left\{
                              \begin{array}{ll}
                                (u_0, 0, \tau_0), & \hbox{$(\xi, \eta)\in\Sigma_1$,} \\[2pt]
                                (u_-, v_{-}, \tau_{-})(\xi, \eta), & \hbox{$(\xi, \eta)\in\Sigma_2$,} \\[2pt]
(\bar{u}, \bar{v}, \bar{\tau})(\theta), & \hbox{$(\xi, \eta)\in\Sigma_3$,} \\[2pt]
(u_i, v_i, \tau_i)(\xi, \eta), & \hbox{$(\xi, \eta)\in\Sigma_i$, $i=4, 5$,} \\[2pt]
                                \text{vacuum}, & \hbox{otherwise}.
                              \end{array}
                            \right.
\end{equation}
Then the $(u, v, \tau)$ defined in (\ref{8220d}) is a self-similar solution to the 2D Riemann problem (\ref{Potential}), (\ref{initial})--(\ref{boundary}) for $\tau_1^e\leq \tau_0<\tau_2^i$ and $0<\theta_w<\alpha_v+\pi$.
This completes the proof of Theorem \ref{thmfs+fs1}.

\subsection{Wave interaction for the case ${\theta}_v<\theta_w<\pi+\sigma_{po}$}
When $\theta_w>\alpha_v+\pi$,
the oncoming flow $(0)$ is deflected through an oblique shock–fan composite wave to a constant state $(u_w, v_w, \tau_w)$ at the corner. The flow near the corner can be represented by (\ref{8230a}).

%=============
\begin{figure}[htbp]
\begin{center}
 \includegraphics[scale=0.55]{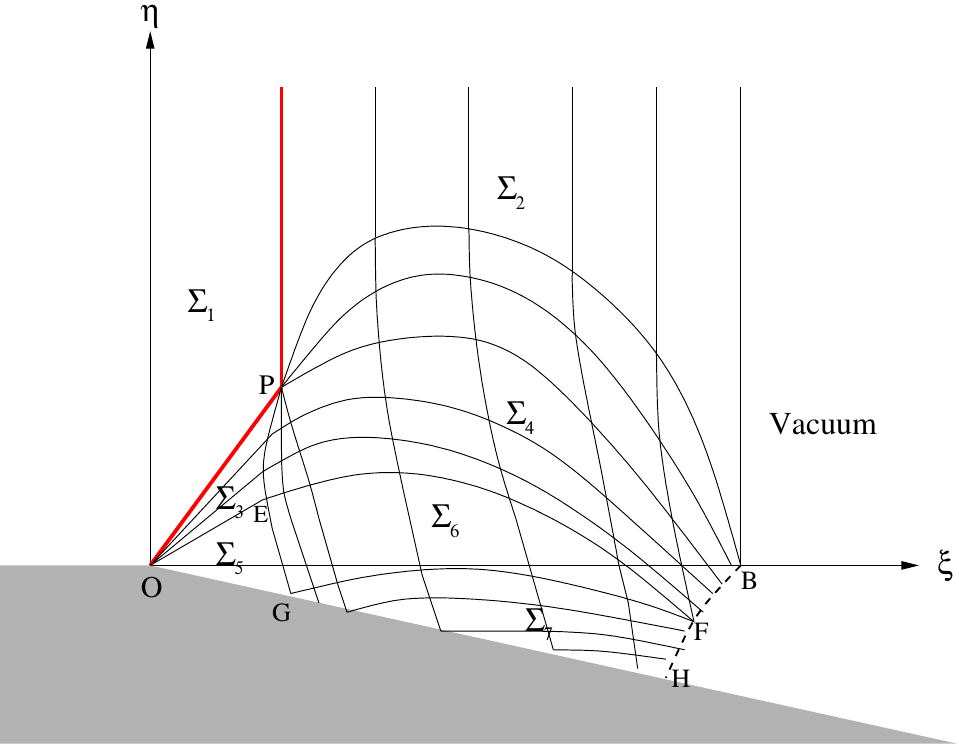}
\caption{\footnotesize Interaction of shock–fan composite waves for $\tau_1^e\leq\tau_0<\tau_2^i$ and  $\theta_w>\alpha_v+\pi$.}
\label{Fig5.4}
\end{center}
\end{figure}
%==========

To construct a global-in-time solution to the 2D Riemann problem (\ref{Potential}), (\ref{initial})--(\ref{boundary})  for $\tau_1^e\leq\tau_0<\tau_2^i$ and  $\theta_w>\alpha_v+\pi$,
we also make the following assumptions:
\begin{description}
%\item[($\overline{H}_{21}$)] $\min\big\{\alpha_2(\mathrm{P})+\frac{\pi}{2},  \phi_{po}-\beta_3(\mathrm{P})\big\}>2\Theta(\tau_{po})$, $\alpha_2(\mathrm{P})-\beta_3(\mathrm{P})<\pi$, $\sigma_2(\mathrm{P})>\sigma_3(\mathrm{P})$;
  \item[($h_{25}$)] $\theta>\beta_3(\mathrm{P})+2\Theta(\bar{\tau}(\theta))$ for  $\alpha_w\leq\theta\leq\phi_{po}$;
\vskip 1pt
%\item[($\overline{H}_{23}$)] $\check{\alpha}_{+}(\tau)>\beta_3(\mathrm{P})+2\Theta(\tau)$ for  $\tau\in[\tau_{po}, \tau_m]$;
 %  \item[($\overline{H}_{24}$)] $\check{\beta}_{-}(\tau)<\alpha_2(\mathrm{P})-2\Theta(\tau)$ for  $\tau\in[\tau_{po}, \tau_m]$.
  \item[($h_{26}$)] $\beta_3(\mathrm{P})+2\Theta(\tau_{po})<\theta_w-\pi<\frac{\pi}{2}+\beta_3(\mathrm{P})$.
\end{description}
The main result is stated as follows.
\begin{thm}\label{thmfs+fs2}
Assume $\tau_1^i<\tau_0<\tau_2^i$ and ${\theta}_v<\theta_w<\pi+\sigma_{po}$. Assume further that assumptions ($h_{21}$) and ($h_{23}$)--($h_{26}$) hold.
Then the 2D Riemann problem (\ref{Potential}), (\ref{initial})--(\ref{boundary})  admits a global-in-time self-similar solution.
\end{thm}

%By the definition of $\sigma_{po}$ we have
%$$
%\frac{{\rm d}}{{\rm d\phi_{po}}}\left(\sigma_{po}-\phi_{po}+2\arctan \left(\frac{\rho_0(1+\sin\phi_{po})}{\rho_{po}\cos\phi_{po}}\right)\right)\Bigg|_{\pi_{po}=\frac{\pi}{2}}=\frac{\rho_{po}}{\rho_0}>0.
%$$
%So, there exists a small $\varepsilon_2>0$ such that when $0<u_0-\frac{\rho_{po}c_{po}}{\rho_0}<\varepsilon_2$ the hold the following inequalities:
%$$
%\phi_{po}-\beta_3(P)>2\Theta(\tau_{po}), \quad \sigma_{po}-\beta_3(P)>\frac{\pi}{2}, \quad
%\sigma_{po}-\phi_{po}+2\arctan \left(\frac{\rho_0(1+\sin\phi_{po})}{\rho_{po}\cos\phi_{po}}\right)<\frac{\pi}{2}.
%$$
%We now assume $0<u_0-\frac{\rho_{po}c_{po}}{\rho_0}<\varepsilon_2$ and $\Theta(\tau_{po})<\frac{\pi}{4}$. Then when $\phi_{po}-\alpha_w>0$ is sufficiently small, we have
%$$
%2\Theta(\tau_{po})<\theta_w-\pi=\arctan\Big(\frac{v_2}{u_w}\Big)-\beta_3(P)<\frac{\pi}{2}
%$$
%and
%$$
%\bar{\alpha}(\theta)>\beta_3(P)+2\Theta(\bar{\tau}(\theta))\quad \text{for}\quad \alpha_w\leq \theta\leq \phi_{po}.
%$$
%This implies that assumptions ($h_{25}$)-($h_{26}$) can be satisfies in some cases.

Let $\varepsilon'$ be a sufficiently small positive constant such that
\begin{equation}\label{8240a}
\alpha_2(\mathrm{P})-\beta_3(\mathrm{P})+3\varepsilon'<\pi\quad \text{and}\quad 2\theta_w-2\beta_3(\mathrm{P})-3\pi+3\varepsilon'<0.
\end{equation}
For any $\tau\geq \tau_{po}$, we define
$$
\Gamma_{22}(\tau)=\Big\{(\alpha, \beta)\mid \bar{\alpha}_{2l}(\tau)<\alpha<\bar{\alpha}_{2r},~\bar{\beta}_{2l}<\beta<\bar{\beta}_{2r}(\tau),~\alpha-\beta>\delta_2\Big\},
$$
where
$$
\begin{aligned}
&\qquad\qquad\bar{\alpha}_{2l}(\tau)=2\Theta(\tau)+\beta_3(\mathrm{P})-\varepsilon', \quad \bar{\beta}_{2l}=\beta_3(\mathrm{P})-\varepsilon', \\& \bar{\alpha}_{2r}=\max\big\{-\beta_3(\mathrm{P})+2(\theta_w-\pi), \alpha_2(\mathrm{P})\big\}+2\varepsilon', \quad \bar{\beta}_{2r}(\tau)=\bar{\alpha}_{2r}-2\Theta(\tau),
\end{aligned}
$$
 and $\delta_2<\delta_*$ is a small positive constant yet to be determined.

By assumption (${h}_{26}$), we have
\begin{equation}\label{113006}
\bar{\alpha}_{2l}(\tau)<2(\theta_w-\pi)-\bar{\beta}_{2l}<\bar{\alpha}_{2r}
\end{equation}
and
\begin{equation}\label{113007}
\bar{\alpha}_{2l}(\tau)<\theta_w-\pi<\bar{\beta}_{2r}(\tau)
\end{equation}
for any $\tau\geq\tau_{po}$.

Let $\widetilde{\mathrm{PE}}$ be the $C_{-}$ cross characteristic curve through $\mathrm{E}$ in the centered simple wave of the oblique shock-fan composite wave.
We consider (\ref{PEC}) with the boundary conditions
\begin{equation}\label{32004}
(u, v, \tau)=\left\{
  \begin{array}{ll}
(u_{-}, v_{-}, \tau_{-})(\xi), & \hbox{$(\xi, \eta)\in\widetilde{\mathrm{PB}}$;} \\[3pt]
    (\bar{u}, \bar{v}, \bar{\tau})(\theta), & \hbox{$(\xi, \eta)\in\widetilde{\mathrm{PE}}$.}
  \end{array}
\right.
\end{equation}
By Lemma \ref{8230b}, the DGP (\ref{PEC}), (\ref{32004}) admits a global continuous and piecewise smooth solution on a domain $\overline{\Sigma_4}$ bounded by $\widetilde{\mathrm{PB}}$, $\widetilde{\mathrm{PE}}$, $\widetilde{\mathrm{EF}}$, and $\widetilde{\mathrm{BF}}$, where $\widetilde{\mathrm{EF}}$ is a forward $C_{+}$ characteristic curve emanating from $\mathrm{E}$ and  $\widetilde{\mathrm{BF}}$ is a vacuum boundary. Moreover, the solution satisfies
\begin{equation}\label{62601}
(\alpha, \beta)\in \Gamma_{22}(\tau)\quad \mbox{in}\quad \overline{\Sigma_4}.
\end{equation}

We draw a forward $C_{-}$ characteristic line from $\mathrm{E}$. %Then the direction of this $C_{-}$ characteristic line is $\beta=\beta(\mathrm{E})$.
Using (\ref{82301}) and (\ref{8}), we have that along this characteristic line,
$\beta\equiv\beta(\mathrm{E})$ and $c\bar{\partial}_{-}(\alpha+\beta)=c\bar{\partial}_{-}\alpha=2\sin^2A>0$.
As in (\ref{1981002}), we have
\begin{equation}\label{8240b}
\beta(\mathrm{E})<\theta_w-\pi.
\end{equation}
By (\ref{62601}) we have
$$
\beta(\mathrm{E})>\bar{\beta}_{2l}=\beta_3(\mathrm{P})-\varepsilon.
$$
Then by (\ref{8240a}) we have
\begin{equation}\label{8240c}
\theta_w-\pi-\beta(\mathrm{E})<\frac{\pi}{2}.
\end{equation}

Therefore, by (\ref{8240b}) and (\ref{8240c}) we know that this $C_{-}$ characteristic line intersects $\partial\mathcal{W}_{+}$ at a point $\mathrm{G}$ before reaching the sonic circle of the state $(u_w, v_w, \tau_w)$.
Let $\overline{\Sigma_5}$ be the triangle domain bounded by $\overline{\mathrm{OE}}$, $\overline{\mathrm{EG}}$, and $\overline{\mathrm{OG}}$.
As in Section 4.3.3, the solution in $\overline{\Sigma_5}$ is the constant state $(u_w, v_w, \tau_w)$.
Moreover, we have $(\alpha, \beta)\in \Gamma_{22}(\tau_w)$ on $\overline{\mathrm{EG}}$.

We denote the solution of the Goursat problem (\ref{PEC}), (\ref{I2}) by $(u, v, \tau)=(u_4, v_4, \tau_4)(\xi, \eta)$, $(\xi, \eta)\in \overline{\Sigma_4}$.
We next consider (\ref{PEC}) with the boundary conditions
\begin{equation}\label{32002}
(u, v, \tau)=\left\{
                       \begin{array}{ll}
                          (u_{4}, v_{4}, \tau_{4})(\xi, \eta) & \hbox{on $\widetilde{\mathrm{EF}}$,} \\[2pt]
                         (u_w, v_w, \tau_w) & \hbox{on $\overline{\mathrm{EG}}$.}
                       \end{array}
                     \right.
\end{equation}
The Goursat problem (\ref{PEC}), (\ref{32002}) admits a simple wave solution with straight $C_{+}$ characteristic lines on the curvilinear triangular domain $\overline{\Sigma_6}$ bounded by $\widetilde{\mathrm{EF}}$, $\overline{\mathrm{EG}}$, and $\widetilde{\mathrm{GF}}$, where $\widetilde{\mathrm{GF}}$ is a forward $C_{+}$ characteristic curve emanating from $\mathrm{G}$. Moreover, the solution satisfies (\ref{62601}).

We denote the simple wave solution of the Goursat problem (\ref{PEC}), (\ref{32002}) by $(u, v, \tau)=(u_6, v_6, \tau_6)(\xi, \eta)$.
We now consider (\ref{PEC}) with the boundary conditions
\begin{equation}\label{32003}
\left\{
               \begin{array}{ll}
                 (u, v, \tau)=(u_{6}, v_{6}, \tau_{6})(\xi,\eta) & \hbox{on $\widetilde{\mathrm{GF}}$;} \\[2pt]
                v= u\tan\theta_w & \hbox{on $\partial\mathcal{W}_{+}$.}
               \end{array}
             \right.
\end{equation}
As in Lemma \ref{lem419},
The mixed boundary value problem (\ref{PEC}), (\ref{32003}) admits a classical solution on the triangular domain $\Sigma_7$ bounded by $\widetilde{\mathrm{GF}}$, $\partial\mathcal{W}_{+}$, and a vacuum boundary $\widetilde{\mathrm{FH}}$.

We then obtain a self-similar solution to the 2D Riemann problem (\ref{Potential}), (\ref{initial})--(\ref{boundary}) for $\tau_1^e\leq \tau_0<\tau_2^i$ and ${\theta}_v<\theta_w<\pi+\sigma_{po}$. This completes the proof of Theorem \ref{thmfs+fs2}.

%%%%%%%%%%%%%%%%%%%%%%%%%%%%%%%%%%%%%%%%%%%%%%%%%%%%%%%%%%%%%%%%%%%%%%%%%%
%%%%%%%%%%%%%%%%%%%%%%%%%%%%%%%%%%%%%%%%%%%%%%%%%%%%%%%%%%%%%%%%%%%%%%%%%%%5
\section{Interaction between two fan–shock–fan composite waves}
In this section, we study the 2D Riemann problem (\ref{Potential}), (\ref{initial})--(\ref{boundary}) for $\tau_0<\tau_1^e$ and $u_0>c_0$.
Away from the corner, the oncoming flow and the vacuum
is connected by a planar fan–shock–fan composite wave, which can be represented by (\ref{pfsf}).
In this section, we only discuss the case $\theta_w<\pi+\alpha_v$.
Then, at the corner the oncoming flow is deflected by an oblique fan–shock–fan composite wave to a vacuum, which can be represented by (\ref{ofsf2}). 
In the subsequent discussion, the  constants $\xi_0$, $\xi_1$, $\xi_2$, $\hat{u}_1$, $\hat{u}_2$, $\sigma_2^e$, $\alpha_e$, $\alpha_w$, $\alpha_v$, as well as the functions $(\hat{u}_l, \hat{v}_l, \hat{\tau}_l)(\xi)$, $(\hat{u}_r, \hat{v}_r, \hat{\tau}_r)(\xi)$, $(\bar{u}_l, \bar{v}_l, \bar{\tau}_l)(\theta)$, $(\bar{u}_r, \bar{v}_r, \bar{\tau}_r)(\theta)$  are all as defined in Section 3.6.%
%see Figure \ref{Fig6.1}.
%The planar fan–shock–fan composite wave and the oblique fan–shock–fan composite wave will interact with each other.

%=============
\begin{figure}[htbp]
\begin{center}
\includegraphics[scale=0.55]{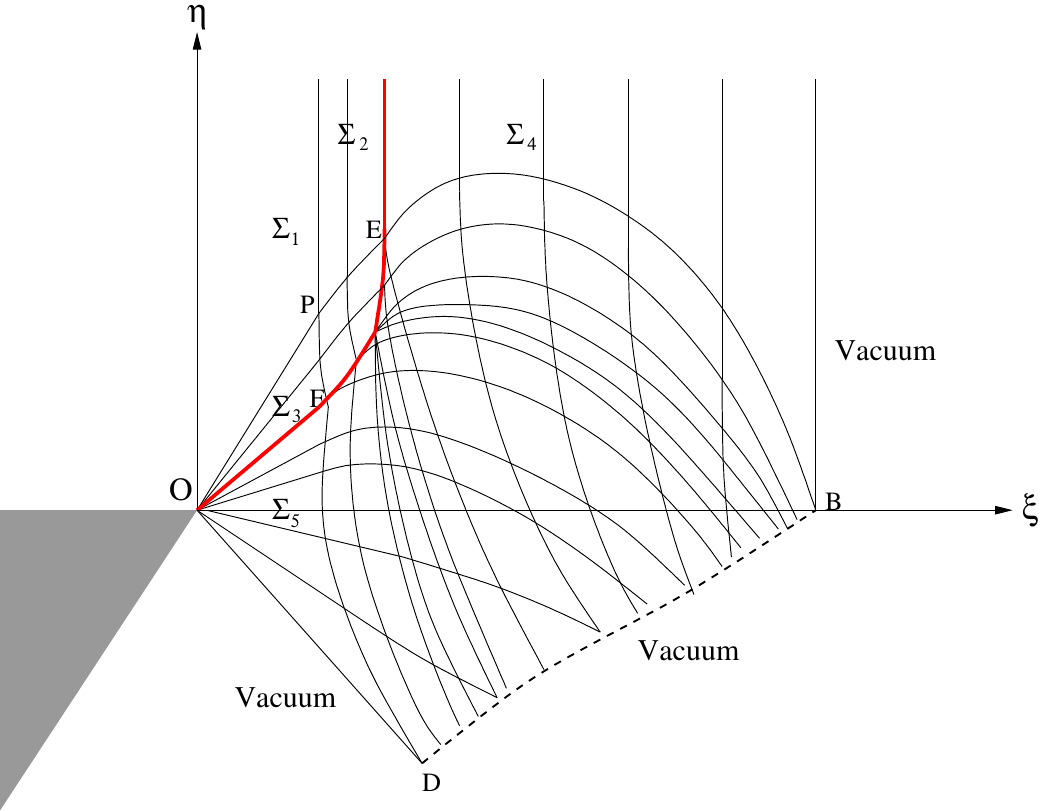}
\caption{\footnotesize Interaction between two fan–shock–fan composite waves.}
\label{Fig6.1}
\end{center}
\end{figure}
%==========

For convenience,  we use the subscripts $4$ and $5$ to denote the flow states of the downstream planar rarefaction wave and the downstream centered simple wave of the double-sonic shocks, respectively; that is,
$$
(u_4, v_4, \tau_4)(\xi, \eta)=(\hat{u}_r, \hat{v}_r, \hat{\tau}_r)(\xi)\quad \text{and} \quad (u_5, v_5, \tau_5)(\xi, \eta)=(\bar{u}_r, \bar{v}_r, \bar{\tau}_r)(\theta).
$$

%\subsection{Characteristic boundary value problem}

%In order to find the flow inside the domain $\Sigma$, we consider (\ref{PEC}) with the boundary condition
%\begin{equation}\label{63001}
%(u, v, \tau)=\left\{
 %              \begin{array}{ll}
  %                (\hat{u}_l, \hat{v}_l, \hat{\tau}_l)(\xi), & \hbox{$(\xi, \eta)\in \widetilde{\mathrm{PE}}$;} \\[2pt]
   %               (\hat{u}_r, \hat{v}_r, \hat{\tau}_r)(\xi), & \hbox{$(\xi, \eta)\in \widetilde{\mathrm{EB}}$;} \\[2pt]
    %           (\bar{u}_l, \bar{v}_l, \bar{\tau}_l)(\theta), & \hbox{$(\xi, \eta)\in \widetilde{\mathrm{PF}}$;} \\[2pt]
     %             (\bar{u}_r, \bar{v}_r, \bar{\tau}_r)(\theta), & \hbox{$(\xi, \eta)\in \widetilde{\mathrm{FD}}$.}
      %         \end{array}
       %      \right.
%\end{equation}
%where the functions $(\hat{u}_l, \hat{v}_l, \hat{\tau}_l)(\xi)$, $(\hat{u}_r, \hat{v}_r, \hat{\tau}_r)(\xi)$,
%$(\bar{u}_l, \bar{v}_l, \bar{\tau}_l)(\theta)$, and $(\bar{u}_r, \bar{v}_r, \bar{\tau}_r)(\theta)$ are defined in Section 1.4.4.
%In what follows, we are going to construct a piecewise smooth solution to the characteristic boundary value problem (\ref{PEC}), (\ref{63001}) piece by piece.

\subsection{Assumptions and main result}
Before stating the main result of this section, we introduce the following constants:
\begin{equation}\label{100902}
\begin{aligned}
&\xi_e=u_0-\frac{\rho_2^ec_2^e}{\rho_1^e}, \quad\phi_{d}=\arcsin\Big(\frac{\rho_2^ec_2^e}{\rho_1^eu_0}\Big),\quad \eta_e =\xi_e\tan\phi_d, \quad u_4^e=u_0-\frac{\rho_2^ec_2^e}{\rho_1^e}+c_2^e,\\&~~\alpha_{4}^e=-\frac{\pi}{2}+2\arctan \left(\frac{\tau_2^e(1+\sin\phi_d)}{\tau_1^e\cos\phi_d}\right),\quad
\sigma_{4}^e=-\frac{\pi}{2}+\arctan \left(\frac{\tau_2^e(1+\sin\phi_d)}{\tau_1^e\cos\phi_d}\right),\\&\quad
\beta_5^e=\phi_d-2\arctan \left(\frac{\tau_2^e(1+\sin\phi_d)}{\tau_1^e\cos\phi_d}\right), \quad
\sigma_5^e=\phi_d-\arctan \left(\frac{\tau_2^e(1+\sin\phi_d)}{\tau_1^e\cos\phi_d}\right),\\&\qquad
u_5^e=c_2^e\sin\phi_d+u_0\cos^2\phi_d, \quad\text{and}\quad v_5^e=u_0\sin\phi_d\cos\phi_d-c_2^e\cos\phi_d.
\end{aligned}
\end{equation}

\begin{rem}
The superscript $e$ represents the corresponding states for the scenario where the specific volume of the oncoming flow is exactly equal to $\tau_1^e$.
\end{rem}

As in (\ref{4905}), we have
$$
\beta_5^e<-\frac{\pi}{2}.
$$

%Set $\mu^1=\cot\beta_5^e$ and $\nu_2=\tan\alpha_4^e$
 We consider (\ref{61902}) with the initial data
\begin{equation}\label{92701a}
(\check{\alpha}_{+},\check{\beta}_{+})(\tau_{2}^e)=(\alpha_4^e, -\frac{\pi}{2}),
\end{equation}
and (\ref{120602}) with the initial data
\begin{equation}\label{92701b}
(\check{\alpha}_{-},\check{\beta}_{-})(\tau_{2}^e)=(\phi_{d}, \beta_5^e).
\end{equation}
As in Section 5.3, when $\sigma_4^e>\sigma_5^e$, there exists a $\tau_m^e>\tau_2^e$ such that both initial value problems admit a solution on $[\tau_2^e, \tau_m^e]$. Moreover,
the solutions satisfy
$$
(\check{\alpha}_{+},\check{\beta}_{+})(\tau_{m}^{e})=(\check{\alpha}_{-},\check{\beta}_{-})(\tau_{m}^e).
$$

We also consider (\ref{102304}) with the initial data
\begin{equation}\label{92703a}
(\bar{q}, \bar{\tau}, \bar{\sigma})(\phi_{d})=(q_{5}^e, \tau_{2}^e, \sigma_{5}^e),
\end{equation}
where $q_{5}^e=\sqrt{(u_{5}^e)^2+(v_{5}^e)^2}$.
There exists a ${\alpha}_v^{e}<\phi_d$ such that
the initial value problem (\ref{102304}), (\ref{92703a}) admits a solution $(\bar{q}, \bar{\tau}, \bar{\sigma})(\theta)$ on $({\alpha}_v^{e},\phi_d)$, and the solution satisfies $\bar{\tau}(\theta)\rightarrow 0$ as $\theta\rightarrow {\alpha}_v^{e}$. %Let $\bar{\alpha}(\theta)=\bar{\sigma}(\theta)$

To construct a global-in-time self-similar solution to the 2D Riemann problem (\ref{Potential}), (\ref{initial})--(\ref{boundary})  for  $\tau_0<\tau_1^e$ and $\theta_w<\pi+\theta_v$,
we impose the following assumptions:
\begin{description}
\item[($h_{31}$)] $\min\big\{\alpha_4^e+\frac{\pi}{2},  \phi_d-\beta_5^e\big\}>2\Theta(\tau_{2}^e)$, $\alpha_4^e-\beta_5^e<\pi$, and $\sigma_4^e>\sigma_5^e$;
  \item[($h_{32}$)] $\theta>\beta_5^e+2\Theta(\bar{\tau}(\theta))$ for  $\theta\in[\alpha_v^{e}, \phi_d]$;
\item[($h_{33}$)] $\check{\alpha}_{+}(\tau)>\beta_5^e+2\Theta(\tau)$ for  $\tau\in[\tau_{2}^e, \tau_m^{e}]$;
   \item[($h_{34}$)] $\check{\beta}_{-}(\tau)<\alpha_4^e-2\Theta(\tau)$ for  $\tau\in[\tau_{2}^e, \tau_m^{e}]$.
\end{description}
As in Section 5.1,  when $\phi_{d}$ is close to $\frac{\pi}{2}$, i.e., $u_0$ is close to $c_1^e$, assumptions ($h_{31}$), ($h_{33}$) and ($h_{34}$) hold.
The main result of this section can be stated as the following theorem.
\begin{thm}\label{fsf/fsf}
Suppose that  $\tau_0<\tau_1^e$, $\theta_w<\pi+\alpha_v$,
and assumptions ($h_{31}$)--($h_{34}$) hold.
Then, for sufficiently small $\tau_1^e-\tau_0>0$, the 2D Riemann problem (\ref{Potential}), (\ref{initial})--(\ref{boundary})  admits a global-in-time piecewise smooth self-similar solution.
\end{thm}

The two composite waves begin to interact with each other at the point $\mathrm{P}(\xi_{_\mathrm{P}}, \xi_{_\mathrm{P}}\tan\alpha_0)$, where $\xi_{_\mathrm{P}}=u_0-c_0$ and $\alpha_0=\arcsin(\frac{c_0}{u_0})$. Through $\mathrm{P}$, we draw a $C_{+}$ cross characteristic curve in the planar composite wave and a $C_{-}$ cross characteristic curve in the oblique composite wave.
The $C_{+}$ cross characteristic curve intersects the double-sonic shock of the planar composite wave at a point $\mathrm{E}(\xi_{_\mathrm{E}}, \eta_{_\mathrm{E}})$ and terminates at point $\mathrm{B}(\xi_2, 0)$.
The $C_{-}$ cross characteristic curve intersects the double-sonic shock of the oblique composite wave at a point $\mathrm{F}(\xi_{_\mathrm{F}}, \eta_{_\mathrm{F}})$ and terminates at a point $\mathrm{D}$ on the tail of the oblique  composite wave.
As in Section 4.1, these two cross characteristic curves can be determined.
To obtain a solution of the interaction between the fan–shock–fan composite waves, we consider (\ref{PEC}) with the boundary conditions
\begin{equation}\label{63001}
(u, v, \tau)=\left\{
               \begin{array}{ll}
                  (\hat{u}_l, \hat{v}_l, \hat{\tau}_l)(\xi), & \hbox{$(\xi, \eta)\in \widetilde{\mathrm{PE}}$;} \\[2pt]
                  (\hat{u}_r, \hat{v}_r, \hat{\tau}_r)(\xi), & \hbox{$(\xi, \eta)\in \widetilde{\mathrm{EB}}$;} \\[2pt]
               (\bar{u}_l, \bar{v}_l, \bar{\tau}_l)(\theta), & \hbox{$(\xi, \eta)\in \widetilde{\mathrm{PF}}$;} \\[2pt]
                  (\bar{u}_r, \bar{v}_r, \bar{\tau}_r)(\theta), & \hbox{$(\xi, \eta)\in \widetilde{\mathrm{FD}}$.}
              \end{array}
             \right.
\end{equation}
We shall show that under the assumptions in Theorem \ref{fsf/fsf}, the boundary value problem (\ref{PEC}), (\ref{63001}) admits a piecewise smooth solution
on the domain  bounded by $\widetilde{\mathrm{PE}}\cup\widetilde{\mathrm{EB}}$, $\widetilde{\mathrm{PF}}\cup\widetilde{\mathrm{FD}}$, and a vacuum boundary connecting point $\mathrm{B}$ and point $\mathrm{D}$. %Then the flow outside $\Sigma$ is known. So, in order to construct a self-similar solution of the 2D Riemann problem , (\ref{initial})--(\ref{boundary}) for $\tau_0<\tau_1^e$ and $\theta_w<\pi+\theta_v$, it only remains to construct the solution in $\Sigma$.

\subsection{Interaction of fan waves}
We first consider
 (\ref{PEC}) with the boundary conditions
\begin{equation}\label{63001A}
(u, v, \tau)=\left\{
               \begin{array}{ll}
                  (\hat{u}_l, \hat{v}_l, \hat{\tau}_l)(\xi), & \hbox{$(\xi, \eta)\in \widetilde{\mathrm{PE}}$;} \\[2pt]
               (\bar{u}_l, \bar{v}_l, \bar{\tau}_l)(\theta), & \hbox{$(\xi, \eta)\in \widetilde{\mathrm{PF}}$.}
               \end{array}
             \right.
\end{equation}
Problem (\ref{PEC}), (\ref{63001A}) is a standard Goursat problem.

\begin{figure}[htbp]
\begin{center}
\includegraphics[scale=0.47]{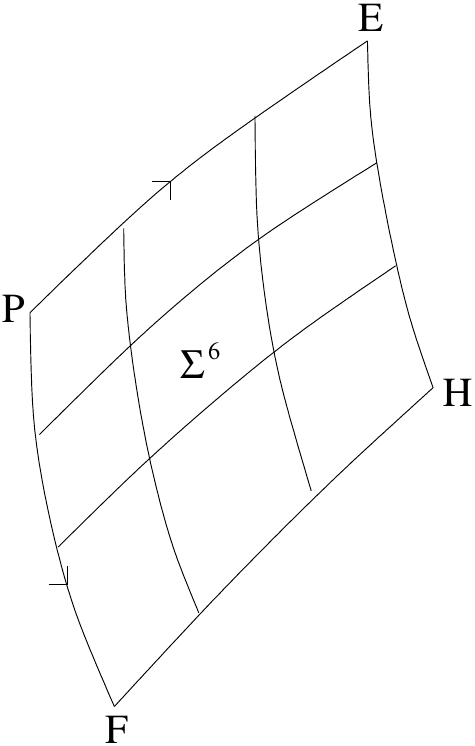}\qquad\qquad\includegraphics[scale=0.47]{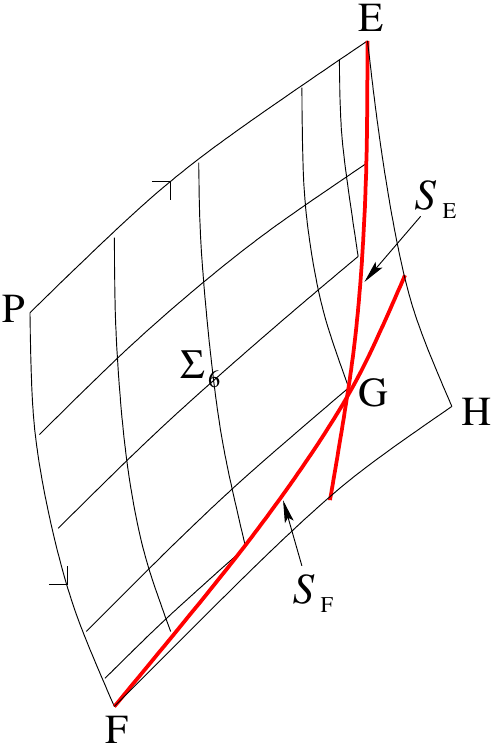}\qquad\qquad \includegraphics[scale=0.47]{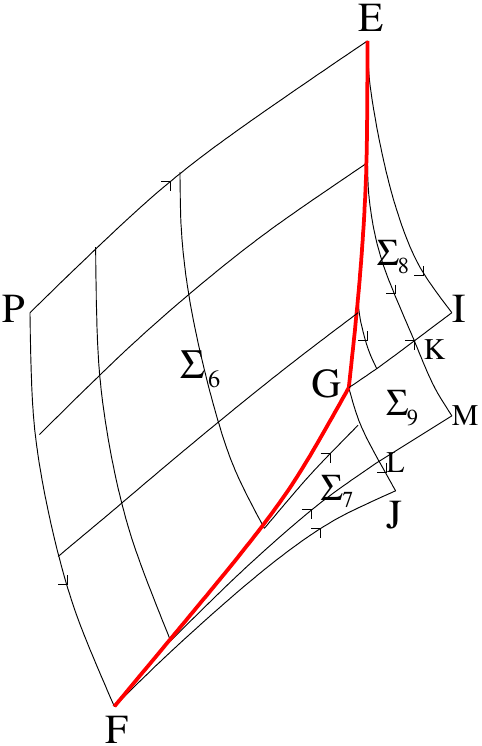}
\caption{\footnotesize Post-sonic shocks in the interaction region.}
\label{Fig6.2}
\end{center}
\end{figure}

\begin{lem}\label{lem7.2}
Suppose that  $\tau_1^e-\tau_0$ is sufficiently small. Then the Goursat problem (\ref{PEC}), (\ref{63001A}) admits a classical solution on a curvilinear quadrilateral domain $\overline{\Sigma^6}$ bounded by $\widetilde{\mathrm{PE}}$, $\widetilde{\mathrm{PF}}$, a forward $C_{+}$ characteristic curve emanating from $\mathrm{E}$, and  a forward $C_{-}$ characteristic curve emanating from $\mathrm{F}$; see Figure \ref{Fig6.2} (left). Moreover, the solution satisfies
\begin{equation}\label{6811}
\Big\|\big(\tau-\tau_1^e, ~\alpha-\phi_{d}, ~\beta+\frac{\pi}{2}, ~\bar{\partial}_{+}\rho+\mathcal{L}_1,~  \bar{\partial}_{-}\rho+\mathcal{L}_2\big)\Big\|_{{0; \Sigma^6}}\rightarrow 0\quad \text{as} \quad \tau_1^e-\tau_0\rightarrow 0,
\end{equation}
where $$\mathcal{L}_1=\frac{2c_1^e\sin (2A_1^e)}{(\tau_1^e)^4p''(\tau_1^e)}, \quad \mathcal{L}_2=
\frac{4(c_1^e)^2\cos^2 A_1^e\cos\phi_{d}}{\xi_e(\tau_1^e)^4p''(\tau_1^e)}, \quad \text{and}\quad A_1^e=\frac{\phi_{d}}{2}+\frac{\pi}{4}.$$
\end{lem}
\begin{proof}
When $\tau_1^e-\tau_0$ is sufficiently small, the existence of $C^1$ solution on $\overline{\Sigma^6}$ follows from Li and Yu \cite{Li-Yu} (Chap. 2). %It can be seen as a local solution of a ``bigger" Gorusat problem.
%We do not give the proof for the existence.
%The estimate (\ref{6811}) can be obtained directly by (\ref{62001})--(\ref{62502}) and the characteristic decompositions (\ref{81104}).
Hence, it only remains to prove the estimate (\ref{6811}).

As in the previous discussion, we have
\begin{equation}\label{93001}
\beta\equiv-\frac{\pi}{2}\quad \text{on}\quad \widetilde{\mathrm{PE}}\cup \widetilde{\mathrm{EB}}.
\end{equation}

As in (\ref{62201}), we have
\begin{equation}
\bar{\partial}_{+}\rho=-\frac{2c\rho^4\sin(2 A)}{p''(\tau)}\quad \text{along}\quad \widetilde{\mathrm{PE}}.
\end{equation}
As in (\ref{62303}), we have
\begin{equation}\label{92602a}
\bar{\partial}_-\rho
=-\frac{4c^2\cos^2 A}{\tau^4p''(\tau)\sqrt{\xi^2+\eta^2}}\quad \text{along}\quad \widetilde{\mathrm{PF}}.
\end{equation}

It is easy to verify that
\begin{equation}\label{92602c}
\alpha(\mathrm{P})\rightarrow \phi_{d}\quad \text{and}\quad   \sqrt{\xi_{_\mathrm{P}}^2+\eta_{_\mathrm{P}}^2}\rightarrow \xi_e \sec\phi_{d}\quad \text{as} \quad \tau_1^e-\tau_0\rightarrow 0.
\end{equation}

Combining (\ref{92602a})--(\ref{92602c}) and using the characteristic decompositions (\ref{81104}), we get (\ref{6811}).
This completes the proof.
\end{proof}

%=============
%==========

For convenience, we use $(u, v, \tau)=(u_6, v_6, \tau_6)(\xi, \eta)$ to denote the solution of the Goursat problem (\ref{PEC}), (\ref{63001A}) on the domain $\overline{\Sigma^6}$.

%\subsection{Post-sonic shocks in the interaction region}
\subsection{Transition from double-sonic to post-sonic}
We shall show that the two double-sonic shocks will propagate into the region $\Sigma^6$ and become post-sonic shocks. We use ${\it S}_{_\mathrm{E}}$ and ${\it S}_{_\mathrm{F}}$ to denote the shocks propagating from point $\mathrm{E}$ and point $\mathrm{F}$, respectively.

We first discuss the shock wave ${\it S}_{_\mathrm{F}}$.
Since ${\it S}_{_\mathrm{F}}$ is a 2-shock, we
set
\begin{equation}\label{62604}
L=U\cos \phi+V\sin \phi, \quad N=U \sin \phi- V\cos \phi,
\end{equation}
 where the variable $\phi$ denotes the inclination angle of ${\it S}_{_\mathrm{F}}$.
By the Rankine-Hugoniot relations we have that along ${\it S}_{_\mathrm{F}}$,
\begin{equation}\label{RH1}
\left\{
  \begin{array}{ll}
    \rho_f N_f=\rho_b N_b=m,  \\[4pt]
L_f=L_b,\\[4pt]
   N_f^2+2h(\tau_f)=N_b^2+2h(\tau_b),
  \end{array}
\right.
\end{equation}
where  subscripts $f$ and $b$  denote the states on the front and back sides of ${\it S}_{_\mathrm{F}}$, respectively.

%=============
%\begin{figure}[htbp]
%\begin{center}
%\includegraphics[scale=0.55]{Shock1.pdf}\qquad\qquad\qquad \qquad\qquad \qquad\includegraphics[scale=0.57]{Shock2.pdf}
%\caption{\footnotesize Post-sonic shocks in the interaction region.}
%\label{Fig6.9}
%\end{center}
%\end{figure}
%==========
%That is,
%We denote the values of $(u, v, \tau, c, \alpha, \beta, \sigma)$  the forward and the backside states of the post-sonic shock $\widetilde{BG}$, respectively.
%We denote by $\phi$ the inclination angle of the shock curve ${\it S}_{_\mathrm{F}}$. % ; see Figure \ref{Fig0} (right).
%Note that $\phi(B)=\pi-\theta$. We shall show that $\widetilde{BG}$ is convex and stays in the region $\Sigma^1$.

By the first and third relations of (\ref{RH1}), we have
\begin{equation}\label{250101}
m^2=-\frac{2h(\tau_f)-2h(\tau_b)}{\tau_f^2-\tau_b^2}.
\end{equation}

We express the normal velocity component $N$ as
$$
N=U\sin\phi-V\cos\phi=q\cos\sigma\sin\phi-q\sin\sigma\cos\phi=q\sin(\phi-\sigma).
$$
Hence, 
\begin{equation}\label{61112}
\phi=\sigma+\arcsin\Big(\frac{N}{q}\Big)\quad \text{on}\quad{\it S}_{_\mathrm{F}}.
\end{equation}

Since $\tau_f(\mathrm{F})=\tau_1^e$ and $\tau_b(\mathrm{F})=\tau_2^e$,
we have 
\begin{equation}\label{250401}
N_b=c_b=q_b\sin A_b\quad \text{and}\quad \alpha_b=\phi=\alpha_e\quad  \text{at}\quad \mathrm{F};
\end{equation}
\begin{equation}\label{102002}
N_f=c_f=q_f\sin A_f\quad \text{and}\quad \alpha_f=\phi=\alpha_e\quad \text{at}\quad \mathrm{F}.
\end{equation}
%where the meaning and the definition of the constant $\alpha_e$ are given in Section 1.4.4.

%We define the directional derivative:
%$$
%\bar{\partial}_s=\cos\phi\partial_{\xi}+\sin\phi\partial_{\eta}.
%$$
%A direct computation yields
%$$
%(\cos\phi, \sin\phi)=\frac{\sin (\beta_b-\phi)}{\sin (\beta_b-\alpha_b)}(\cos\alpha_b, \sin\alpha_b)-
%\frac{\sin (\alpha_b-\phi)}{\sin (\beta_b-\alpha_b)}(\cos\beta_b, \sin\beta_b).
%$$
%Thus, we have
%\begin{equation}\label{250109}
%\bar{\partial}_{s}=\frac{\sin (\phi-\beta_b)}{\sin (2A_b)}\bar{\partial}_{+}+\frac{\sin (\alpha_b-\phi)}{\sin %(2A_b)}\bar{\partial}_{-}.
%\end{equation}

From (\ref{250101}), we have
\begin{equation}\label{250104}
\begin{aligned}
2m \bar{\partial}_sm ~=~&\frac{2\tau_f}{\tau_f^2-\tau_b^2}\left(-p'(\tau_f)+\frac{2h(\tau_f)-2h(\tau_b)}{\tau_f^2-\tau_b^2}\right)
\bar{\partial}_s\tau_f\\&
-\frac{2\tau_b}{\tau_f^2-\tau_b^2}\left(-p'(\tau_b)+\frac{2h(\tau_f)-2h(\tau_b)}{\tau_f^2-\tau_b^2}\right)
\bar{\partial}_s\tau_b\quad \text{along}\quad {\it S}_{_\mathrm{F}}.
\end{aligned}
\end{equation}
Combining this with $\tau_f(\mathrm{F})=\tau_1^e$ and $\tau_b(\mathrm{F})=\tau_2^e$,  we have
\begin{equation}\label{250402}
\bar{\partial}_sm=0\quad \text{at}\quad \mathrm{F}.
\end{equation}

From (\ref{U}) and (\ref{62604}) we have
$$
m=\rho  N =\rho  q  (\cos\sigma \sin\phi-\sin\sigma \cos\phi).
$$
Thus,
\begin{equation}\label{4250103}
\bar{\partial}_s m=\frac{N }{q }\partial_{s}(\rho q )-\rho L \bar{\partial}_s\sigma +\rho L \bar{\partial}_s \phi\quad \text{along}\quad {\it S}_{_\mathrm{F}}.
\end{equation}

From Bernoulli's law (\ref{5802}) we have
$$
q \bar{\partial}_s q +\tau  p'(\tau )\bar{\partial}_s \tau +L  =0\quad \text{along}\quad {\it S}_{_\mathrm{F}}.
$$
This yields
$$
\bar{\partial}_s (\rho q )=\frac{1}{q }\big(q ^2+\tau ^2 p'(\tau )\big)\bar{\partial}_s\rho -\frac{\rho L }{q }\quad \text{along}\quad {\it S}_{_\mathrm{F}}.
$$
Hence, we get
\begin{equation}\label{250102}
\bar{\partial}_s \phi=\frac{\bar{\partial}_s m}{\rho L }-\frac{N }{\rho L q ^2}\big(q ^2+\tau ^2 p'(\tau )\big)\bar{\partial}_s\rho +\frac{N }{q ^2}+\bar{\partial}_s\sigma \quad \text{along}\quad {\it S}_{_\mathrm{F}}.
\end{equation}
%From the Bernoulli's law (\ref{5802}), we also have
%\begin{equation}\label{42803a}
%\bar{\partial}_s q =\frac{\tau ^3 p'(\tau )\bar{\partial}_s \rho }{q }-\frac{L }{q }=-\tau q \sin^2 A \bar{\partial}_s \rho -\frac{L }{q }\quad \text{along}\quad {\it S}_{_\mathrm{F}}.
%\end{equation}
From (\ref{250102}) we have
\begin{equation}\label{42804a}
\bar{\partial}_s \sigma =\bar{\partial}_s \phi-\frac{\bar{\partial}_s m}{\rho L }+\frac{N }{\rho L q ^2}\big(q ^2+\tau ^2 p'(\tau )\big)\bar{\partial}_s\rho -\frac{N }{q ^2}\quad \text{along}\quad {\it S}_{_\mathrm{F}}.
\end{equation}

A direct computation yields
\begin{equation}\label{50701}
\bar{\partial}_su=\bar{\partial}_s(U+\xi)=\bar{\partial}_s(q\cos\sigma)+\cos\phi=\cos\sigma\bar{\partial}_sq-q\sin\sigma\bar{\partial}_s\sigma+\cos\phi
\end{equation}
and
\begin{equation}\label{50702}
\bar{\partial}_sv=\bar{\partial}_s(V+\eta)=\bar{\partial}_s(q\sin\sigma)+\sin\phi=\sin\sigma\bar{\partial}_sq
+q\cos\sigma\bar{\partial}_s\sigma+\sin\phi.
\end{equation}

%We now assume that shocks emanating from $\mathrm{D}$ and $\mathrm{D}$ are post-sonic.
%We shall show that the post-sonic shocks emanating from $\mathrm{D}$ and $\mathrm{D}$ will stay in $\Sigma^1$ until they intersect with each other at some point $\mathrm{G}$ on the $\xi$-axis; see Figure \ref{Fig13}.
 %We denote by $\widetilde{\mathrm{DG}}$ and $\widetilde{\mathrm{DG}}$ the post-sonic shocks emanating from $\mathrm{D}$ and $\mathrm{D}$, respectively.
%In this section we discuss the existence of the post-sonic shock curves $\widetilde{\mathrm{DG}}$ and $\widetilde{\mathrm{DG}}$. In the next section we shall show that

%We shall construct the post-sonic shock from $B$ in detail.
%$\widetilde{BG}$ in $\Sigma^1$.
%The post-sonic shock from $D$ can be obtained by the symmetry.

We assume a priori that the shock ${\it S}_{_\mathrm{F}}$ is post-sonic and lies in the interior of $\Sigma^6$.
Consequently, the states on the front side of the shock are given by
$$
(u_f, v_f, \tau_f)=(u_6, v_6, \tau_6)(\xi, \eta)\quad \text{on}\quad  {\it S}_{_\mathrm{F}}.
$$
Furthermore, we assume that
\begin{equation}\label{72202}
\tau_1^e<\tau_f<\tau_2^i\quad \text{on}\quad  {\it S}_{_\mathrm{F}}\setminus \mathrm{F}.
\end{equation}
Then by Proposition \ref{51003}, we have
\begin{equation}\label{62603}
\tau_b=\psi_{po}(\tau_f)\quad \text{and}\quad
m^2=-p'(\psi_{po}(\tau_f))>-p'(\tau_f)\quad \text{on}\quad  {\it S}_{_\mathrm{F}}\setminus \mathrm{F}.
\end{equation}
In addition, 
\begin{equation}\label{6810}
N_b=c_b=\sqrt{-\tau_b^2p'(\tau_b)}\quad \text{and}\quad N_f>c_f=\sqrt{-\tau_f^2p'(\tau_f)}\quad \text{along}\quad {\it S}_{_\mathrm{F}}\setminus \mathrm{F}.
\end{equation}

%We still denote by $\phi$ the inclination angle of the shock curve $ {\it S}_{_\mathrm{F}}$.
By (\ref{6810}) and (\ref{61112}), we have
\begin{equation}\label{61113}
\phi=\sigma_b+A_b=\alpha_b\quad\text{and}\quad \phi>\sigma_f+A_f=\alpha_f  \quad \text{on}\quad  {\it S}_{_\mathrm{F}}\setminus \mathrm{F}.
\end{equation}
%It is easily to see
%\begin{equation}\label{61114}
%\phi(B)=\pi+\theta.
%\end{equation}

%Since $\tau_f(\mathrm{F})=\tau_f^e$, we have $N_f(\mathrm{F})=c_f(\mathrm{F})$, and consequently $\phi=\alpha_f$ at $\mathrm{F}$. This implies
%\begin{equation}\label{61114}
%\phi(\mathrm{F})=\alpha_f(\mathrm{F})=\pi+\theta.
%\end{equation}

We differentiate $m^2=-p'(\psi_{po}(\tau_f))$ along the shock $ {\it S}_{_\mathrm{F}}$ to obtain
\begin{equation}\label{61105}
\bar{\partial}_s m=\frac{-p''(\psi_{po}(\tau_f))\psi_{po}'(\tau_f)\bar{\partial}_s \tau_f}{2m }\quad \text{along}\quad {\it S}_{_\mathrm{F}}.
\end{equation}

%As in (\ref{4250103}), we have
%\begin{equation}\label{42802}
%\bar{\partial}_s m=\frac{N_f}{q_f}\partial_{s}(\rho_fq_f)-\rho_fL_f\bar{\partial}_s\sigma_f+\rho_fL_f\bar{\partial}_s \phi\quad \text{on}\quad  {\it S}_{_\mathrm{F}}.
%\end{equation}

As in (\ref{250102}), we have
\begin{equation}\label{61115}
\bar{\partial}_s \phi=\frac{\bar{\partial}_s m}{\rho_fL_f}-\frac{N_f}{\rho_fL_fq_f^2}\big(q_f^2+\tau_f^2 p'(\tau_f)\big)\bar{\partial}_s\rho_f+\frac{N_f}{q_f^2}+\bar{\partial}_s\sigma_f\quad \text{along}\quad {\it S}_{_\mathrm{F}}.
\end{equation}
%This immediately implies that the curved shock wave $\widetilde{BG}$ is convex.

We now establish the existence of the post-sonic shock curve $S_{_\mathrm{F}}$.
Let $\eta=\eta_{+}(\xi)$ for $\xi_{_\mathrm{F}}\leq \xi\leq \xi_{_\mathrm{H}}$ denote the $C_{+}$ characteristic curve $\widetilde{\mathrm{FH}}$, and let $\eta=\eta_{s}(\xi)$ for $\xi\geq \xi_{\mathrm{F}}$ denote the post-sonic shock curve $S_{\mathrm{F}}$.
These curves satisfy the following differential equations and initial conditions:
\begin{equation}\label{61102}
\begin{aligned}
&\eta_{s}'(\xi)=\tan\big(\phi(\xi, \eta_s(\xi))\big), \quad \eta_{+}'(\xi)=\tan\big(\alpha_f(\xi, \eta_+(\xi))\big), \\&
\eta_{s}(\xi_{_\mathrm{F}})=\eta_{+}(\xi_{_\mathrm{F}})=\eta_{_\mathrm{F}}, \quad \eta_{s}'(\xi_{_\mathrm{F}})=\eta_{+}'(\xi_{_\mathrm{F}})=\tan\alpha_e.
\end{aligned}
\end{equation}

To obtain $\xi_s(\eta)$, we only need to solve (\ref{61115}) with the initial data
\begin{equation}\label{61114}
\phi(\mathrm{F})=\alpha_f(\mathrm{F})=\alpha_e.
\end{equation}
The initial value problem (\ref{61115}), (\ref{61114}) relies on the solution in $\Sigma^6$. Therefore, to establish the existence of the post-sonic shock curve $ {\it S}_{_\mathrm{F}}$, we also need to prove
\begin{equation}\label{61101}
\eta_{s}(\xi)>\eta_{+}(\xi), \quad \xi>\xi_{_\mathrm{F}}.
\end{equation}

From $\tan(\alpha_f(\eta_{+}(\xi),\xi))=\eta_{+}'(\xi)$, we have
\begin{equation}\label{61103}
\eta_{+}''(\xi)=\frac{\bar{\partial}_{+}\alpha_6(\xi, \eta_{+}(\xi))}{\cos^3(\alpha_6(\xi, \xi_{+}(\eta)))}.
\end{equation}
%where $\bar{\partial}_{+}\alpha_f(\xi, \eta_{+}(\xi))=\cos(\xi, \eta_{+}(\xi))\partial_{\xi}\alpha_6(\xi, %\eta_{+}(\xi))+
%\sin(\xi, \eta_{+}(\xi))\partial_{\eta}\alpha_6(\xi, \eta_{+}(\xi))$.
From $\tan(\phi(\xi, \eta_{s}(\xi)))=\eta_{s}'(\xi)$, we have
\begin{equation}\label{61104}
\eta_{s}''(\xi)=\frac{\bar{\partial}_{s}\phi(\xi, \eta_{s}(\xi))}{\cos^3(\phi(\xi, \eta_{s}(\xi)))}.
\end{equation}

%In what follows, we are going to $\bar{\partial}_{+}\alpha_f$ and $\bar{\partial}_{s}\phi$ at the point $B$.

%Using (\ref{61105}) and recalling $\psi_{po}'(\tau_f^e)=0$, we have
%\begin{equation}\label{82001}
%\bar{\partial}_s m=0\quad \text{at}\quad \mathrm{D}.
%\end{equation}
From (\ref{61115}), (\ref{7}), (\ref{10}), (\ref{61114}), (\ref{102002}), and (\ref{250402}), we have that at the point $\mathrm{F}$,
\begin{equation}\label{61108}
\begin{aligned}
\bar{\partial}_s \phi&=\frac{\bar{\partial}_s m}{\rho_fL_f}-\frac{N_f}{\rho_fL_fq_f^2}(q_f^2+\tau_f^2 p'(\tau_f))\bar{\partial}_s\rho_f+\frac{N_f}{q_f^2}+\bar{\partial}_s\sigma_f\\&=
-\frac{N_f}{\rho_fL_fq_f^2}(q_f^2+\tau_f^2 p'(\tau_f))\bar{\partial}_{+}\rho_6+\frac{N_f}{q_f^2}+\bar{\partial}_{+}\Big(\frac{\alpha_6+\beta_6}{2}\Big)
\\&=-\frac{\sin 2A_f}{2\rho_f}\bar{\partial}_{+}\rho_6+\frac{N_f}{q_f^2}+\frac{\bar{\partial}_{+}\alpha_6}{2}-\frac{p''(\tau_f)\tan A_f}{4c_f^2\rho_f^4}\bar{\partial}_{+}\rho_6-\frac{\sin^2 A_f}{c_f}
\\&=-\frac{\sin 2A_f}{2\rho_f}\bar{\partial}_{+}\rho_6+\frac{\bar{\partial}_{+}\alpha_6}{2}-\frac{p''(\tau_f)\tan A_f}{4c_f^2\rho_f^4}\bar{\partial}_{+}\rho_6\\&>-\frac{\sin 2A_f}{2\rho_f}\bar{\partial}_{+}\rho_6+\frac{\bar{\partial}_{+}\alpha_6}{2}
\end{aligned}
\end{equation}
and
$$
\begin{aligned}
\bar{\partial}_+ \alpha_6&=-\frac{p''(\tau_f)}{4c_f^2\rho_f^4}\left(-\frac{4p'(\tau_f)+\tau_f p''(\tau_f)}{\tau_f p''(\tau_f)}-\tan^2A_f\right)\sin 2A_f\bar{\partial}_{+}\rho_6\\&<
-\frac{p''(\tau_f)}{4c_f^2\rho_f^4}\left(-\frac{4p'(\tau_f)}{\tau_f p''(\tau_f)}\right)\sin 2A_f\bar{\partial}_{+}\rho_6\\&=-\frac{\sin 2A_f}{\rho_f}\bar{\partial}_{+}\rho_6.
\end{aligned}
$$
This implies
\begin{equation}\label{61106}
\bar{\partial}_s \phi(\xi_{_\mathrm{F}}, \eta_{_\mathrm{F}})>\bar{\partial}_+ \alpha_6(\xi_{_\mathrm{F}}, \eta_{_\mathrm{F}}).
\end{equation}

Combining (\ref{61102}) and (\ref{61106}), there exists a sufficiently small $\varepsilon>0$ such that $\eta_{s}(\xi)>\eta_{+}(\xi)$ for $\xi_{_\mathrm{F}}<\xi<\xi_{_\mathrm{F}}+\varepsilon$.
Furthermore,  inequality (\ref{61113}) implies that
\begin{equation}\label{62603a}
\phi(\xi, \eta_{s}(\xi))>\alpha_f(\xi, \eta_{s}(\xi))
\end{equation}
provided that $(\xi, \eta_{s}(\xi))$ lies in the interior region of $\Sigma^6$.
Thus, a continuity argument yields (\ref{61101}).

\vskip 8pt

We now establish the existence of the post-sonic shock curve ${\it S}_{_\mathrm{E}}$.
We still use subscripts $f$ and $b$ to denote the states on the front and back sides of ${\it S}_{_\mathrm{E}}$,  and use $\phi$ to denote the inclination angle ${\it S}_{_\mathrm{E}}$.
Since ${\it S}_{_\mathrm{E}}$ is a 1-shock,
we set
\begin{equation}\label{100701a}
L=U\cos \phi+V\sin \phi, \quad N=-U \sin \phi+V\cos \phi.
\end{equation}
%Then we still have the shock relations (\ref{RH1}) on ${\it S}_{_\mathrm{E}}$.

As in (\ref{250102}), we have
\begin{equation}\label{61115a}
\bar{\partial}_s \phi=-\frac{\bar{\partial}_s m}{\rho_fL_f}+\frac{N_f}{\rho_fL_fq_f^2}(q_f^2+\tau_f^2 p'(\tau_f))\bar{\partial}_s\rho_f-\frac{N_f}{q_f^2}+\bar{\partial}_s\sigma_f\quad \text{along}\quad {\it S}_{_\mathrm{E}}.
\end{equation}
We denote by $\xi=\xi_{-}(\eta)$, $\eta\leq \eta_{_\mathrm{E}}$ the $C_{-}$ characteristic curve $\widetilde{\mathrm{EH}}$ and by
$\xi=\xi_{s}(\eta)$, $\eta\leq \eta_{_\mathrm{E}}$ the post-sonic shock curve $ {\it S}_{_\mathrm{E}}$.
We have
\begin{equation}\label{61102a}
\begin{aligned}
&\xi_{s}'(\eta)=\cot(\phi(\xi_s(\eta), \eta)), \quad \xi_{-}'(\eta)=\cot(\beta_f(\xi_{-}(\eta), \eta)), \quad\\&
\qquad \xi_{s}(\eta_{_\mathrm{E}})=\xi_{-}(\eta_{_\mathrm{E}}), \quad  \xi_{s}'(\eta_{_\mathrm{E}})=\xi_{-}'(\eta_{_\mathrm{E}})=0.
\end{aligned}
\end{equation}

To obtain $\xi_s(\eta)$, we consider (\ref{61115a}) with the initial data
\begin{equation}\label{61114a}
\phi(\mathrm{E})=\beta_f(\mathrm{E})=-\frac{\pi}{2}.
\end{equation}

%For the post-sonic shock $\widetilde{BG}$, the front side state of it is $(u_f, v_f, \tau_f)$ and the backside state of it can be determined by the front side state.
%Actually, in order to obtain $\xi_s(\eta)$, we only need to consider (\ref{61115}) with data
%\begin{equation}\label{61114}
%\phi(\mathrm{F})=\alpha_f(\mathrm{F})=\alpha_e.
%\end{equation}
%Since the the initial value problem (\ref{61115}, \ref{61114}) use the data in $\Sigma^1$, in order to prove the existence of $ {\it S}_{_\mathrm{F}}$ we need to prove
%\begin{equation}\label{61101}
%\eta_{s}(\xi)>\eta_{+}(\xi), \quad \xi_{_\mathrm{F}}<\xi<\xi_{_\mathrm{G}}.
%\end{equation}

From $\cot(\beta_6(\xi_{-}(\eta),\eta)))=\xi_{-}'(\eta)$, we have
\begin{equation}\label{61103a}
\xi_{-}''(\eta)=-\frac{\bar{\partial}_{-}\beta_6(\xi_{-}(\eta),\eta))}{\cos^3(\beta_6(\xi_{-}(\eta),\eta)))}.
\end{equation}
%where $\bar{\partial}_{+}\alpha_f(\xi, \eta_{+}(\xi))=\cos(\xi, \eta_{+}(\xi))\partial_{\xi}\alpha_6(\xi, \eta_{+}(\xi))+
%\sin(\xi, \eta_{+}(\xi))\partial_{\eta}\alpha_6(\xi, \eta_{+}(\xi))$.
From $\cot(\phi(\xi_{s}(\eta),\eta)))=\xi_{s}'(\eta)$, we have
\begin{equation}\label{61104a}
\xi_{s}''(\eta)=-\frac{\bar{\partial}_{s}\phi(\xi_{s}(\eta),\eta))}{\sin^3(\phi(\xi_{s}(\eta),\eta)))}.
\end{equation}

Since $\tau_f(\mathrm{E})=\tau_1^e$ and $\tau_b(\mathrm{E})=\tau_2^e$,
we have
\begin{equation}
N_b=c_b=q_b\sin A_b\quad \text{and}\quad \beta_b=\phi=-\frac{\pi}{2}\quad  \text{at}\quad \mathrm{E};
\end{equation}

We  assume a priori that
\begin{equation}\label{72202aa}
\tau_1^e<\tau_f<\tau_2^i\quad \text{on}\quad  {\it S}_{_\mathrm{E}}\setminus \mathrm{E}.
\end{equation}
As in (\ref{61108}), we have that at the point $\mathrm{E}$,
\begin{equation}\label{61108a}
\begin{aligned}
\bar{\partial}_s \phi&=-\frac{\bar{\partial}_s m}{\rho_fL_f}+\frac{N_f}{\rho_fL_fq_f^2}(q_f^2+\tau_f^2 p'(\tau_f))\bar{\partial}_s\rho_f-\frac{N_f}{q_f^2}+\bar{\partial}_s\sigma_f\\&=
\frac{N_f}{\rho_fL_fq_f^2}(q_f^2+\tau_f^2 p'(\tau_f))\bar{\partial}_{-}\rho_6+\frac{N_f}{q_f^2}+\bar{\partial}_{-}\Big(\frac{\alpha_6+\beta_6}{2}\Big)
\\&=\frac{\sin (2A_f)}{2\rho_f}\bar{\partial}_{-}\rho_6-\frac{N_f}{q_f^2}+\frac{\bar{\partial}_{-}\beta_6}{2}+\frac{p''(\tau_f)\tan A_f}{4c_f^2\rho_f^4}\bar{\partial}_{-}\rho_6+\frac{\sin^2 A_f}{c_f}
\\&=\frac{\sin (2A_f)}{2\rho_f}\bar{\partial}_{-}\rho_6+\frac{\bar{\partial}_{-}\beta_6}{2}+\frac{p''(\tau_f)\tan A_f}{4c_f^2\rho_f^4}\bar{\partial}_{-}\rho_6\\&<\frac{\sin (2A_f)}{2\rho_f}\bar{\partial}_{-}\rho_6+\frac{\bar{\partial}_{-}\beta_6}{2}
\end{aligned}
\end{equation}
and
$$
\begin{aligned}
\bar{\partial}_{-} \beta_6&=\frac{p''(\tau_f)}{4c_f^2\rho_f^4}\left(-\frac{4p'(\tau_f)+\tau_f p''(\tau_f)}{\tau_f p''(\tau_f)}-\tan^2A_f\right)\sin 2A_f\bar{\partial}_{-}\rho_6\\&>
\frac{p''(\tau_f)}{4c_f^2\rho_f^4}\left(-\frac{4p'(\tau_f)}{\tau_f p''(\tau_f)}\right)\sin (2A_f)\bar{\partial}_{+}\rho_6\\&=\frac{\sin (2A_f)}{\rho_f}\bar{\partial}_{-}\rho_6.
\end{aligned}
$$
This implies
\begin{equation}\label{61106a}
\bar{\partial}_s \phi<\bar{\partial}_- \beta_6\quad \text{at}\quad \mathrm{E}.
\end{equation}
Then, for sufficiently small $\varepsilon>0$,  $\xi_{s}(\eta)<\xi_{-}(\eta)$ for $\eta_{_\mathrm{E}}-\varepsilon<\eta<\eta_{_\mathrm{E}}$.
Furthermore, as in (\ref{62603a}) we have
$$
\phi(\xi_{s}(\eta),\eta))<\beta_f(\xi_{s}(\eta),\eta))
$$
provided that $(\xi_{s}(\eta),\eta))$ lies in the interior region of $\Sigma^6$.

A direct computation yields
$$
(\cos\phi, \sin\phi)=\frac{\sin (\beta_f-\phi)}{\sin (\beta_f-\alpha_f)}(\cos\alpha_f, \sin\alpha_f)-
\frac{\sin (\alpha_f-\phi)}{\sin (\beta_f-\alpha_f)}(\cos\beta_f, \sin\beta_f).
$$
Thus, we have
\begin{equation}\label{82003}
\bar{\partial}_{s}\tau_f=\frac{\sin (\beta_6-\phi)}{\sin (\beta_6-\alpha_6)}\bar{\partial}_{+}\tau_6-\frac{\sin (\alpha_6-\phi)}{\sin (\beta_6-\alpha_6)}\bar{\partial}_{-}\tau_6.
\end{equation}

From (\ref{61115}), (\ref{61114}), (\ref{61115a}), and (\ref{61114a}) we have
\begin{equation}\label{82004}
\|\phi+\frac{\pi}{2}\|_{0;  {\it S}_{_\mathrm{E}}}\rightarrow 0\quad \text{and}\quad
\|\phi-\phi_{d}\|_{0;  {\it S}_{_\mathrm{F}}}\rightarrow 0 \quad \text{as}\quad \tau_f^e-\tau_0\rightarrow 0.
\end{equation}
Therefore, by (\ref{6811}), (\ref{82003}) and (\ref{82004}), we deduce that
\begin{equation}\label{82008}
\bar{\partial}_{s}\tau_f>0\quad \text{along}\quad  {\it S}_{_\mathrm{F}}\cup {\it S}_{_\mathrm{E}},
\end{equation}
provided that $\tau_1^e-\tau_0$ is sufficiently small.
Combining this with $\tau_f(\mathrm{E})=\tau_f(\mathrm{F})=\tau_1^e$ and (\ref{6811}), we have that 
(\ref{72202}) and (\ref{72202aa}) hold for sufficiently small $\tau_1^e-\tau_0$.
This completes the proof for the existence of the post-sonic shocks $ {\it S}_{_\mathrm{E}}$ and $ {\it S}_{_\mathrm{F}}$.

The states on the backside of $ {\it S}_{_\mathrm{E}}$  and $ {\it S}_{_\mathrm{F}}$, denoted by $(u_b, v_b, \tau_b)(\xi, \eta)$, are then determined by (\ref{RH1}), (\ref{62603}), and (\ref{6810}).
Moreover, by (\ref{82004}) and (\ref{82008}), 
the two post-sonic shocks intersect at a point $\mathrm{G}(\xi_{_{\mathrm{G}}}, \eta_{_\mathrm{G}})$ within $\Sigma^6$.
In the subsequent discussion, we denote by $\widehat{\mathrm{FG}}$ and  $\widehat{\mathrm{EG}}$ the post-sonic shock segments  propagating from $\mathrm{E}$ and $\mathrm{F}$, respectively.

%Following the idea of Li and Zheng \cite{Li3},
%we shall use the method of hodograph transformation to construct a solution near the backside of the post-sonic shock curve $ {\it S}_{_\mathrm{F}}$. In order to establish the global one-to-one inversion of the hodograph transformation, we need to establish some important monotonicity conditions about the inversion. The following estimates  are crucial to deduce the monotonicity conditions.
\subsection{Structural conditions on the post-sonic shocks}
We shall establish some structural conditions on the post-sonic shocks $\widehat{\mathrm{FG}}$ and  $\widehat{\mathrm{EG}}$, which are crucial for ensuring the existence of the downstream flow.

\begin{lem}\label{61302}
Assume that $\tau_1^e-\tau_0$ is sufficiently small. Then the following estimates hold:
\begin{equation}\label{61111}
\cos\alpha_b \bar{\partial}_s u_b+\sin\alpha_b \bar{\partial}_s v_b>0\quad \text{and}\quad \cos\beta_b \bar{\partial}_s u_b+\sin\beta_b \bar{\partial}_s v_b<0 \quad \text{along}\quad  \widehat{\mathrm{FG}};
\end{equation}
\begin{equation}\label{61305}
\cos\alpha_b \bar{\partial}_s u_b+\sin\alpha_b \bar{\partial}_s v_b<0\quad \text{and}\quad \cos\beta_b \bar{\partial}_s u_b+\sin\beta_b \bar{\partial}_s v_b>0\quad \text{along}\quad  \widehat{\mathrm{EG}}.
\end{equation}
\end{lem}
\begin{proof}
%We only prove (\ref{61111}). The proof for (\ref{61305}) will be similar.
%The proof proceeds in two steps

Using (\ref{7}), (\ref{10}), (\ref{6811}), (\ref{102002}), (\ref{250402}), and (\ref{61115}), we obtain
\begin{equation}\label{250407}
\begin{aligned}
\bar{\partial}_s \phi&=\frac{\bar{\partial}_s m}{\rho_fL_f}-\frac{N_f}{\rho_fL_fq_f^2}(q_f^2+\tau_f^2 p'(\tau_f))\bar{\partial}_s\rho_f+\frac{N_f}{q_f^2}+\bar{\partial}_s\sigma_f\\&=
-\frac{N_f}{\rho_fL_fq_f^2}(q_f^2-c_f^2)\bar{\partial}_{+}\rho_6+\frac{N_f}{q_f^2}
+\bar{\partial}_{+}\Big(\frac{\alpha_6+\beta_6}{2}\Big)
\\&=-\frac{\sin 2A_6}{2\rho_6}\bar{\partial}_{+}\rho_6-\frac{p''(\tau_6)}{4c_6^2\rho_6^4}
\Big(\underbrace{\big(\omega(\tau_6)-\tan^2A_6\big)\sin A_6\cos A_6+\tan A_6}_{>0}\Big)\bar{\partial}_{+}\rho_6\\&>-\frac{\sin 2A_6}{2\rho_6}\bar{\partial}_{+}\rho_6>0\quad\text{at}\quad \mathrm{F}.
\end{aligned}
\end{equation}
Using (\ref{6811}), (\ref{250407}), and (\ref{250402}) we obtain that when $\tau_1^e-\tau_0$ is sufficiently small,
\begin{equation}\label{61109}
\bar{\partial}_s \phi>\frac{\tau_1^e\sin(2A_1^e)}{4}\mathcal{L}_1\quad \text{along}\quad  \widehat{\mathrm{FG}}.
\end{equation}
Similarly, when $\tau_1^e-\tau_0$ is sufficiently small,
\begin{equation}\label{61109a}
\bar{\partial}_s \phi<-\frac{\tau_1^e\sin(2A_1^e)}{4}\mathcal{L}_2\quad \text{along}\quad  \widehat{\mathrm{EG}}.
\end{equation}
%This immediately implies that when $\tau_1^e-\tau_0$ is sufficiently small $ \widehat{\mathrm{FG}}$ is convex.

Let $U_b=u_b-\xi$ and $V_b=v_b-\eta$. Then by (\ref{61113}) we have
\begin{equation}\label{42801}
\begin{aligned}
&\cos\alpha_b \bar{\partial}_s u_b+\sin\alpha_b \bar{\partial}_s v_b\\~=~&
\cos\alpha_b \bar{\partial}_s (U_b+\xi)+\sin\alpha_b \bar{\partial}_s (V_b+\eta)\\~=~&
\cos\alpha_b \bar{\partial}_s (q_b\cos\sigma_b)+\sin\alpha_b \bar{\partial}_s (q_b\sin\sigma_b)+\cos\alpha_b\cos\phi+\sin\alpha_b\sin\phi\\~=~&
\cos A_b \bar{\partial}_s q_b+q_b\sin A_b \bar{\partial}_s \sigma_b+1\quad \text{along}\quad \widehat{\mathrm{FG}}.
\end{aligned}
\end{equation}

Using Bernoulli's law (\ref{5802}) and recalling the first equality in (\ref{6810}), we have
\begin{equation}\label{42803}
\bar{\partial}_s q_b=\frac{\tau_b^3 p'(\tau_b)\bar{\partial}_s \rho_b}{q_b}-\frac{L_b}{q_b}=-\tau_bq_b\sin^2 A_b\bar{\partial}_s \rho_b-\cos A_b\quad \text{along}\quad  \widehat{\mathrm{FG}}.
\end{equation}

From (\ref{42804a}), we have
\begin{equation}\label{42804}
\bar{\partial}_s \sigma_b=\bar{\partial}_s \phi-\frac{\bar{\partial}_s m}{\rho_bL_b}+\frac{N_b}{\rho_bL_bq_b^2}(q_b^2+\tau_b^2 p'(\tau_b))\bar{\partial}_s\rho_b-\frac{N_b}{q_b^2}\quad \text{along}\quad  \widehat{\mathrm{FG}}.
\end{equation}

Inserting (\ref{42803}) and  (\ref{42804}) into  (\ref{42801}) and recalling the first equality in (\ref{6810}), we get
$$
\begin{aligned}
&\cos\alpha_b \bar{\partial}_s u_b+\sin\alpha_b \bar{\partial}_s v_b\\=~&
\cos A_b\left(-\tau_bq_b\sin^2 A_b\bar{\partial}_s \rho_b-\cos A_b\right)
+q_b\sin A_b\left(\bar{\partial}_s \phi-\frac{\bar{\partial}_s m}{\rho_bL_b}+\frac{\sin A_b\cos A_b}{\rho_b}\bar{\partial}_s\rho_b-\frac{N_b}{q_b^2}\right)+1\\=~&
q_b\sin A_b \left(\bar{\partial}_s \phi-\frac{\bar{\partial}_s m }{\rho_bL_b}\right) \quad \text{along}\quad \widehat{\mathrm{FG}}.
\end{aligned}
$$
Therefore, by (\ref{250402}) and (\ref{61109}) we have that when $\tau_1^e-\tau_0$ is sufficiently small, $$\cos\alpha_b \bar{\partial}_s u_b+\sin\alpha_b \bar{\partial}_s v_b>0\quad  \mbox{along}\quad  \widehat{\mathrm{FG}}.$$

 From (\ref{61109}), we have
\begin{equation}\label{82006}
\begin{aligned}
&\cos\beta_b \bar{\partial}_s u_b+\sin\beta_b \bar{\partial}_s v_b\\~=~&
\cos\beta_b \bar{\partial}_s (U_b+\xi)+\sin\beta_b \bar{\partial}_s (V_b+\eta)\\~=~&
\cos\beta_b \bar{\partial}_s (q_b\cos\sigma_b)+\sin\beta_b \bar{\partial}_s (q_b\sin\sigma_b)+\cos\beta_b\cos\phi+\sin\beta_b\sin\phi\\~=~&
\cos A_b \bar{\partial}_s q_b-q_b\sin A_b \bar{\partial}_s \sigma_b+\cos\beta_b\cos\phi+\sin\beta_b\sin\phi\\~=~&
\cos A_b \bar{\partial}_s q_b-q_b\sin A_b \bar{\partial}_s \sigma_b+\cos(2A_b)\quad \text{along}\quad \widehat{\mathrm{FG}}.
\end{aligned}
\end{equation}
Inserting (\ref{42803}) and (\ref{42804}) into (\ref{82006}), we get
$$
\begin{aligned}
&\cos\beta_b \bar{\partial}_s u_b+\sin\beta_b \bar{\partial}_s v_b\\=~&
\cos A_b\left(-\tau_bq_b\sin^2 A_b\bar{\partial}_s \rho_b-\cos A_b\right)
-q_b\sin A_b\left(\bar{\partial}_s \phi-\frac{\bar{\partial}_s m}{\rho_bL_b}+\frac{\sin A_b\cos A_b}{\rho_b}\bar{\partial}_s\rho_b-\frac{N_b}{q_b^2}\right)\\&+\cos(2A_b)\\=~&
-q_b\sin A_b \left(\bar{\partial}_s \phi-\frac{\bar{\partial}_s m }{\rho_bL_b}\right)-
2\tau_bq_b\cos A_b\sin^2 A_b\bar{\partial}_s \rho_b\quad \text{along}\quad \widehat{\mathrm{FG}}.
\end{aligned}
$$

Using $\tau_b=\psi_{po}(\tau_f)$ and (\ref{82008}), we obtain that
when $\tau_1^e-\tau_0$ is sufficiently small,
\begin{equation}\label{61308}
\bar{\partial}_s \rho_b=-\frac{1}{\tau_b^2}\psi_{po}'(\tau_f)\bar{\partial}_s \tau_f>0 \quad \text{along}\quad \widehat{\mathrm{FG}}.
\end{equation}
Therefore, when $\tau_1^e-\tau_0$ is sufficiently small, $ \cos\beta_b \bar{\partial}_s u_b+\sin\beta_b \bar{\partial}_s v_b<0$ along $\widehat{\mathrm{FG}}$.
This completes the proof for  (\ref{61111}).
The estimate (\ref{61111}) also admits a clear geometric interpretation, as illustrated in Figure \ref{Fig6.3} (left).

We can use (\ref{61109a}) to establish (\ref{61305}).  Since the proof is analogous to that of (\ref{61111}). We omit the details.
This completes the proof of the lemma.
\end{proof}

\begin{lem}\label{lem43}
We have the following convergence results as $\tau_0\rightarrow \tau_1^e$:
\begin{equation}\label{72203}
\lim\limits_{\tau_0\rightarrow \tau_1^e}\Big\|\big(\alpha_b-\alpha_4^e,~ \beta_b+\frac{\pi}{2},~ ~\tau_b-\tau_2^e\big)\Big\|_{0; \widehat{\mathrm{EG}}}~=~ 0
\end{equation}
and
\begin{equation}
\lim\limits_{\tau_0\rightarrow \tau_1^e}\Big\|\big(\alpha_b-\phi_{d},~ \beta_b-\beta_5^e,~ ~\tau_b-\tau_2^e\big)\Big\|_{0; \widehat{\mathrm{FG}}}~=~ 0.
\end{equation}
%where the constants $u_*$ and $\sigma_*$ are defined in (\ref{62903}).
\end{lem}
\begin{proof}
It is straightforward to verify that
$$
\mathrm{E}\rightarrow (\xi_e, \eta_e)\quad\text{and}\quad
(u_4, v_4)(\mathrm{E})\rightarrow (u_4^e, 0)\quad \text{as}\quad \tau_0\rightarrow\tau_1^e.
$$
Consequently, we have $(\alpha_4, \beta_4, \tau_4)(\mathrm{E})\rightarrow (\alpha_4^e, -\frac{\pi}{2}, \tau_2^e)$ as $\tau_0\rightarrow\tau_1^e$.
From (\ref{6811}), (\ref{61113}), (\ref{61105}), and (\ref{61115}), we deduce that $(\bar{\partial}_s\tau_b, \bar{\partial}_s\alpha_b, \bar{\partial}_s\beta_b)$ is uniformly bounded on $\widehat{\mathrm{EG}}$ with respect to the small parameter $\tau_1^e-\tau_0$.
Therefore, the limit relation (\ref{72203}) follows immediately.
The proof for the other is analogous, and we omit the details.
\end{proof}

 %So, we call this type of characteristic boundary value problem the singular characteristic boundary value problem.

\subsection{Flow downstream of the post-sonic shocks}
\subsubsection{The flow downstream of the post-sonic shock $\widehat{\mathrm{FG}}$}
%From (\ref{61113}) we know that $\widehat{\mathrm{FG}}$ and $C_{+}$ characteristic have the same direction at each point on $\widehat{\mathrm{FG}}$. While, by the second inequality in (\ref{61111}) and the first equation of (\ref{form}) we know that  $\widehat{\mathrm{FG}}$ is not a characteristic.
In what follows, we use the hodograph transformation method to find the flow downstream of the post-sonic shock $\widehat{\mathrm{FG}}$.
We shall also show that $\widehat{\mathrm{FG}}$ is an envelope of the $C_{+}$ characteristic curves of the flow downstream of it, and that the  flow downstream of $\widehat{\mathrm{FG}}$  satisfies $\bar{\partial}_{-}\rho=-\infty$ on  $\widehat{\mathrm{FG}}$.

By Lemmas \ref{61302} and \ref{lem43}, we know that when $\tau_1^e-\tau_0$ is sufficiently small, the map
$$(u, v)=(u_b, v_b)(\xi, \eta), \quad (\xi, \eta)\in \widehat{\mathrm{FG}},$$
has an inverse denoted by
$$
(\xi, \eta)=(\xi_b, \eta_b)(u, v), \quad(u, v)\in \widehat{\mathrm{F_7G_7}},
$$
where
$\widehat{\mathrm{F_7G_7}}=\big\{(u, v)|(u, v)=(u_b, v_b)(\xi, \eta), (\xi, \eta)\in \widehat{\mathrm{FG}}\big\}$, $\mathrm{F}_7=(u_b, v_b)(\mathrm{F})$, and $\mathrm{G}_7=(u_b, v_b)(\mathrm{G})$.

We now consider (\ref{52305}) with the initial data
\begin{equation}\label{6401}
(\alpha, \beta, \tau)=(\alpha_h, \beta_h, \tau_h)(u, v)\quad \text{on}\quad \widehat{\mathrm{F_7 G_7}},
\end{equation}
where $(\alpha_h, \beta_h, \tau_h)(u, v)=(\alpha_b, \beta_b, \tau_b)\big(\xi_b(u, v), \eta_b(u,v)\big)$.

%=============
\begin{figure}[htbp]
\begin{center}
\includegraphics[scale=0.38]{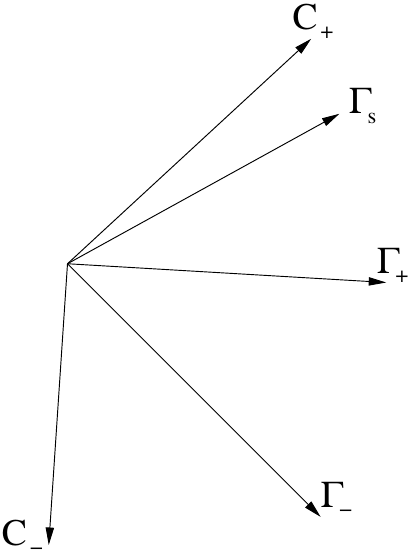}\qquad \qquad \includegraphics[scale=0.45]{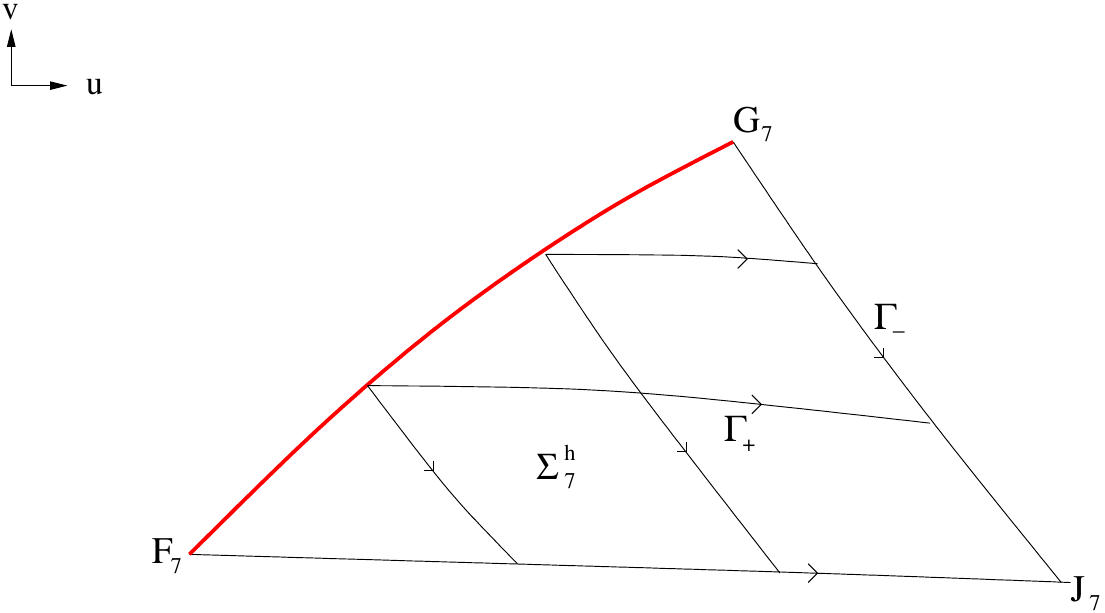}
\caption{\footnotesize Flow downstream of the post-sonic shock $\widehat{\mathrm{FG}}$ in the hodograph plane.}
\label{Fig6.3}
\end{center}
\end{figure}
%==========

From Lemma \ref{61302}, we see that
the $\Gamma_{\pm}$ characteristic directions $(-\sin\beta, \cos\beta)$ and $(\sin\alpha, -\cos\alpha)$ on
$\widehat{\mathrm{F_7G_7}}$ point to the same side of $\widehat{\mathrm{F_7G_7}}$; see Figure \ref{Fig6.3} (right).
Thus, the curve $\widehat{\mathrm{F_7G_7}}$ is a space-like non-characteristic curve,
and
the
problem (\ref{52305}), (\ref{6401}) is actually a Cauchy problem.

%For $\delta>0$ we let
%$$
%\Gamma_{\delta}=\big\{(u, v)\big| (u+\delta, v)\in \widehat{\mathrm{F_7 G_7}}\big\}.
%$$
%\begin{lem}
%There exists a small $\delta>0$ such that the Cauchy problem (\ref{52305}, \ref{6401}) admits a local solution in a domain $\Sigma_7^h(\delta)$ closed by $\Gamma_{\delta}$, a forward $\Gamma_{+}$ characteristic curve emanating from $\mathrm{F}_7$, and a forward $\Gamma_{-}$ characteristic curve emanating from $\mathrm{G}_7$; see Figure \ref{Fig15} (left).
%\end{lem}
%\begin{proof}
%The local existence can be obtained by the classical characteristic method; we omit the details.
%\end{proof}

% Next, we will extend the local solution to the whole determinate region.

\begin{lem}
(Uniform $C^1$ norm estimate) Consider the Cauchy problem (\ref{52305}), (\ref{6401}). Assume that $\tau_1^e-\tau_0$ is sufficiently small. Assume also that the problem admits a $C^1$ solution on some domain $\Sigma_{loc}$ determined by the characteristic method. Then, the $C^1$ norm of $(\alpha, \beta, \tau)$ has a uniform bound, which depends only on the $C^1$ norm of the Cauchy data (\ref{6401}).  That is, there is a constant $\mathcal{D}>0$, depending only on the Cauchy data, but not on $\Sigma_{loc}$, such that
$$
\|(\alpha, \beta, \tau)\|_{1; \Sigma_{loc}}\leq \mathcal{D}.
$$
\end{lem}
\begin{proof}
System (\ref{52305}) is  a linearly degenerate hyperbolic system. Thus, the derivatives of solution will not blow up unless $A=0$ or $A=\frac{\pi}{2}$.  Integrating (\ref{52305}) along the forward $\Gamma_{\pm}$ characteristic curves emanating from $\widehat{\mathrm{F_7G_7}}$ and  recalling (\ref{72203}), we know that
for any small $\varepsilon>0$, there exists a $\delta>0$ independent of $\Sigma_{loc}$ such that when $0<\tau_1^e-\tau_0<\delta$, the solution satisfies $$
\big\|(\alpha-\phi_{d},~ \beta-\beta_5^e,~\tau-\tau_2^e)\big\|_{0;\Sigma_{loc}}<\varepsilon.
$$
In fact, from (\ref{52305}) one can see that  $\|(\nabla\alpha,~ \nabla\beta,~\nabla\tau)\|_{0; \Sigma_{loc}}$ is independent of the derivatives of the Cauchy data (\ref{6401}).
Thus, one can obtained the existence of $\delta$.
The derivatives of the solution can be controlled by  (\ref{52305}), (\ref{52502}), and the  $C^1$ norm of the Cauchy data (\ref{6401}).
This completes the proof.
\end{proof}

By the classical extension method, we have the following existence result.
\begin{lem}
 Assume that $\tau_1^e-\tau_0$ is sufficiently small. Then the Cauchy problem (\ref{52305}), (\ref{6401}) admits a global $C^1$ solution denoted by $(\alpha, \beta, \tau)=(\alpha_7^h, \beta_7^h, \tau_7^h)(u,v)$ on the curvilinear triangular domain $\overline{\Sigma_7^h}$ bounded by $\widehat{\mathrm{F_7 G_7}}$, $\widetilde{\mathrm{F_7 J_7}}$, and $\widetilde{\mathrm{G_7 J_7}}$, where $\widetilde{\mathrm{F_7 J_7}}$ is a forward $\Gamma_{+}$ characteristic curve emanating from $\mathrm{F}_7$, and $\widetilde{\mathrm{G_7 J_7}}$ is a forward $\Gamma_{-}$ characteristic curve emanating from $\mathrm{G}_7$; see Figure \ref{Fig6.3} (right). Moreover, the solution satisfies
\begin{equation}\label{72205}
\big\|(\alpha-\phi_{d},~ \beta-\beta_5^e,~\tau-\tau_2^e)\big\|_{0; \Sigma_7^h}\rightarrow0\quad \text{as}\quad \tau_1^e-\tau_0\rightarrow 0.
\end{equation}
\end{lem}

We next show that the hodograph transformation is globally one-to-one.

\begin{lem}\label{102502}
For the $C^1$ solution of  the Cauchy problem (\ref{52305}), (\ref{6401}),
we have
\begin{equation}\label{61316}
\mathcal{Z}_{+}<0\quad \text{and}\quad \mathcal{Z}_{-}>0\quad \text{in}\quad \overline{\Sigma_7^h}\setminus\widehat{\mathrm{F_7G_7}},
\end{equation}
where the variables $\mathcal{Z}_{\pm}$ are defined in (\ref{62701}).
\end{lem}
\begin{proof}
We introduce the directional derive
$$
\hat{\partial}_{s}=\Gamma_s^1\partial_u+\Gamma_s^2\partial_v\quad \text{along}\quad \widehat{\mathrm{F_7G_7}},
$$
where $(\Gamma_s^1, \Gamma_s^2)=(\bar{\partial}_s u_b, \bar{\partial}_s v_b)$.
Then, by (\ref{317}) we have
\begin{equation}\label{61309}
\hat{\partial}_{s}=\frac{(\Gamma_s^1\cos\alpha+\Gamma_s^2\sin\alpha)\hat{\partial}_{+}}{\sin 2A}+
\frac{(\Gamma_s^1\cos\beta+\Gamma_s^2\sin\beta)\hat{\partial}_{-}}{\sin 2A}.
\end{equation}

Since $\phi=\alpha_b$ on $\widehat{\mathrm{FG}}$, we have
$$
\hat{\partial}_{s}\eta=(\tan\alpha_{7}^{h})\hat{\partial}_{s}\xi \quad \text{on}\quad \widehat{\mathrm{F_7G_7}}.
$$
Thus, by (\ref{61309}) we obtain
$$
\begin{aligned}
&(\Gamma_s^1\cos\alpha_{7}^{h}+\Gamma_s^2\sin\alpha_{7}^{h})\hat{\partial}_{+}\eta+
(\Gamma_s^1\cos\beta_{7}^{h}+\Gamma_s^2\sin\beta_{7}^{h})\hat{\partial}_{-}\eta\\[2pt]~=~&\Big(
(\Gamma_s^1\cos\alpha_{7}^{h}+\Gamma_s^2\sin\alpha_{7}^{h})\hat{\partial}_{+}\xi+
(\Gamma_s^1\cos\beta_{7}^{h}+\Gamma_s^2\sin\beta_{7}^{h})\hat{\partial}_{-}\xi\Big)\tan\alpha_{7}^{h}\quad \text{on}\quad \widehat{\mathrm{F_7G_7}}.
\end{aligned}
$$
Inserting (\ref{38}) into this, we get
$$
(\Gamma_s^1\cos\beta_{7}^{h}+\Gamma_s^2\sin\beta_{7}^{h})(\tan\beta_{7}^{h}-\tan\alpha_{7}^{h})\hat{\partial}_{-}\xi=0  \quad \text{on}\quad \widehat{\mathrm{F_7G_7}}.
$$
Thus, by the second inequality in (\ref{61111}) we have
\begin{equation}\label{61310}
\hat{\partial}_{-}\xi=0 \quad \text{on}\quad \widehat{\mathrm{F_7G_7}}.
\end{equation}
Combining this with (\ref{61312}), (\ref{61319}), and (\ref{62701}), we obtain
\begin{equation}\label{82201}
\mathcal{Z}_{-}=0\quad \text{on}\quad \widehat{\mathrm{F_7G_7}}.
\end{equation}

A direct computation yields
$$
\hat{\partial}_{s}\xi=\xi_u\bar{\partial}_s u+\xi_v\bar{\partial}_s v=\bar{\partial}_s\xi=\cos\alpha_{7}^{h}\quad \text{on}\quad \widehat{\mathrm{F_7G_7}}.
$$
Hence, by (\ref{61312}), (\ref{61311}), (\ref{61111}), (\ref{61309}), and (\ref{61310}), we have
\begin{equation}\label{61314}
\mathcal{Z}_{+}=-\frac{2\sin^2A_{7}^{h}\sin (2A_{7}^{h})}{c_{7}^{h}(\Gamma_s^1\cos\alpha_{7}^{h}+\Gamma_s^2\sin\alpha_{7}^{h})}<0\quad \text{on}\quad \widehat{\mathrm{F_7G_7}},
\end{equation}
where $A_{7}^{h}=\frac{\alpha_{7}^{h}-\beta_{7}^{h}}{2}$.

From (\ref{72205}), we deduce that when $\tau_1^e-\tau_0$ is sufficiently small,
$$
\frac{\tau_{7}^{h} p''(\tau_{7}^{h})}{4c(\tau_{7}^{h})p'(\tau_{7}^{h})}<0
\quad \text{on}\quad \Sigma_7^h.
$$
Therefore,
 integrating (\ref{322}) along the forward $\Gamma_{\pm}$ characteristic curves emanating from $\widehat{\mathrm{F_7G_7}}$ and recalling (\ref{82201}) and (\ref{61314}), we obtain (\ref{61316}).
This completes the proof.
\end{proof}

%From (\ref{61312})--(\ref{61319}), (\ref{72205}) and (\ref{61316}) we have that when $\tau_1^e-\tau_0$ is small,
%\begin{equation}\label{61321}
%\hat{\partial}_{+}\xi=-\frac{c\cos\alpha}{2\sin^2A}\mathcal{Z}_{+}<0\quad \text{and}\quad
%\hat{\partial}_{-}\xi=\frac{c\cos\beta}{2\sin^2A}\mathcal{Z}_{-}<0 \quad \text{in}\quad \Sigma_7^h\setminus\widehat{\mathrm{F_7G_7}}.
%\end{equation}
%Combining this with (\ref{38}) we have
%\begin{equation}\label{61322}
%\hat{\partial}_{+}\eta<0\quad \text{and}\quad
%\hat{\partial}_{-}\eta>0 \quad \text{in}\quad \Sigma_7^h\setminus\widehat{\mathrm{F_7G_7}}.
%\end{equation}
%For any two points in $\Sigma_7^h\setminus\widehat{\mathrm{F_7G_7}}$, there exist two characteristic curves with different type connecting the two points. By (\ref{61321}) and (\ref{61322}) we know that either $\xi$ or $\eta$ is monotone along the connecting path. Thus, no two points from $(u, v)$ maps to one point in the $(\xi, \eta)$-plane.
We define the functions
$$
\xi_7^{h}(u, v)=u-\frac{c(\tau_7^{h}(u, v))\cos(\sigma_7^{h}(u, v))}{\sin (A_7^{h}(u, v))}\quad\text{and}\quad \eta_7^{h}(u, v)=
v-\frac{c(\tau_7^{h}(u, v))\sin(\sigma_7^{h}(u, v))}{\sin (A_7^{h}(u, v))}
$$
for $(u, v)\in \overline{\Sigma_7^h}$.
By the sign-preserving property (\ref{61316}), we can use Theorem \ref{onetoone1} and Remark \ref{92505aa} to obtain that the hodograph transformation can be inverted to the physical plane, i.e.,
 the mapping
$$
(\xi, \eta)=\big(\xi_7^{h}, \eta_7^{h}\big)(u, v),\quad (u, v)\in \overline{\Sigma_7^h}
$$
has an inverse mapping denoted by
$$
(u, v)=(u_7, v_7)(\xi, \eta), \quad  (\xi, \eta)\in \overline{\Sigma_7},
$$
where
$$
\overline{\Sigma_7}=\Big\{(\xi, \eta) \mid  (\xi, \eta)=(\xi_7^{h}(u, v), \eta_7^{h}(u, v)), (u, v)\in \overline{\Sigma_7^h} \Big\}.
$$

Let
$\tau_7(\xi, \eta)=\tau_7^h(u_7(\xi, \eta), v_7(\xi, \eta)).$
We have $(u_7, v_7, \tau_7)(\xi, \eta)\in C^{1}(\overline{\Sigma_7}\setminus\widehat{\mathrm{FG}})\cap C^0(\overline{\Sigma_7})$.
Let
\begin{equation}\label{82202}
(u, v, \tau)=(u_7, v_7, \tau_7)(\xi, \eta), \quad (\xi, \eta)\in \overline{\Sigma_7}.
\end{equation}
Then,
the function $(u, v, \tau)$ defined in (\ref{82202}) is a solution of (\ref{PEC}) in $\overline{\Sigma_7}\setminus\widehat{\mathrm{FG}}$, and satisfies 
$$(u, v, \tau)=(u_b, v_b, c_b)\quad  \mbox{on} \quad \widehat{\mathrm{FG}}.$$
The domain $\overline{\Sigma_7}$ is bounded by $\widehat{\mathrm{FG}}$, $\widetilde{\mathrm{FJ}}$, and $\widetilde{\mathrm{GJ}}$, where $\widetilde{\mathrm{FJ}}$ is a forward $C_{+}$ characteristic curve emanating from $\mathrm{F}$ and $\widehat{\mathrm{GJ}}$ is a forward $C_{-}$ characteristic curve emanating from $\mathrm{G}$; see Figure \ref{Fig6.2} (right).
From (\ref{72205}), we also have that the solution satisfies
\begin{equation}\label{72801}
\big\|(\alpha-\phi_{d},~ \beta-\beta_5^e,~\tau-\tau_2^e)\big\|_{0; \Sigma_7}\rightarrow0
\quad \text{as}\quad \tau_1^e-\tau_0\rightarrow 0.
\end{equation}
From (\ref{6807}), (\ref{6808}), and (\ref{61316}) we have
\begin{equation}\label{102509}
\bar{\partial}_{+}\rho\mid_{\widetilde{\mathrm{FJ}}}~<~0\quad \text{and}\quad \bar{\partial}_{-}\rho\mid_{\widetilde{\mathrm{GJ}}}~<~0.
\end{equation}

From (\ref{6808}) and (\ref{82201}), we deduce that the solution also satisfies 
\begin{equation}\label{26101501}
\bar{\partial}_{-}\rho=-\infty\quad  \text{on} \quad \widehat{\mathrm{FG}}.
\end{equation}
Moreover, by (\ref{322}), (\ref{26102402}), (\ref{6808}), (\ref{82201}), and (\ref{61314}), we also have that along each $C_{-}$ characteristic curve emanating from $\widehat{\mathrm{FG}}$,
\begin{equation}\label{26102403}
(\bar{\partial}_{-}\rho)\cdot\sqrt{\ell}=O(1)\quad \mbox{as}\quad \ell \rightarrow 0,
\end{equation}
where $\ell$ denotes the arc length of $C_{-}$ characteristic curve. 

%\subsubsection{\bf Post-sonic shock through $B$}
%Let $\widetilde{BG}$ be defined by a post-sonic shock.

%We consider (\ref{PEC}) with data
%\begin{equation}\label{63003}
%(u, v, \tau)=(u_b, v_b, \tau_v)(\xi, \eta), \quad (\xi, \eta)\in\widetilde{BG}.
%\end{equation}

%The singular Cauchy problem (\ref{PEC}), \ref{63003} admits a classical solution in a triangle domain $\Sigma_8$ bounded by $\widetilde{\mathrm{BG}}$, $\widetilde{\mathrm{BI}}$, and $\widetilde{\mathrm{GI}}$, where $\widetilde{\mathrm{GI}}$ is a $C_{+}$ characteristic curve emanating from $\mathrm{G}$ and $\widetilde{\mathrm{BI}}$ is a $C_{-}$ characteristic curve emanating from $B$; see Figure \ref{Fig15}(right).
%From (\ref{72205}) we also have that the solution satisfies
%\begin{equation}\label{72801}
%\big\|\big(\alpha_b-\alpha_3(B),~ \beta_b+\frac{\pi}{2},~\tau-\tau_2^e\big)\big\|_{0; \Sigma_7}\rightarrow0
%\quad \text{as}\quad \tau_1^e-\tau_0\rightarrow 0.
%\end{equation}
%From (\ref{6807}), (\ref{6808}), and (\ref{61316}) we have
%\begin{equation}\label{102509}
%\bar{\partial}_{+}\rho\mid_{\widetilde{\mathrm{DJ}}}~<~0.
%\end{equation}
%From (\ref{6808}) and (\ref{82201}) we can see that the solution also satisfies $\bar{\partial}_{-}\rho=-\infty$ on $\widetilde{\mathrm{DJ}}$.

%\subsection{Singular characteristic boundary value problem}

\subsubsection{The flow downstream of the post-sonic shock $\widehat{\mathrm{EG}}$}
%Similarly, we can construct the flow downstream the post-sonic shock $\widehat{\mathrm{EG}}$.

When $\tau_1^e-\tau_0$ is sufficiently small, the mapping
$(u, v)=(u_b, v_b)(\xi, \eta)$, $(\xi, \eta)\in \widehat{\mathrm{EG}}$
has an inverse denoted by
$$
(\xi, \eta)=(\xi_b, \eta_b)(u, v), \quad(u, v)\in \widehat{\mathrm{E_8G_8}},
$$
where
$\widehat{\mathrm{E_8G_8}}=\big\{(u, v)|(u, v)=(u_b, v_b)(\xi, \eta), (\xi, \eta)\in \widehat{\mathrm{EG}}\big\}$, $\mathrm{E}_8=(u_b, v_b)(\mathrm{E})$, and $\mathrm{G}_8=(u_b, v_b)(\mathrm{G})$.
We now consider (\ref{52305}) with data
\begin{equation}\label{6401a}
(\alpha, \beta, \tau)=(\alpha_h, \beta_h, \tau_h)(u, v)\quad \text{on}\quad \widehat{\mathrm{E_8 G_8}},
\end{equation}
where $(\alpha_h, \beta_h, \tau_h)(u, v)=(\alpha_b, \beta_b, \tau_b)\big(\xi_b(u, v), \eta_b(u,v)\big)$.

\begin{lem}
 Assume that $\tau_1^e-\tau_0$ is sufficiently small. Then, the Cauchy problem (\ref{52305}), (\ref{6401a}) admits a global $C^1$ solution denoted by $(\alpha, \beta, \tau)=(\alpha_8^h, \beta_8^h, \tau_8^h)(u,v)$ on a curvilinear triangular domain $\overline{\Sigma_8^h}$ bounded by $\widehat{\mathrm{E_8 G_8}}$, $\widetilde{\mathrm{E_8 I_8}}$, and $\widetilde{\mathrm{G_8 I_8}}$, where $\widetilde{\mathrm{E_8 I_8}}$ is a forward $\Gamma_{-}$ characteristic curve emanating from $\mathrm{E}_8$ and $\widetilde{\mathrm{G_8 I_8}}$ is a forward $\Gamma_{+}$ characteristic curve emanating from $\mathrm{G}_8$. Moreover, the solution satisfies
\begin{equation}\label{72205a}
\Big\|(\alpha-\alpha_4^e,~ \beta+\frac{\pi}{2},~\tau-\tau_2^e)\Big\|_{0, \Sigma_8^h}\rightarrow0\quad \text{as}\quad \tau_1^e-\tau_0\rightarrow 0.
\end{equation}
\end{lem}

We define the functions
$$
\xi_8^{h}(u, v)=u-\frac{c(\tau_8^{h}(u, v))\cos(\sigma_8^{h}(u, v))}{\sin (A_8^{h}(u, v))}\quad\text{and}\quad\eta_8^{h}(u, v)=v-\frac{c(\tau_8^{h}(u, v))\sin(\sigma_8^{h}(u, v))}{\sin (A_8^{h}(u, v))}
$$
for $(u, v)\in\overline{\Sigma_8^h}$.
Then the mapping
$$
(\xi, \eta)=\big(\xi_8^{h}, \eta_8^{h}\big)(u, v),\quad (u, v)\in \overline{\Sigma_8^h},
$$
has an inverse denoted by
$$
(u, v)=(u_8, v_8)(\xi, \eta), \quad  (\xi, \eta)\in \overline{\Sigma_8},
$$
where
$$
\overline{\Sigma_8}=\Big\{(\xi, \eta) \mid  (\xi, \eta)=(\xi_8^{h}(u, v), \eta_8^{h}(u, v)), (u, v)\in \overline{\Sigma_8^h} \Big\}.
$$

Let
$\tau_8(\xi, \eta)=\tau_8^h(u_8(\xi, \eta), v_8(\xi, \eta)).$
We have $(u_8, v_8, \tau_8)(\xi, \eta)\in C^{1}(\overline{\Sigma_8}\setminus\widehat{\mathrm{EG}})\cap C^0(\overline{\Sigma_8})$.
Let
\begin{equation}\label{82202A}
(u, v, \tau)=(u_8, v_8, \tau_8)(\xi, \eta), \quad (\xi, \eta)\in \overline{\Sigma_8}.
\end{equation}
Then,
the $(u, v, \tau)$ defined in (\ref{82202}) is a solution of (\ref{PEC}) on $\overline{\Sigma_8}\setminus\widehat{\mathrm{FG}}$, and this solution satisfies 
$$(u, v, \tau)=(u_b, v_b, c_b)\quad  \mbox{on}\quad  \widehat{\mathrm{EG}}.$$
The domain $\overline{\Sigma_8}$ is bounded by $\widehat{\mathrm{EG}}$, $\widetilde{\mathrm{EI}}$, and $\widetilde{\mathrm{GI}}$, where $\widetilde{\mathrm{EI}}$ is a forward $C_{-}$ characteristic curve emanating from $\mathrm{E}$ and $\widetilde{\mathrm{GI}}$ is a forward $C_{+}$ characteristic curve emanating from $\mathrm{G}$; see Figure \ref{Fig6.2} (right).
From (\ref{72205a}), we also have that the solution satisfies
\begin{equation}\label{72801A}
\Big\|(\alpha-\alpha_4^e,~ \beta+\frac{\pi}{2},~\tau-\tau_2^e)\Big\|_{0; \Sigma_8}\rightarrow0
\quad \text{as}\quad \tau_1^e-\tau_0\rightarrow 0.
\end{equation}
From (\ref{6807}), (\ref{6808}), and (\ref{61316}), we have
\begin{equation}\label{102509A}
\bar{\partial}_{+}\rho\mid_{\widetilde{\mathrm{GI}}}~<~0\quad \text{and}\quad \bar{\partial}_{-}\rho\mid_{\widetilde{\mathrm{EI}}}~<~0.
\end{equation}

From (\ref{6808}) and (\ref{82201}), we deduce that the solution also satisfies 
\begin{equation}\label{26101502}
\bar{\partial}_{+}\rho=-\infty\quad  \text{on}\quad  \widehat{\mathrm{EG}}.
\end{equation}
Moreover, as in (\ref{26102403}), along each $C_{+}$ characteristic curve emanating from $\widehat{\mathrm{EG}}$,
\begin{equation}\label{26102404}
(\bar{\partial}_{+}\rho)\cdot\sqrt{\ell}=O(1)\quad \mbox{as}\quad \ell \rightarrow 0,
\end{equation}
where $\ell$ denotes the arc length of $C_{+}$ characteristic curve. 

\subsection{Global solvability of the fan-shock-fan composite wave interaction problem}
We next establish the existence of the solution in the remaining part of the fan-shock-fan composite wave interaction region.
To control the hyperbolicity of the system, we first construct an invariant region of $(\alpha, \beta)$.
Let $\tau_{_\mathrm{G}}=\tau_8(\mathrm{G})=\tau_7(\mathrm{G})=\psi_{po}(\tau_6(\mathrm{G}))$.
For $\tau\geq \tau_{_\mathrm{G}}$, we define
\begin{equation}
\Gamma_{3}(\tau)=\Big\{(\alpha, \beta)\mid \alpha_{3l}(\tau)<\alpha<\alpha_{3r},~\beta_{3l}<\beta<\beta_{3r}(\tau),~\alpha-\beta>\delta\Big\},
\end{equation}
where
\begin{equation}\label{102504a}
\begin{aligned}
&\alpha_{3l}(\tau)=\beta_5^e-\varepsilon+2\Theta(\tau), \quad \beta_{3l}=\beta_5^e-\varepsilon,\\&  \alpha_{3r}=\alpha_4^e+\varepsilon, \quad \beta_{3r}(s)=\alpha_4^e+\varepsilon-2\Theta(\tau),
\end{aligned}
\end{equation}
 $\varepsilon=\frac{1}{4}(\pi-\alpha_4^e+\beta_5^e)$, and $\delta$ is a small positive constant.

\begin{lem}
Assume that assumptions ($h_{31}$)--($h_{32}$) hold.
Then, when $\tau_1^e-\tau_0$ is sufficiently small, we have
\begin{equation}\label{92803a}
(\alpha, \beta)(\xi, \eta)\in \Gamma_{3}(\tau(\xi, \eta))\quad \text{for} \quad (\xi, \eta)\in\widetilde{\mathrm{EB}}\cup\widetilde{\mathrm{FD}}\cup
\widetilde{\mathrm{GI}}\cup \widetilde{\mathrm{GJ}} \cup \widetilde{\mathrm{EI}}\cup \widetilde{\mathrm{FG}}.
\end{equation}
\end{lem}
\begin{proof}
It can be verified that $\alpha_4(E)\rightarrow \alpha_4^e$ as $\tau_1^e-\tau_0\rightarrow 0$. Thus, when $\tau_1^e-\tau_0$ is sufficiently small, we have $\alpha_4(\mathrm{E})+\frac{\pi}{2}>2\Theta(\tau_2^e)$. Thus, as in Lemma \ref{lem3.1}, the $C_{+}$ cross characteristic $\widetilde{\mathrm{EB}}$ exists, and $2\Theta(\tau)-\frac{\pi}{2}<\alpha\leq \alpha_4(\mathrm{E})$ on $\widetilde{\mathrm{EB}}$.
This implies that when $\tau_1^e-\tau_0$ is sufficiently small, $(\alpha, \beta)\in \Gamma_3(\tau)$ on $\widetilde{\mathrm{EB}}$.
Furthermore,
we have $\alpha_5(\mathrm{F})=\alpha_e\rightarrow \phi_d$ and $\beta_5(\mathrm{F})\rightarrow \beta_5^e$ as $\tau_1^e-\tau_0\rightarrow 0$. Thus, when $\tau_1^e-\tau_0$ is sufficiently small,
$\alpha_5(\mathrm{F})-\beta_5(\mathrm{F})>2\Theta(\tau_{2}^e)$.
Then, as in Lemma \ref{lem3.3}, the $C_{-}$ cross characteristic $\widetilde{\mathrm{FD}}$ exists, and there exists a $\delta_0>0$ such that $$\delta_0<A\leq \frac{\alpha_{5}(\mathrm{E})-\beta_5(\mathrm{F})}{2}\quad \mbox{on}\quad  \widetilde{\mathrm{FD}}.$$

It can be verified that $\alpha_e\rightarrow \phi_d$ and $\alpha_v\rightarrow \alpha_v^e$ as $\tau_1^e-\tau_0\rightarrow 0$.
Hence, when $\tau_1^e-\tau_0$ is sufficiently small,
$\theta>\beta_5(E)+2\Theta(\bar{\tau}_r(\theta))$ for  $\theta\in[\alpha_v, \alpha_e]$.
Then, when $\tau_1^e-\tau_0$ is sufficiently small, as in Lemma \ref{112602}, we obtain $(\alpha, \beta)\in \Gamma_{3}(\tau)$ on $\widetilde{\mathrm{FD}}$.

From (\ref{72205a}) and (\ref{72801A}), we immediately have $(\alpha, \beta)\in \Gamma_3(\tau)$ on $\widetilde{\mathrm{GI}}\cup \widetilde{\mathrm{GJ}} \cup \widetilde{\mathrm{EI}}\cup \widetilde{\mathrm{FG}}$, provided that $\tau_1^e-\tau_0$ is sufficiently small.
This completes the proof.
\end{proof}

First, we consider (\ref{PEC}) with the boundary conditions
\begin{equation}\label{63011a}
(u, v, \tau)=\left\{
               \begin{array}{ll}
                 (u_7, v_7, \tau_7)(\xi, \eta), & \hbox{$(\xi, \eta)\in \widetilde{\mathrm{GI}}$;} \\[2pt]
                 (u_8, v_8, \tau_8)(\xi, \eta), & \hbox{$(\xi, \eta)\in \widetilde{\mathrm{GJ}}$.}
               \end{array}
             \right.
\end{equation}
Problem (\ref{PEC}), (\ref{63011a}) is a discontinuous Goursat problem.

By (\ref{72801}), (\ref{72801A}), and assumption ($h_{31}$), we have that when $\tau_0$ is close to $\tau_1^e$,
$$
\sigma_{8}(\mathrm{G})>\sigma_7(\mathrm{G}).
$$
Thus, there exist $C_{+}$ and $C_{-}$ centered waves emanating from the point $\mathrm{G}$.

There exists a $q_{m}\in (q_{_{\mathrm{G}}}, q_{_{lim}})$ such that
$$
\int_{q_{_\mathrm{G}}}^{q_{m}}\frac{\sqrt{q^2-c^2(\grave{\tau}(q))}}{qc(\grave{\tau}(q))}{\rm d}q=\frac{\sigma_8(\mathrm{G})
-\sigma_7(\mathrm{G})}{2},
$$
where $q_{_\mathrm{G}}=\sqrt{(u_8(\mathrm{G})-\xi_{_\mathrm{G}})^2+(v_8(\mathrm{G})-\eta_{_\mathrm{G}})^2}$,  $q_{_{lim}}=\sqrt{q_{_\mathrm{G}}^2+2h(\tau_{_\mathrm{G}})-h_{\infty}}$, and the function $\grave{\tau}(q)$ is determined by Bernoulli's law
$q^2+2h(\grave{\tau}(q))=q_{_\mathrm{G}}^2+2h(\tau_{_\mathrm{G}})$.

We introduce the following constants:
$$
\begin{aligned}
&\sigma_m=\frac{\sigma_8(\mathrm{G})+\sigma_7(\mathrm{P})}{2}, \quad \tau_m=\grave{\tau}(q_m), \quad A_m=\arcsin\Big(\frac{c(\tau_m)}{q_m}\Big), \\ &\qquad \qquad\qquad~ \alpha_m=\sigma_m+A_m, \quad \text{and}\quad \beta_m=\sigma_m-A_m.
\end{aligned}
$$
We also define
\begin{equation}\label{121405a}
\nu_1=\tan\alpha_m\quad \text{and} \quad \mu_2=\cot\beta_m.
\end{equation}

 We solve (\ref{61902}) with the initial data
\begin{equation}\label{61803}
(\check{\alpha}_{+},\check{\beta}_{+})(\tau_{_\mathrm{G}})=(\alpha_8(\mathrm{G}), \beta_8(\mathrm{G})),
\end{equation}
and (\ref{120602}) with the initial data
\begin{equation}\label{120603a}
(\check{\alpha}_{-},\check{\beta}_{-})(\tau_{_\mathrm{G}})=(\alpha_7(\mathrm{G}), \beta_7(\mathrm{G})).
\end{equation}
Both initial value problems admit admit a solution $(\check{\alpha}_{-},\check{\beta}_{-})(\tau)$ on $[\tau_{_\mathrm{G}}, \tau_m]$. Moreover,
$$
(\check{\alpha}_{+},\check{\beta}_{+})(\tau_m)=(\check{\alpha}_{-},\check{\beta}_{-})(\tau_m)=(\alpha_m, \beta_m).
$$

By (\ref{72801}) and (\ref{72801A}), we deduce that
under  assumptions ($h_{31}$)-($h_{34}$),  when $\tau_1^e-\tau_0$ is sufficiently small, the following estimates hold:
\begin{description}
%\item[($\overline{H}_{61}$)] $\min\big\{\alpha_2(\mathrm{P})+\frac{\pi}{2},  \phi_{po}-\beta_3(\mathrm{P})\big\}>2\Theta(\tau_{po})$, $\alpha_2(\mathrm{P})-\beta_3(\mathrm{P})<\pi$, $\sigma_2(\mathrm{P})>\sigma_3(\mathrm{P})$;
  %\item[($\overline{H}_{62}$)] $\bar{\alpha}(\theta)>\beta_3(\mathrm{P})+2\Theta(\bar{\tau}(\theta))$ for  $\theta\in[\alpha_v, \alpha_0]$;
\item[($\overline{h}_{h3}$)] $\check{\alpha}_{+}(\tau)>\beta_5^e+2\Theta(\tau)$ for  $\tau\in[\tau_{_\mathrm{G}}, \tau_m]$;
   \item[($\overline{h}_{h4}$)] $\check{\beta}_{-}(\tau)<\alpha_4^e-2\Theta(\tau)$ for  $\tau\in[\tau_{_\mathrm{G}}, \tau_m]$.
\end{description}
Therefore,  as in Lemma \ref{dgploc}, we have the following result.
\begin{lem}
Let $\mathrm{K}$ be a point on $\widetilde{\mathrm{GI}}$, and let $\mathrm{L}$ be a point on $\widetilde{\mathrm{GJ}}$. When $\mathrm{K}$ and $\mathrm{L}$ are close to $\mathrm{G}$,
the DGP (\ref{PEC}), (\ref{63011a}) admits a continuous and piecewise $C^1$ solution on a curvilinear quadrilateral domain $\overline{\Sigma_9}$ bounded by  $\widetilde{\mathrm{GL}}$, $\widetilde{\mathrm{GK}}$, $\widetilde{\mathrm{LM}}$, and $\widetilde{\mathrm{KM}}$, where $\widetilde{\mathrm{LM}}$ is a forward $C_{+}$ characteristic curve emanating from $\mathrm{L}$ and $\widetilde{\mathrm{KM}}$ is a forward $C_{-}$ characteristic curve emanating from $\mathrm{K}$. Moreover, the solution satisfies the following estimates:
$$
(\alpha, \beta)\in \Gamma_3(\tau)\quad \text{in}\quad \overline{\Sigma_9}; \quad \bar{\partial}_{+}\rho<0\quad \text{on}\quad \widetilde{\mathrm{LM}}; \quad
\bar{\partial}_{-}\rho<0\quad \text{on}\quad \widetilde{\mathrm{KM}}.
$$
\end{lem}

%=============
\begin{figure}[htbp]
\begin{center}
\includegraphics[scale=0.2]{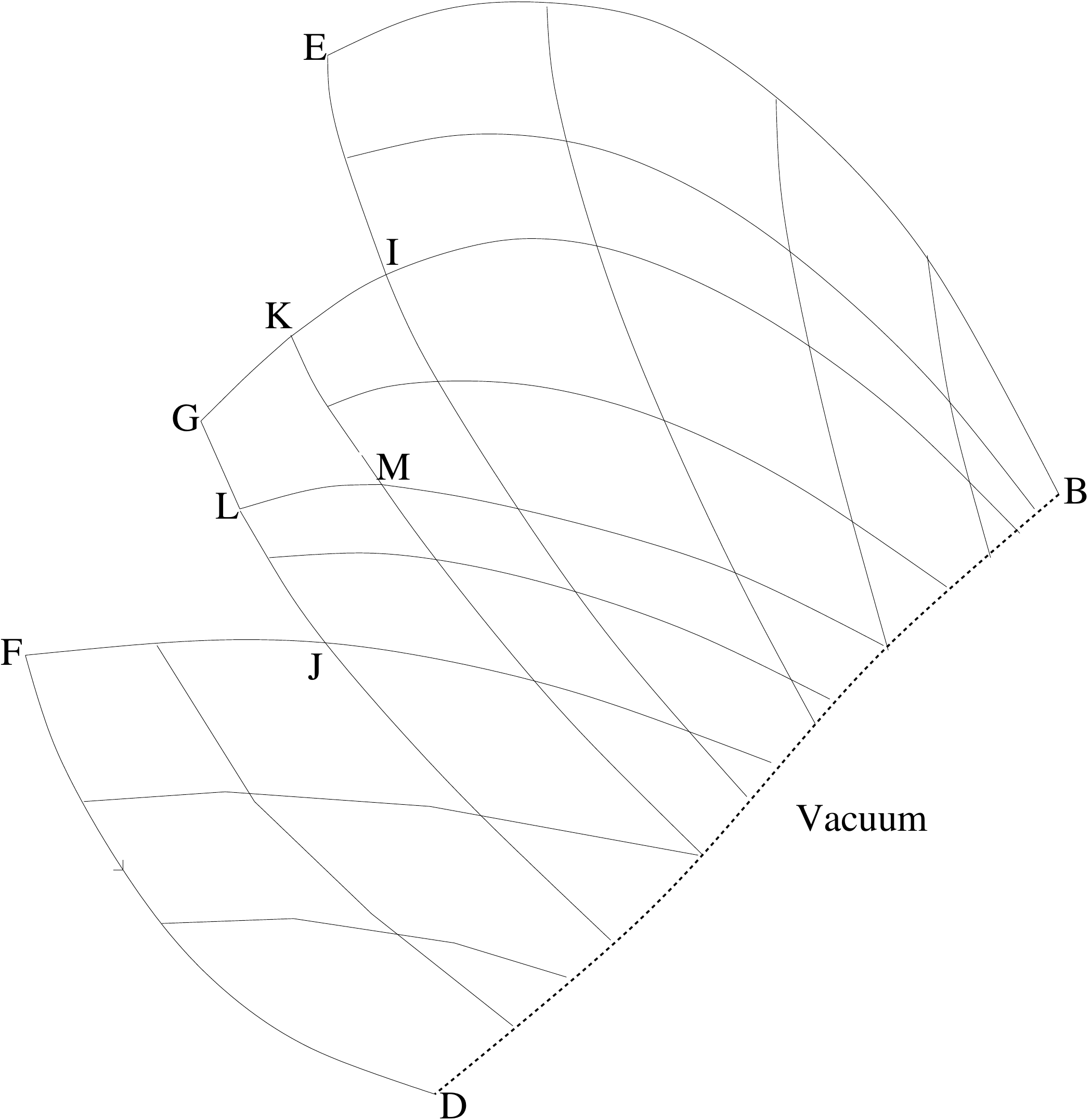}
\caption{\footnotesize Flow in the remaining part of the wave interaction region.}
\label{Fig6.4}
\end{center}
\end{figure}
%==========

For convenience, we use $(u, v, \tau)=(u_9, v_9, \tau_9)(\xi, \eta)$, $(\xi,\eta)\in\overline{\Sigma_9}\setminus\{\mathrm{G}\}$ to denote the local solution of the DGP (\ref{PEC}), (\ref{63011a}).
Second, we consider (\ref{PEC}) with the boundary conditions
\begin{equation}\label{63011aa}
(u, v, \tau)=\left\{
               \begin{array}{ll}
               (\hat{u}_r, \hat{v}_r, \hat{\tau}_r)(\xi), & \hbox{$(\xi, \eta)\in \widetilde{{\mathrm{EB}}}$;}\\[2pt]
                (\bar{u}_r, \bar{v}_r, \bar{\tau}_r)(\theta), & \hbox{$(\xi, \eta)\in \widetilde{\mathrm{FD}}$;}\\[2pt]
                 (u_7, v_7, \tau_7)(\xi, \eta), & \hbox{$(\xi, \eta)\in \widetilde{\mathrm{LJ}}\cup \widetilde{\mathrm{FJ}}$;}\\[2pt]
(u_8, v_8, \tau_8)(\xi, \eta), & \hbox{$(\xi, \eta)\in \widetilde{\mathrm{KI}}\cup\widetilde{\mathrm{EI}}$;}\\[2pt]
(u_9, v_9, \tau_9)(\xi, \eta), & \hbox{$(\xi, \eta)\in \widetilde{\mathrm{LM}}\cup\widetilde{\mathrm{KM}} $.}
               \end{array}
             \right.
\end{equation}
Problem (\ref{PEC}), (\ref{63011aa}) is a standard characteristic boundary value problem.
Let
$$
\Pi_{+}=\widetilde{\mathrm{EB}}\cup\widetilde{\mathrm{KI}}\cup \widetilde{\mathrm{LM}}\cup\widetilde{\mathrm{FJ}}  \quad \text{and}\quad \Pi_{-}= \widetilde{\mathrm{FD}}\cup \widetilde{\mathrm{LJ}}\cup\widetilde{\mathrm{EI}}\cup \widetilde{\mathrm{KM}}
$$
be the $C_{+}$ and $C_{-}$ characteristic boundaries, respectively.

From the previous discussion, we have that the boundary data satisfy the hyperbolicity estimate
$$
(\alpha, \beta)\in \Gamma_3(\tau)\quad \text{on}\quad \Pi_{+}\cup \Pi_{-},
$$
and the monotonicity condition
$$
 \bar{\partial}_{\pm}\rho<0\quad \text{on}\quad \Pi_{\pm}.
$$
Thus, using the previous method, we obtain that the problem (\ref{PEC}), (\ref{63011aa})
 admits a global classical solution on a domain $\Sigma_{10}$ bounded by $\Pi_{\pm}$ and a vacuum boundary $\widetilde{\mathrm{BD}}$ connecting $\mathrm{B}$ and $\mathrm{D}$.
For convenience, we denote the solution of the problem (\ref{PEC}), (\ref{63011aa}) by $(u, v, \tau)=(u_{10}, v_{10}, \tau_{10})(\xi, \eta)$, $(\xi, \eta)\in \Sigma_{10}$.

We define
\begin{equation}\label{91601ee}
(u, v, \tau)=
(u_i, v_i, \tau_i)(\xi, \eta), \quad (\xi,\eta)\in \Sigma_i \quad i= 6,~ 7,~ 8,~ 9, ~10.
\end{equation}
Then the function $(u, v, \tau)$ defined in (\ref{91601ee}) is a piecewise smooth solution to the BVP (\ref{PEC}), (\ref{63001}). This completes the proof of Theorem \ref{fsf/fsf}.

%\subsection{$\alpha_v<\theta_w-\pi<\sigma_2^e$}

\section{Interaction between a shock–fan composite wave and a shock wave}
%\subsection{Elementary waves}
In this section, we study the 2D Riemann problem (\ref{Potential}), (\ref{initial})--(\ref{boundary}) for $\tau_1^e\leq\tau_0<\tau_2^i$ and $u_0>\frac{\rho_{po}c_{po}}{\rho_0}$. % and $\pi+\sigma_{po}<\theta_w<\pi$, where the constant $\sigma_{po}$ is defined in (\ref{102408}).
In the subsequent discussion of this section, the constants $u_{w}$, $v_{w}$, $\tau_{w}$, $\tau_{po}$, $\phi_w$, $\xi_{_\mathrm{B}}$, $\xi_{_\mathrm{P}}$, as well as the functions $u_{-}(\xi)$, $v_{-}(\xi)$, $\tau_{-}(\xi)$ are as defined in Section 3.5.
Away from the corner, the oncoming flow and the vacuum
are connected by a planar shock–fan composite wave which can be represented by (\ref{71501c}). We assume that at the corner the oncoming flow $(0)$ is deflected through an oblique shock wave to a constant state $(u_{w}, v_{w}, \tau_{w})$  where $\tau_w\in(\tau_1^i, \tau_2^e)$. The inclination angle of the oblique shock is $\phi_{w}$, and the flow near the corner is represented by (\ref{41801}).

%

%where the function $(u_{-}, v_{-}, \tau_{-})(\xi)$ is defined in Section 1.4.2.

%Let
%$$
%\begin{aligned}
%\xi_{_\mathrm{P}}=u_0-\frac{\rho_{po}c_{po}}{\rho_0}\quad \text{and} \quad \eta_{_\mathrm{P}}=\xi_{_\mathrm{P}}\tan\phi_{w}.
%\end{aligned}
%$$

For convenience, we use the subscripts $2$ and $3$ to denote the corresponding flow states of the planar rarefaction fan wave and the constant state $(u_{w}, v_{w}, \tau_{w})$, respectively; that is,
$$
(u_2, v_2, \tau_2)(\xi, \eta)=(u_{-}, v_{-}, \tau_{-})(\xi)\quad \text{and} \quad (u_3, v_3, \tau_3)(\xi, \eta)=(u_{w}, v_{w}, \tau_{w}).
$$
%We assume $\tau_1^i\leq\tau_0<\tau_2^i$, $u_0>\frac{\rho_{po}c_{po}}{\rho_0}$ and $\pi+\sigma_{po}<\theta_w<\pi$. Then the oncoming supersonic flow is connected to the vacuum by a planar shock–fan wave far away from the origin, and is deflected by an oblique shock wave to a vacuum or a constant state near the corner.
%We denote by $(u_{w}, v_{w}, \tau_{w})$ the state of the flow downstream of the oblique shock and $\phi_{w}$ the inclination angle of the oblique shock.
We divide the discussion into the following two cases: (i) $\tau_{w}\geq \tau_2^i$; (ii) $\tau_{w}<\tau_2^i$.

 %(1). (2).
%When $\tau_{w}\geq \tau_2^i$ the interaction generates two centered waves emanating from $\mathrm{P}$; see Figure \ref{FS+S} (left). When $\tau_{w}<\tau_2^i$ the interaction generates a centered wave and a shock emanating from $\mathrm{P}$; see Figure \ref{FS+S} (right).

\subsection{Wave interaction for the case $\tau_{w}\geq \tau_2^i$}
\begin{figure}[htbp]
\begin{center}
  \includegraphics[scale=0.45]{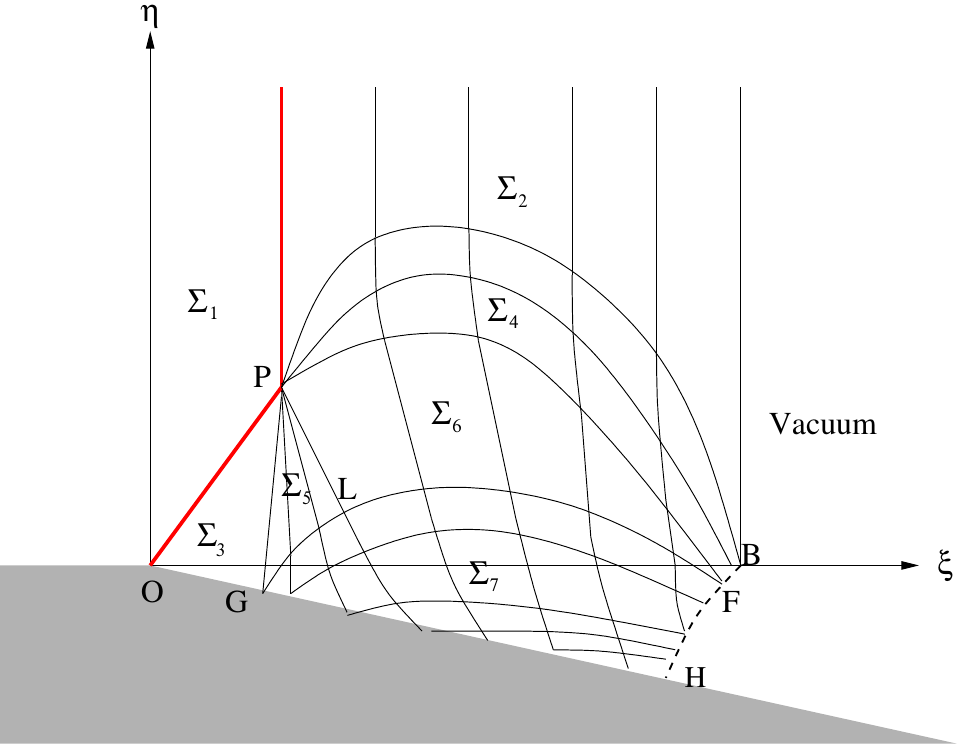} \qquad\quad\includegraphics[scale=0.45]{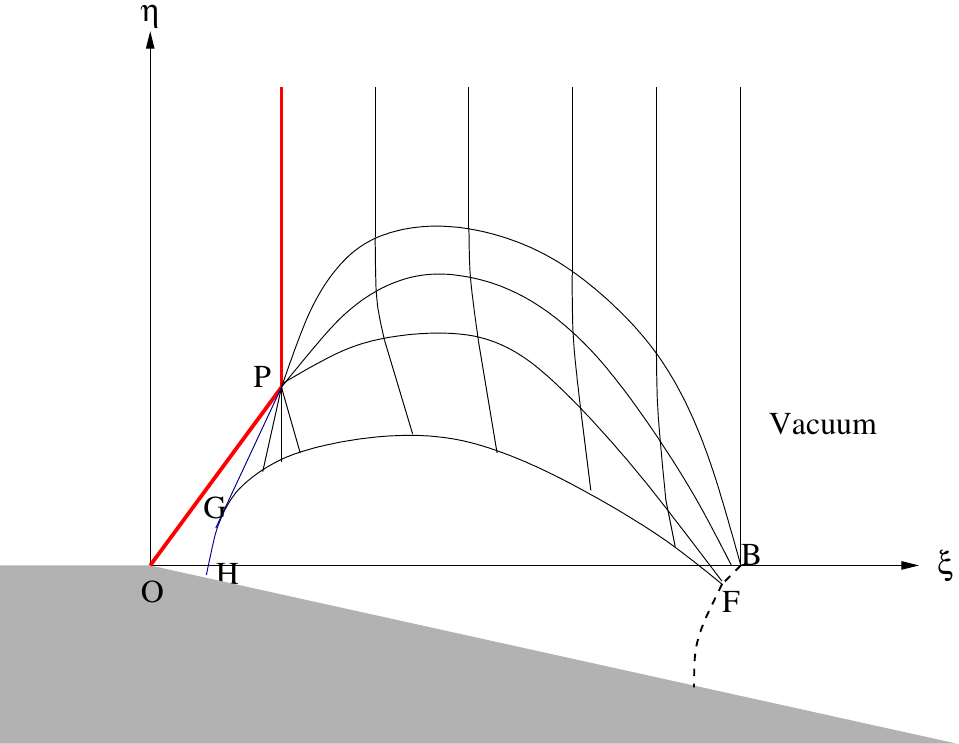}
\caption{\footnotesize Interaction between a planar shock–fan composite wave and an oblique shock wave for $\tau_{w}\geq \tau_2^i$.  }
\label{Fig7.1}
\end{center}
\end{figure}

We first assume
\begin{description}
\item[($h_{40}$)] $(u_w-\xi_{_\mathrm{P}})^2+(v_w-\eta_{_\mathrm{P}})^2>c^2(\tau_w)$,
\end{description}
where $\xi_{_\mathrm{P}}=u_0-\frac{\rho_{po}c_{po}}{\rho_0}$ and $\eta_{_\mathrm{P}}=\xi_{_\mathrm{P}}\tan\phi_{w}$.
Then the two elementary waves begin to interact at the point $\mathrm{P}(\xi_{_\mathrm{P}}, \eta_{_\mathrm{P}})$.
Through $\mathrm{P}$, we draw a $C_{+}$ cross characteristic curve of the planar fan wave.
Assume $\alpha_2(\mathrm{P})+\frac{\pi}{2}>2\Theta(\tau_{po})$. Then the $C_{+}$ cross characteristic curve  can cross the whole planar fan wave and terminates at the point $\mathrm{B}(\xi_{_\mathrm{B}}, 0)$. Moreover, along this $C_{+}$ cross characteristic curve $\widetilde{\mathrm{PB}}$, we have
\begin{equation}\label{7201a}
2\Theta(\tau)-\frac{\pi}{2}<\alpha\leq \alpha_2(\mathrm{P}).
\end{equation}
Through $\mathrm{P}$, we draw a straight $C_{-}$ characteristic line of the constant state $(u_{w}, v_{w}, \tau_{w})$. The characteristic angle of this $C_{-}$ straight characteristic is $\beta_3(\mathrm{P})$. For convenience, we denote this straight $C_{-}$ characteristic line by $C_{-}^{\mathrm{P}}$; see Figure \ref{Fig7.1}.
The following two cases may occur:
\begin{description}
  \item[(i)] If $\theta_w-\pi-\beta_3(\mathrm{P})<\frac{\pi}{2}$, then $C_{-}^{\mathrm{P}}$ intersects $\partial\mathcal{W}_{+}$ at a point $\mathrm{G}$ before reaching the sonic circle $(\xi-u_{w})^2+(\eta-v_{w})^2=c_{w}^2$; see Figure \ref{Fig7.1} (left).
Moreover,
\begin{equation}\label{8250d}
\alpha_3(\mathrm{P})\leq \alpha\leq 2(\theta_w-\pi)-\beta_3(\mathrm{P})\quad \text{on}\quad \overline{\mathrm{PG}}.
\end{equation}
  \item[(ii)] If $\theta_w-\pi-\beta_3(\mathrm{P})>\frac{\pi}{2}$, then
$C_{-}^{\mathrm{P}}$ is tangent to the sonic circle of the constant state $(u_w, v_w, \tau_w)$ at a point
 $\mathrm{G}$ before reaching the wall $\partial{\mathcal{W}}_{+}$; see Figure \ref{Fig7.1} (right).
\end{description}

Assume that
$$
\sigma_2(\mathrm{P})-\sigma_3(\mathrm{P})-\int_{q_3(\mathrm{P})}
^{q_2(\mathrm{P})}\frac{\sqrt{q^2-c^2(\grave{\tau}(q))}}{qc(\grave{\tau}(q))}{\rm d}q>0,
$$
where
%$q_2(\mathrm{P})=\sqrt{(u_{-}(\xi_{_\mathrm{P}})-\xi_{_\mathrm{P}})^2+\eta_{_\mathrm{P}}^2}$, $q_3(\mathrm{P})=\sqrt{(u_w-\xi_{_\mathrm{P}})^2+(v_w-\eta_{_\mathrm{P}})^2}$, and
 the function $\grave{\tau}(q)$ is determined by
$$
\frac{1}{2}q^2+\int_{\tau_{po}}^{\grave{\tau}(q)}\tau p'(\tau)~{\rm d}\tau=\frac{1}{2}q_2^2(\mathrm{P}).
$$
Then, as in Section 5.3,
there exist two centered waves of different types emanating from $\mathrm{P}$. To determine the principal parts of the centered waves, we consider (\ref{principal})--(\ref{principal1}) with the initial data
\begin{equation}\label{62802a}
(\widetilde{u}_{+},\widetilde{v}_{+}, \widetilde{\rho}_{+})(\nu_2)=(u_{-}(\xi_{_P}), 0, \rho_{po}),
\end{equation}
and (\ref{82306})--(\ref{82307}) with the initial data
\begin{equation}\label{62803a}
(\widetilde{u}_{-},\widetilde{v}_{-}, \widetilde{\rho}_{-})(\mu_1)=(u_{w}, v_{w}, \rho_{w}),
\end{equation}
where $\nu_2=\tan(\alpha_2(P))$ and $\mu_1=\cot(\beta_3(\mathrm{P}))$.

Assume further that
\begin{equation}\label{8260d}
2\int_{q_2(\mathrm{P})}^{q_{_{lim}}}\frac{\sqrt{q^2-c^2(\grave{\tau}(q))}}{qc(\grave{\tau}(q))}{\rm d}q>
\sigma_2(\mathrm{P})-\sigma_3(\mathrm{P})-\int_{q_3(\mathrm{P})}
^{q_2(\mathrm{P})}\frac{\sqrt{q^2-c^2(\grave{\tau}(q))}}{qc(\grave{\tau}(q))}{\rm d}q,
\end{equation}
where  $q_{_{lim}}=\sqrt{q_2^2(\mathrm{P})-2\int_{\tau_{po}}^{+\infty}\tau p'(\tau)~{\rm d}\tau}$.
Then the nonlinear system
$$
\left\{
  \begin{array}{ll}
    \displaystyle\widetilde{\sigma}+ \int_{q_{2}(\mathrm{P})}^{\widetilde{q}}\frac{\sqrt{q^2-c^2(\grave{\tau}(q))}}{qc(\grave{\tau}(q))}{\rm d}q=\sigma_2(\mathrm{P}),\\[12pt]
   \displaystyle \widetilde{\sigma}- \int_{q_{2}(\mathrm{P})}^{\widetilde{q}}
   \frac{\sqrt{q^2-c^2(\grave{\tau}(q))}}{qc(\grave{\tau}(q))}{\rm d}q=\sigma_3(\mathrm{P})+\int_{q_3(\mathrm{P})}^{q_2(\mathrm{P})}\frac{\sqrt{q^2-c^2(\grave{\tau}(q))}}{qc(\grave{\tau}(q))}{\rm d}q
  \end{array}
\right.
$$
 admits a unique solution $(\widetilde{\sigma}, \widetilde{q})=(\sigma_m, q_m)$, with $\widetilde{q}>q_2(\mathrm{P})$.

%=============
%==========

Let
 $$\nu_1=\tan(\sigma_m+A_m) \quad \text{and}\quad \mu_2=\cot(\sigma_m-A_m),$$
where $A_m=\arcsin\Big(\frac{c(\tau_m)}{q_m}\Big)$ and $\tau_m=\grave{\tau}(q_m)$.
 As in Lemma \ref{8250c}, the initial value problem (\ref{principal})--(\ref{principal1}), (\ref{62802a}) admits a solution on $[\nu_1, \nu_2]$, and the initial value problem (\ref{82306})--(\ref{82307}), (\ref{62803a}) admits a solution on $[\mu_2, \mu_1]$, The solutions satisfy
$
(\widetilde{u}_{-},\widetilde{v}_{-}, \widetilde{\rho}_{-})(\mu_2)=(\widetilde{u}_{+},\widetilde{v}_{+}, \widetilde{\rho}_{+})(\nu_1).
$
Moreover,
$\widetilde{\tau}_{+}'(\nu)<0$ for $\nu_1\leq\nu\leq \nu_2$;
$\widetilde{\tau}_{-}'(\mu)<0$ for $\mu_2\leq\nu\leq \mu_1$.

Let $\mu=\widetilde{\mu}_{-}(\tau)$, $\tau\in[\tau_{w}, \tau_m]$ be the inverse function of $\tau=\widetilde{\tau}_{-}(\mu)$,  $\mu\in[\mu_2, \mu_1]$,
and let
$\nu=\widetilde{\nu}_{+}(\tau)$, $\tau\in[\tau_{po}, \tau_m]$ be the inverse function of $\tau=\widetilde{\tau}_{+}(\nu)$, $\nu\in [\nu_1, \nu_2]$.
%at the function $\tau=\widetilde{\tau}_{-}(\mu)$ has an inverse function $\mu=\widetilde{\mu}_{-}(\tau)$.
Let  $$(\check{\alpha}_{-}(\tau), \check{\beta}_{-}(\tau))=\big(\widetilde{\alpha}_{-}(\widetilde{\mu}_{-}(\tau)), \widetilde{\beta}_{-}(\widetilde{\mu}_{-}(\tau))\big), \quad \tau\in[\tau_{w}, \tau_m]$$ and
  $$(\check{\alpha}_{+}(\tau), \check{\beta}_{+}(\tau))=\big(\widetilde{\alpha}_{+}(\widetilde{\nu}_{+}(\tau)), \widetilde{\beta}_{+}(\widetilde{\nu}_{+}(\tau))\big), \quad \tau\in[\tau_{po}, \tau_m].$$

We first consider (\ref{PEC})
with the boundary conditions
\begin{equation}\label{93001c}
\left\{
               \begin{array}{ll}
                 (u, v, \tau)=(u_{-}, v_{-}, \tau_{-})(\xi), & \hbox{$(\xi, \eta)\in\widetilde{\mathrm{PB}}$;} \\[2pt]
                 (u, v, \tau)=(u_w, v_w, \tau_w), & \hbox{$(\xi, \eta)\in\overline{\mathrm{PG}}$.}\\
                % v=u\tan(\theta_w-\pi), & \hbox{$(\xi, \eta)\in\mathcal{W}_{+}$;}
               \end{array}
             \right.
\end{equation}
To obtain a global classical solution to the BVP (\ref{PEC}), (\ref{93001c}),
we impose the following assumptions:
\begin{description}
%\item[($h_{41}$)]
%$\alpha_2(\mathrm{P})-\beta_3(\mathrm{P})<\pi$,  $\alpha_3(\mathrm{P})-\beta_3(\mathrm{P})>2\Theta(\tau_w)$,  $\alpha_2(\mathrm{P})+\frac{\pi}{2}>2\Theta(\tau_{po})$;
%\vskip 2pt
%\item[($h_{42}$)] $\sigma_2(\mathrm{P})>\sigma_3(\mathrm{P})+{\displaystyle\int_{q_3(\mathrm{P})}
%^{q_2(\mathrm{P})}}\frac{\sqrt{q^2-c^2(\grave{\tau}(q))}}{qc(\grave{\tau}(q))}{\rm d}q$;
%\vskip 2pt
\item[($h_{41}$)] $\check{\alpha}_{+}(\tau)+\frac{\pi}{2}>2\Theta(\tau)$ for  $\tau\in[\tau_{po}, \tau_m]$;
\item[($h_{42}$)] $\sigma_2(\mathrm{P})>\sigma_3(\mathrm{P})+{\displaystyle\int_{q_3(\mathrm{P})}
^{q_2(\mathrm{P})}}\frac{\sqrt{q^2-c^2(\grave{\tau}(q))}}{qc(\grave{\tau}(q))}{\rm d}q$.
%\vskip 4pt
 %  \item[($h_{44}$)] $\check{\beta}_{-}(\tau)<\alpha_2(\mathrm{P})-2\Theta(\tau)$ for  $\tau\in[\tau_{w}, \tau_m]$.
 % \item[($h_{45}$)] $\beta_3(\mathrm{P})+2\Theta(\tau_{po})<\theta_w-\pi<\frac{\pi}{2}+\beta_3(\mathrm{P})$.
\end{description}

One main result of this section is stated in the following theorem, which solves the interaction between the planar shock–fan composite wave with the oblique shock for $\tau_{w}<\tau_2^i$.
\begin{thm}\label{1001ba}
Under assumptions ($h_{40}$)--($h_{42}$), the BVP (\ref{PEC}), (\ref{93001c}) admits a piecewise smooth solution in the domain bounded by $\widetilde{\mathrm{PB}}$,  $\overline{\mathrm{PG}}$, $\widetilde{\mathrm{GF}}$, and $\widetilde{\mathrm{BF}}$, where $\widetilde{\mathrm{GF}}$ is a forward $C_{+}$ characteristic curve emanating from $\mathrm{G}$ and $\widetilde{\mathrm{BF}}$ is a vacuum boundary.
\end{thm}

Let $\varepsilon_1$ be a small positive constant such that
$$
\alpha_2(\mathrm{P})+\frac{\pi}{2}+4\varepsilon_1<\pi.
$$
For $\tau\geq \tau_{po}$, we define
$$
\Gamma_{41}(\tau)=\Big\{(\alpha, \beta)\mid {\alpha}_{4l}(\tau)<\alpha<{\alpha}_{4r},~{\beta}_{4l}<\beta<{\beta}_{4r}(\tau),~\alpha-\beta>\delta_1\Big\},
$$
where
$$
\begin{aligned}
{\alpha}_{4l}(\tau)=2\Theta(\tau)+\beta_{4l}, \quad {\beta}_{4l}= -\frac{\pi}{2}-\varepsilon_1, \quad {\alpha}_{4r}=\alpha_2(\mathrm{P})+\varepsilon_1, \quad {\beta}_{4r}(\tau)={\alpha}_{4r}-2\Theta(\tau),
\end{aligned}
$$
and $\delta_1$ is a small positive constant.
%From assumption (${h}_{45}$), one has
%\begin{equation}
%{\alpha}_{4l}(\tau)<2(\theta_w-\pi)-{\beta}_{4l}<{\alpha}_{4r}\quad
%\text{and}\quad
%{\alpha}_{4l}(\tau)<\theta_w-\pi<{\beta}_{4r}(\tau)
%\end{equation}
%for any $\tau\geq \tau_w$.

As in Section 4,
under assumption ($h_{41}$),
 we have the following estimates:
\begin{description}
  \item[$(e1)$]$(\alpha, \beta)\in \Gamma_{41}(\tau)$ on $\widetilde{\mathrm{PB}}$;
  %\item[$(e2)$]$(\alpha, \beta)\in \Gamma_{41}(\tau)$ on $\overline{\mathrm{PG}}$;
  \item[$(e2)$]$(\check{\alpha}_{+},\check{\beta}_{+})(\tau)\in \Gamma_{41}(\tau)$ for $\tau\in [\tau_{po}, \tau_m]$.
 % \item[$(e3)$]$(\check{\alpha}_{-},\check{\beta}_{-})(\tau)\in \Gamma_{41}(\tau)$ for $\tau\in [\tau_{w}, \tau_m]$.
 \end{description}

%In what follows, we are going to construct a piecewise smooth solution to the BVP (\ref{PEC}), (\ref{93001a}) piece by piece.

Applying the hodograph transformation method and using ($e_1$) and  ($e_2$),  we deduce that there exists a  $C_{+}$ centered wave emanating from $\mathrm{P}$ on a curvilinear triangular domain $\Sigma_4$ bounded by $\widetilde{\mathrm{PB}}$, $\widetilde{\mathrm{PF}}$, and a vacuum boundary $\widetilde{\mathrm{BF}}$, where $\widetilde{\mathrm{PF}}$ is a $C_{+}$ characteristic curve emanating from $\mathrm{P}$ with the slope $\nu_1$ at $\mathrm{P}$. The principal part of the centered wave is $(\widetilde{u}_{+},\widetilde{v}_{+}, \widetilde{\rho}_{+})(\nu)$, $\nu_1\leq \nu\leq \nu_2$. Moreover, the solution satisfies
$$
\bar{\partial}_{+}\rho<0\quad \text{and}\quad (\alpha, \beta)\in \Gamma_{41}(\tau)\quad \text{in}\quad \overline{\Sigma_4}\setminus\{\mathrm{P}\}.
$$
For convenience, we denote this $C_{+}$ centered wave solution  by $(u, v, \tau)=(u_4, v_4, \tau_4)(\xi, \eta)$, $(\xi, \eta)\in \overline{\Sigma_4}\setminus\{\mathrm{P}\}$.

%\subsubsection{\bf Centered simple wave emanating from $P$}
By Lemma \ref{simple}, there exists a
centered simple wave with straight $C_{-}$ characteristics emanating from $\mathrm{P}$. The centered simple wave flow region $\Sigma_5$ is bounded by $\overline{\mathrm{PG}}$, $\overline{\mathrm{PL}}$, and $\widetilde{\mathrm{GL}}$, where $\overline{\mathrm{PL}}$ is a straight $C_{-}$ characteristic line emanating from $\mathrm{P}$ with the slope $\mu_2^{-1}$ at $\mathrm{P}$ and $\widetilde{\mathrm{GL}}$ is a $C_{+}$ cross characteristic curve from $\mathrm{G}$. Moreover,  this centered simple wave satisfies
$$
\bar{\partial}_{+}\rho<0\quad  \text{in}\quad \overline{\Sigma_5}\setminus\{\mathrm{P}\}.
$$
We denote the centered simple wave by $(u, v, \tau)=(u_5, v_5, \tau_5)(\xi, \eta)$, $(\xi, \eta)\in \overline{\Sigma_5}\setminus\{\mathrm{P}\}$.

%\subsubsection{\bf Non-planar simple wave emanating from $\widetilde{PF}$}
We then consider (\ref{PEC}) with the boundary conditions
\begin{equation}\label{62805}
\left\{
  \begin{array}{ll}
(u, v, \tau)=(u_{4}, v_{4}, \tau_{4})(\xi, \eta), & \hbox{$(\xi, \eta)\in\widetilde{\mathrm{PF}}$;} \\[2pt]
    (u, v, \tau)=(u_{5}, v_{5}, \tau_{5})(\xi, \eta), & \hbox{$(\xi, \eta)\in\overline{\mathrm{PL}}$.}
% (u, v, \tau)=  (u_{w}, v_{w}, \tau_{w}), & \hbox{$(\xi, \eta)\in\widetilde{\mathrm{EG}}$;} \\[3pt]
   % v=u\tan(\theta_w-\pi), & \hbox{$(\xi, \eta)\in\mathcal{W}_{-}$.}
  \end{array}
\right.
\end{equation}
From Lemma \ref{simple},
the problem  (\ref{PEC}), (\ref{62805}) admits a simple wave solution on a triangular region $\overline{\Sigma_6}$ bounded by $\overline{\mathrm{PF}}$, $\overline{\mathrm{PL}}$, and $\overline{\mathrm{LF}}$, where $\overline{\mathrm{LF}}$ is a $C_{+}$ characteristic curve emanating from $\mathrm{L}$. Moreover, the solution satisfies
$$
\bar{\partial}_{+}\rho<0\quad  \text{in}\quad \overline{\Sigma_6}.
$$
We denote the simple wave by $(u, v, \tau)=(u_6, v_6, \tau_6)(\xi, \eta)$, $(\xi, \eta)\in \overline{\Sigma_6}$.
This completes the proof of Theorem \ref{1001ba}.
Therefore, when $\tau_w>\tau_2^i$, the interaction between the shock–fan composite wave and the oblique shock wave generates a non-planar simple wave emanating from the interaction region. %This simple wave will be reflected by the solid wall.

%\subsubsection{Reflection of the non-planar simple wave at the ramp}
If $\theta_w-\pi-\beta_3(\mathrm{P})<\frac{\pi}{2}$,
then the non-planar simple wave emanating from the interaction region will be reflected at the ramp. To construct a solution of the wave reflection problem,
we consider  (\ref{PEC})
with the boundary conditions
\begin{equation}\label{91501d}
\left\{
  \begin{array}{ll}
  (u, v, \tau)=(u_{-}, v_{-}, \tau_{-})(\xi), & \hbox{$(\xi, \eta)\in\widetilde{\mathrm{PB}}$,} \\[2pt]
                 (u, v, \tau)=(u_w, v_w, \tau_w), & \hbox{$(\xi, \eta)\in\overline{\mathrm{PG}}$,}\\[2pt]
    %(u, v, \tau)=(u_{b}, v_{b}, \tau_{b}), & \hbox{$(\xi, \eta)\in\overline{\mathrm{PG}}$.}
% (u, v, \tau)=  (u_{w}, v_{w}, \tau_{w}), & \hbox{$(\xi, \eta)\in\widetilde{\mathrm{EG}}$;} \\[3pt]
    v=u\tan(\theta_w-\pi), & \hbox{$(\xi, \eta)\in\partial\mathcal{W}_{+}$.}
  \end{array}
\right.
\end{equation}

To control the hyperbolicity of the system (\ref{PEC}), we make the following assumptions:
\begin{description}
\item[($h_{43}$)]
$\max\big\{\alpha_2(\mathrm{P}), \alpha_3(\mathrm{P}), 2(\theta_w-\pi)-\beta_3(\mathrm{P}), 2(\theta_w-\pi)+\frac{\pi}{2}\big\}-\min\big\{\beta_3(\mathrm{P}), -\frac{\pi}{2}\big\}<\pi$;
\item[($h_{44}$)] $\check{\alpha}_{+}(\tau)-\min\big\{\beta_3(\mathrm{P}), -\frac{\pi}{2}\big\}>2\Theta(\tau)$ for  $\tau\in[\tau_{po}, \tau_m]$;
   \item[($h_{45}$)] $\check{\beta}_{-}(\tau)<\max\big\{\alpha_2(\mathrm{P}), \alpha_3(\mathrm{P}), 2(\theta_w-\pi)-\beta_3(\mathrm{P}), 2(\theta_w-\pi)+\frac{\pi}{2}\big\}-2\Theta(\tau)$ for  $\tau\in[\tau_{w}, \tau_m]$;
  \item[($h_{46}$)] $\min\big\{\beta_3(\mathrm{P}), -\frac{\pi}{2}\big\}+2\Theta(\tau_{po})<\theta_w-\pi<\frac{\pi}{2}+\min\big\{\beta_3(\mathrm{P}), -\frac{\pi}{2}\big\}$.
\end{description}

Let $\varepsilon_2$ be a small positive constant such that
$$
\max\Big\{\alpha_2(\mathrm{P}), \alpha_3(\mathrm{P}), 2(\theta_w-\pi)-\beta_3(\mathrm{P}), 2(\theta_w-\pi)+\frac{\pi}{2}\Big\}-\min\Big\{\beta_3(\mathrm{P}), -\frac{\pi}{2}\Big\}+4\varepsilon_2<\pi.
$$
For $\tau\geq \tau_w$, we define
$$
\Gamma_{42}(\tau)=\Big\{(\alpha, \beta)\mid \bar{\alpha}_{4l}(\tau)<\alpha<\bar{\alpha}_{4r},~\bar{\beta}_{4l}<\beta<\bar{\beta}_{4r}(\tau),~\alpha-\beta>\delta_2\Big\},
$$
where
$$
\begin{aligned}
&\qquad \bar{\alpha}_{4r}=\max\Big\{\alpha_2(\mathrm{P}), \alpha_3(\mathrm{P}), 2(\theta_w-\pi)-\beta_3(\mathrm{P}), 2(\theta_w-\pi)+\frac{\pi}{2}\Big\}+\varepsilon_2, \\& \bar{\alpha}_{4l}(\tau)=2\Theta(\tau)+\beta_{4l}, \quad \bar{\beta}_{4l}=\min\Big\{\beta_3(\mathrm{P}), -\frac{\pi}{2}\Big\}-\varepsilon_2,\quad \bar{\beta}_{4r}(\tau)=\bar{\alpha}_{4r}-2\Theta(\tau),
\end{aligned}
$$
and $\delta_2$ is a small positive constant.
From assumption ($h_{46}$), we have
$$
\bar{\alpha}_{4l}(\tau)<2(\theta_w-\pi)-\bar{\beta}_{4l}<\bar{\alpha}_{4r}\quad
\text{and}\quad
\bar{\alpha}_{4l}(\tau)<\theta_w-\pi<\bar{\beta}_{4r}(\tau)
$$
for any $\tau>\tau_d$.

%As in the previous discussion, we have that the solution of the  BVP (\ref{PEC}), (\ref{93001c}) satisfies
%\begin{equation}
%\bar{\partial}_{+}\rho<0\quad \text{and}\quad (\alpha, \beta)\in \Gamma_{42}(\tau)\quad \text{on}\quad \widetilde{\mathrm{GF}}.
%\end{equation}

As in Lemma \ref{lem419},
  %For convenience, we denote the solution of the problem (\ref{PEC}), (\ref{91501a}) by
 %$(u, v, \tau)=(u_7, v_7, \tau_7)(\xi, \eta)$, $(\xi, \eta)\in \overline{\Sigma_7}$.
under assumptions ($h_{40}$) and ($h_{42}$)--($h_{46}$),  the mixed boundary value problem (\ref{PEC}), (\ref{91501d}) admits a classical solution on a curvilinear
quadrilateral domain  bounded by $\widetilde{\mathrm{PB}}$, $\overline{\mathrm{PG}}$,  $\partial\mathcal{W}_{+}$, and a vacuum boundary $\widetilde{\mathrm{BH}}$.
Moreover, the solution satisfies $\bar{\partial}_{\pm}\rho<0$ and $(\alpha, \beta)\in \Gamma_{42}(\tau)$.

%Therefore, we obtain a global piecewise smooth solution to the boundary value problem (\ref{PEC}), (\ref{93001a}).
%This completes the proof of Theorem \ref{1001b}.

\begin{rem}
From Remark \ref{26102201}, we see that if the oblique shock is merely a post-sonic shock, then 
$$
\theta_w-\pi=\arctan\left(\frac{v_w}{u_w}\right)>\frac{\pi}{2}+\beta_3(\mathrm{P}).
$$
This implies that case (ii) can still occur in some cases even when the oblique shock is a transonic shock.
%We guess that in this case, the point $\mathrm{G}$ and the ramp are connected by an arc $\wideparen{\mathrm{GH}}$ of the sonic circle $(\xi-u_{w})^2+(\eta-v_{w})^2=c_{w}^2$, as shown in Figure \ref{Fig7.1} (right),  and the flow across the arc may be subsonic. 
The present approach does not apply to this problem, and even the reflection configuration cannot be readily conjectured; we leave it open.
%Numerical simulations can be used to conjecture the Mach-like reflection configurations and to guide subsequent rigorous analysis.
\end{rem}

\subsection{Wave interaction for the case $\tau_{w}<\tau_2^i$}
%\subsubsection{Main theorem}
If $\tau_{w}<\tau_2^i$,
the interaction generates a new straight post-sonic 1-shock propagating from the point $\mathrm{P}$.
The state of the flow upstream of the shock is $(u_{w}, v_{w}, \tau_{w})$ and the specific volume of the flow downstream of the shock is $\tau_d=\psi_{po}(\tau_{w})$. We use $\phi_{_\mathrm{P}}$ to denote the inclination angle and $(u_d, v_d)$ the downstream flow velocity of the post-sonic shock.
%The ray $\xi=\xi_{_P}+r\cos\phi_{w}$, $\eta=\eta_{_P}+r\sin\phi_{w}$ ($r>0$) intersects $\partial{\mathcal{W}}_{+}$ at some point $G$.

Let
\begin{equation}\label{8260a}
\begin{aligned}
&(U_d, V_d)=(u_d-\xi_{_\mathrm{P}}, v_d-\eta_{_\mathrm{P}}), \quad (U_{w}, V_{w})=(u_{w}-\xi_{_\mathrm{P}}, v_{w}-\eta_{_\mathrm{P}}), \quad
\\& \quad N_d=V_d\cos\phi_{_\mathrm{P}}-U_d\sin\phi_{_\mathrm{P}},\quad L_d=U_d\cos\phi_{_\mathrm{P}}+V_d\sin\phi_{_\mathrm{P}},\\& N_{w}=V_{w}\cos\phi_{_\mathrm{P}}-U_{w}\sin\phi_{_\mathrm{P}},\quad\text{and} \quad L_{w}=U_{w}\cos\phi_{_\mathrm{P}}+V_{w}\sin\phi_{_\mathrm{P}}.
\end{aligned}
\end{equation}
Since the shock is post-sonic, we have
$
 N_d=c_d=\sqrt{-\tau_d^2p'(\tau_d)}.
$.
Consequently, 
%$$
%\rho_dN_d=\rho_{w}N_{w},\quad \rho_{w}=\tau_{w}^{-1}, \quad N_d=c_d=\sqrt{-\tau_d^2p'(\tau_d)}.
%$$
%Thus, we have
\begin{equation}\label{8260b}
N_{w}=\frac{\rho_d c_d}{\rho_w}, \quad L_w=L_d=\sqrt{q_3^2(\mathrm{P})-N_w^2}, \quad \text{and} \quad \phi_{_\mathrm{P}}=\sigma_3(\mathrm{P})-\arcsin\left(\frac{N_w}{q_3(\mathrm{P})}\right).
\end{equation}
From (\ref{8260a}) and (\ref{8260b}), we obtain $u_d$ and $v_d$. Hence, the oblique shock front emanating from $\mathrm{P}$ and the downstream state are determined.

%=============
\begin{figure}[htbp]
\begin{center}
\includegraphics[scale=0.45]{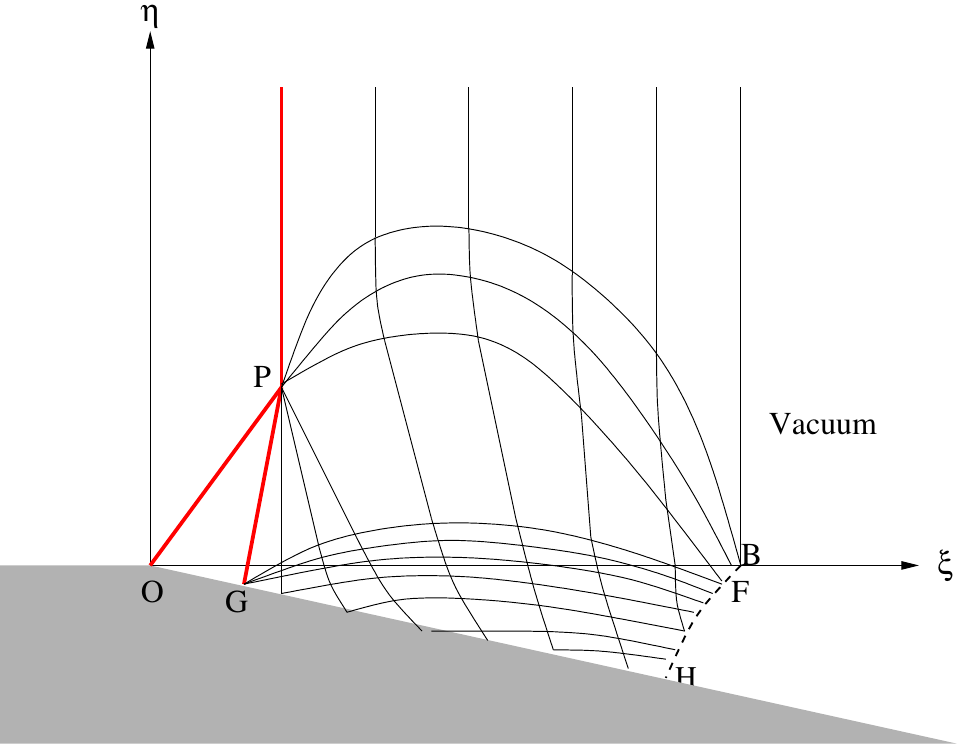}\quad\qquad\includegraphics[scale=0.45]{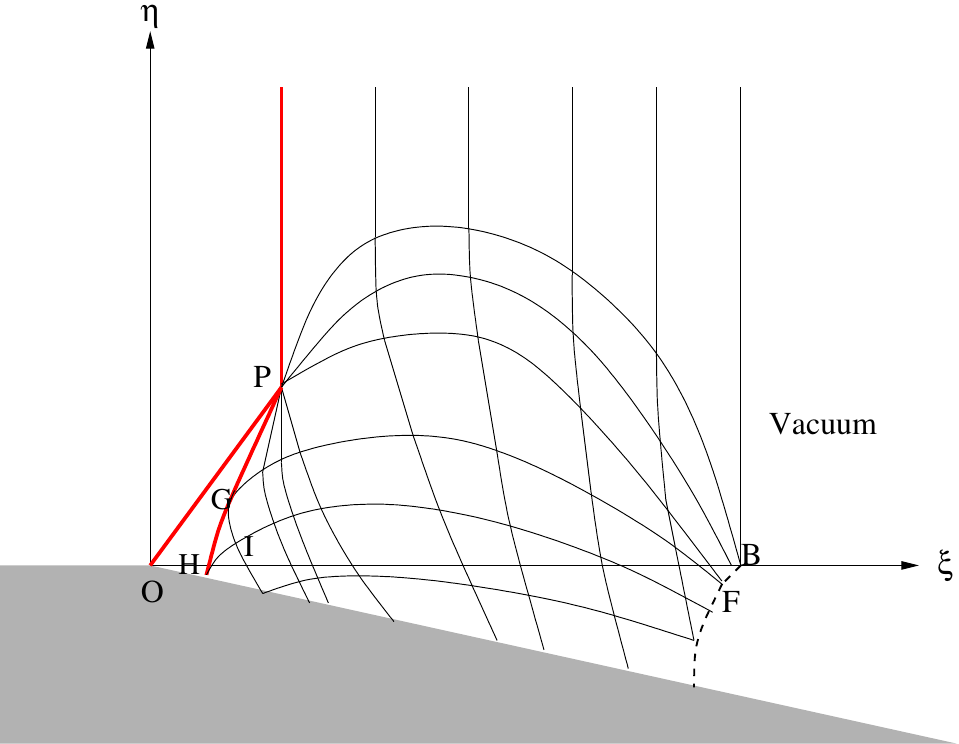}
\caption{\footnotesize Interaction between a planar shock–fan composite wave and an oblique shock wave for $\tau_{w}<\tau_2^i$. Left: Regular-like reflection; right: Mach-like reflection.}
\label{Fig7.2}
\end{center}
\end{figure}
%==========

%We have
%$$
%\beta_3(P)>\phi_{_\mathrm{P}}=\beta_5(P).
%$$

Let $q_d=\sqrt{(u_d-\xi_{_\mathrm{P}})^2+(v_d-\eta_{_\mathrm{P}})^2}$ and $\sigma_d=\arctan\big(\frac{v_d-\eta_{_\mathrm{P}}}{u_d-\xi_{_\mathrm{P}}}\big)$.
When
$$
\sigma_2(\mathrm{P})
-\sigma_d-\int_{q_2(\mathrm{P})}^{q_d}\frac{\sqrt{q^2-c^2(\grave{\tau}(q))}}{qc(\grave{\tau}(q))}{\rm d}q>0,
$$
there exists a centered simple waves with straight $C_{-}$ characteristics and a $C_{+}$ type centered wave emanating from the point $\mathrm{P}$. To determine the principal parts of the centered waves, we consider (\ref{principal})--(\ref{principal1}) with  initial data (\ref{62802a})
and (\ref{82306})--(\ref{82307}) with  initial data
\begin{equation}\label{62810a}
(\widetilde{u}_{-},\widetilde{v}_{-}, \widetilde{\rho}_{-})(\bar{\mu}_1)=(u_{d}, v_{d}, \rho_{d}),
\end{equation}
where $\bar{\mu}_1=\cot\phi_{_\mathrm{P}}$.

Assume that
$$
2\int_{q_{2}(\mathrm{P})}^{q_{_{lim}}}\frac{\sqrt{q^2-c^2(\grave{\tau}(q))}}{qc(\grave{\tau}(q))}{\rm d}q>\sigma_2(\mathrm{P})
-\sigma_d-\int_{q_2(\mathrm{P})}^{q_d}\frac{\sqrt{q^2-c^2(\grave{\tau}(q))}}{qc(\grave{\tau}(q))}{\rm d}q.
$$
Then the system
$$
\left\{
  \begin{array}{ll}
    \displaystyle\widetilde{\sigma}+ \int_{q_{2}(\mathrm{P})}^{\widetilde{q}}\frac{\sqrt{q^2-c^2(\grave{\tau}(q))}}{qc(\grave{\tau}(q))}{\rm d}q=\sigma_2(\mathrm{P}),\\[12pt]
   \displaystyle \widetilde{\sigma}- \int_{q_d}^{\widetilde{q}}\frac{\sqrt{q^2-c^2(\grave{\tau}(q))}}{qc(\grave{\tau}(q))}{\rm d}q
   =\sigma_d+\int_{q_3(\mathrm{P})}^{q_d}\frac{\sqrt{q^2-c^2(\grave{\tau}(q))}}{qc(\grave{\tau}(q))}{\rm d}q
  \end{array}
\right.
$$
 admits a unique solution $(\widetilde{\sigma}, \widetilde{q})=(\bar{\sigma}_m, \bar{q}_m)$ satisfying $\widetilde{q}>q_2(\mathrm{P})$.

% We define the following constants:
%$$
%\begin{aligned}
%&\sigma_m=\frac{\sigma_2(\mathrm{P})+\sigma_d}{2}+\frac{1}{2}\int_{q_d}^{q_2(\mathrm{P})}\frac{\sqrt{q^2-c^2(\grave{\tau}(q))}}{qc(\grave{\tau}(q))}{\rm d}q, \quad \tau_m=\grave{\tau}(q_m),  \\ &A_m=\arcsin\Big(\frac{c(\tau_m)}{q_m}\Big), \quad \alpha_m=\sigma_m+A_m, \quad \text{and}\quad \beta_m=\sigma_m-A_m.
%\end{aligned}
%$$
We define
$$\bar{\nu}_1=\tan(\bar{\sigma}_m+\bar{A}_m) \quad \text{and}\quad \bar{\mu}_2=\cot(\bar{\sigma}_m-\bar{A}_m),$$
where $\bar{A}_m=\arcsin\big(\frac{c(\bar{\tau}_m)}{\bar{q}_m}\big)$ and $\bar{\tau}_m=\grave{\tau}(\bar{q}_m)$.

The initial value problem (\ref{principal})--(\ref{principal1}), (\ref{62802a}) admits a solution on $[\bar{\nu}_1, \nu_2]$, and the initial value problem (\ref{82306})--(\ref{82307}), (\ref{62810a}) admits a solution on $[\bar{\mu}_2, \bar{\mu}_1]$.
Moreover, $
(\widetilde{u}_{-},\widetilde{v}_{-}, \widetilde{\rho}_{-})(\bar{\mu}_2)=(\widetilde{u}_{+},\widetilde{v}_{+}, \widetilde{\rho}_{+})(\bar{\nu}_1)
$;
$\widetilde{\tau}_{+}'(\nu)<0$ for $\bar{\nu}_1\leq\nu\leq \nu_2$;
$\widetilde{\tau}_{-}'(\mu)<0$ for $\bar{\mu}_2\leq\mu\leq \bar{\mu}_1$.
Let $\mu=\widetilde{\mu}_{-}(\tau)$, $\tau\in[\tau_{d}, \bar{\tau}_m]$ be the inverse function of $\tau=\widetilde{\tau}_{-}(\mu)$,  $\mu\in [\bar{\mu}_2, \bar{\mu}_1]$ and
$\nu=\widetilde{\nu}_{+}(\tau)$, $\tau\in[\tau_{po}, \bar{\tau}_m]$ the inverse function of $\tau=\widetilde{\tau}_{+}(\nu)$, $\nu\in[\bar{\nu}_1, \nu_2]$.
%at the function $\tau=\widetilde{\tau}_{-}(\mu)$ has an inverse function $\mu=\widetilde{\mu}_{-}(\tau)$.
Let  $(\check{\alpha}_{-}(\tau), \check{\beta}_{-}(\tau))=\big(\widetilde{\alpha}_{-}(\widetilde{\mu}_{-}(\tau)), \widetilde{\beta}_{-}(\widetilde{\mu}_{-}(\tau))\big)$ and
  $(\check{\alpha}_{+}(\tau), \check{\beta}_{+}(\tau))=\big(\widetilde{\alpha}_{+}(\widetilde{\nu}_{+}(\tau)), \widetilde{\beta}_{+}(\widetilde{\nu}_{+}(\tau))\big)$.

If $(\theta_w-\pi)-\phi_{_\mathrm{P}}<\frac{\pi}{2}$, then
the straight post-sonic shock propagating from $\mathrm{P}$ intersects $\partial{\mathcal{W}}_{+}$ at a point $\mathrm{G}$,
and
the state on the backside of the post-sonic shock is uniform hyperbolic; see Figure \ref{Fig7.2} (left).
If $(\theta_w-\pi)-\phi_{_\mathrm{P}}>\frac{\pi}{2}$, the straight post-sonic shock line is tangent to the sonic circle  $(\xi-u_d)^2+(\eta-v_d)^2=c^2(\tau_d)$ at a point
 $\mathrm{G}$ before reaching $\partial{\mathcal{W}}_{+}$; see Figure \ref{Fig7.2} (right).

We first consider (\ref{PEC})
with the boundary conditions
\begin{equation}\label{93001d}
\left\{
               \begin{array}{ll}
                 (u, v, \tau)=(u_{-}, v_{-}, \tau_{-})(\xi), & \hbox{$(\xi, \eta)\in\widetilde{\mathrm{PB}}$;} \\[2pt]
                 (u, v, \tau)=(u_d, v_d, \tau_d), & \hbox{$(\xi, \eta)\in\overline{\mathrm{PG}}$.}\\
                % v=u\tan(\theta_w-\pi), & \hbox{$(\xi, \eta)\in\mathcal{W}_{+}$;}
               \end{array}
             \right.
\end{equation}

To obtain a global classical solution to the BVP (\ref{PEC}), (\ref{93001d}),
we impose the following assumptions:
\begin{description}
%\item[($h_{41}$)]
%$\alpha_2(\mathrm{P})-\beta_3(\mathrm{P})<\pi$,  $\alpha_3(\mathrm{P})-\beta_3(\mathrm{P})>2\Theta(\tau_w)$,  $\alpha_2(\mathrm{P})+\frac{\pi}{2}>2\Theta(\tau_{po})$;
%\vskip 2pt
%\item[($h_{42}$)] $\sigma_2(\mathrm{P})>\sigma_3(\mathrm{P})+{\displaystyle\int_{q_3(\mathrm{P})}
%^{q_2(\mathrm{P})}}\frac{\sqrt{q^2-c^2(\grave{\tau}(q))}}{qc(\grave{\tau}(q))}{\rm d}q$;
%\vskip 2pt
\item[($h_{41}'$)] $\check{\alpha}_{+}(\tau)>\frac{\pi}{2}+2\Theta(\tau)$ for  $\tau\in[\tau_{po}, \tau_m]$;
\item[($h_{42}'$)] $\sigma_2(\mathrm{P})>\sigma_d+{\displaystyle\int_{q_2(\mathrm{P})}^{q_d}}\frac{\sqrt{q^2-c^2(\grave{\tau}(q))}}{qc(\grave{\tau}(q))}{\rm d}q$.
%\vskip 4pt
 %  \item[($h_{44}$)] $\check{\beta}_{-}(\tau)<\alpha_2(\mathrm{P})-2\Theta(\tau)$ for  $\tau\in[\tau_{w}, \tau_m]$.
 % \item[($h_{45}$)] $\beta_3(\mathrm{P})+2\Theta(\tau_{po})<\theta_w-\pi<\frac{\pi}{2}+\beta_3(\mathrm{P})$.
\end{description}
As in Theorem \ref{1001ba}, we have the following theorem, which solves the interaction of the planar shock–fan composite wave with the oblique shock for $\tau_{w}<\tau_2^i$.
\begin{thm}\label{1001baa}
Under assumptions ($h_{40}$), ($h_{41}'$)--($h_{42}'$), the boundary value problem (\ref{PEC}), (\ref{93001d}) admits a piecewise smooth solution in the curvilinear quadrilateral domain bounded by $\wideparen{\mathrm{PB}}$,  $\overline{\mathrm{PG}}$, $\widetilde{\mathrm{GF}}$, and $\widetilde{\mathrm{BF}}$, where $\widetilde{\mathrm{GF}}$ is a forward $C_{+}$ characteristic curve emanating from $\mathrm{G}$ and $\wideparen{\mathrm{BF}}$ is a vacuum boundary.
\end{thm}

When $\tau_w<\tau_2^i$, the interaction generates a non-planar shock–fan composite wave propagating from the interaction region. This shock–fan composite wave  will be reflected at the ramp.
For the needs of the subsequent discussion, we use subscript $4$ to denote the solution of the BVP (\ref{PEC}), (\ref{93001d}).

%\subsubsection{Reflection of the shock–fan composite wave at the ramp}
If $(\theta_w-\pi)-\phi_{_\mathrm{P}}<\frac{\pi}{2}$,
then we consider (\ref{PEC}) with the boundary conditions
\begin{equation}\label{8270a}
\left\{
  \begin{array}{ll}
 (u, v, \tau)=(u_{-}, v_{-}, \tau_{-})(\xi), & \hbox{$(\xi, \eta)\in\widetilde{\mathrm{PB}}$,} \\[2pt]
                 (u, v, \tau)=(u_d, v_d, \tau_d), & \hbox{$(\xi, \eta)\in\overline{\mathrm{PG}}$,}\\[2pt]
% (u, v, \tau)=  (u_{w}, v_{w}, \tau_{w}), & \hbox{$(\xi, \eta)\in\widetilde{\mathrm{EG}}$;} \\[3pt]
    v=u\tan(\theta_w-\pi), & \hbox{$(\xi, \eta)\in\partial\mathcal{W}_{+}$.}
  \end{array}
\right.
\end{equation}

Since $\phi_{_\mathrm{P}}=\arctan\big(\frac{v_d-v_w}{u_d-u_w}\big)$ and $\theta_w-\pi=\arctan\frac{v_w}{u_w}$, from $(\theta_w-\pi)-\phi_{_\mathrm{P}}<\frac{\pi}{2}$ we have
\begin{equation}\label{91601b}
\sigma_4(\mathrm{G})=\arctan\bigg(\frac{v_d-\eta_{_\mathrm{G}}}{u_d-\xi_{_\mathrm{G}}}\bigg)>\theta_w-\pi.
\end{equation}
Hence,
the boundary data (\ref{8270a}) is discontinuous at the point $\mathrm{G}$. From (\ref{91601b}), we see that there exists a $C_{+}$ centered wave emanating from the point $\mathrm{G}$.
To obtain the principal part of the centered wave, we consider (\ref{principal})--(\ref{principal1}) with the initial data
\begin{equation}\label{8270b}
(\widetilde{u}_{+},\widetilde{v}_{+}, \widetilde{\tau}_{+})(\hat{\nu}_2)=(u_d, v_d, \tau_d),
\end{equation}
where $\hat{\nu}_2=\tan\alpha_4(\mathrm{G})$.

There exists a $\hat{\nu}_1<\hat{\nu}_2$ such that the problem (\ref{principal})--(\ref{principal1}), (\ref{8270b}) admits a solution in $[\hat{\nu}_1, \hat{\nu}_2]$, and the solution satisfies
$\arctan\left(\frac{\tilde{v}_{+}(\hat{\nu}_1)}{\tilde{u}_{+}(\hat{\nu}_1)}\right)=\theta_w-\pi$ and $\widetilde{\tau}_{+}'(\nu)<0$ for $\hat{\nu}_1\leq\nu\leq \hat{\nu}_2$. Let
$\nu=\widetilde{\nu}_{+}(\tau)$, $\tau\in[\tau_{d}, \hat{\tau}_m]$ the inverse function of $\tau=\widetilde{\tau}_{+}(\nu)$, $\nu\in [\hat{\nu}_1, \hat{\nu}_2]$. Let $(\breve{\alpha}_{+}(\tau), \breve{\beta}_{+}(\tau))=\big(\widetilde{\alpha}_{+}(\widetilde{\nu}_{+}(\tau)), \widetilde{\beta}_{+}(\widetilde{\nu}_{+}(\tau))\big)$, $\tau\in [\tau_d, \hat{\tau}_m]$.

To control the hyperbolicity of the system for the BVP (\ref{Potential}), (\ref{8270a}),
we impose the following assumptions:
\begin{description}
\item[($h_{43}'$)] $\max\big\{\alpha_2(\mathrm{P}), \alpha_4(\mathrm{P}), 2(\theta-\pi)-\phi_{_\mathrm{P}}, 2(\theta-\pi)+\frac{\pi}{2}\big\}-\min\big\{\phi_{_\mathrm{P}}, -\frac{\pi}{2}\big\}<\pi$;
%\item[($h_{42}'$)] $\sigma_2(\mathrm{P})>\sigma_d-\int_{q_2(P)}^{q_d}\frac{\sqrt{q^2-c^2(\grave{\tau}(q))}}{qc(\grave{\tau}(q))}{\rm d}q$;
\item[($h_{44}'$)] $\check{\alpha}_{+}(\tau)>\min\big\{\phi_{_\mathrm{P}}, -\frac{\pi}{2}\big\}+2\Theta(\tau)$ for  $\tau\in[\tau_{po}, {\tau}_m]$;
   \item[($h_{45}'$)] $\check{\beta}_{-}(\tau)<\big\{\alpha_2(\mathrm{P}), \alpha_4(\mathrm{P}), 2(\theta-\pi)-\phi_{_\mathrm{P}}, 2(\theta-\pi)+\frac{\pi}{2}\big\}-2\Theta(\tau)$ for  $\tau\in[\tau_{d}, {\tau}_m]$;
  \item[($h_{46}'$)] $\min\big\{\phi_{_\mathrm{P}}, -\frac{\pi}{2}\big\}+2\Theta(\tau_{po})<\theta_w-\pi<\frac{\pi}{2}+\min\big\{\phi_{_\mathrm{P}}, -\frac{\pi}{2}\big\}$;
  \item[($h_{47}'$)] $\breve{\alpha}_{+}(\tau)>\min\big\{\phi_{_\mathrm{P}}, -\frac{\pi}{2}\big\}+2\Theta(\tau)$ for  $\tau\in[\tau_{d}, \hat{\tau}_m]$.
\end{description}
For $\tau\geq \tau_d$, we define 
$$
\Gamma_{43}(\tau)=\Big\{(\alpha, \beta)\mid \tilde{\alpha}_{4l}(\tau)<\alpha<\tilde{\alpha}_{4r},~\tilde{\beta}_{4l}<\beta<\tilde{\beta}_{4r}(\tau),~
\alpha-\beta>\delta_3\Big\},
$$
where
$$
\begin{aligned}
&\qquad\quad\tilde{\alpha}_{4r}=\max\Big\{\alpha_2(\mathrm{P}), \alpha_4(\mathrm{P}), 2(\theta-\pi)-\phi_{_\mathrm{P}}, 2(\theta-\pi)+\frac{\pi}{2}\Big\}+\varepsilon_3, \\&\tilde{\beta}_{4l}=\min\Big\{\phi_{_\mathrm{P}}, -\frac{\pi}{2}\Big\}-\epsilon_3\quad, \tilde{\alpha}_{4l}(\tau)=2\Theta(\tau)+\tilde{\beta}_{4l}, \quad \tilde{\beta}_{4r}(s)=\tilde{\alpha}_{4r}-2\Theta(\tau),
\end{aligned}
$$
and $\varepsilon_3$ and $\delta_3$ are two small positive constants.
From assumption ($h_{46}'$), one has
$$
\tilde{\alpha}_{4l}(\tau)<2(\theta_w-\pi)-\tilde{\beta}_{4l}<\tilde{\alpha}_{4r}\quad
\text{and}\quad
\tilde{\alpha}_{4l}(\tau)<\theta_w-\pi<\tilde{\beta}_{4r}(\tau)
$$
for any $\tau>\tau_d$.

%Assume $\tau_1^i\leq\tau_0<\tau_2^i$, $u_0>\frac{\rho_{po}c_{po}}{\rho_0}$, $\pi+\sigma_{po}<\theta_w<\pi$, and $\tau_0\leq \tau_{w}<\tau_2^i$.
As in Lemma \ref{lem419},
under assumptions ($h_{42}'$)--($h_{46}'$),
 the  BVP (\ref{Potential}), (\ref{8270a})  admits a piecewise smooth solution in a curvilinear
quadrilateral domain  bounded by $\widetilde{\mathrm{PB}}$, $\overline{\mathrm{PG}}$,  $\mathcal{W}_{-}$, and a vacuum boundary $\widetilde{\mathrm{BH}}$.
Moreover, the solution satisfies $\bar{\partial}_{\pm}\rho<0$ and $(\alpha, \beta)\in \Gamma_{43}(\tau)$.

\begin{rem}
If $(\theta_w-\pi)-\phi_{_\mathrm{P}}>\frac{\pi}{2}$, we can construct a circular post-sonic shock emanating from $\mathrm{G}$, and this shock intersect the ramp a point $\mathrm{H}$. The state on its front side is $(u_w, v_w, \tau_w)$, and the states $(u_b, v_b, \tau_b)(\xi,\eta)$ on its back side can be determined by Rankine-Hugoniot conditions.
Moreover, we have $(u_b-\xi)^2+(v_b-\xi)^2=-\tau_b^2p'(\tau_b)$ on $\widehat{\mathrm{GH}}$.
We consider (\ref{PEC}) with the data 
\begin{equation}\label{26102401}
  (u, v, \tau)=(u_b, v_b, \tau_b)(\xi, \eta), \quad (\xi, \eta)\in \widehat{\mathrm{GH}}.
\end{equation}
%The problem (\ref{PEC}), (\ref{26102401}) is a degenerate Cauchy problem, since the system (\ref{PEC}) is degenerate hyperbolic on 
%$\widehat{\mathrm{GH}}$. 
By the result of Chen and Lai \cite{ChenLai}, the Cauchy problem (\ref{PEC}), (\ref{26102401}) admits a sonic-supersonic solution\footnote{This sonic-supersonic flow resembles the flow downstream of a spherical Chapman-Jouguet detonation wave; see e.g., Landau and Lifschitz \cite{Lan} (Chapter 14).} on a curvilinear triangular domain bounded by $\widehat{\mathrm{GH}}$, $\widetilde{\mathrm{GI}}$, and $\widetilde{\mathrm{HI}}$,  where $\widetilde{\mathrm{GI}}$ is  a forward $C_{-}$ characteristic curve emanating from $G$ and $\widetilde{\mathrm{HI}}$ is a forward $C_{+}$ characteristic curve emanating from $H$; see Figure \ref{Fig7.2} (right). Moreover, the solution satisfies $\bar{\partial}_{\pm}\rho<0$ in this triangular domain. 
To construct the solution, we need to solve some standard Goursat problems. 
However, constructing pentagonal invariant regions of $(\alpha, \beta)$ to control the hyperbolicity of system (\ref{PEC}) does not work. 
Whether there exists another method to ensure the hyperbolicity remains an open question.
\end{rem}

\section{Interaction between  shock–fan and  fan–shock composite waves}
We consider the 2D Riemann problem (\ref{Potential}), (\ref{initial})--(\ref{boundary}) for $\tau_1^e\leq\tau_0<\tau_1^i$ and $u_0>\frac{\rho_{po}c_{po}}{\rho_0}$.
Then, away from the corner, the oncoming flow $(0)$ and the vacuum
are connected by a planar shock–fan composite wave which can be represented by (\ref{71501c}). In this section, we assume that at the corner the oncoming flow is deflected by an oblique fan–shock composite wave to a constant state $(u_{w}, v_{w}, \tau_{w})$,  where $\tau_w\in(\tau_1^i, \tau_2^e)$. The flow near the corner can be represented by (\ref{2662601}); see Figure \ref{Fig7.1}.
In the subsequent discussion of this section, the constants $\alpha_0$, $\tau_{po}$, $c_{po}$, $u_{w}$, $v_{w}$, $\tau_{w}$, $\phi_{pr}$, $\xi_{_\mathrm{B}}$, as well as the functions $(u_{-}, v_{-}, \tau_{-})(\xi)$, $(\bar{u}, \bar{v}, \bar{\tau})(\theta)$ are as defined in Section 3.5. % and $\pi+\sigma_{po}<\theta_w<\pi$, where the constant $\sigma_{po}$ is defined in (\ref{102408}).
We use the subscripts $2$ and $3$ to denote the corresponding flow states of the planar fan wave and the centered simple wave, respectively; that is,
$$
(u_2, v_2, \tau_2)(\xi, \eta)=(u_{-}, v_{-}, \tau_{-})(\xi)\quad \text{and} \quad (u_3, v_3, \tau_3)(\xi, \eta)=(\bar{u}, \bar{v}, \bar{\tau})(\theta).
$$

%where the function $(u_{-}, v_{-}, \tau_{-})(\xi)$ is defined in Section 1.4.2.

\subsection{Post-sonic shock free boundary problem}
The two composite waves begin to interact  at the point $\mathrm{P}(\xi_{_\mathrm{P}}, \eta_{_\mathrm{P}})$, where
$$
\begin{aligned}
\xi_{_\mathrm{P}}=u_0-\frac{\rho_{po}c_{po}}{\rho_0}\quad \text{and} \quad \eta_{_\mathrm{P}}=\xi_{_\mathrm{P}}\tan\alpha_0.
\end{aligned}
$$
Through  $\mathrm{P}$, we draw a forward $C_{+}$ cross characteristic curve of the planar fan wave. As in Section 4, this $C_{+}$ cross characteristic curve can cross the entire planar fan wave and terminates at point $\mathrm{B}(\xi_{_\mathrm{B}}, 0)$.

The planar shock will propagate into the centered simple wave region from $\mathrm{P}$. We consider (\ref{PEC}) with the following boundary conditions:
\begin{equation}\label{free11}
(u, v, \tau)=(u_{-}, v_{-}, \tau_{-})(\xi)\quad \text{on}\quad \widetilde{\mathrm{PB}};
\end{equation}
\begin{equation}\label{free12}
(\rho  U -\rho_3U_3)(U -U_3)+  (\rho  V -\rho_3V_3)(V -V_3)=0\quad \text{on}\quad \xi=\psi(\eta), \quad \eta\leq \eta_{_\mathrm{P}};
\end{equation}
\begin{equation}\label{free13}
\frac{{\rm d}\psi(\eta)}{{\rm d}\xi}=-\frac{V-V_3}{U-U_3}, \quad \psi(\eta_{_\mathrm{P}})=\xi_{_\mathrm{P}},
\end{equation}
where $\xi=\psi(\eta)$, $\eta\leq \eta_{_\mathrm{P}}$ represents the shock propagating from $\mathrm{P}$, and $(U_3, V_3)=(u_3(\xi, \eta)-\xi, v_3(\xi, \eta)-\eta)$.

The main result of this section is the following theorem, which characterizes the propagation of a post-sonic shock in a central simple wave of a different type.
\begin{thm}\label{thm42201}
Assume that the amplitude of the fan wave of the oblique fan–shock composite wave is sufficiently small, i.e.,  $\bar{\tau}(\phi_{pr})-\tau_0$ is sufficiently small.
The shock free boundary problem (\ref{PEC}), (\ref{free11})--(\ref{free13}) admits a piecewise smooth solution on a curvilinear quadrilateral domain bounded by $\widehat{\mathrm{PE}}$, $\widetilde{\mathrm{PB}}$, $\widetilde{\mathrm{EF}}$, and $\wideparen{\mathrm{BF}}$. Here,  $\widehat{\mathrm{PE}}$ represents the shock propagating from $\mathrm{P}$, point $\mathrm{E}$ lies on the oblique pre-sonic shock line $\eta=\xi\tan\phi_{pr}$, $\widetilde{\mathrm{EF}}$ is a forward $C_{+}$ characteristic curve emanating from $\mathrm{E}$, and $\wideparen{\mathrm{BF}}$ is a vacuum boundary.
Moreover, the shock $\widehat{\mathrm{PE}}$ is post-sonic.
\end{thm}

%=============
\begin{figure}[htbp]
\begin{center}
\includegraphics[scale=0.45]{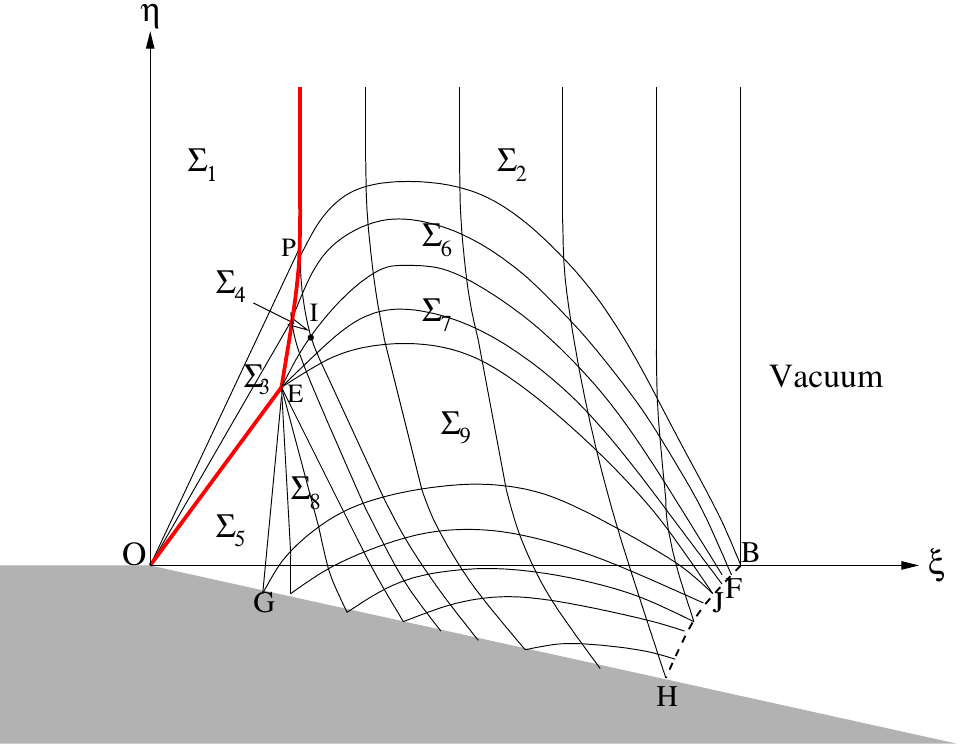}\qquad\quad\includegraphics[scale=0.45]{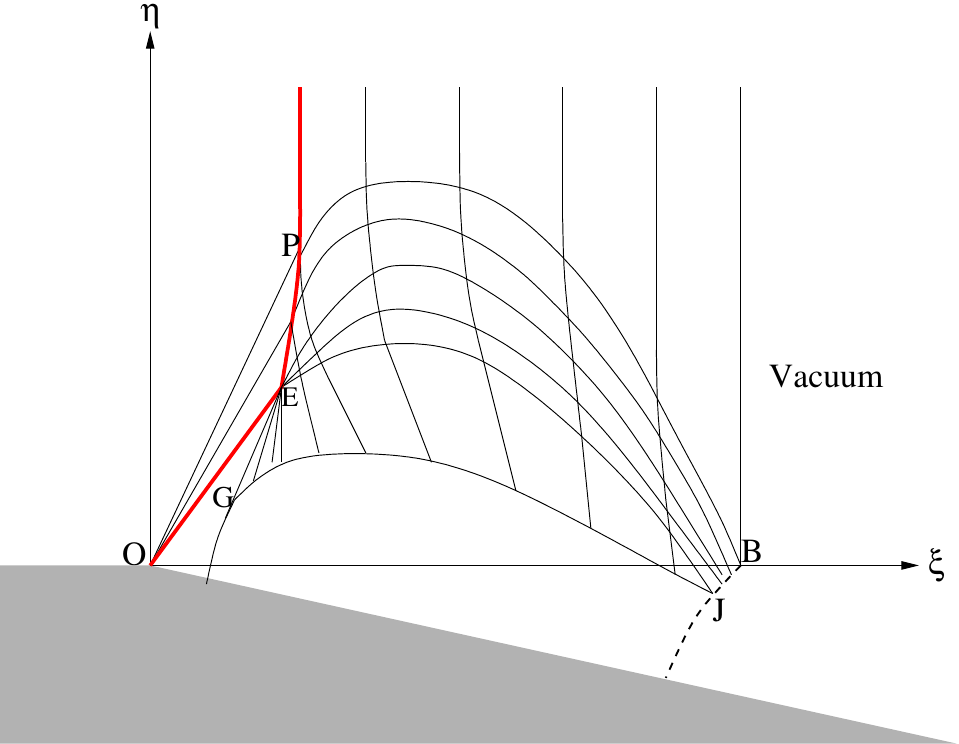}
\caption{\footnotesize Interaction between shock–fan and fan–shock composite waves for $\tau_w>\tau_2^i$.}
\label{Fig8.1}
\end{center}
\end{figure}
%==========

%In what follows we are going to construct a piecewise sooth solution to the free boundary problem (\ref{PEC}), (\ref{free11})--(\ref{free13}).
%\subsubsection{\bf Invariant region of $(\alpha, \beta)$}
%To start off with, we define the following invariant region:

 %For convenience, we denote the shock by  ${\it S}_{_\mathrm{P}}$.
We first construct the post-sonic shock propagating from $\mathrm{P}$, which is initially denoted by ${\it S}_{_\mathrm{P}}$.
We still use the subscripts $f$ and $b$ to denote the states on the front and back sides of ${\it S}_{_\mathrm{P}}$.
Then we have $(u_f, v_f, \tau_f)=(u_3, v_3, \tau_3)(\xi, \eta)$ and $\tau_b=\tau_{po}(\tau_f)$.
We also define
$$
L=U\cos \phi+V\sin \phi\quad \text{and}\quad N=-U \sin \phi+V\cos \phi,
$$
where the variable $\phi$ denotes the inclination angle of ${\it S}_{_\mathrm{P}}$.
Along the shock ${\it S}_{_\mathrm{P}}$, the Rankine-Hugoniot conditions (\ref{RH1}) hold.

As in (\ref{61115a}), we have
\begin{equation}\label{83001a}
\bar{\partial}_s \phi=-\frac{\bar{\partial}_s m}{\rho_fL_f}+\frac{N_{f}}{\rho_{f}L_{f}q_{f}^2}(q_{f}^2+\tau_{f}^2 p'(\tau_{f}))\bar{\partial}_s\rho_{f}-\frac{N_{f}}{q_{f}^2}+\bar{\partial}_s\sigma_{f}\quad\text{along}\quad {\it S}_{_\mathrm{P}}.
\end{equation}
From (\ref{250104}), we have
$$
\begin{aligned}
\bar{\partial}_sm &~=~\frac{\tau_f}{\rho_fN_f(\tau_f^2-\tau_b^2)}\bigg(-p'(\tau_f)+\frac{2h(\tau_f)-2h(\tau_b)}{\tau_f^2-\tau_b^2}\bigg)
\bar{\partial}_s\tau_f
\\&~=~\frac{\tau_f(N_f^2-c_f^2)}{\rho_fN_f(\tau_f^2-\tau_b^2)}
\bar{\partial}_s\rho_f\quad \text{along}\quad {\it S}_{_\mathrm{P}}.
\end{aligned}
$$
Hence, 
\begin{equation}\label{100903}
\bar{\partial}_s \phi=-\frac{\tau_f(N_f^2-c_f^2)}{\rho_f^2N_fL_f(\tau_f^2-\tau_b^2)}
\bar{\partial}_s\rho_f+\frac{N_{f}\cos^2 A_f}{\rho_{f}L_{f}}\bar{\partial}_s\rho_{f}-\frac{N_{f}}{q_{f}^2}+\bar{\partial}_s\sigma_{f}\quad\text{along}\quad {\it S}_{_\mathrm{P}}.
\end{equation}
Therefore, the post-sonic shock front ${\it S}_{_\mathrm{P}}$ can be determined by solving (\ref{100903}) with initial data $\phi(P)=-\frac{\pi}{2}$. The states on the backside of the post-sonic shock can be determined by
$$
u_b=\xi+L_b\cos\phi-N_b\sin\phi,\quad  v_b=\eta+L_b\sin\phi+N_b\cos\phi,
$$
where
$L_b=L_f$ and  $N_b=\frac{\tau_bN_f}{\tau_f}$.

The post-sonic shock front ${\it S}_{_\mathrm{P}}$ and the pre-sonic shock of the oblique fan–shock composite wave intersect at point $\mathrm{E}$. We now denote the post-sonic shock curve propagating from the point $\mathrm{P}$ by $\widehat{\mathrm{PE}}$.

\begin{lem}\label{2662602}
We have the following estimates:
\begin{equation}\label{100702a}
\cos\alpha_b \bar{\partial}_s u_b+\sin\alpha_b \bar{\partial}_s v_b<0\quad \text{and}\quad \cos\beta_b \bar{\partial}_s u_b+\sin\beta_b \bar{\partial}_s v_b>0\quad \text{along}\quad  \widehat{\mathrm{PE}}.
\end{equation}
\end{lem}
\begin{proof}
We emphasize that
the proof of this lemma is different to that of Lemma \ref{61302},
since $\bar{\partial}_sm\neq0$ at the starting point of the post-sonic shock.

We decompose $\bar{\partial}_{s}$ in two characteristic directions
$$
\bar{\partial}_{s}=t_{+}\bar{\partial}_{+}+t_{-}\bar{\partial}_{-},
$$
where
$$
t_{+}=\frac{\sin(\phi-\beta_f)}{\sin(2A_f)}\quad \text{and}\quad t_{-}=\frac{\sin(\alpha_f -\phi)}{\sin(2A_f)}.
$$
Then, by (\ref{82301}), (\ref{8}), and (\ref{83001a}), we have
\begin{equation}\label{61108f}
\begin{aligned}
\bar{\partial}_s \phi+\frac{\bar{\partial}_s m}{\rho_fL_f}&=\frac{N_{f}}{\rho_{f}L_{f}q_{f}^2}\big(q_{f}^2+\tau_{f}^2 p'(\tau_{f})\big)\bar{\partial}_s\rho_{f}-\frac{N_{f}}{q_{f}^2}+\bar{\partial}_s\sigma_{f}\\&=
\frac{N_{f}\cos^2 A}{\rho_{f}L_{f}}t_{-}\bar{\partial}_{-}\rho_{3}-\frac{N_{f}}{q_{f}^2}+
t_{-}\bar{\partial}_{-}\Big(\frac{\alpha_{3}+\beta_{3}}{2}\Big)-\frac{t_{+}\sin A_f}{c_f}
\\&=
\Big(\frac{\sin(2A_f)}{2\rho_f}+\frac{N_{f}\cos^2 A}{\rho_{f}L_{f}}\Big)t_{-}\bar{\partial}_{-}\rho_{3}-\frac{N_{f}}{q_{f}^2}+t_{-}\frac{\sin^2 A_f}{c_f}-\frac{t_{+}\sin A_f}{c_f}
\\&=
\Big(\frac{\sin(2A_f)}{2\rho_f}+\frac{N_{f}\cos^2 A}{\rho_{f}L_{f}}\Big)t_{-}\bar{\partial}_{-}\rho_{3}-\frac{\sin(\sigma_f-\phi)}{q_{f}}+\frac{\sin(\alpha_f-\phi)}{2q_f\cos A_f}+\frac{\sin(\beta_f-\phi)}{2q_f\cos A_f}
\\&=
\Big(\frac{\sin(2A_f)}{2\rho_f}+\frac{N_{f}\cos^2 A}{\rho_{f}L_{f}}\Big)\bar{\partial}_{s}\rho_{f}\quad \text{along}\quad  \widehat{\mathrm{PE}}.
\end{aligned}
\end{equation}

As in (\ref{62303}), we have
$$
\bar{\partial}_{s}\rho_f=t_{-}\bar{\partial}_{-}\rho_{3}<0\quad \text{along}\quad  \widehat{\mathrm{PE}}.
$$
Thus, we obtain
\begin{equation}
\begin{aligned}
\bar{\partial}_s \phi+\frac{\bar{\partial}_s m}{\rho_bL_b}&=
\bar{\partial}_s \phi+\frac{\bar{\partial}_s m}{\rho_fL_f}-\left(\frac{1}{\rho_fL_f}-\frac{1}{\rho_bL_b}\right)\bar{\partial}_s m\\
~&=~\Big(\frac{\sin(2A_f)}{2\rho_f}+\frac{N_{f}\cos^2 A}{\rho_{f}L_{f}}\Big)\bar{\partial}_{s}\rho_{f}-\frac{\tau_f(N_f^2-c_f^2)}{\rho_fN_fL_f(\tau_f+\tau_b)}
\bar{\partial}_s\rho_f\\
~&<~\Big(\frac{\sin(2A_f)}{2\rho_f}+\frac{N_{f}\cos^2 A}{\rho_{f}L_{f}}\Big)\bar{\partial}_{s}\rho_{f}-\frac{N_f^2-c_f^2}{\rho_fN_fL_f}
\bar{\partial}_s\rho_f\\&=
\left(\frac{\sin(2A_f)}{2\rho_f}+\frac{(q_f^2-N_f^2)\sin^2 A_f}{\rho_fL_fN_f}\right)\bar{\partial}_{s}\rho_{f}<0\quad \text{along}\quad \widehat{\mathrm{PE}}.
\end{aligned}
\end{equation}

Therefore,
as in the proof of Lemma \ref{61302}, we have that along $\widehat{\mathrm{PE}}$,
$$
\begin{aligned}
&\cos\beta_b \bar{\partial}_s u_b+\sin\beta_b \bar{\partial}_s v_b\\~=~&
\cos A_b \bar{\partial}_s q_b-q_b\sin A_b \bar{\partial}_s \sigma_b+1
\\=~&
\cos A_b\left(-\tau_bq_b\sin^2 A_b\bar{\partial}_s \rho_b-\cos A_b\right)
-q_b\sin A_b\left(\bar{\partial}_s \phi+\frac{\bar{\partial}_s m}{\rho_bL_b}-\frac{\sin 2A_b}{2\rho_b}\bar{\partial}_s\rho_b+\frac{N_b}{q_b^2}\right)+1\\=~&
-q_b\sin A_b \left(\bar{\partial}_s \phi+\frac{\bar{\partial}_s m }{\rho_bL_b}\right)>0
\end{aligned}
$$
and
\begin{equation}\label{26101001}
\begin{aligned}
&\cos\alpha_b \bar{\partial}_s u_b+\sin\alpha_b \bar{\partial}_s v_b\\~=~&
\cos A_b \bar{\partial}_s q_b+q_b\sin A_b \bar{\partial}_s \sigma_b+\cos (2A_b)
\\=~&
\cos (2A_b)-\cos A_b\left(\tau_bq_b\sin^2 A_b\bar{\partial}_s \rho_b+\cos A_b\right)
+q_b\sin A_b\left(\bar{\partial}_s \phi+\frac{\bar{\partial}_s m}{\rho_bL_b}-\frac{\sin 2A_b}{2\rho_b}\bar{\partial}_s\rho_b+\frac{N_b}{q_b^2}\right)\\=~&
q_b\sin A_b \left(\bar{\partial}_s \phi+\frac{\bar{\partial}_s m }{\rho_bL_b}\right) -
2\tau_bq_b\cos A_b\sin^2 A_b\bar{\partial}_s \rho_b<0.
\end{aligned}
\end{equation}

This completes the proof of the lemma.
\end{proof}

We next consider (\ref{PEC}) with the initial data
\begin{equation}\label{100703a}
(u, v, \tau)=(u_b, v_b, \tau_b)(\xi, \eta), \quad (\xi, \eta)\in \widehat{\mathrm{PE}}.
\end{equation}
As in Section 6.5, applying the hodograph transformation method, we obtain that when $\bar{\tau}(\phi_{pr})-\tau_0$ is sufficiently small,
the singular Cauchy problem (\ref{PEC}), (\ref{100703a}) admits a solution on a curvilinear triangular domain $\overline{\Sigma_4}$ bounded by $\widehat{\mathrm{PE}}$, $\widetilde{\mathrm{PI}}$, and $\widetilde{\mathrm{EI}}$, where $\widetilde{\mathrm{PI}}$ is a forward $C_{-}$ characteristic curve emanating from $\mathrm{P}$ and $\widetilde{\mathrm{EI}}$ is a forward $C_{+}$ characteristic curve emanating from $\mathrm{E}$.
Moreover, the solution satisfies the following estimates:
\begin{equation}\label{42201}
\Big|(\alpha, \beta)-\big(\alpha_2(\mathrm{P}), -\frac{\pi}{2}\big)\Big|<\varepsilon_0\quad \text{on}\quad \overline{\Sigma_4};\quad
\bar{\partial}_{-}\rho<0\quad \text{along}\quad \widetilde{\mathrm{PI}}; \quad \bar{\partial}_{+}\rho<0\quad \text{along}\quad \widetilde{\mathrm{EI}},
\end{equation}
where $\varepsilon_0=\frac{1}{4}\big(\frac{\pi}{2}-\alpha_2(\mathrm{P})\big)$.

For convenience, we use $(u, v, \tau)=(u_4, v_4, \tau_4)(\xi, \eta)$ to denote the solution of the singular Cauchy problem.
We next consider (\ref{PEC}) with the boundary conditions
\begin{equation}\label{100704a}
(u, v, \tau)=\left\{
               \begin{array}{ll}
               (u_{-}, v_{-}, \tau_{-})(\xi), & \hbox{$(\xi, \eta)\in \widetilde{\mathrm{PB}}$;} \\[2pt]
                 (u_4, v_4, \tau_4)(\xi, \eta), & \hbox{$(\xi, \eta)\in \widetilde{\mathrm{PI}}$.}
               \end{array}
             \right.
\end{equation}

For $\tau>\tau_{po}$,
we define
\begin{equation}
\Gamma_{50}(\tau)=\Big\{(\alpha, \beta)\mid \alpha_{5l}(\tau)<\alpha<\alpha_{5r},~\beta_{5l}<\beta<\beta_{5r}(\tau),~\alpha-\beta>\delta_0\Big\},
\end{equation}
where
$$
 \alpha_{5r}=\alpha_2(\mathrm{P})+\varepsilon_0,\quad \beta_{5r}=\alpha_{5r}-2\Theta(\tau), \quad \beta_{5l}=-\frac{\pi}{2}-\varepsilon_0, \quad\alpha_{5l}=\beta_{5l}+2\Theta(\tau),
$$
 and $\delta_0$ is a small positive constant.

Assume
$\alpha_2(\mathrm{P})+\frac{\pi}{2}>2\Theta(\tau_{po})$.
As in Section 4.1,  we have
\begin{equation}\label{42203}
(\alpha, \beta)\in \Gamma_{50}(\tau)\quad \text{and}\quad \bar{\partial}_{+}\rho<0\quad \text{on}\quad \widetilde{\mathrm{PB}}.
\end{equation}
Based on (\ref{42201}) and (\ref{42203}),
the Goursat problem (\ref{PEC}), (\ref{100704a}) admits a classical solution in a domain $\overline{\Sigma_6}$ bounded by $\widetilde{\mathrm{PB}}$, $\widetilde{\mathrm{PI}}$, $\widetilde{\mathrm{IF}}$, and $\wideparen{\mathrm{BF}}$, where
$\widetilde{\mathrm{IF}}$ is a forward $C_{+}$ characteristic curve emanating from $\mathrm{I}$ and $\wideparen{\mathrm{BF}}$ is a vacuum boundary. Moreover, the solution satisfies
\begin{equation}\label{42204}
(\alpha, \beta)\in \Gamma_{50}(\tau)\quad \mbox{and} \quad \bar{\partial}_{+}\rho<0\quad \text{on}\quad \widetilde{\mathrm{IF}}.
\end{equation}
We use $(u, v, \tau)=(u_6, v_6, \tau_6)(\xi, \eta)$ to denote the solution of the  Goursat problem (\ref{PEC}), (\ref{100704a}).

We define
\begin{equation}\label{42205}
(u, v, \tau)=\left\{
               \begin{array}{ll}
                 (u_4, v_4, \tau_4)(\xi, \eta), & \hbox{$(\xi, \eta)\in\overline{\Sigma_4}$;} \\[2pt]
                 (u_6, v_6, \tau_6)(\xi, \eta), & \hbox{$(\xi, \eta)\in\overline{\Sigma_6}$.}
               \end{array}
             \right.
\end{equation}
Then the function $(u, v, \tau)$ defined in (\ref{42205}) is a solution of the free boundary problem (\ref{PEC}), (\ref{free11})--(\ref{free13}).
This completes the proof of Theorem \ref{thm42201}.

%\subsubsection{Nonexistence of transonic shock solution}

\subsection{Wave interaction for the case $\tau_{w}\geq \tau_2^i$}
For convenience, we use the subscript $5$ to denote the constant state $(u_{w}, v_{w}, \tau_{w})$; that is,
$
 (u_5, v_5, \tau_5)(\xi, \eta)=(u_{w}, v_{w}, \tau_{w}).
$
The post-sonic shock $\widehat{\mathrm{PE}}$ and the pre-sonic shock $\overline{\mathrm{OE}}$  interact at the point $\mathrm{E}$.
Assume that
$$
\sigma_4(\mathrm{E})-\sigma_5(\mathrm{E})-\int_{q_5(\mathrm{E})}
^{q_4(\mathrm{E})}\frac{\sqrt{q^2-c^2(\grave{\tau}(q))}}{qc(\grave{\tau}(q))}{\rm d}q>0,
$$
where the function  $\grave{\tau}(q)$ is determined by
$$
\frac{1}{2}q^2+\int_{\tau_{w}}^{\grave{\tau}(q)}\tau p'(\tau)~{\rm d}\tau=\frac{1}{2}q_5^2(\mathrm{E}).
$$
Then
the interaction of these two shocks generates two centered waves emanating from $\mathrm{E}$. To determine the principal parts of the centered waves, we consider (\ref{principal})--(\ref{principal1}) with  initial data
\begin{equation}\label{62802ar}
(\widetilde{u}_{+},\widetilde{v}_{+}, \widetilde{\rho}_{+})(\nu_2)=(u_4, v_4, \tau_4)(\mathrm{E}),
\end{equation}
and (\ref{82306})--(\ref{82307}) with  initial data
\begin{equation}\label{62803ar}
(\widetilde{u}_{-},\widetilde{v}_{-}, \widetilde{\rho}_{-})(\mu_1)=(u_{w}, v_{w}, \rho_{w}),
\end{equation}
where $\nu_2=\tan(\alpha_4(\mathrm{E}))$ and $\mu_1=\cot(\beta_5(\mathrm{E}))$.

%where $q_2(\mathrm{P})=\sqrt{(u_{-}(\xi_{_\mathrm{P}})-\xi_{_\mathrm{P}})^2+\eta_{_\mathrm{P}}^2}$.

Assume further that
\begin{equation}\label{8260dr}
2\int_{q_4(\mathrm{E})}^{q_{_{lim}}}\frac{\sqrt{q^2-c^2(\grave{\tau}(q))}}{qc(\grave{\tau}(q))}{\rm d}q>
\sigma_4(\mathrm{E})-\sigma_5(\mathrm{E})-\int_{q_5(\mathrm{E})}
^{q_4(\mathrm{E})}\frac{\sqrt{q^2-c^2(\grave{\tau}(q))}}{qc(\grave{\tau}(q))}{\rm d}q.
\end{equation}
%here $q_3(\mathrm{P})=\sqrt{(u_{w}-\xi_{_\mathrm{P}})^2+(v_{w}-\eta_{_\mathrm{P}})^2}$ and $q_{_{lim}}=\sqrt{q_2^2(\mathrm{P})-2\int_{\tau_{po}}^{+\infty}\tau p'(\tau)~{\rm d}\tau}$.
Then the equations
$$
\left\{
  \begin{array}{ll}
    \displaystyle\widetilde{\sigma}+ \int_{q_{4}(\mathrm{E})}^{\widetilde{q}}\frac{\sqrt{q^2-c^2(\grave{\tau}(q))}}{qc(\grave{\tau}(q))}{\rm d}q=\sigma_4(\mathrm{E}),\\[12pt]
   \displaystyle \widetilde{\sigma}- \int_{q_{4}(\mathrm{E})}^{\widetilde{q}}
   \frac{\sqrt{q^2-c^2(\grave{\tau}(q))}}{qc(\grave{\tau}(q))}{\rm d}q=\sigma_5(\mathrm{E})+\int_{q_5(\mathrm{E})}^{q_4(\mathrm{E})}\frac{\sqrt{q^2-c^2(\grave{\tau}(q))}}{qc(\grave{\tau}(q))}{\rm d}q
  \end{array}
\right.
$$
 admits a unique solution $(\widetilde{\sigma}, \widetilde{q})=(\sigma_m, q_m)$ with $\widetilde{q}>q_4(\mathrm{E})$.

Let
 $$\nu_1=\tan(\sigma_m+A_m) \quad \text{and}\quad \mu_2=\cot(\sigma_m-A_m),$$
where $A_m=\arcsin\Big(\frac{c(\tau_m)}{q_m}\Big)$ and $\tau_m=\grave{\tau}(q_m)$.
 As in Lemma \ref{8250c}, the initial value problem (\ref{principal})--(\ref{principal1}), (\ref{62802a}) admits a solution on $[\nu_1, \nu_2]$ and the initial value problem (\ref{82306})--(\ref{82307}), (\ref{62803a}) admits a solution on $[\mu_2, \mu_1]$, and the solutions satisfy
$
(\widetilde{u}_{-},\widetilde{v}_{-}, \widetilde{\rho}_{-})(\mu_2)=(\widetilde{u}_{+},\widetilde{v}_{+}, \widetilde{\rho}_{+})(\nu_1).
$
Moreover,
$\widetilde{\tau}_{+}'(\nu)<0$ for $\nu_1\leq\nu\leq \nu_2$;
$\widetilde{\tau}_{-}'(\mu)<0$ for $\mu_2\leq\nu\leq \mu_1$.

Let $\mu=\widetilde{\mu}_{-}(\tau)$, $\tau\in[\tau_{w}, \tau_m]$ be the inverse function of $\tau=\widetilde{\tau}_{-}(\mu)$,  $\mu\in[\mu_2, \mu_1]$,
and let
$\nu=\widetilde{\nu}_{+}(\tau)$, $\tau\in[\tau_{po}, \tau_m]$ be the inverse function of $\tau=\widetilde{\tau}_{+}(\nu)$, $\nu\in [\nu_1, \nu_2]$.
%at the function $\tau=\widetilde{\tau}_{-}(\mu)$ has an inverse function $\mu=\widetilde{\mu}_{-}(\tau)$.
Let  $$(\check{\alpha}_{-}(\tau), \check{\beta}_{-}(\tau))=\big(\widetilde{\alpha}_{-}(\widetilde{\mu}_{-}(\tau)), \widetilde{\beta}_{-}(\widetilde{\mu}_{-}(\tau))\big), \quad \tau\in[\tau_{w}, \tau_m]$$ and
  $$(\check{\alpha}_{+}(\tau), \check{\beta}_{+}(\tau))=\big(\widetilde{\alpha}_{+}(\widetilde{\nu}_{+}(\tau)), \widetilde{\beta}_{+}(\widetilde{\nu}_{+}(\tau))\big), \quad \tau\in[\tau_{po}, \tau_m].$$

Through $\mathrm{E}$, we draw a forward straight $C_{-}$ characteristic line with the characteristic angle $\beta_5(\mathrm{E})$.
This straight line intersects the ramp $\partial{\mathcal{W}}_{+}$ or the sonic circle $(\xi-u_w)^2+(\eta-v_w)^2=c_w^2$ at a point $\mathrm{G}$.
We then consider
 (\ref{PEC}) with the boundary conditions
\begin{equation}\label{100705a}
\left\{
               \begin{array}{ll}
               (u, v, \tau)=(u_4, v_4, \tau_4)(\xi, \eta), & \hbox{$(\xi, \eta)\in \widetilde{\mathrm{EI}}$;} \\[2pt]
                 (u, v, \tau)=(u_6, v_6, \tau_6)(\xi, \eta), & \hbox{$(\xi, \eta)\in \widetilde{\mathrm{IF}}$;}\\[2pt]
(u, v, \tau)=(u_w, v_w, \tau_w), & \hbox{$(\xi, \eta)\in\overline{\mathrm{EG}}$.}
                 %v=u\tan(\theta_w-\pi), & \hbox{$(\xi, \eta)\in\mathcal{W}_{+}$;}
               \end{array}
             \right.
\end{equation}

To obtain a global classical solution to the BVP (\ref{PEC}), (\ref{100705a}), we impose the following assumptions
\begin{description}
%\item[($h_{41}$)]
%$\alpha_2(\mathrm{P})-\beta_3(\mathrm{P})<\pi$,  $\alpha_3(\mathrm{P})-\beta_3(\mathrm{P})>2\Theta(\tau_w)$,  $\alpha_2(\mathrm{P})+\frac{\pi}{2}>2\Theta(\tau_{po})$;
%\vskip 2pt
%\item[($h_{42}$)] $\sigma_2(\mathrm{P})>\sigma_3(\mathrm{P})+{\displaystyle\int_{q_3(\mathrm{P})}
%^{q_2(\mathrm{P})}}\frac{\sqrt{q^2-c^2(\grave{\tau}(q))}}{qc(\grave{\tau}(q))}{\rm d}q$;
%\vskip 2pt
\item[($h_{51}$)] $\check{\alpha}_{+}(\tau)>\frac{\pi}{2}+2\Theta(\tau)$ for  $\tau\in[\tau_{po}, \tau_m]$;
\item[($h_{52}$)] $\sigma_4(\mathrm{E})>\sigma_5(\mathrm{E})+{\displaystyle\int_{q_5(\mathrm{E})}
^{q_4(\mathrm{E})}}\frac{\sqrt{q^2-c^2(\grave{\tau}(q))}}{qc(\grave{\tau}(q))}{\rm d}q$.
%\vskip 4pt
 %  \item[($h_{44}$)] $\check{\beta}_{-}(\tau)<\alpha_2(\mathrm{P})-2\Theta(\tau)$ for  $\tau\in[\tau_{w}, \tau_m]$.
 % \item[($h_{45}$)] $\beta_3(\mathrm{P})+2\Theta(\tau_{po})<\theta_w-\pi<\frac{\pi}{2}+\beta_3(\mathrm{P})$.
\end{description}
As in Section 7.1, we obtain that
under assumptions ($h_{51}$) and ($h_{52}$), the boundary value problem (\ref{PEC}), (\ref{100705a}) admits a classical solution in a domain bounded by
$\widetilde{\mathrm{EF}}$, $\overline{\mathrm{EG}}$, $\widetilde{\mathrm{GJ}}$, and $\widetilde{\mathrm{FJ}}$, where $\widetilde{\mathrm{GJ}}$ is a forward $C_{+}$ characteristic curve emanating from $\mathrm{G}$
and $\widetilde{\mathrm{FJ}}$ is a vacuum boundary; see Figure \ref{Fig8.1}.
%\begin{proof}
%The proof is similar to that of Lemma, we omit the details.
%\end{proof}
Therefore, when $\tau_{w}>\tau_2^i$, under assumptions ($h_{51}$) and ($h_{52}$), the interaction between the shock–fan and fan–shock composite waves generates a non-planar simple wave  emanating from the interaction region. 

\subsection{Wave interaction for the case $\tau_{w}<\tau_2^i$}
If $\tau_{w}<\tau_2^i$,
there is a straight post-sonic 1-shock propagating from the point $\mathrm{E}$.
The state of the flow upstream of the shock is $(u_{w}, v_{w}, \tau_{w})$ and the specific volume of the flow downstream of the shock is $\tau_d=\psi_{po}(\tau_{w})$. We denote by $\phi_{_\mathrm{E}}$ the inclination angle and $(u_d, v_d)$ the downstream flow velocity of the oblique post-sonic shock.
%The ray $\xi=\xi_{_P}+r\cos\phi_{w}$, $\eta=\eta_{_P}+r\sin\phi_{w}$ ($r>0$) intersects $\partial{\mathcal{W}}_{+}$ at some point $G$.

Let
\begin{equation}\label{8260ar}
\begin{aligned}
&(U_d, V_d)=(u_d-\xi_{_\mathrm{P}}, v_d-\eta_{_\mathrm{P}}), \quad (U_{w}, V_{w})=(u_{w}-\xi_{_\mathrm{P}}, v_{w}-\eta_{_\mathrm{E}}), \quad
\\& \qquad N_d=V_d\cos\phi_{_\mathrm{E}}-U_d\sin\phi_{_\mathrm{E}},\quad L_d=U_d\cos\phi_{_\mathrm{E}}+V_d\sin\phi_{_\mathrm{E}},\\&\quad N_{w}=V_{w}\cos\phi_{_\mathrm{E}}-U_{w}\sin\phi_{_\mathrm{E}},\quad L_{w}=U_{w}\cos\phi_{_\mathrm{E}}+V_{w}\sin\phi_{_\mathrm{E}}.
\end{aligned}
\end{equation}
Since the shock is post-sonic, we have
$
 N_d=c_d=\sqrt{-\tau_d^2p'(\tau_d)}
$.
Hence,
%$$
%\rho_dN_d=\rho_{w}N_{w},\quad \rho_{w}=\tau_{w}^{-1}, \quad N_d=c_d=\sqrt{-\tau_d^2p'(\tau_d)}.
%$$
%Thus, we have
\begin{equation}\label{8260br}
N_{w}=\frac{\rho_d c_d}{\rho_w}, \quad L_w=L_d=\sqrt{q_5^2(\mathrm{E})-N_w^2}, \quad \text{and} \quad \phi_{_\mathrm{E}}=\sigma_5(\mathrm{E})-\arcsin\left(\frac{N_w}{q_5(\mathrm{E})}\right).
\end{equation}
From (\ref{8260ar}) and (\ref{8260br}), we obtain $u_d$ and $v_d$. Hence, the oblique shock front emanating from $\mathrm{E}$ and the state on the back side of it are determined.

%We have
%$$
%\beta_3(P)>\phi_{_\mathrm{P}}=\beta_5(P).
%$$

Let $q_d=\sqrt{(u_d-\xi_{_\mathrm{E}})^2+(v_d-\eta_{_\mathrm{E}})^2}$ and $\sigma_d=\arctan\big(\frac{v_d-\eta_{_\mathrm{E}}}{u_d-\xi_{_\mathrm{E}}}\big)$.
When
$$
\sigma_4(\mathrm{E})
-\sigma_d-\int_{q_5(\mathrm{E})}^{q_d}\frac{\sqrt{q^2-c^2(\grave{\tau}(q))}}{qc(\grave{\tau}(q))}{\rm d}q>0,
$$
there exists a centered simple wave with straight $C_{-}$ characteristics and a $C_{+}$ type centered wave emanating from the point $\mathrm{E}$.

To determine the principal parts of the centered waves, we consider (\ref{principal})--(\ref{principal1}) with  initial data (\ref{62802a}),
and (\ref{82306})--(\ref{82307}) with  initial data
\begin{equation}\label{62810ad}
(\widetilde{u}_{-},\widetilde{v}_{-}, \widetilde{\rho}_{-})(\bar{\mu}_1)=(u_{d}, v_{d}, \rho_{d}),
\end{equation}
where $\bar{\mu}_1=\cot\phi_{_\mathrm{E}}$.

%=============
\begin{figure}[htbp]
\begin{center}
\includegraphics[scale=0.45]{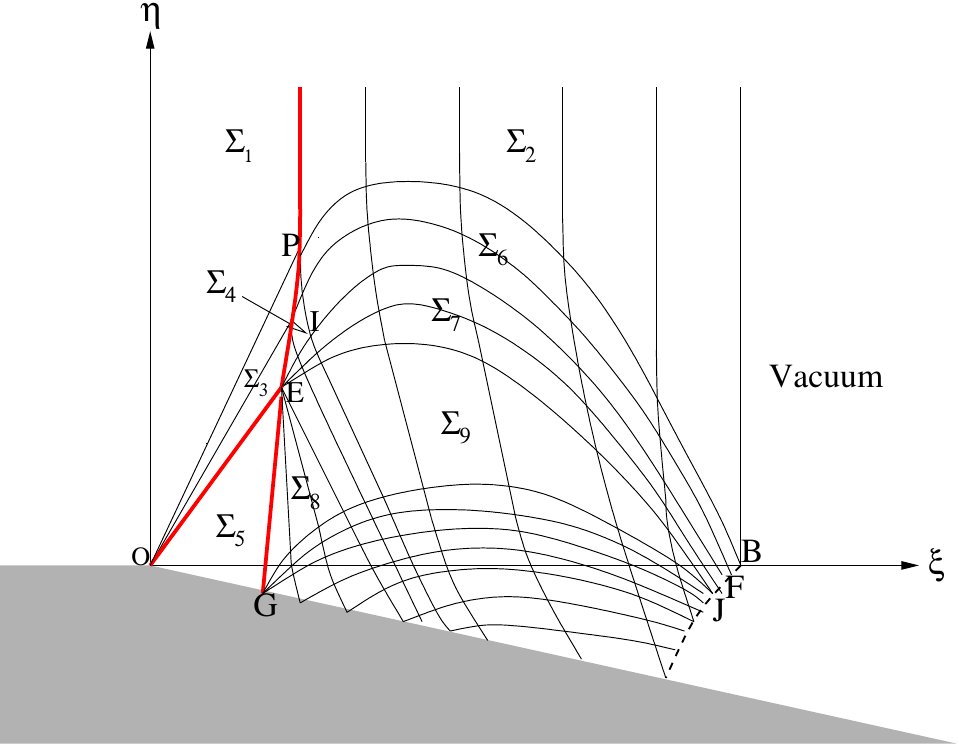}\qquad\quad \includegraphics[scale=0.45]{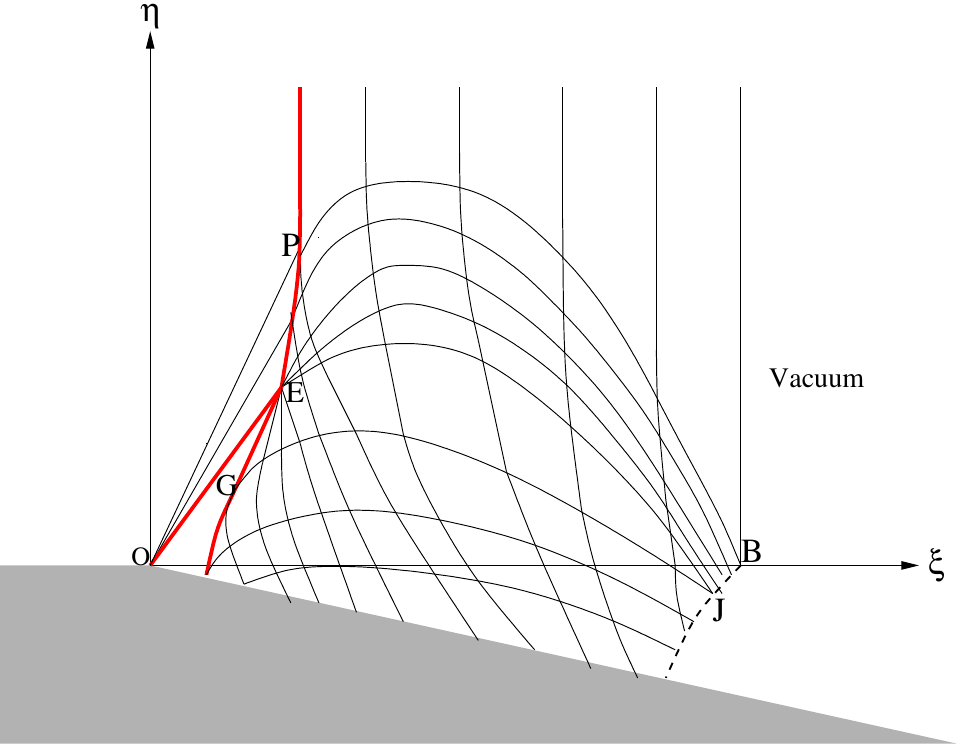}
\caption{\footnotesize Interaction between shock–fan and fan–shock composite waves for $\tau_w<\tau_2^i$.}
\label{Fig8.2}
\end{center}
\end{figure}
%==========

Assume further that
$$
2\int_{q_{4}(E)}^{q_{_{lim}}}\frac{\sqrt{q^2-c^2(\grave{\tau}(q))}}{qc(\grave{\tau}(q))}{\rm d}q>\sigma_4(\mathrm{E})
-\sigma_d-\int_{q_5(\mathrm{E})}^{q_d}\frac{\sqrt{q^2-c^2(\grave{\tau}(q))}}{qc(\grave{\tau}(q))}{\rm d}q.
$$
Then, the equations
$$
\left\{
  \begin{array}{ll}
    \displaystyle\widetilde{\sigma}+ \int_{q_{4}(\mathrm{E})}^{\widetilde{q}}\frac{\sqrt{q^2-c^2(\grave{\tau}(q))}}{qc(\grave{\tau}(q))}{\rm d}q=\sigma_4(\mathrm{E}),\\[12pt]
   \displaystyle \widetilde{\sigma}- \int_{q_d}^{\widetilde{q}}\frac{\sqrt{q^2-c^2(\grave{\tau}(q))}}{qc(\grave{\tau}(q))}{\rm d}q
   =\sigma_d+\int_{q_5(\mathrm{E})}^{q_d}\frac{\sqrt{q^2-c^2(\grave{\tau}(q))}}{qc(\grave{\tau}(q))}{\rm d}q
  \end{array}
\right.
$$
 admits a unique solution $(\widetilde{\sigma}, \widetilde{q})=(\bar{\sigma}_m, \bar{q}_m)$ with $\widetilde{q}>q_4(\mathrm{E})$.

We define
$$\bar{\nu}_1=\tan(\bar{\sigma}_m+\bar{A}_m) \quad \text{and}\quad \bar{\mu}_2=\cot(\bar{\sigma}_m-\bar{A}_m),$$
where $\bar{A}_m=\arcsin\big(\frac{c(\bar{\tau}_m)}{\bar{q}_m}\big)$ and $\bar{\tau}_m=\grave{\tau}(\bar{q}_m)$.

The initial value problem (\ref{principal})--(\ref{principal1}), (\ref{62802a}) admits a solution on $[\bar{\nu}_1, \nu_2]$, and the initial value problem (\ref{82306})--(\ref{82307}), (\ref{62810ad}) admits a solution on $[\bar{\mu}_2, \bar{\mu}_1]$, and the solutions satisfy
$
(\widetilde{u}_{-},\widetilde{v}_{-}, \widetilde{\rho}_{-})(\bar{\mu}_2)=(\widetilde{u}_{+},\widetilde{v}_{+}, \widetilde{\rho}_{+})(\bar{\nu}_1).
$
Moreover,
$\widetilde{\tau}_{+}'(\nu)<0$ for $\bar{\nu}_1\leq\nu\leq \nu_2$;
$\widetilde{\tau}_{-}'(\mu)<0$ for $\bar{\mu}_2\leq\mu\leq \bar{\mu}_1$.
Let $\mu=\widetilde{\mu}_{-}(\tau)$, $\tau\in[\tau_{d}, \bar{\tau}_m]$ be the inverse function of $\tau=\widetilde{\tau}_{-}(\mu)$,  $\mu\in [\bar{\mu}_2, \bar{\mu}_1]$ and
$\nu=\widetilde{\nu}_{+}(\tau)$, $\tau\in[\tau_{po}, \bar{\tau}_m]$ the inverse function of $\tau=\widetilde{\tau}_{+}(\nu)$, $\nu[\bar{\nu}_1, \nu_2]$.
%at the function $\tau=\widetilde{\tau}_{-}(\mu)$ has an inverse function $\mu=\widetilde{\mu}_{-}(\tau)$.
Let  $(\check{\alpha}_{-}(\tau), \check{\beta}_{-}(\tau))=\big(\widetilde{\alpha}_{-}(\widetilde{\mu}_{-}(\tau)), \widetilde{\beta}_{-}(\widetilde{\mu}_{-}(\tau))\big)$ and
  $(\check{\alpha}_{+}(\tau), \check{\beta}_{+}(\tau))=\big(\widetilde{\alpha}_{+}(\widetilde{\nu}_{+}(\tau)), \widetilde{\beta}_{+}(\widetilde{\nu}_{+}(\tau))\big)$.

This straight post-sonic shock intersects $\partial{\mathcal{W}}_{+}$ or the sonic circle $(\xi-u_d)^2+(\eta-v_d)^2=c_d^2$ at a point $\mathrm{G}$. We consider  (\ref{PEC})
with the boundary conditions
\begin{equation}\label{10021ar}
\left\{
               \begin{array}{ll}
                   (u, v, \tau)=(u_4, v_4, \tau_4)(\xi, \eta), & \hbox{$(\xi, \eta)\in \widetilde{\mathrm{EI}}$;} \\[2pt]
                 (u, v, \tau)=(u_6, v_6, \tau_6)(\xi, \eta), & \hbox{$(\xi, \eta)\in \widetilde{\mathrm{IF}}$;}\\[2pt]
                 (u, v, \tau)=(u_d, v_d, \tau_d), & \hbox{$(\xi, \eta)\in\overline{\mathrm{EG}}$.}
               \end{array}
             \right.
\end{equation}

To obtain a global classical solution to the BVP (\ref{PEC}), (\ref{10021ar}),
we impose the following assumptions:
\begin{description}
%\item[($h_{41}$)]
%$\alpha_2(\mathrm{P})-\beta_3(\mathrm{P})<\pi$,  $\alpha_3(\mathrm{P})-\beta_3(\mathrm{P})>2\Theta(\tau_w)$,  $\alpha_2(\mathrm{P})+\frac{\pi}{2}>2\Theta(\tau_{po})$;
%\vskip 2pt
%\item[($h_{42}$)] $\sigma_2(\mathrm{P})>\sigma_3(\mathrm{P})+{\displaystyle\int_{q_3(\mathrm{P})}
%^{q_2(\mathrm{P})}}\frac{\sqrt{q^2-c^2(\grave{\tau}(q))}}{qc(\grave{\tau}(q))}{\rm d}q$;
%\vskip 2pt
\item[($h_{51}'$)] $\check{\alpha}_{+}(\tau)>\frac{\pi}{2}+2\Theta(\tau)$ for  $\tau\in[\tau_{po}, \bar{\tau}_m]$;
\vskip 2pt
\item[($h_{52}'$)] $\displaystyle\sigma_4(\mathrm{E})
>\sigma_d+\int_{q_4(\mathrm{E})}^{q_d}{\frac{\sqrt{q^2-c^2(\grave{\tau}(q))}}{qc(\grave{\tau}(q))}}{\rm d}q$.
%\vskip 4pt
 %  \item[($h_{44}$)] $\check{\beta}_{-}(\tau)<\alpha_2(\mathrm{P})-2\Theta(\tau)$ for  $\tau\in[\tau_{w}, \tau_m]$.
 % \item[($h_{45}$)] $\beta_3(\mathrm{P})+2\Theta(\tau_{po})<\theta_w-\pi<\frac{\pi}{2}+\beta_3(\mathrm{P})$.
\end{description}
As in Theorem \ref{1001ba}, we have that
under assumptions ($h_{51}'$) and ($h_{52}'$),
the BVP (\ref{PEC}), (\ref{10021ar}) admits a classical solution on a domain bounded by
$\widetilde{\mathrm{EF}}$, $\overline{\mathrm{EG}}$, $\widetilde{\mathrm{GJ}}$, and $\wideparen{\mathrm{FJ}}$, where $\widetilde{\mathrm{GJ}}$ is a forward $C_{+}$ characteristic curve emanating from $\mathrm{G}$ and
$\wideparen{\mathrm{FJ}}$ is a vacuum boundary; see Figure \ref{Fig8.2}.
Therefore, when $\tau_{w}>\tau_2^i$, under assumptions ($h_{51}'$) and ($h_{52}'$), the interaction of the shock–fan and fan–shock composite waves generates a non-planar shock–fan composite wave emanating from the interaction region. 

%\subsection{Wave reflection by the ramp}
%%%%%%%%%%%%%%%%%%%%%%%%%%%%%%%%%%%%%%%%%%%%%%%%%%%%%%%%%%%%%%%%%%%%%%%%%%%%%%%%%%%%%%%%%%%%%%%%%

%%%%%%%%%%%%%%%%%%%%%%%%%%%%%%%%%%%%%%%%%%%%%%%%%%%%%%%%%%%%%%%%%%%%%%%%%%%%%%%%%%%%%%%%%%%%%%%%%

\section{Interaction between fan–shock–fan and fan–shock composite waves}

%\subsection{Elementary waves}
In this section, we study the 2D Riemann problem (\ref{Potential}), (\ref{initial})--(\ref{boundary}) for $\tau_0<\tau_1^e$ and $u_0>c_0$. % and $\pi+\sigma_{po}<\theta_w<\pi$, where the constant $\sigma_{po}$ is defined in (\ref{102408}).
Then, far away from the corner, the oncoming flow $(0)$ and the vacuum
are connected by a planar fan–shock–fan composite wave which can be represented by (\ref{pfsf}). We assume further that near the corner the oncoming flow is deflected by an oblique fan–shock composite wave to a constant state $(u_{w}, v_{w}, \tau_{w})$, which can be represented by (\ref{100901}).
%In this section, we discuss the interaction of
%In what follows we are going to discuss the interaction of the fan–shock–fan and fan–shock composite waves. % to a constant state $(u_{w}, v_{w}, \tau_{w})$
%We shall divide the discussion into the following two cases: (i) $\tau_{w}\geq \tau_2^i$; (ii) $\tau_{w}<\tau_2^i$.
In the subsequent discussion, the  constants $\xi_0$, $\xi_1$, $\xi_2$, $\hat{u}_1$, $\hat{u}_2$, $\sigma_2^e$, $\phi_{pr}$, as well as and the functions $(\hat{u}_l, \hat{v}_l, \hat{\tau}_l)(\xi)$, $(\hat{u}_r, \hat{v}_r, \hat{\tau}_r)(\xi)$, $(\bar{u}, \bar{v}, \bar{\tau})(\theta)$  are as defined in Section 3.6.

\subsection{Interaction of fan waves}
The two composite waves begin to interact at the point  $\mathrm{P}(\xi_{0}, \eta_{_\mathrm{P}})$,
 where $\eta_{_\mathrm{P}}=\xi_{0}\tan\alpha_0$.
Through $\mathrm{P}$, we draw a $C_{+}$ cross characteristic curve in the planar fan–shock–fan composite wave and a $C_{-}$ cross characteristic curve in the oblique fan–shock composite wave.
The $C_{+}$ cross characteristic curve intersects the double-sonic shock of the planar fan–shock–fan composite wave at a point $\mathrm{E}(\xi_{_\mathrm{E}}, \eta_{_\mathrm{E}})$ and terminates at a point $\mathrm{B}$;
the $C_{-}$ cross characteristic curve intersects the pre-sonic shock of the oblique fan–shock composite wave at a point $\mathrm{F}(\xi_{_\mathrm{F}}, \eta_{_\mathrm{F}})$. %Through the point $\mathrm{D}$ we draw a straight $C_{-}$ characteristic line of the constant state $(u_{w}, v_{w}, \tau_{w})$. This straight characteristic line intersects the ramp at a point $J$.

%We consider (\ref{PEC}) with the boundary conditions
%\begin{equation}\label{100905}
%(u, v, \tau)=\left\{
 %              \begin{array}{ll}
  %                (\hat{u}_l, \hat{v}_l, \hat{\tau}_l)(\xi), & \hbox{$(\xi, \eta)\in \widetilde{\mathrm{PE}}$;} \\[2pt]
   %               (\hat{u}_r, \hat{v}_r, \hat{\tau}_r)(\xi), & \hbox{$(\xi, \eta)\in \widetilde{\mathrm{EB}}$;} \\[2pt]
    %           (\bar{u}, \bar{v}, \bar{\tau})(\theta), & \hbox{$(\xi, \eta)\in \widetilde{\mathrm{PD}}$;} \\[2pt]
     %            (u_{w}, v_{w}, \tau_{w}), & \hbox{$(\xi, \eta)\in \widetilde{\mathrm{DJ}}$,}
      %         \end{array}
       %      \right.
%\end{equation}

%The main result can be sated as the following theorem
%\begin{thm}
%Assume $\tau_0-\tau_1^e$, $u_0-c_0$ and $\alpha_0-\alpha_w$ are sufficiently small. Then the boundary value problem
%(\ref{PEC}), (\ref{100905}) admits a piecewise solution in a region bounded by $\widetilde{\mathrm{PE}}\cup\widetilde{\mathrm{EB}}$, $ \widetilde{\mathrm{PD}}\cup\widetilde{\mathrm{DJ}}$, $\partial\mathcal{W}_{+}$, and a vacuum boundary $\widetilde{BH}$; see Figure \ref{Fig9.1}.
%\end{thm}

%\subsubsection{\bf Interaction of fan waves}

We first consider
 (\ref{PEC}) with the boundary conditions
\begin{equation}\label{100906}
(u, v, \tau)=\left\{
               \begin{array}{ll}
                  (\hat{u}_l, \hat{v}_l, \hat{\tau}_l)(\xi), & \hbox{$(\xi, \eta)\in \widetilde{\mathrm{PE}}$;} \\[2pt]
               (\bar{u}, \bar{v}, \bar{\tau})(\theta), & \hbox{$(\xi, \eta)\in \widetilde{\mathrm{PF}}$.}
               \end{array}
             \right.
\end{equation}
%Problem (\ref{PEC}), (\ref{100906}) is a standard Goursat problem.

\begin{lem}
Assume that $\tau_1^e-\tau_0$ is sufficiently small. Then, the Goursat problem (\ref{PEC}), (\ref{100906}) admits a classical solution on a curvilinear quadrilateral domain $\overline{\Sigma^6}$ bounded by $\widetilde{\mathrm{PE}}$, $\widetilde{\mathrm{PF}}$, $\widetilde{\mathrm{ES}}$, and $\widetilde{\mathrm{FS}}$, where $\widetilde{\mathrm{FS}}$ is a forward $C_{+}$ characteristic curve emanating from $\mathrm{F}$ and $\widetilde{\mathrm{ES}}$ is a forward $C_{-}$ characteristic curve emanating from $\mathrm{E}$. %Moreover, the solution satisfies
%\begin{equation}
%\Big\|\big(\tau-\tau_1^e, ~\alpha-\phi_{d}, ~\beta+\frac{\pi}{2}, ~\bar{\partial}_{+}\rho+\mathcal{L}_1,~  \bar{\partial}_{-}\rho+\mathcal{L}_2\big)\Big\|_{{0; \Sigma^6}}\rightarrow 0\quad \text{as} \quad \tau_1^e-\tau_0\rightarrow 0\quad \text{and}\quad \alpha_0-\alpha_w\rightarrow 0.
%\end{equation}
\end{lem}
\begin{proof}
Using the smallness of $\tau_1^e-\tau_0$ one can establish a priori estimates for the hyperbolicity of the system (\ref{PEC}) and $C^1$ norm of the solution. Then the global solution can be obtained by the method of characteristics.
\end{proof}

%=============
\begin{figure}[htbp]
\begin{center}
\includegraphics[scale=0.39]{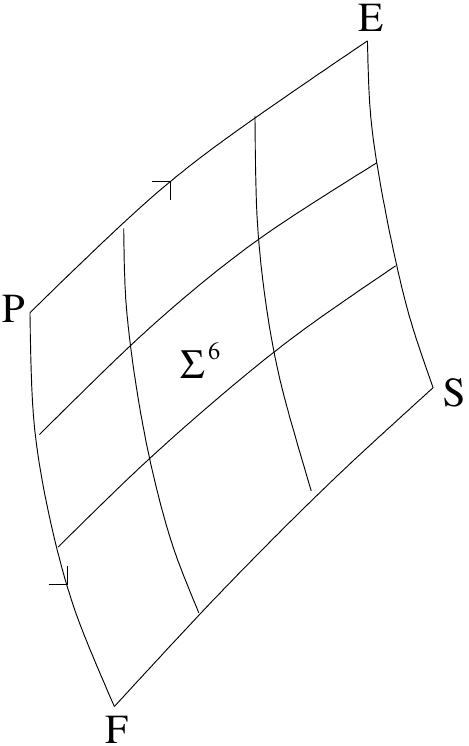}\qquad\qquad\qquad\includegraphics[scale=0.39]{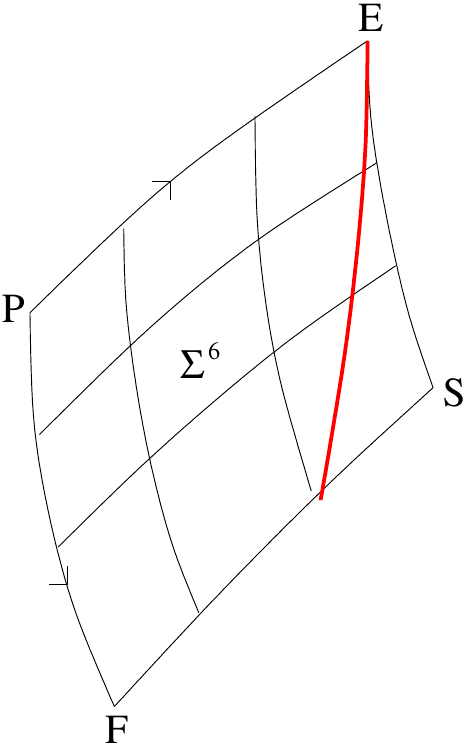}\qquad\qquad
\qquad\includegraphics[scale=0.39]{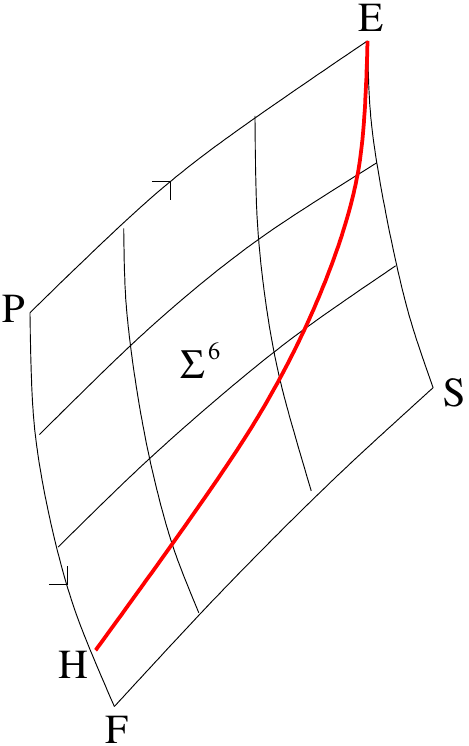}
\caption{\footnotesize Interaction region of the fan waves.}
\label{Fig9.0}
\end{center}
\end{figure}
%==========

For convenience, we use $(u, v, \tau)=(u_6, v_6, \tau_6)(\xi, \eta)$ to denote the solution to the Goursat problem (\ref{PEC}), (\ref{100906}).

\subsection{Post-sonic shock propagating from point E}
As in Section 6, the double-sonic shock of the fan–shock–fan composite wave propagates into the interior region $\Sigma^6$ from the point $\mathrm{E}$ and become a post-sonic shock.
For convenience, we denote this shock by $S_{_\mathrm{E}}$. The shock front $S_{_\mathrm{E}}$ can be determined applying the method described in Section 6.3. However,
unlike the interaction of the fan–shock–fan composite waves,
there are the following two cases:
\begin{itemize}
  \item there is a shock propagating from the point $\mathrm{F}$, and this shock intersects $S_{_\mathrm{E}}$;
  \item $S_{_\mathrm{E}}$ intersects the characteristic curve $\widetilde{\mathrm{PF}}$ at a point, and then propagates from there into the centered simple wave region of the fan–shock composite wave.
\end{itemize}

%Let $\tau_2$ be the same as that defined in (\ref{def12}).
%Then we have $\tau_2=\psi_{pr}(\bar{\tau}(\phi_{pr}))<\tau_2^e$.
By the result of Section 8.3, we see that for a fixed $\tau_0<\tau_1^e$,  for sufficiently small  $\bar{\tau}(\phi_{pr})-\tau_0>0$, the shock $S_{_\mathrm{E}}$
 stays  in $\Sigma^6$ until it
intersects  $\widetilde{\mathrm{FS}}$; see Figure \ref{Fig9.0} (mid). We next show that the second case can occur in some cases.

We use the subscripts $f$ and $d$ to denote the states on the front and back sides of $S_{_\mathrm{E}}$, respectively.
We set
$L=U\cos \phi+V\sin \phi$ and  $N=-U\sin \phi+V\cos \phi$,
where the variable $\phi$ denotes the inclination angle of ${\it S}_{_{\mathrm{E}}}$.
\begin{thm}\label{111503}
%Let the value of $\bar{\tau}(\phi_{pr})>\tau_1^e$ be fixed.
Let the amplitude of the oblique fan–shock composite wave be fixed.
When $\tau_1^e-\tau_0$ is sufficiently small, the shock curve $S_{_\mathrm{E}}$ intersects   $\widetilde{\mathrm{PF}}$ at a point $\mathrm{H}$; see Figure \ref{Fig9.0} (right). Moreover, $\mathrm{H}$ approaches $\mathrm{P}$ as $\tau_1^e-\tau_0\rightarrow 0$.
\end{thm}
\begin{proof}
The shock curve $S_{_\mathrm{E}}$ can be represented by a function
$\xi=\xi_{s}(\eta)$, which satisfies $$\xi_{s}'(\eta)=\cot(\phi(\xi_s(\eta), \eta)), \quad \xi_{s}(\eta_{_\mathrm{E}})=\xi_{_\mathrm{E}}.$$
We describe the characteristic line $\widetilde{\mathrm{PD}}$  by a function $\xi=\xi_{-}(\eta)$, which satisfies
$$\xi_{-}'(\eta)=\cot(\beta_6(\xi_{-}(\eta), \eta)), \quad \xi_{-}'(\eta_{_\mathrm{P}})=\xi_0.$$
It can be verified that $\eta_{_\mathrm{E}}-\eta_{_\mathrm{P}}\rightarrow 0$ and $\xi_{_\mathrm{E}}-\xi_{0}\rightarrow 0$ as $\tau_1^e-\tau_0\rightarrow 0$.

For $\tau'\in (\tau_1^e, \bar{\tau}(\phi_{pr}))$, we let $\mathrm{F}'$ be the point on $\widetilde{\mathrm{PF}}$ such that $\tau(\mathrm{F}')=\tau'$.
From $\mathrm{F}'$,
 we draw a forward $C_{+}$ characteristic curve of the Goursat problem (\ref{PEC}), (\ref{100906}). This characteristic curve intersects $\widetilde{\mathrm{ES}}$ at a point $\mathrm{S}'$.
We denote by $\Sigma'$ the curvilinear quadrilateral domain bounded by $\widetilde{{\mathrm{PE}}}$, $\widetilde{{\mathrm{PF}'}}$, $\widetilde{{\mathrm{ES}'}}$, and $\widetilde{{\mathrm{F'S}'}}$.

As in Lemma \ref{lem7.2}, we have
\begin{equation}\label{110407}
\Big\|\big(\tau-\tau_1^e, ~\alpha_6-\phi_{d}, ~\beta_6+\frac{\pi}{2}, ~  \bar{\partial}_{-}\rho+\mathcal{L}_2\big)\Big\|_{{0; \Sigma'}}\rightarrow 0\quad \text{as} \quad \tau_1^e-\tau_0\rightarrow 0\quad~ \text{and}\quad ~ \tau'-\tau_1^e\rightarrow 0,
\end{equation}
where the constants $\mathcal{L}_2$  and $\phi_d$ are as defined in Lemma \ref{lem7.2}.
%as $\tau_1^e-\tau_0\rightarrow 0$ and $\tau_2^e-\tau_2'\rightarrow 0$.

Using (\ref{8}), (\ref{110407}), and (\ref{61108a}), there exists a small $\delta>0$, such that when $\tau_1^e-\tau_0<\delta$ and $\tau'-\tau_1^e=\delta$, we have
$$
\bar{\partial}_{-}\beta_6>-\frac{\sin (2A_1^e)}{\rho_1^e}\mathcal{L}_2+\frac{\tan A_1^e p''(\tau_1^e)}{4(c_1^e)^2(\rho_1^e)^4}\mathcal{L}_2\quad \text{along}\quad \widetilde{\mathrm{PF}'},
$$
and
$$
\begin{aligned}
\bar{\partial}_{s}\phi&= -\frac{\bar{\partial}_{s} m}{\rho_fL_f}+\frac{N_{f}}{\rho_{f}L_{f}q_{f}^2}(q_{f}^2+\tau_{f}^2 p'(\tau_{f}))\bar{\partial}_{s}\rho_{f}+\bar{\partial}_{s}\sigma_{f}\\&=-\frac{\bar{\partial}_{_\Gamma} m}{\rho_fL_f}+
\frac{N_f}{\rho_fL_fq_f^2}(q_f^2+\tau_f^2 p'(\tau_f))\bar{\partial}_{-}\rho_6+\frac{N_f}{q_f^2}+\bar{\partial}_{-}\Big(\frac{\alpha_6+\beta_6}{2}\Big)
\\&=-\frac{\bar{\partial}_{s} m}{\rho_fL_f}+\frac{\sin (2A_f)}{2\rho_f}\bar{\partial}_{-}\rho_6-\frac{N_f}{q_f^2}+\frac{\bar{\partial}_{-}\beta_6}{2}+\frac{p''(\tau_f)\tan A_f}{4c_f^2\rho_f^4}\bar{\partial}_{-}\rho_6+\frac{\sin^2 A_f}{c_f}
\\&=-\frac{\bar{\partial}_{s} m}{\rho_fL_f}+\frac{\sin (2A_f)}{2\rho_f}\bar{\partial}_{-}\rho_6+\frac{\bar{\partial}_{-}\beta_6}{2}+\frac{p''(\tau_f)\tan A_f}{4c_f^2\rho_f^4}\bar{\partial}_{-}\rho_6\\&<-\frac{\bar{\partial}_{_\Gamma} m}{\rho_fL_f}+\frac{\sin (2A_f)}{2\rho_f}\bar{\partial}_{-}\rho_6+\frac{\bar{\partial}_{-}\beta_6}{2}
\\&<-\frac{\sin (2A_1^e)}{\rho_1^e}\mathcal{L}_2+\frac{\tan A_1^e p''(\tau_1^e)}{8(c_1^e)^2(\rho_1^e)^4}\mathcal{L}_2\quad \text{along}\quad  S_{_\mathrm{E}},
\end{aligned}
$$
where $(u_f, v_f, \tau_f)(\xi, \eta)=(u_6, v_6, \tau_6)(\xi,\eta)$.

Therefore, for any fixed  $\tau'\in (\tau_1^e, \bar{\tau}(\phi_{pr}))$, when $\tau_1^e-\tau_0$ is sufficiently small, there exists an $\eta_{_\mathrm{H}}\in (\eta_{_{\mathrm{F}'}}, \eta_{_\mathrm{P}})$ such that
$\xi_{s}(\eta_{_\mathrm{H}})=\xi_{-}(\eta_{_\mathrm{H}})$, and $\xi_{s}(\eta)>\xi_{-}(\eta)$ for $\eta_{_\mathrm{H}}<\eta<\eta_{_\mathrm{P}}$.
Moreover, $\eta_{_\mathrm{H}}-\eta_{_\mathrm{P}}\rightarrow 0$ as $\tau_1^e-\tau_0\rightarrow 0$.
This completes the proof of the theorem.
\end{proof}

If $S_{_\mathrm{E}}$ intersects $\widetilde{\mathrm{PF}}$ at a point $\mathrm{H}$,
 then the post-sonic shock $S_{_\mathrm{E}}$ propagates from this point into the simple wave region of the oblique fan–shock composite wave  and intersects the pre-sonic shock at some point $\mathrm{G}$. For convenience, we denote this post-sonic shock by $\widehat{\mathrm{HG}}$.
As in (\ref{100903}), the shock front $\widehat{\mathrm{HG}}$ can be determined by
\begin{equation}\label{260901}
\bar{\partial}_s \phi=-\frac{\tau_f(N_f^2-c_f^2)}{\rho_f^2N_fL_f(\tau_f^2-\tau_b^2)}
\bar{\partial}_s\rho_f+\frac{N_{f}\cos^2 A_f}{\rho_{f}L_{f}}\bar{\partial}_s\rho_{f}-\frac{N_{f}}{q_{f}^2}+\bar{\partial}_s\sigma_{f},
\end{equation}
where $(u_f, v_f, \tau_f)(\xi,\eta)=(\bar{u}, \bar{v}, \bar{\tau})(\theta)$.
Then, the states on the backside of $\widehat{\mathrm{HG}}$ can be determined by the Rankine-Hugoniot conditions and the post-sonic condition.
 For convenience, we denote the shock propagating from $\mathrm{E}$ by $\widehat{\mathrm{EG}}$; that is, $\widehat{\mathrm{EG}}=\widehat{\mathrm{EH}}\cup \widehat{\mathrm{HG}}$. %We use $\Sigma_5$
 Then, the flow in the domain bounded by $\widetilde{\mathrm{PE}}$, $\widetilde{\mathrm{PH}}$, and $\widehat{\mathrm{EH}}$ can be represented by $(u, v, \tau)=(u_6, v_6, \tau_6)(\xi,\eta)$.
%Then the flow in $\Omega_5$ is
We next establish the flow downstream of the shock $\widehat{\mathrm{EG}}$.

%=============
\begin{figure}[htbp]
\begin{center}
\includegraphics[scale=0.37]{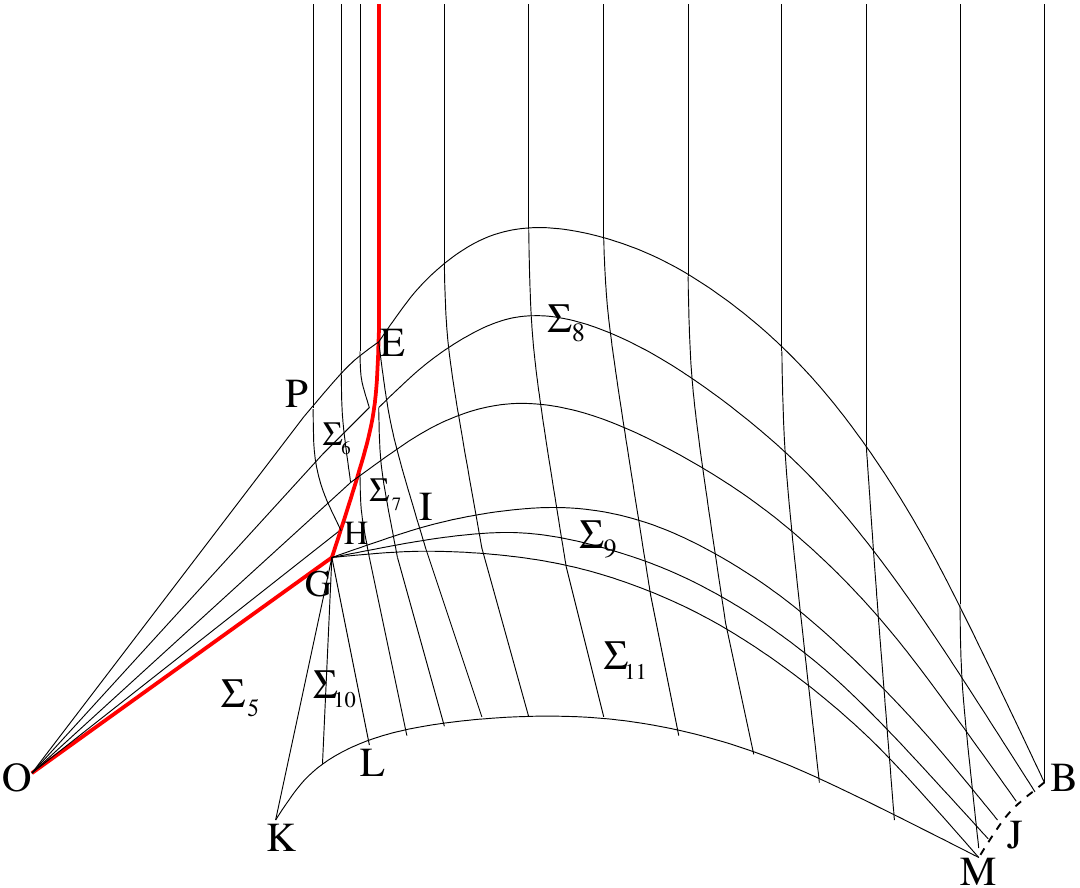}\qquad\qquad\includegraphics[scale=0.37]{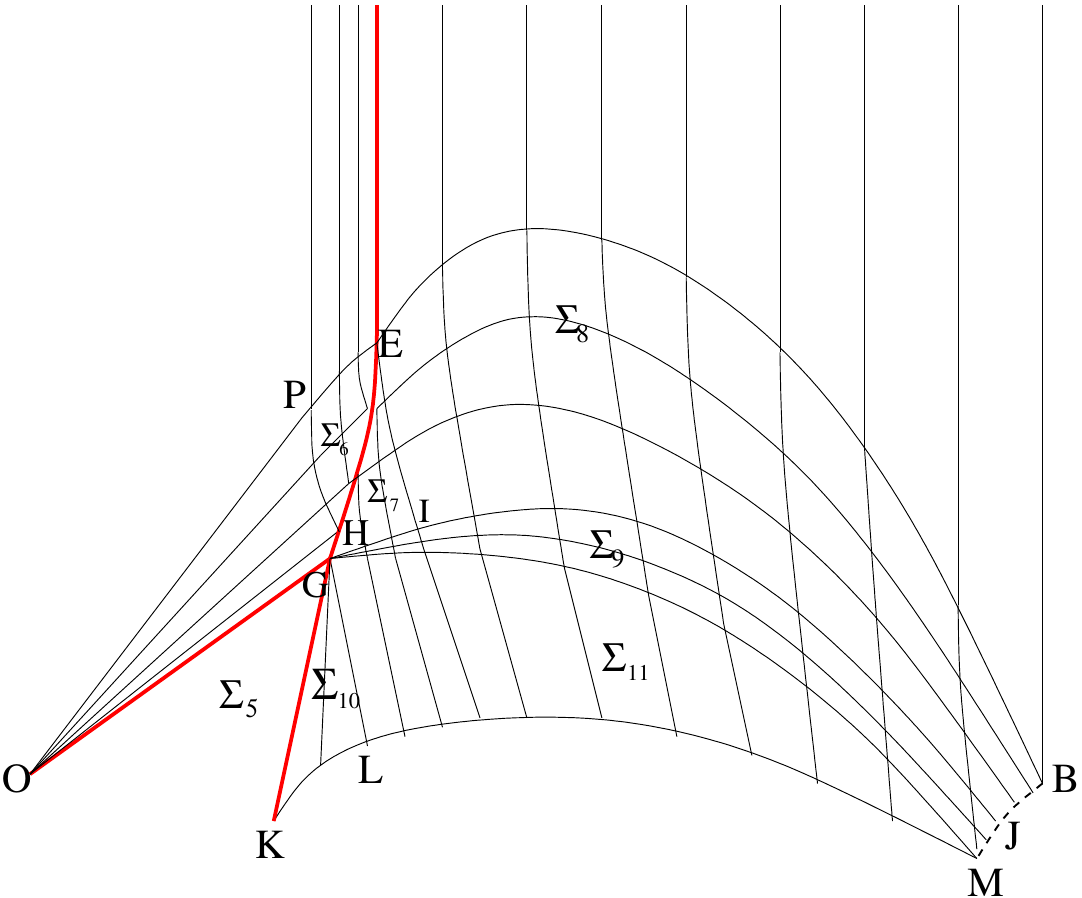}
\caption{\footnotesize Interaction between a planar fan–shock–fan composite wave and an oblique fan–shock composite wave.}
\label{Fig9.1a}
\end{center}
\end{figure}
%==========

As in Lemmas \ref{61302} and \ref{2662602}, we have the following structural conditions.
\begin{lem}\label{121101}
(Structural conditions)
Assume that $\tau_1^e-\tau_0$ is sufficiently small. Then, we have
\begin{equation}\label{120901a}
\cos\alpha_b \bar{\partial}_s u_b+\sin\alpha_b \bar{\partial}_s v_b<0\quad \text{and}\quad \cos\beta_b \bar{\partial}_s u_b+\sin\beta_b \bar{\partial}_s v_b>0\quad \text{along}\quad \widehat{\mathrm{EG}}.
\end{equation}
%where $r_{pr}=r(u_d, v_d)$, $s_d=s(u_d, v_d)$, and $(u_d, v_d)(x,y)$ represent the backside  states on $\widetilde{EG}$.
\end{lem}
%\begin{proof}
%As in Lemma \ref{110411}, when  $\tau_1^e-\tau_0$ is sufficiently small the desired estimates hold along the section $\widetilde{FH}$. As in Lemma \ref{120801a}, the desired estimates hold along $\widetilde{HG}$.
%\end{proof}

Using Lemma \ref{121101},
we can solve a
singular Cauchy problem for (\ref{PEC})
to obtain the flow downstream of the post-sonic shock $\widehat{\mathrm{EG}}$. This flow is defined
on a curvilinear triangular domain $\overline{\Sigma}_7$ bounded by $\widetilde{\mathrm{EI}}$, $\widehat{\mathrm{EG}}$, and $\widetilde{\mathrm{GI}}$, where $\widetilde{\mathrm{EI}}$ is a forward $C_{-}$
characteristic curve emanating from $\mathrm{E}$ and $\widetilde{\mathrm{GI}}$ is a forward $C_{+}$ characteristic curve emanating from $\mathrm{G}$. Moreover, the solution satisfies
$$
\bar{\partial}_{\pm}\rho<0 \quad \text{on}\quad \overline{\Sigma}_7\setminus\widehat{\mathrm{EG}}.
$$
We use the subscript $7$ to denote this flow. That is, $(u, v, \tau)=(u_7, v_7, \tau_7)(\xi,\eta)$, $(\xi,\eta)\in \overline{\Sigma}_7$.

Next, we solve a standard Gorusat problem for (\ref{PEC}) to obtain the flow in a curvilinear quadrilateral domain ${\Sigma}_8$ bounded by $\widetilde{\mathrm{EB}}$,  $\widetilde{\mathrm{EI}}$, $\widetilde{\mathrm{IJ}}$ and $\wideparen{\mathrm{BJ}}$, where $\widetilde{\mathrm{IJ}}$ is a forward $C_{+}$ characteristic line emanating from $\mathrm{I}$ and $\wideparen{\mathrm{BJ}}$ is a vacuum boundary. Moreover, the solution satisfies $\bar{\partial}_{\pm}\rho<0$ in $\overline{\Sigma}_8\setminus\wideparen{\mathrm{BJ}}$.
We use subscript $8$ to denote the flow in domain $\overline{\Sigma}_8\setminus\wideparen{\mathrm{BJ}}$. That is, $(u, v, \tau)=(u_8, v_8, \tau_8)(\xi,\eta)$, $(\xi,\eta)\in \overline{\Sigma}_8\setminus\wideparen{\mathrm{BJ}}$.

When $\tau_w\geq\tau_2^i$,
we draw a forward straight $C_{-}$ characteristic from the point $\mathrm{G}$ with the characteristic angle $\beta_5(\mathrm{G})$.
This characteristic line intersects $\partial{\mathcal{W}}_{+}$ or the sonic circle $(\xi-u_w)^2+(\eta-v_w)^2=c_w^2$ at a point $\mathrm{K}$.
%Let $K$ be an arbitrary point on this straight $C_{-}$ characteristic line.
We then consider (\ref{PEC}) with the boundary conditions
\begin{equation}\label{121102}
(u, v, \tau)=\left\{
               \begin{array}{ll}
                  (u_7, v_7, \tau_7)(\xi,\eta), & \hbox{$(\xi, \eta)\in\widetilde{\mathrm{GI}}$;} \\[2pt]
 (u_8, v_8, \tau_8)(\xi,\eta), & \hbox{$(\xi, \eta)\in\widetilde{\mathrm{IJ}}$;}  \\[2pt]
               (u_w, v_w, \tau_w), & \hbox{$(\xi, \eta)\in\overline{\mathrm{GK}}$.}
               \end{array}
             \right.
\end{equation}
Under some assumptions, the problem (\ref{PEC}), (\ref{121102})
admits a solution on a curved quadrilateral domain bounded by $\widetilde{\mathrm{GI}}\cup\widetilde{\mathrm{IJ}}$, $\overline{\mathrm{GK}}$, $C_{-}^{\mathrm{L}}$, and $C_{+}^{\mathrm{K}}$.
This solution consists of a $C_{+}$-type centered wave, a $C_{-}$-type centered simple wave,  and a non-planar $C_{-}$-type simple wave.
The centered wave flow region $\Sigma_9$ is bounded by $\widetilde{\mathrm{GI}}\cup\widetilde{\mathrm{IJ}}$, $\widetilde{\mathrm{GM}}$, and $\wideparen{\mathrm{JM}}$, where $\widetilde{\mathrm{GM}}$ is a forward $C_{+}$ characteristic line emanating from $\mathrm{G}$ and $\wideparen{\mathrm{JM}}$ is a vacuum boundary. The straight $C_{-}$ characteristic lines of the simple wave emanate from $\mathrm{G}$ and $\widetilde{\mathrm{GM}}$.

There is a slight difference in the discussion between the case $\tau_w<\tau_2^i$ and the case $\tau_w\geq\tau_2^i$.
When $\tau_w<\tau_2^i$, we need only replace the straight $C_{-}$ characteristic line emanating from $\mathrm{G}$ by a straight post-sonic 2-shock.
The state on the front side of the post-shock is $(u_w, v_w, \tau_w)$; the back-side state and the inclination angle of the post-sonic shock can be determined by the Rankine-Hugoniot conditions.

Therefore, we have the following conclusion about the interaction between the fan–shock–fan composite wave and the fan–shock composite wave.
\begin{itemize}
  \item If $\tau_w\geq\tau_2^i$, the interaction generates a simple wave; see Figure \ref{Fig9.1a} (left);
  \item if $\tau_w<\tau_2^i$, the interaction   generates a shock–fan composite wave; see Figure \ref{Fig9.1a} (right).
\end{itemize}

\subsection{Transonic shock propagating from point F}
%\subsubsection{Transonic shock free boundary problem}
Assume that $\tau_1^e-\tau_0$ is sufficiently small. Then, by Lemma \ref{111503}, we know that
there is a small $\epsilon>0$, depending only on  $\tau_1^e-\tau_0$,  such that
 the pre-sonic shock curve $S_{_\mathrm{E}}$ determined by (\ref{260901}) intersects the characteristic curve $\widetilde{\mathrm{FS}}$ if and only if $0<\bar{\tau}(\phi_{pr})-\tau_1^e<\epsilon$.
Moreover, $\epsilon\rightarrow0$ as $\tau_1^e-\tau_0\rightarrow 0$.
%In what follows, we are going to discuss the case of $0<\tau_2^e-\tau_2<\epsilon$.
We shall show that when $0<\bar{\tau}(\phi_{pr})-\tau_1^e<\epsilon$, the pre-sonic shock of the oblique fan–shock composite wave will propagate into the domain $\Omega^6$ from the point $\mathrm{F}$ and become a transonic shock.

From $\mathrm{F}$, we draw a straight $C_{-}$ characteristic line of the constant state $(u_w, v_w, \tau_w)$. This characteristic line intersects the wall $\partial\mathcal{W}_{+}$ or the sonic circle $(\xi-u_w)^2+(\eta-v_w)^2=c_w^2$ at a point $\mathrm{K}$.
%Let $K$ be an arbitrary point on this straight characteristic line.
To obtain the shock propagating from  $\mathrm{F}$ and the flow downstream of it, we consider (\ref{PEC}) with the following boundary conditions:
\begin{equation}\label{102904}
(u, v, \tau)=(u_w, v_w, \tau_w)\quad \text{on}~ \overline{\mathrm{FK}};
\end{equation}
\begin{equation}\label{102905}
(\rho  U -\rho_6U_6)(U -U_6)+  (\rho  V -\rho_6V_6)(V -V_6)=0\quad \text{on}\quad \eta=\eta_{t}(\xi),
\end{equation}
\begin{equation}\label{102906}
\left\{
  \begin{array}{ll}
    \displaystyle\frac{{\rm d}\eta_{t}}{{\rm d}\xi}=-\frac{U(\xi, \eta_t)-U_6(\xi, \eta_t)}{V(\xi, \eta_t)-V_6(\xi, \eta_t)}, & \hbox{$\xi>\xi_{_\mathrm{F}}$;} \\[8pt]
    \eta_t(\xi_{_\mathrm{F}})=\eta_{_\mathrm{F}}.
  \end{array}
\right.
\end{equation}
Here, $\eta=\eta_t(\xi)$, $\xi>\xi_{_\mathrm{F}}$ represents the shock propagating from $\mathrm{F}$.
%The problem (\ref{PEC}), (\ref{102904})--(\ref{102906}) is a shock free boundary problem.

%=============
\begin{figure}[htbp]
\begin{center}
\includegraphics[scale=0.37]{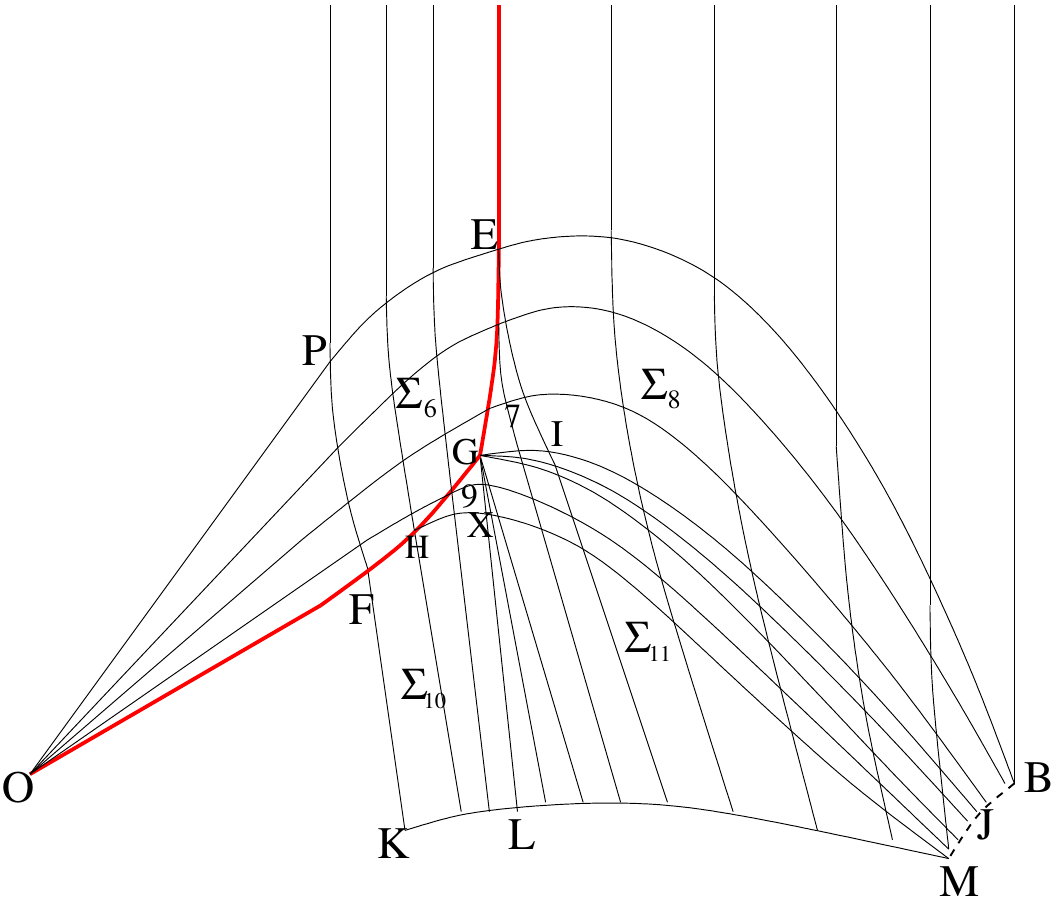}\qquad\quad \includegraphics[scale=0.37]{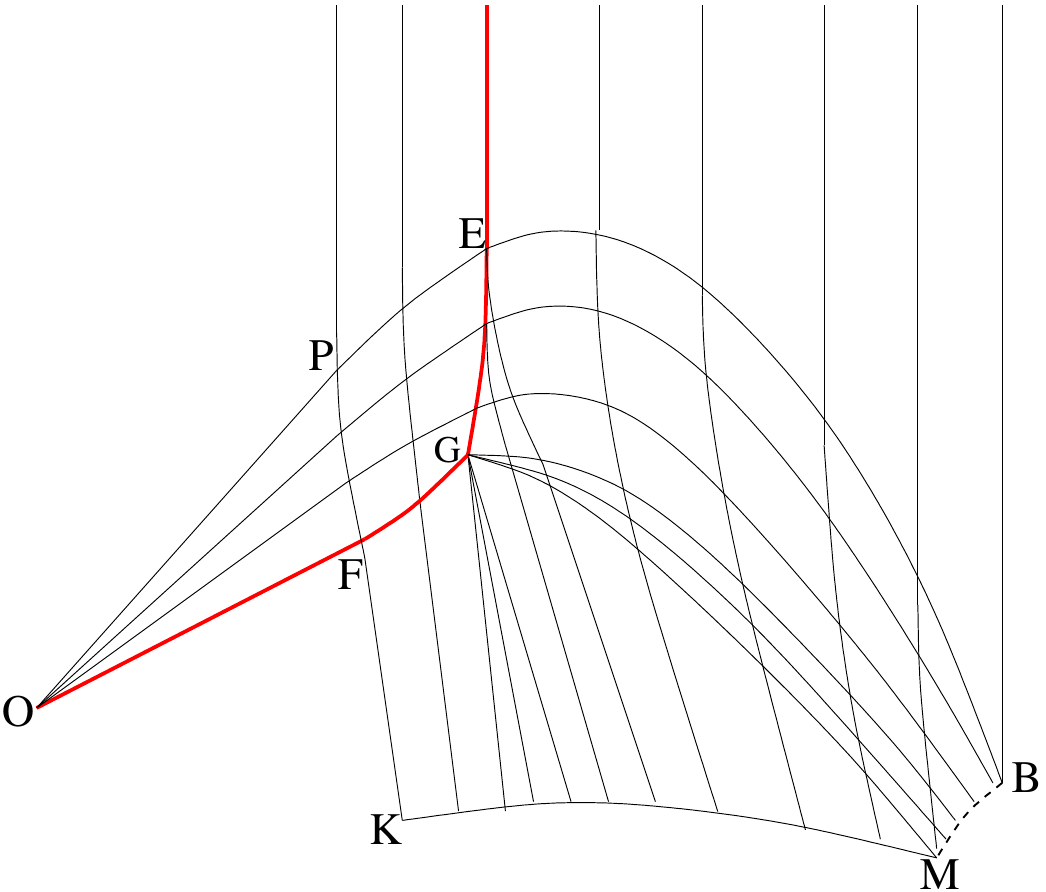}
\caption{\footnotesize Interaction between a planar shock–fan composite wave and an oblique fan–shock composite wave.}
\label{Fig9.2}
\end{center}
\end{figure}
%==========

\begin{thm}
The free boundary problem (\ref{PEC}), (\ref{102904})--(\ref{102906}) admits a (local) simple wave solution with straight $C_{-}$ characteristic lines.
Moreover, the simple wave solution satisfies $\bar{\partial}_{+}\rho<0$.
\end{thm}
\begin{proof}
We assume a priori that the free boundary problem admits a transonic shock solution. Then, for any point on the transonic shock front, the backward $C_{+}$ characteristic line emanating from this point goes back to $\overline{\mathrm{FK}}$.
By the characteristic decompositions (\ref{81104}), the solution of this shock free boundary problem must be a simple wave if it exists a transonic shock solution.

For convenience, we use $S_{_\mathrm{F}}$ to denoted the transonic shock propagating from $\mathrm{F}$.
Let $\phi$ be the inclination angle of $S_{_\mathrm{F}}$.
We set
$L=U\cos \phi+V\sin \phi$, $N=U\sin \phi- V\cos \phi$,
$L_6=U_6\cos \phi+V_6\sin \phi$, $N_6=U_6\sin \phi- V_6\cos \phi$, $m=\rho_6 N_6$ on $S_{_\mathrm{F}}$.

As in (\ref{250102}), we have
\begin{equation}\label{103002}
\bar{\partial}_s \phi=\frac{\bar{\partial}_s m}{\rho_fL_f}-\frac{N_f\cos^2 A_f}{\rho_fL_f}\bar{\partial}_s\rho_f+\frac{N_f}{q_f^2}+\bar{\partial}_s\Big(\frac{\alpha_f+\beta_f}{2}\Big)
\end{equation}
and
$$
\bar{\partial}_s \phi=\frac{\bar{\partial}_s m}{\rho_bL_b}-\frac{N_b\cos^2 A_b}{\rho_bL_b}\bar{\partial}_s\rho_b+\frac{N_b}{q_b^2}+\bar{\partial}_s\Big(\frac{\alpha_b+\beta_b}{2}\Big)\quad \text{along}\quad {\it S}_{_\mathrm{F}}.
$$
%along ${\it S}_{_\mathrm{F}}$.
Hence, we obtain
$$
\begin{aligned}
&(\tau_f-\tau_b)\frac{\bar{\partial}_s m}{L_f}-\frac{N_f\cos^2 A_f}{\rho_fL_f}\bar{\partial}_s\rho_f+\frac{N_f}{q_f^2}+\bar{\partial}_s\Big(\frac{\alpha_f+\beta_f}{2}\Big)
\\=&-\frac{N_b\cos^2 A_b}{\rho_bL_b}\bar{\partial}_s\rho_b+\frac{N_b}{q_b^2}+\bar{\partial}_s\Big(\frac{\alpha_b+\beta_b}{2}\Big)\quad \text{along}\quad {\it S}_{_\mathrm{F}}.
\end{aligned}
$$
Inserting (\ref{250104}) into this relation, we get
\begin{equation}\label{260902}
\begin{aligned}
&\frac{\tau_f(N_f^2-c_f^2)}{\rho_fL_fN_f(\tau_f+\tau_b)}\bar{\partial}_s\rho_f-\frac{N_f\cos^2 A_f}{\rho_fL_f}\bar{\partial}_s\rho_f+\frac{N_f}{q_f^2}+\bar{\partial}_s\Big(\frac{\alpha_f+\beta_f}{2}\Big)
\\~=~&\frac{\tau_b(N_b^2-c_b^2)}{\rho_bL_bN_b(\tau_f+\tau_b)}\bar{\partial}_s\rho_b-\frac{N_b\cos^2 A_b}{\rho_bL_b}\bar{\partial}_s\rho_b+\frac{N_b}{q_b^2}+\bar{\partial}_s\Big(\frac{\alpha_b+\beta_b}{2}\Big)\quad \text{along}\quad {\it S}_{_\mathrm{F}}.
\end{aligned}
\end{equation}

Since the downstream flow is a simple wave, by (\ref{7}) and (\ref{10}) we have
\begin{equation}\label{260903}
\begin{aligned}
\bar{\partial}_s\Big(\frac{\alpha_b+\beta_b}{2}\Big)&=
\frac{\sin(\phi-\beta_b)}{\sin(2A_b)}\bar{\partial}_+\Big(\frac{\alpha_b+\beta_b}{2}\Big)
\\&=\frac{\sin(\phi-\beta_b)}{\sin(2A_b)}\left(-\frac{\sin(2A_b)}{2\rho_b}\bar{\partial}_+\rho_b-\frac{\sin A_b}{q_b}\right)\\&=-\frac{\sin(2A_b)}{2\rho_b}\bar{\partial}_{s}\rho_b-\frac{\sin(\phi-\beta_b)}{2q_b\cos A_b}
\quad \text{along}\quad {\it S}_{_\mathrm{F}}.
\end{aligned}
\end{equation}

Combining (\ref{260902}) and (\ref{260903}), we get
\begin{equation}\label{103004}
\begin{aligned}
&\left(\frac{\tau_b(N_b^2-c_b^2)}{\rho_bL_bN_b(\tau_f+\tau_b)}-\frac{N_b\cos^2 A_b}{\rho_bL_b}-\frac{\sin (2A_b)}{2\rho_b}\right)\bar{\partial}_s\rho_b\\=&~
t_{+}^{f}\left(\frac{\tau_f(N_f^2-c_f^2)}{\rho_fL_fN_f(\tau_f+\tau_b)}-\frac{N_f\cos^2 A_f}{\rho_fL_f}-\frac{\sin (2A_f)}{2\rho_f}\right)\bar{\partial}_+\rho_f\\&~+~t_{-}^{f}\left(\frac{\tau_f(N_f^2-c_f^2)}
{\rho_fL_fN_f(\tau_f+\tau_b)}-\frac{N_f\cos^2 A_f}{\rho_fL_f}+\frac{\sin (2A_f)}{2\rho_f}\right)\bar{\partial}_{-}\rho_f\\&~+~\frac{N_f}{q_f^2}+(t_{-}^f-t_{+}^f)\frac{\sin^2 A_f}{c_f}+\frac{\sin(\phi-\beta_b)}{2q_b\cos A_b}\quad \text{along}\quad {\it S}_{_\mathrm{F}},
\end{aligned}
\end{equation}
where
$$
t_{+}^{f}=\frac{\sin(\phi-\beta_6)}{\sin(2A_6)}\quad \text{and}\quad t_{-}^{f}=\frac{\sin(\alpha_6 -\phi)}{\sin(2A_6)}.
$$

Since the flow downstream of $S_{_\mathrm{F}}$ is a simple wave,  $S_{_\mathrm{F}}$ can be determined by the density of the flow on the backside of it. Therefore, to obtain the transonic shock $S_{_\mathrm{F}}$ and the states on the backside of it, we consider (\ref{103002}) and (\ref{103004}) with the initial data
\begin{equation}\label{103005}
(\phi, \rho)(\xi_{_\mathrm{F}}, \eta_{_\mathrm{F}})=(\phi_{pr}, \rho_w).
\end{equation}
The problem (\ref{103002}), (\ref{103004}), and (\ref{103005}) is an initial value problem for a system of  ordinary differential equations. Since this initial value problem also usese the flow states in the domain $\Omega^6$, we need to establish an a priori estimate about the position relation between $S_{_\mathrm{F}}$ and $\widetilde{\mathrm{FS}}$.
We denote by $\eta=\eta_{+}(\xi)$, $\xi_{_{\mathrm{F}}}\leq \xi\leq \xi_{_{\mathrm{S}}}$ the $C_{+}$ characteristic curve $\widetilde{{\mathrm{FS}}}$.

We first establish an a priori estimate about the transonic shock front $S_{_\mathrm{F}}$.
%Since the shock is pre-sonic at the point $F$.
From (\ref{103004}) and $t_{-}^f(\mathrm{F})=0$, we have
\begin{equation}
\bar{\partial}_{s}\rho_b<0\quad \text{at}\quad  \mathrm{F}.
\end{equation}
Then, by (\ref{250104}),
we obtain
\begin{equation}
m\bar{\partial}_{s}m =\frac{\tau_b^3}{\bar{\tau}^2(\phi_{pr})-\tau_b^2}\left(-p'(\tau_b)+\frac{2h(\bar{\tau}(\phi_{pr}))-2h(\tau_b)}{\bar{\tau}^2(\phi_{pr})-\tau_b^2}\right)
\bar{\partial}_{s}\rho_b>0 \quad \text{at}\quad \mathrm{F}.
\end{equation}
%where $\tau_*=\psi_{pr}^{-1}(\tau_w)$.

As in (\ref{61108}), we have that at the point $\mathrm{F}$,
\begin{equation}
\begin{aligned}
\bar{\partial}_{s} \phi&=\frac{\bar{\partial}_{s} m}{\rho_6L_6}-\frac{N_6}{\rho_6L_6q_6^2}\big(q_6^2+\tau_6^2 p'(\tau_6)\big)\bar{\partial}_{s}\rho_6+\bar{\partial}_{s}\sigma_6+\frac{N_6}{q_6^2}\\&=\frac{\bar{\partial}_{s} m}{\rho_6L_6}-\frac{\sin (2A_6)}{2\rho_6}\bar{\partial}_{+}\rho_{6}+\frac{\bar{\partial}_{+}\alpha_{6}}{2}-\frac{p''(\tau_6)\tan A_6}{4c_6^2\rho_6^4}\bar{\partial}_{+}\rho_{6}\\&>-\frac{\sin (2A_6)}{2\rho_6}\bar{\partial}_{+}\rho_{6}+\frac{\bar{\partial}_{+}\alpha_{6}}{2}
\end{aligned}
\end{equation}
and
\begin{equation}
\bar{\partial}_+ \alpha_{6}<-\frac{\sin(2A_6)}{\rho_6}\bar{\partial}_{+}\rho_{6}.
\end{equation}
Hence, we have
\begin{equation}\label{103009}
\bar{\partial}_{s} \phi(\xi_{_\mathrm{F}}, \eta_{_\mathrm{F}})>\bar{\partial}_{+}\alpha_{6}(\xi_{_\mathrm{F}}, \eta_{_\mathrm{F}}).
\end{equation}
Then, there exists a sufficiently small $\varepsilon>0$ such that
$\phi(\xi, \eta_t(\xi))>\alpha_6(\xi, \eta_t(\xi))$ for $\xi_{_{\mathrm{F}}}<\xi<\xi_{_{\mathrm{F}}}+\varepsilon$.
This implies
%$y_{s}(x)>y_{+}(x)$ for $y_{_{F}}<x<x_{_{F}}+\varepsilon$.
$$
\eta_{t}(\xi)>\eta_{+}(\xi), \quad \xi_{_{\mathrm{F}}}<\xi<\xi_{_{\mathrm{F}}}+\varepsilon,
$$
and that the transonic shock $S_{_\mathrm{F}}$ exists.

For convenience, we denote the solution to (\ref{103002}), (\ref{103004}), and (\ref{103005}) by
$(\phi, \rho)=(\phi, \rho_b)(\xi, \eta)$, $(\xi, \eta)\in S_{_\mathrm{F}}$.
The states on the backside of $S_{_\mathrm{F}}$ can  be determined by
$$
u_b=L_b\cos\phi-N_b\sin\phi+\xi,\quad  v_b=L_b\sin\phi+N_b\cos\phi+\eta,
$$
where
$L_b=L_6$ and  $N_b=\frac{N_6}{\rho_b\tau_f}$.

The simple wave and its straight characteristic lines can then be determined by
the backside states of the transonic shock $S_{_\mathrm{F}}$.
This completes the proof of the theorem.
\end{proof}

Since $\epsilon\rightarrow0$ as $\tau_1^e-\tau_0\rightarrow 0$,
the radius of the domain $\Omega^6$ approaches $0$ as $\tau_1^e-\tau_0\rightarrow 0$ and $0<\bar{\tau}(\phi_{pr})-\tau_1^e<\epsilon$. Consequently, we obtain
\begin{equation}\label{103011}
\Big\|\big(\tau_6-\tau_1^e, ~\alpha_6-\phi_{d}, ~\beta_6+\frac{\pi}{2}, ~\bar{\partial}_{+}\rho_6+\mathcal{L}_1,~  \bar{\partial}_{-}\rho_6+\mathcal{L}_2\big)\Big\|_{{0; \Sigma^6}}\rightarrow 0
\end{equation}
as $0<\bar{\tau}(\phi_{pr})-\tau_1^e<\epsilon$ and $\tau_1^e-\tau_0\rightarrow 0$, where the constants $\mathcal{L}_1$, $\mathcal{L}_2$, and $\phi_d$ are as defined in Lemma \ref{lem7.2}.
Then, when  $\tau_1^e-\tau_0$ is sufficiently small and $0<\bar{\tau}(\phi_{pr})-\tau_1^e<\epsilon$, $\phi-\phi_{d}$ is small along $S_{_\mathrm{F}}$, and hence
\begin{equation}\label{2601}
\bar{\partial}_{s}\rho_{6}=t_{+}^f\bar{\partial}_{+}\rho_{6}+t_{-}^f\bar{\partial}_{-}\rho_{6}<0\quad \text{along}\quad S_{_\mathrm{F}}.
\end{equation}
Combining this with $\tau_{6}(\mathrm{F})=\bar{\tau}(\phi_{pr})>\tau_1^e$, we deduce that $S_{_\mathrm{F}}$ does not intersect the characteristic curve $\widetilde{\mathrm{PE}}$.

Moreover, by (\ref{103004}) and (\ref{103011}), we also have that when $\tau_1^e-\tau_0$ is sufficiently small and  $0<\bar{\tau}(\phi_{pr})-\tau_1^e<\epsilon$, the solution of the free boundary problem (\ref{PEC}), (\ref{102904})--(\ref{102906}) satisfies
\begin{equation}\label{12201}
\bar{\partial}_{s}\rho<0\quad \text{along}\quad S_{_\mathrm{F}}.
\end{equation}

We next assert that if $\phi-\alpha_b<0$ along  $S_{_\mathrm{F}}$, then $S_{_\mathrm{F}}$  does not intersect the characteristic curve $\widetilde{\mathrm{FS}}$.
%We next show that $S_{_F}$ does not intersect the characteristic line $\widetilde{FS}$.
Suppose that $S_{_\mathrm{F}}$ intersects $\widetilde{\mathrm{FS}}$ at a point. Then, there exists a $\xi'>\xi_{_F}$ such that
$$
\phi(\xi', \eta_t(\xi'))=\alpha_{6}(\xi', \eta_t(\xi')),\quad \text{and} \quad \phi(\xi, \eta_t(\xi))>\alpha_{6}(\xi, \eta_t(\xi))\quad \text{for}\quad \xi_{_F}<\xi<\xi'.
$$
However, by (\ref{2601}) and (\ref{12201}), we have
$$
\tau_b(\xi', \eta_t(\xi'))>\tau_w=\psi_{pr}(\tau_{6}(F))>\psi_{pr}(\tau_{6}(\xi', \eta_t(\xi'))),
$$
and hence $\phi(\xi', \eta_t(\xi'))>\alpha_{6}(\xi', \eta_t(\xi'))$.
%as in (\ref{103009}) we have
%$$
%\bar{\partial}_{_\Gamma} \phi(\xi_{0}, \eta_{t}(\xi_0))>\bar{\partial}_{+}\alpha_{3}(\xi_{0}, \eta_{t}(\xi_0)).
%$$
This leads to a contradiction.
This completes the proof of the assertion. %As described in the proof of Lemma \ref{91011}, we see that
% $S_{_F}$ does not intersect $\widetilde{PE}$ and $\widetilde{PF}$.

\subsection{Transition from transonic to post-sonic}
There are  two cases about the shock of the free boundary problem (\ref{PEC}), (\ref{102904})--(\ref{102906}):
\begin{itemize}            
  \item The shock $S_{_\mathrm{F}}$ stays in $\Omega^6$ and remains transonic until it intersects $S_{_\mathrm{E}}$ at a point $\mathrm{G}$, see Figure \ref{Fig9.2} (right);
  \item The shock $S_{_\mathrm{F}}$ stays in  $\Omega^6$ and remains transonic until there is a point $\mathrm{H}(\xi_{_\mathrm{H}}, \eta_t(\xi_{_\mathrm{H}}))$ on $S_{_\mathrm{F}}$ such that $\phi(\mathrm{H})-\alpha_b(\mathrm{H})=0$. Moreover, the transonic shock does not intersect $S_{_\mathrm{E}}$, see Figure \ref{Fig9.2} (left).
\end{itemize}

We show that the second subcase can happen in some cases.
\begin{thm}
(Transition from transonic to  post-sonic)
For a fixed small $\tau_1^e-\tau_0$,
when $\bar{\tau}(\phi_{pr})-\tau_1^e$ is sufficiently small, there exists a $\xi_{_\mathrm{H}}>\xi_{_\mathrm{F}}$  such that $(\phi-\alpha_b)(\xi_{_\mathrm{H}}, \eta_t(\xi_{_\mathrm{H}}))=0$, and $(\phi-\alpha_b)(\xi, \eta_t(x))<0$ for $\xi_{_\mathrm{F}}<\xi<\xi_{_\mathrm{H}}$.
\end{thm}
\begin{proof}
A direct computation yields that along ${\it S}_{_\mathrm{F}}$,
$$
0>\frac{\tau_b(N_b^2-c_b^2)}{\rho_bL_bN_b(\tau_f+\tau_b)}-\frac{N_b\cos^2 A_b}{\rho_bL_b}-\frac{\sin (2A_b)}{2\rho_b}>
-\frac{\sin (2A_b)}{2\rho_b}
$$
and
$$
\frac{\tau_f(N_f^2-c_f^2)}{\rho_fL_fN_f(\tau_f+\tau_b)}-\frac{N_f\cos^2 A_f}{\rho_fL_f}-\frac{\sin (2A_f)}{2\rho_f}<
-\frac{\sin (2A_f)}{2\rho_f}.
$$
Hence, by (\ref{103004}), we know that
when $\tau_1^e-\tau_0$ is sufficiently small and $0<\bar{\tau}(\phi_{pr})-\tau_1^e<\epsilon$,
\begin{equation}
\bar{\partial}_{s}\rho_b<-\frac{\tau_f\sin (2A_f)}{2\tau_b\sin (2A_b)}\mathcal{L}_1\quad \text{along}\quad S_{_\mathrm{F}}.
\end{equation}
From Proposition \ref{51003} and (\ref{2601}), we have
$$
\bar{\partial}_{s}\big(\rho_b-\psi_{po}(\tau_{6})\big)
=\bar{\partial}_{s}\rho_b+\psi_{po}'(\tau_{6})\bar{\partial}_{s}\tau_{6}<
\bar{\partial}_{s}\rho_b\quad \text{along}\quad S_{_\mathrm{F}}.
$$
We also have
$$\tau_b(\mathrm{F})-\psi_{po}(\tau_{6}(\mathrm{F}))=\tau_w-\psi_{po}(\bar{\tau}(\phi_{pr}))\rightarrow 0\quad \text{as}\quad \bar{\tau}(\phi_{pr})-\tau_1^e\rightarrow 0.$$
Hence,
 for fixed $\tau_0<\tau_1^e$, when $\bar{\tau}(\phi_{pr})-\tau_1^e$ is sufficiently small, there exists a $\xi_{_\mathrm{H}}>\xi_{_\mathrm{F}}$ such that
$\tau_b(\xi_{_\mathrm{H}}, \eta_t(\xi_{_\mathrm{H}}))-\psi_{po}(\tau_{6}(\xi_{_\mathrm{H}}, \eta_t(\xi_{_\mathrm{H}})))=0$. Moreover, $\xi_{_\mathrm{H}}\rightarrow \xi_{_\mathrm{F}}$ as $\tau_2^e-\tau_2\rightarrow 0$.

This completes the proof of the theorem.
\end{proof}

If the second case happens, we do not use (\ref{103002}) and (\ref{103004}) to look for the shock propagating from $\mathrm{H}$, since this would violate Liu's extended entropy condition. Instead, we look for a post-sonic shock.
As in Section 6, we can find a post-sonic shock propagated from $\mathrm{H}$, and this post-sonic shock stays in $\Omega^6$ until it intersects ${\it S}_{_\mathrm{E}}$ at a point $\mathrm{G}$.

%Therefore, for both cases, the shock waves $S_{_E}$ and $S_{_F}$ intersect at a point $G$ inside the domain $\Omega_{6}$; see Fig. \ref{Fig5.8.3}.
For convenience, we denote by $\wideparen{\mathrm{EG}}$ and $\wideparen{\mathrm{FG}}$ the shocks propagated from $\mathrm{E}$ and $\mathrm{F}$, respectively. The shock $\wideparen{\mathrm{EG}}$ is post-sonic. For the first case, the shock $\wideparen{\mathrm{FG}}$ is transonic.   For the second case, the shock $\wideparen{\mathrm{FG}}$ consists of a transonic part $\wideparen{\mathrm{FH}}$ and a post-sonic part $\wideparen{\mathrm{HG}}$.
In what follows we only discuss the second case, since the first case can be discussed similarly.

We first construct the flow downstream of the post-sonic shock $\widehat{\mathrm{EG}}$.
We solve a singular Cauchy problem to find the flow behind the post-sonic shock $\widehat{\mathrm{EG}}$, i.e., the flow on a triangular domain $\overline{\Omega}_7$ bounded by $\widehat{\mathrm{EG}}$, $\widetilde{\mathrm{EI}}$, and $\widetilde{\mathrm{GI}}$, where
$\widetilde{\mathrm{GI}}$ is a forward $C_{+}$ characteristic line emanating from $\mathrm{G}$ and
$\widetilde{\mathrm{EI}}$ is a forward $C_{-}$ characteristic line emanating from $\mathrm{E}$. We next solve a standard Goursat problem to obtain the flow  on a curvilinear quadrilateral domain $\overline{\Omega}_8$ bounded by $\widetilde{\mathrm{EB}}$,  $\widetilde{\mathrm{EI}}$, $\widetilde{\mathrm{IJ}}$ and $\widetilde{\mathrm{BJ}}$, where
 $\widetilde{\mathrm{IJ}}$ is a forward $C_{+}$ characteristic curve emanating from $\mathrm{I}$.

We next construct the flow downstream of the shock $\widehat{\mathrm{FG}}$.
There is a simple wave with straight $C_{-}$ characteristic lines emanating from the transonic shock $\widehat{\mathrm{FH}}$.
We solve a singular Cauchy problem to find the flow downstream the post-sonic shock $\widehat{\mathrm{HG}}$, i.e., the flow on a curvilinear triangular domain $\overline{\Omega}_{9}$ bounded by $\widehat{\mathrm{HG}}$, $\widetilde{\mathrm{HX}}$, and $\widetilde{\mathrm{GX}}$, where
$\widetilde{\mathrm{HX}}$ is a forward $C_{+}$ characteristic line emanating from $\mathrm{H}$ and
$\widetilde{\mathrm{GX}}$ is a forward $C_{-}$ characteristic line emanating from $\mathrm{G}$.
Moreover, when
$\tau_1^e-\tau_0$ is sufficiently small and  $0<\bar{\tau}(\phi_{pr})-\tau_1^e<\epsilon$, we have $\tau>\tau_2^i$ and $\bar{\partial}_{+}\tau>0$ on $\widetilde{\mathrm{HX}}$.
We then solve a simple wave problem to obtain a simple wave with straight $C_{-}$ characteristic lines emanating from $\widetilde{\mathrm{HX}}$.
%When $\tau_1^e-\tau_0$ is sufficiently small and  $0<\tau_2^e-\tau_2<\epsilon$, we have $\tau_w>\tau_2^i$.
The straight characteristic lines emanating from $\widetilde{\mathrm{FH}}$ and $\widetilde{\mathrm{HX}}$ do not intersect with each other.
Therefore, we establish the flow in a domain $\overline{\Omega}_{10}$ bounded by $\widehat{\mathrm{FH}}$, $\overline{\mathrm{HX}}$, $\overline{\mathrm{FK}}$, $\overline{\mathrm{XL}}$, and $\widetilde{\mathrm{KL}}$, where $\overline{\mathrm{XL}}$ is a forward straight $C_{-}$ characteristic line emanating from $\mathrm{X}$ and $\widetilde{\mathrm{KL}}$ is a forward $C_{+}$ characteristic line emanating from $\mathrm{K}$.

For convenience, we use the subscript $i$ to denote the flow states in domains $\overline{\Omega}_{i}$, i.e., $(u, v, \tau)=(u_i, v_i, \tau_i)(\xi, \eta)$, $(\xi, \eta)\in \overline{\Omega}_i$, $i=7,~8,~9,~10$.
Finally, we consider (\ref{PEC}) with the boundary conditions
\begin{equation}\label{121105}
(u, v, \tau)=\left\{
               \begin{array}{ll}
                  (u_7, v_7, \tau_7)(\xi,\eta), & \hbox{$(\xi, \eta)\in\widetilde{\mathrm{GI}}$;} \\[2pt]
(u_8, v_8, \tau_8)(\xi,\eta), & \hbox{$(\xi, \eta)\in\widetilde{\mathrm{IJ}}$;} \\[2pt]
               (u_9, v_9, \tau_9)(\xi,\eta), & \hbox{$(\xi, \eta)\in\widetilde{\mathrm{GX}}$;}\\[2pt]
   (u_{10}, v_{10}, \tau_{10})(\xi,\eta), & \hbox{$(\xi, \eta)\in\overline{\mathrm{XL}}$.}
               \end{array}
             \right.
\end{equation}
%The DGP (\ref{PEC}), (\ref{121105})
%admits a classical solution on a curvilinear quadrilateral domain $\Sigma_{11}$ bounded by $\widetilde{\mathrm{GI}}\cup\widetilde{\mathrm{IJ}}$, $\widetilde{\mathrm{GX}}\cup\overline{\mathrm{XL}}$, $\wideparen{\mathrm{JM}}$, and $\widetilde{\mathrm{LM}}$, where $\wideparen{\mathrm{JM}}$ is a vacuum boundary and
%$\widetilde{\mathrm{LM}}$ is a $C_{+}$ characteristic curve emanating from $\mathrm{L}$.
The discussion for the boundary value problem (\ref{PEC}), (\ref{121105}) is similar to that of the previous section, and we omit the details. This solves the interaction between the shock-fan-shock composite wave and the fan-shock composite wave.

%%%%%%%%%%%%%%%%%%%%%%%%%%%%%%%%%%%%%%%%%%%%%%%%%%%%%%%%%%%%%%%%%%%%%%%%%%%%%%%%%%%%%%%

\section{Interaction between a fan–shock–fan composite wave and a fan wave}

In this section, we consider the 2D Riemann problem (\ref{Potential}), (\ref{initial})--(\ref{boundary}) for $\tau_0<\tau_1^e$ and $u_0>c_0$.
Then, far away from the corner, the oncoming flow and the vacuum
are connected by a planar fan–shock–fan composite wave (\ref{pfsf}). We assume that near the corner, the oncoming flow is deflected by a centered simple wave (\ref{820a}) to a constant state $(u_w, v_w, \tau_w)$.

\begin{figure}[htbp]
\begin{center}
\includegraphics[scale=0.4]{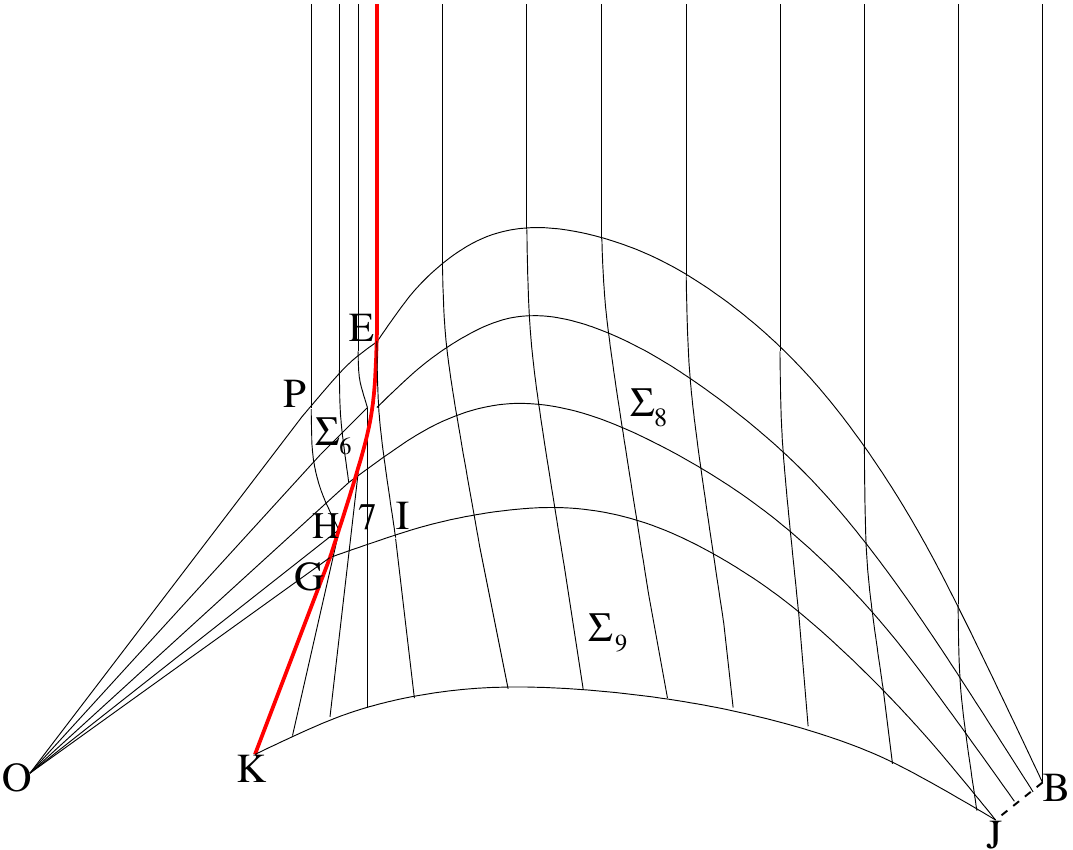}\qquad\quad\includegraphics[scale=0.4]{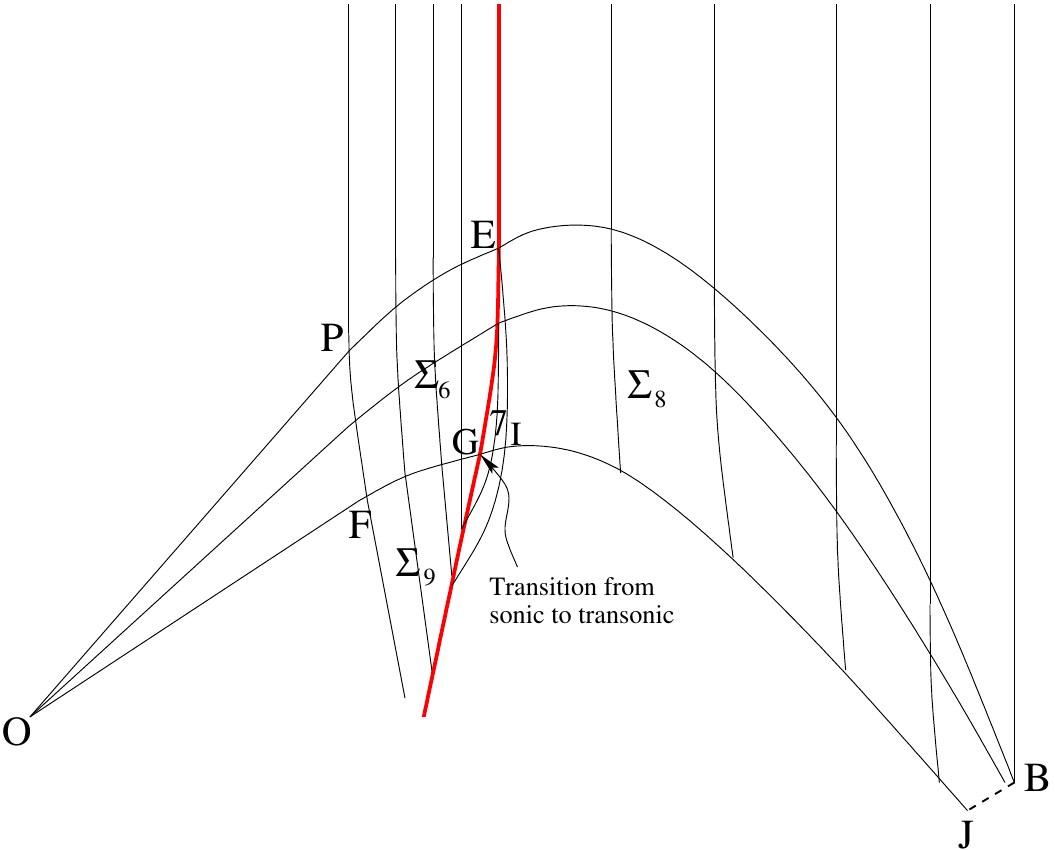}
\caption{\footnotesize Interaction between a fan–shock–fan composite wave and a fan wave.}
\label{Fig101}
\end{center}
\end{figure}

As shown in Figure \ref{Fig101}, the two composite waves begin to interact at point $\mathrm{P}$.
Through $\mathrm{P}$, we draw a $C_{+}$ cross characteristic curve in the planar fan–shock–fan composite wave and a $C_{-}$ cross characteristic curve in the oblique fan–shock composite wave.
The $C_{+}$ cross characteristic curve intersects the double-sonic shock of the planar fan–shock–fan composite wave at $\mathrm{E}$ and terminates at $\mathrm{B}$;
the $C_{-}$ cross characteristic curve intersects the pre-sonic shock of the oblique fan–shock composite wave at $\mathrm{F}$. 

We first consider a standard Goursat problem for (\ref{PEC}) with $\widetilde{\mathrm{PE}}$ and $\widetilde{\mathrm{PF}}$ as the characteristic boundaries. Assume that $\tau_1^e-\tau_0$ is sufficiently small. Then, as in Section 9, this Goursat problem admits a classical solution on a curvilinear quadrilateral domain ${\Sigma^6}$ bounded by $\widetilde{\mathrm{PE}}$, $\widetilde{\mathrm{PF}}$, $\widetilde{\mathrm{ES}}$, and $\widetilde{\mathrm{FS}}$, where $\widetilde{\mathrm{FS}}$ is a forward $C_{+}$ characteristic curve emanating  from $\mathrm{F}$ and $\widetilde{\mathrm{ES}}$ is a forward $C_{-}$ characteristic curve  emanating  from $\mathrm{E}$.
Moreover, the solution satisfies 
\begin{equation}\label{92601a}
  \bar{\partial}_{\pm}\rho<0\quad \text{and}\quad \tau<\tau_1^i.
\end{equation}
The double-sonic shock then propagates into the domain $\Sigma^6$ and becomes a post-sonic shock. 
As in Section 9, we have the following two cases:
\begin{enumerate}
  \item the post-sonic shock intersects $\widetilde{\mathrm{PF}}$ at a point $\mathrm{H}$; see Figure \ref{Fig101} (left);
  \item the post-sonic shock intersects $\widetilde{\mathrm{FS}}$ at a point $\mathrm{G}$; see Figure \ref{Fig101} (right).
\end{enumerate}

%=============
%==========

For the first case, 
 the post-sonic shock propagates into the centered simple region from $\mathrm{H}$ and intersects the tail characteristic of the centered simple wave at a point $\mathrm{G}$. By solving a degenerate Cauchy problem, we can find the flow downstream of the post-sonic shock $\widehat{\mathrm{EG}}$, i.e., the flow in a triangular domain $\Sigma_7$ bounded by $\widehat{\mathrm{EG}}$, $\widetilde{\mathrm{EI}}$, and $\widetilde{\mathrm{GI}}$, where $\widetilde{\mathrm{EI}}$ is a forward $C_{-}$ characteristic curve emanating from $\mathrm{\mathrm{E}}$ and $\widetilde{\mathrm{GI}}$ is a forward $C_{+}$ characteristic curve emanating from $\mathrm{G}$. We then solve a Goursat problem to obtain the flow in a domain $\Sigma_8$ bounded by the characteristic curves $\widetilde{\mathrm{EB}}$, $\widetilde{\mathrm{EI}}$, $\widetilde{\mathrm{IJ}}$ and  $\wideparen{\mathrm{BJ}}$, where $\widetilde{\mathrm{IJ}}$ is a forward $C_{+}$ characteristic curve emanating from $\mathrm{I}$ and $\wideparen{\mathrm{BJ}}$ is a vacuum boundary. From $\mathrm{G}$, we draw a straight post-sonic shock, which intersects the wall or the sonic circle $(\xi-u_w)^2+(\eta-v_w)^2=c_w^2$ at a point $\mathrm{\mathrm{K}}$. We then solve a simple wave problem to obtain a simple wave in a domain $\Sigma_9$ bounded by $\widetilde{\mathrm{GI}}$, $\widetilde{\mathrm{IJ}}$, $\overline{\mathrm{GK}}$, and $\widetilde{\mathrm{KJ}}$, where $\widetilde{\mathrm{KJ}}$ is a forward $C_{+}$ characteristic curve emanating from $\mathrm{K}$.
Therefore, when $(\tau_1^e-\tau_0)/(\tau_w-\tau_1^e)$ is sufficiently small, the interaction between the fan–shock–fan composite wave and the centered simple wave generates a shock–fan composite wave. 
 Since this discussion is similar to that in Section 9, we omit the details.

For the second case, 
the flow downstream of the post-sonic shock $\widehat{\mathrm{EG}}$ is similar to the first case. 
For the needs of the following discussion, we use $(u, v, \tau)=(u_7, v_7, \tau_7)(\xi, \eta)$ to denote the flow in the triangular domain $\Sigma_7$ bounded by $\widehat{\mathrm{EG}}$, $\widetilde{\mathrm{EI}}$, and $\widetilde{\mathrm{GI}}$.
By Lemma \ref{simple},  there is a simple wave with straight $C_{-}$ characteristic lines emanating from $\widetilde{\mathrm{FG}}$.
For convenience, we use the subscript $9$ to denote this simple wave. Then the shock propagating from $\mathrm{G}$ will interact with this simple wave. 
  We will show
 that the shock transitions to a transonic shock from $\mathrm{G}$. 
This will be done in two steps. We first show that the shock propagating from $\mathrm{G}$ cannot be a post-sonic shock. Then
we solve a free boundary problem to determine the transonic shock propagating from $\mathrm{G}$.

 \subsection{Impossibility of post-sonic shock propagating from point $\mathrm{G}$}
We use $\phi$ to denote the inclination angle of the shock propagating from $\mathrm{G}$, and the subscripts $f$ and $b$ to denote the states on its front side and backside, respectively. 
Hence, we have $(u_f, v_f, \tau_f)(\xi, \eta)=(u_9, v_9, \tau_9)(\xi,\eta)$ along the shock. 

Along the shock, we define
$$
t_{+}^{f}=\frac{\sin(\phi-\beta_f)}{\sin(2A_f)}\quad \text{and}\quad t_{-}^{f}=\frac{\sin(\alpha_f -\phi)}{\sin(2A_f)}.
$$
From (\ref{92601a}) we have
\begin{equation}\label{92601bb}
\bar{\partial}_{s}\rho_f=t_{+}^{f}\bar{\partial}_{+}\rho_9>0\quad \text{at}\quad  {\mathrm{G}},
\end{equation}
since $\phi<\beta_f$ at $\mathrm{G}$.

Using (\ref{7})--(\ref{82301}), (\ref{111306}), and (\ref{61108f}), we have
\begin{equation}\label{92001a}
\begin{aligned}
\bar{\partial}_s \phi+\frac{\bar{\partial}_s m}{\rho_fL_f}&=\frac{N_{f}}{\rho_{f}L_{f}q_{f}^2}\big(q_{f}^2+\tau_{f}^2 p'(\tau_{f})\big)\bar{\partial}_s\rho_{f}-\frac{N_{f}}{q_{f}^2}+\bar{\partial}_s\sigma_{f}\\&=
\frac{N_{f}\cos^2 A_f}{\rho_{f}L_{f}}t_{+}^{f}\bar{\partial}_{+}\rho_{9}-\frac{N_{f}}{q_{f}^2}+
t_{+}^{f}\bar{\partial}_{+}\Big(\frac{\alpha_{9}+\beta_{9}}{2}\Big)+\frac{t_{-}^{f}\sin^2 A_f}{c_f}
\\&=
\Big(-\frac{\sin(2A_f)}{2\rho_f}+\frac{N_{f}\cos^2 A_f}{\rho_{f}L_{f}}\Big)t_{+}^{f}\bar{\partial}_{+}\rho_{9}-\frac{N_{f}}{q_{f}^2}-t_{+}^{f}\frac{\sin^2 A_f}{c_f}+\frac{t_{-}^{f}\sin^2 A_f}{c_f}
\\&=
\Big(-\frac{\sin(2A_f)}{2\rho_f}+\frac{N_{f}\cos^2 A_f}{\rho_{f}L_{f}}\Big)t_{+}^{f}\bar{\partial}_{+}\rho_{9}
\\&=
\Big(-\frac{\sin(2A_f)}{2\rho_f}+\frac{N_{f}\cos^2 A_f}{\rho_{f}L_{f}}\Big)\bar{\partial}_{s}\rho_{f}\quad \text{at}\quad  {\mathrm{G}}.
\end{aligned}
\end{equation}

From (\ref{250104}) we have
\begin{equation}\label{202609302}
\begin{aligned}
\bar{\partial}_sm &~=~\frac{\tau_f}{\rho_fN_f(\tau_f^2-\tau_b^2)}\bigg(-p'(\tau_f)+\frac{2h(\tau_f)-2h(\tau_b)}{\tau_f^2-\tau_b^2}\bigg)
\bar{\partial}_s\tau_f
\\&~=~\frac{\tau_f(N_f^2-c_f^2)}{\rho_fN_f(\tau_f^2-\tau_b^2)}
\bar{\partial}_s\rho_f<0\quad \text{at}\quad {\mathrm{G}}.
\end{aligned}
\end{equation}
Inserting this into (\ref{92001a}) and using (\ref{92601bb}), we get
\begin{equation}\label{26101201}
\begin{aligned}
\bar{\partial}_s \phi+\frac{\bar{\partial}_s m}{\rho_bL_b}&=
\bar{\partial}_s \phi+\frac{\bar{\partial}_s m}{\rho_fL_f}-\left(\frac{1}{\rho_fL_f}-\frac{1}{\rho_bL_b}\right)\bar{\partial}_s m\\
~&=~\Big(-\frac{\sin(2A_f)}{2\rho_f}+\frac{N_{f}\cos^2 A_f}{\rho_{f}L_{f}}\Big)\bar{\partial}_{s}\rho_{f}-\frac{\tau_f(N_f^2-c_f^2)}{\rho_fN_fL_f(\tau_f+\tau_b)
}
\bar{\partial}_s\rho_f\\
~&=~\left(-\frac{\sin(2A_f)}{2\rho_f}+\frac{N_{f}\cos^2 A_f}{\rho_{f}L_{f}}-\frac{N_f^2-c_f^2}{2\rho_fN_fL_f}
+\frac{(\tau_b-\tau_f)(N_f^2-c_f^2)}{\rho_fN_fL_f(\tau_f+\tau_b)}\right)\bar{\partial}_s\rho_f\\&=\left(\frac{(\sin\delta_f\cos A_f-\cos\delta_f\sin A_f)^2}{\rho_f\sin(2\delta_f)}+\frac{(\tau_b-\tau_f)(N_f^2-c_f^2)}{\rho_fN_fL_f(\tau_f+\tau_b)}\right)\bar{\partial}_{s}\rho_{f}\\&>0\quad \text{at}\quad {\mathrm{G}},
\end{aligned}
\end{equation}
where $\delta_f=\arcsin\big(\frac{N_f}{q_f}\big)>A_f$.
Thus, we have
\begin{equation}\label{26100101}
\cos\beta_b \bar{\partial}_s u_b+\sin\beta_b \bar{\partial}_s v_b=~
-q_b\sin A_b \left(\bar{\partial}_s \phi+\frac{\bar{\partial}_s m }{\rho_bL_b}\right)<0
\quad \text{at}\quad {\mathrm{G}}.
\end{equation}
Suppose the shock propagating from $\mathrm{G}$ is a post-sonic shock. Then by (\ref{5906}) and (\ref{92601bb}) we have
$$
\bar{\partial}_{s}\rho_b<0\quad \text{at}\quad  {\mathrm{G}}.
$$
Thus, we have
\begin{equation}\label{26100102}
\begin{aligned}
&\cos\alpha_b \bar{\partial}_s u_b+\sin\alpha_b \bar{\partial}_s v_b\\=~&
q_b\sin A_b \left(\bar{\partial}_s \phi+\frac{\bar{\partial}_s m }{\rho_bL_b}\right) -
2\tau_bq_b\cos A_b\sin^2 A_b\bar{\partial}_s \rho_b>0\quad \text{at}\quad {\mathrm{G}}.
\end{aligned}
\end{equation}

The relations (\ref{26100101}) and (\ref{26100102}) imply that the shock propagating from $\mathrm{G}$ can not be a post-sonic, since it makes the corresponding degenerate Cauchy problem do not have a solution, i.e., there does not exist a downstream flow.
%We next derive  an a priori estimate for $\beta_b-\phi$ along the shock from $\mathrm{G}$, provided it exists and the downstream flow between it and $\widetilde{\mathrm{GI}}$ is smooth.

%From $\bar{\partial}_{-}u+\lambda_{+}\bar{\partial}_{-}v=0$ and $\phi=\beta_b$ at $\mathrm{G}$, we have
%$$
%\begin{aligned}
%\cos\alpha_b \bar{\partial}_s u_b+\sin\alpha_b \bar{\partial}_s v_b=~
%c_b \left(\bar{\partial}_s \phi+\frac{\bar{\partial}_s m }{\rho_bL_b}\right) -
%\tau_bc_b\sin(2 A_b)\bar{\partial}_s \rho_b=0\quad \text{at}\quad {\mathrm{G}}.
%\end{aligned}
%$$
%Thus, we obtain
%\begin{equation}\label{92601cc}
%\bar{\partial}_s \rho_b=\frac{\rho_b}{\sin (2A_b)}\left(\bar{\partial}_s \phi+\frac{\bar{\partial}_s m }{\rho_bL_b}\right)>0\quad \text{at}\quad {\mathrm{G}}.
%\end{equation}
%From (\ref{92601bb}) and (\ref{92601cc}) we see that the states on both sides of the shock propagating from $G$ can not be connected by a post-sonic shock. 

%From (\ref{8}), (\ref{202609302}), and (\ref{92601cc}), we have
%\begin{equation}\label{202609301}
%\begin{aligned}
 % \bar{\partial}_{s}(\beta_b-\phi)&=  \bar{\partial}_{-}\beta-  \bar{\partial}_{s}\phi
%\\&=\left(\frac{\sin(2A_b)}{\rho_b}-\frac{p''(\tau_b)\tan A_b}{2c_b^2\rho_b^4}\right) \bar{\partial}_{s}\rho_b-  \bar{\partial}_{s}\phi\\&=\frac{\bar{\partial}_sm}{\rho_bL_b}-\frac{p''(\tau_b)\tan A_b}{2c_b^2\rho_b^4}\bar{\partial}_{s}\rho_b<0
%\quad \text{at}\quad {\mathrm{G}}.
%\end{aligned}
%\end{equation}
%This guarantees that Lax's entropy condition holds along the shock.

 \subsection{Transition from post-sonic to transonic at point $\mathrm{G}$}
In the following discussion, we shall look for a transonic shock propagating from $\mathrm{G}$ by solving a free boundary problem. 
For convenience, we apply a coordinate rotation\footnote{It is easy to see that $\rho$ and $\bar{\partial}_{\pm}(\alpha, \beta, \rho)$ are invariant under coordinate rotations.} so that 
$$
\sigma_7(\mathrm{G})=\frac{\pi}{2}.
$$
Under this coordinate rotation,  $$\frac{\pi}{2}<\alpha_7(\mathrm{G})=\frac{\pi}{2}+A_7(\mathrm{G})<\pi\quad  \text{and}\quad 0<\beta_7(\mathrm{G})=\frac{\pi}{2}-A_7(\mathrm{G})<\frac{\pi}{2}.$$
Let 
$$
\phi_*=\arctan\left(\frac{u_{7}(\mathrm{G})-u_{6}(\mathrm{G})}{v_{6}(\mathrm{G})-v_{7}(\mathrm{G})}\right).
$$
Since the shock is post-sonic at $\mathrm{G}$, we have $\phi_*=\beta_7(\mathrm{G})$. We still use $(\xi, \eta)$ to denote the coordinates after this rotation. 

%=============
\begin{figure}[htbp]
\begin{center}
\includegraphics[scale=0.55]{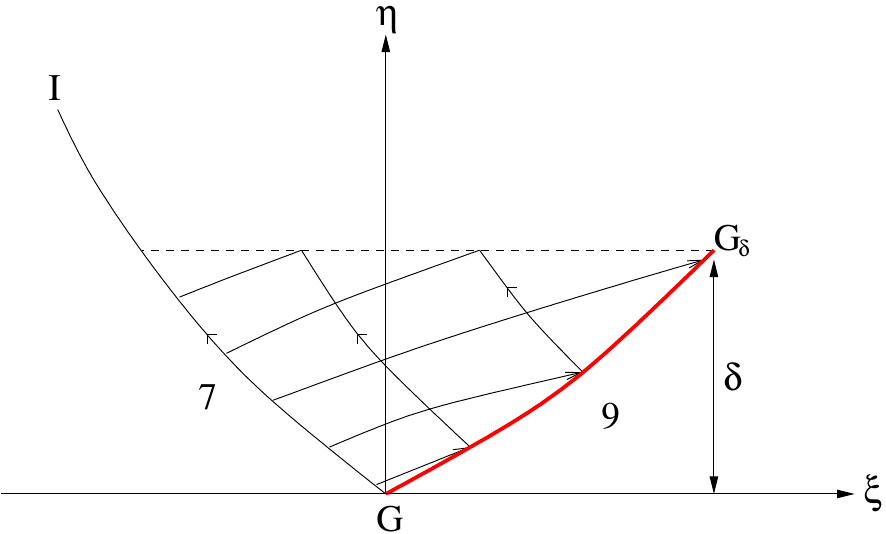}
\caption{\footnotesize  Local solution to the singulary free boundary problem.}
\label{STT}
\end{center}
\end{figure}
%==========

To obtain the transonic shock propagating from $\mathrm{G}$, 
we consider (\ref{PEC}) with the following boundary conditions:
\begin{equation}\label{102904a}
(u, v, \rho)=(u_7, v_7, \rho_7)(\xi, \eta)\quad \text{on}~ \widetilde{\mathrm{GI}};
\end{equation}
\begin{equation}\label{102905a}
(\rho  U -\rho_9U_9)(U -U_9)+  (\rho  V -\rho_9V_9)(V -V_9)=0\quad \text{on}\quad \xi=\psi(\eta);
\end{equation}
\begin{equation}\label{102906a}
\left\{
  \begin{array}{ll}
    \displaystyle\frac{{\rm d}\psi(\eta)}{{\rm d}\eta}=-\frac{V(\psi(\eta), \eta)-V_9(\psi(\eta), \eta)}{U(\psi(\eta), \eta)-U_9(\psi(\eta), \eta)}, & \hbox{$\eta>\eta_{_\mathrm{G}}$,} \\[8pt]
    \psi(\eta_{_\mathrm{G}})=\xi_{_\mathrm{G}}.
  \end{array}
\right.
\end{equation}
Here, $\xi=\psi(\eta)$, $\eta>\eta_{_\mathrm{G}}$ represents the shock propagating from $\mathrm{G}$.
The problem (\ref{PEC}), (\ref{102904a})--(\ref{102906a}) is a free boundary problem.

For convenience, we represent the $C_{+}$ characteristic curve $\widetilde{\mathrm{GI}}$ by the function $\eta=\varphi(\xi)$, where $\xi\in[\xi_{_\mathrm{I}}, \xi_{_\mathrm{G}}]$. Thus, $\varphi(\xi)\in C^1[\xi_{_\mathrm{I}}, \xi_{_\mathrm{G}}]\cap
C^2(\xi_{_\mathrm{I}}, \xi_{_\mathrm{G}}]$, and
$$
\alpha_7(\xi, \varphi(\xi))=\pi+\arctan\varphi'(\xi)\quad \text{and} \quad \bar{\partial}_{+}\alpha_7(\xi, \varphi(\xi))=\frac{\varphi''(\xi)\cos(\alpha_7(\xi, \varphi(\xi)))}{1+[\varphi'(\xi)]^2}\quad \text{along}\quad \widetilde{\mathrm{GI}}.
$$
The main result of this section is stated as follows. 

\begin{thm}\label{thm10.1}
(Transition from post-sonic to transonic)  There exists a small $\delta>0$ such that the free boundary problem  (\ref{PEC}), (\ref{102904a})--(\ref{102906a}) admits a Lipschitz-continuous solution on the domain
 $\Sigma_{10}(\delta)=\big\{(\xi, \eta)~|~\varphi^{-1}(\eta)<\xi<\psi(\eta), \eta_{_G}<\eta<\eta_{_G}+\delta\big\}$,
where $\varphi^{-1}$ represents the inverse function of $\varphi$. Moreover, the solution satisfies 
$\phi>\beta$ along $\xi=\psi(\eta)$, $\eta_{_\mathrm{G}}<\eta<\eta_{_\mathrm{G}}+\delta$. 
\end{thm}
The following two sonic singularity difficulties arise in the free boundary problem (\ref{PEC}), (\ref{102904a})--(\ref{102906a}).
\begin{itemize}
  \item The $C_{-}$ characteristic direction and the shock direction coincide at $\mathrm{G}$, i.e., $\phi=\beta$ at $\mathrm{G}$;
  \item For the boundary data (\ref{102904a}), $\bar{\partial}_{+}\rho_7=-\infty$ at $\mathrm{G}$.
\end{itemize}
Thus, 
the result by Li and Yu \cite{Li-Yu} on typical free boundary value problems for first-order quasilinear hyperbolic systems does not apply here. 
To overcome the difficulties arising from the sonic singularity at $\mathrm{G}$, we consider the following regularized problem.

For any given small positive constant $\varepsilon$,
let 
$$
\rho_{\varepsilon}=\rho_7(\mathrm{G})+\varepsilon,
$$
%where the function $\eta(\varepsilon)$, with $\eta(\varepsilon)>0$ and $\eta(\varepsilon)\rightarrow 0$ as $\varepsilon\rightarrow 0$, will be specified later, 
and let the state $(u_\varepsilon, v_\varepsilon)$ be defined by the Rankine-Hugoniot conditions
\begin{equation}\label{26100401}
\left\{
  \begin{array}{ll}
    (\rho_{\varepsilon}U_{\varepsilon} -\rho_9(\mathrm{G})U_9(\mathrm{G}))(U_{\varepsilon} -U_9(\mathrm{G}))+
(\rho_{\varepsilon}V_{\varepsilon} -\rho_9(\mathrm{G})V_9(\mathrm{G}))(V_{\varepsilon}-V_9(\mathrm{G}))=0, \\[4pt]
\displaystyle\frac{1}{2}q_{\varepsilon}^2+h(\tau_{\varepsilon})=\frac{1}{2}q_9^2(\mathrm{G})+h(\tau_9(\mathrm{G})),
  \end{array}
\right.
\end{equation}
where 
$$\tau_{\varepsilon}=\rho_{\varepsilon}^{-1}, \quad (U_{\varepsilon}, V_{\varepsilon})=(u_\varepsilon, v_\varepsilon)-(\xi_{_\mathrm{G}}, \eta_{_\mathrm{G}}),\quad
(U_{9}, V_{9})(\mathrm{G})=(u_{9}, v_{9})(\mathrm{G})-(\xi_{_\mathrm{G}}, \eta_{_\mathrm{G}}),$$ $$q_{\varepsilon}^2=U_{\varepsilon}^2+V_{\varepsilon}^2, \quad \text{and} \quad q_9^2(\mathrm{G})=U_{9}^2(\mathrm{G})+V_{9}^2(\mathrm{G}).$$

Let 
$$c_{\varepsilon}=c(\rho_{\varepsilon}),\quad \alpha_\varepsilon=\arctan\left(\frac{V_\varepsilon}{U_\varepsilon}\right)+\arcsin\left(\frac{c_\varepsilon}{q_{\varepsilon}}\right)$$ 
$$
\beta_\varepsilon=\arctan\left(\frac{V_\varepsilon}{U_\varepsilon}\right)-\arcsin\left(\frac{c_\varepsilon}{q_{\varepsilon}}\right),\quad \text{and} \quad 
\phi_{\varepsilon}=\arctan\left(-\frac{u_{9}(\mathrm{G})-u_\varepsilon}{v_{9}(\mathrm{G})-v_\varepsilon}\right).
$$
It is easy to check that
$$
(\alpha_\varepsilon, \beta_{\varepsilon})\rightarrow \big(\alpha_7(\mathrm{G}), \beta_7(\mathrm{G})\big)\quad \text{as}\quad \varepsilon\rightarrow 0.
$$
From (\ref{26100401}) we also have that for sufficiently small $\varepsilon>0$,
\begin{equation}\label{26100402}
\beta_\varepsilon<\phi_\varepsilon.
\end{equation}

Let $\varphi_{\varepsilon}(\xi)\in C^2[\xi_{_\mathrm{I}}, \xi_{_\mathrm{G}}]$  satisfy 
$$
\varphi_{\varepsilon}(\xi_{_\mathrm{G}})=\varphi(\xi_{_\mathrm{G}}),\quad
\varphi_{\varepsilon}'(\xi_{_\mathrm{G}})=\tan\alpha_{\varepsilon},\quad\text{and}\quad
\varphi_{\varepsilon}''(\xi)=\left\{
                               \begin{array}{ll}
                                 \varphi''(\xi), & \hbox{$\xi_{_\mathrm{I}}\leq \xi\leq\xi_{_\mathrm{G}}-\varepsilon $;} \\[4pt]
                                 \varphi''(\xi_{_\mathrm{G}}), & \hbox{$\xi_{_\mathrm{G}}-\varepsilon \leq \xi\leq \xi_{_\mathrm{G}}$.}
                               \end{array}
                             \right.
$$
We denote by $\widetilde{\mathrm{GI}_{\varepsilon}}$ the curve $\eta=\varphi_{\varepsilon}(\xi)$, where $\xi_{_\mathrm{I}}\leq\xi\leq \xi_{_\mathrm{G}}$. %We have $\lim\limits_{\varepsilon\rightarrow 0}\varphi_{\varepsilon}(\xi)= \varphi(\xi)$ in $C^1[\xi_{_\mathrm{I}}, \xi_{_\mathrm{G}}]$.

Let the function $(u_7^{\varepsilon}, v_7^{\varepsilon}, \rho_7^{\varepsilon})(\xi, \eta)$ be defined on $\widetilde{\mathrm{GI}_{\varepsilon}}$ such that $(u_7^{\varepsilon}, v_7^{\varepsilon}, \rho_7^{\varepsilon})(\mathrm{G})=(u_\varepsilon, v_\varepsilon, \rho_\varepsilon)$, 
$\widetilde{\mathrm{GI}_{\varepsilon}}$ is a $C_{+}$ characteristic curve
and the corresponding characteristic equations hold along it.
The function $(u_7^{\varepsilon}, v_7^{\varepsilon}, \rho_7^{\varepsilon})(\xi, \eta)$ can be determined by solving an initial value problem for a system of ``nonlinear ordinary differential equations" given by (\ref{10}), (\ref{11}), (\ref{72804}).
We have that for any small $\varepsilon'>0$, 
\begin{equation}\label{26101202}
 \left\|(\bar{\partial}_{+}\rho_7)(\xi, \varphi(\xi))-(\bar{\partial}_{+}\rho_7^{\varepsilon})(\xi, \varphi_{\varepsilon}(\xi))\right\|_{0; [\xi_{_\mathrm{I}}, \xi_{_\mathrm{G}}-\varepsilon']}~\rightarrow~ 0\quad \mbox{as}\quad \varepsilon\rightarrow 0.
\end{equation}

%Obviously, we have 
%\begin{equation}
%\|(u_7, v_7, \rho_7, \bar{\partial}_{+}\rho_7)-(u_7^{\varepsilon}, v_7^{\varepsilon}, \rho_7^{\varepsilon}, \bar{\partial}_{+}\rho_7^{\varepsilon})\|_{\widetilde{GI}}\rightarrow 0\quad\text{and}\quad 
%\beta_7^{\varepsilon}(G)-\phi_{\varepsilon}(G)\rightarrow0\quad\text{as}\quad \varepsilon\rightarrow 0.
%\end{equation}
We then replace the boundary condition (\ref{102904a}) by
\begin{equation}\label{102904an}
(u, v, \rho)=\big(u_7^{\varepsilon}, v_7^{\varepsilon}, \rho_7^{\varepsilon}\big)(\xi, \eta)\quad \text{on}\quad \widetilde{\mathrm{GI}_{\varepsilon}},
\end{equation}
and consider (\ref{PEC}) with the boundary conditions (\ref{102905a}),  (\ref{102906a}) and (\ref{102904an}).
Since $\beta_7^{\varepsilon}(\mathrm{G})=\beta_{\varepsilon}<\arcctg(\psi'(\eta_{_\mathrm{G}}))=\phi_{\varepsilon}$, 
the free boundary problem (\ref{PEC}), (\ref{102905a}),  (\ref{102906a}), (\ref{102904an}) is a typical free boundary problem and admits a local solution; see  Li and Yu \cite{Li-Yu}.

For the needs of the following discussion, we introduce the following constants:
$$
\begin{aligned}
&\qquad\qquad \qquad N_*=q_6(\mathrm{G})\sin(\sigma_{6}(\mathrm{G})-\phi_*), \quad L_*=q_6(\mathrm{G})\cos(\sigma_{6}(\mathrm{G})-\phi_*),
\\[3pt]&
\zeta_1=\Bigg\{\left(-\frac{\sin(2A_6)}{2\rho_6}+\frac{N_{*}\cos^2 A_6}{\rho_{6}L_{*}}-\frac{\tau_6(N_*^2-c_6^2)}{\rho_6N_*L_*(\tau_6+\psi_{po}(\tau_6))
}\right)\frac{\rho_7\sin(\phi_*-\beta_6)\bar{\partial}_{+}\rho_6}{\sin(\alpha_7-\phi_*)\sin (2A_6)}\Bigg\}\Bigg|_{\mathrm{G}}, 
\\[3pt]&
\zeta_2=
\frac{(N_*^2-c_6^2)\sin(\phi_*-\beta_6)\bar{\partial}_{+}\rho_6}{\rho_6^2N_*(\tau_6^2-\psi_{po}^2(\tau_6))\sin (2A_6)}\Bigg|_{\mathrm{G}}, \quad
\zeta_3=\Bigg\{\frac{p''(\tau_7)\tan A_7\sin(\alpha_7-\phi_*)\zeta_1}{2c_7^2\rho_7^4\sin(2A_7)}-\frac{\zeta_2}{\rho_7 L_*}\Bigg\}\Bigg|_{\mathrm{G}}.
\end{aligned}
$$
We have
$$
\zeta_1>0, \quad \zeta_2<0, \quad \mbox{and} \quad \zeta_3>0.
$$
We need to mention that $(\alpha_6, \beta_6, \rho_6, \bar{\partial}_{+}\rho_6)(\mathrm{G})=(\alpha_9, \beta_9, \rho_9, \bar{\partial}_{+}\rho_9)(\mathrm{G})$.
%$$
%\zeta_4=1+\frac{\tan(A_7(\mathrm{G}))}{c_7(\mathrm{G})}, \quad
%$$

\begin{lem}\label{lem10.2}
There is a small $\delta>0$ independent of $\varepsilon$ such that when $\varepsilon$ is sufficiently small, the free boundary problem (\ref{PEC}), (\ref{102905a}),  (\ref{102906a}), (\ref{102904an}) admits a local solution on a curvilinear triangular domain $\Sigma_{10}^{\varepsilon}(\delta)=\big\{(\xi, \eta)\mid \varphi_{\varepsilon}^{-1}(\eta)<\xi<\psi(\eta), \eta_{_\mathrm{G}}<\eta<\eta_{_\mathrm{G}}+\delta\big\}$. Moreover, there exist a small positive constant $\zeta_4$ and a positive constant $\zeta_5$, which are independent of $\delta$ and $\varepsilon$, such that the solution satisfies the following estimates:
\begin{equation}\label{2692801}
\frac{1}{2}\zeta_1<\bar{\partial}_{-}\rho<2\zeta_1\quad \text{and}\quad
-\frac{1}{\zeta_4d(\xi, \eta)}<\bar{\partial}_{+}\rho<0\quad \text{on}\quad 
\overline{\Sigma_{10}^{\varepsilon}(\delta)},
\end{equation}
\begin{equation}
\zeta_5<\bar{\partial}_{s}(\phi-\beta)<2\zeta_5\quad \text{along}\quad \widehat{\mathrm{GG}_{\delta}}.
\end{equation}
Here, $\widehat{\mathrm{GG}_{\delta}}$ represents the shock curve $\xi=\psi(\eta)$ ($\eta_{_\mathrm{G}}<\eta<\eta_{_\mathrm{G}}+\delta$), $\phi=\arcctg \left(\psi'(\eta)\right)$, $d(\xi, \eta)=\sqrt{(\xi-\xi_{_\mathrm{G}})^2+(\eta-\eta_{_\mathrm{G}})^2}$, and $\bar{\partial}_{s}=\cos\phi\partial_{\xi}+\sin\phi\partial_{\eta}$.
\end{lem}
\begin{proof}
Let
$$
\nu_1(\varepsilon)=\tau_{\varepsilon}^2\left(p'(\tau_{\varepsilon})
-\frac{2h(\tau_9(\mathrm{G}))-2h(\tau_{\varepsilon})}{\tau_9^2(\mathrm{G})-\tau_{\varepsilon}^2}\right).
$$
Then,
\begin{equation}\label{26101002}
N^2-c^2=\nu_1(\varepsilon)\quad \text{at}\quad \mathrm{G}.
\end{equation}
By computation, we have $$\nu_1(\varepsilon)\varepsilon^{-1}=O(1) \quad\text{as}\quad  \varepsilon\rightarrow 0,$$
where $O(1)$ denotes a bounded quantity.

A direct computation yields
$$
\begin{aligned}
N^2-c^2&=q^2(\sin^2(\sigma-\phi)-\sin^2(\sigma-\beta))\\
&=q^2[\sin(\sigma-\phi)+\sin(\sigma-\beta)][\sin(\sigma-\phi)-\sin(\sigma-\beta)]\quad \text{along}\quad \widehat{\mathrm{GG}_{\delta}}.
\end{aligned}
$$
Thus, we have
\begin{equation}
\left(\frac{\phi-\beta}{N^2-c^2}\right)(\mathrm{G})={O}(1)\quad \text{as}\quad \varepsilon\rightarrow 0,
\end{equation}

Let 
$$
\nu_2(\varepsilon)=(\bar{\partial}_{+}\rho_{7}^{\varepsilon})(\mathrm{G}).
$$
As in (\ref{26102404}), we have
\begin{equation}\label{26110301}
\nu_2(\varepsilon)\sqrt{\varepsilon}={O}(1)\quad \text{as}\quad \varepsilon\rightarrow 0.
\end{equation}

%Therefore, there exists a function $\eta(\varepsilon)$, with $\eta(\varepsilon)>0$ and $\eta(\varepsilon)\rightarrow 0$ as $\varepsilon\rightarrow 0$, such that 
%\begin{equation}\label{26101003}
%\nu_1^2(\varepsilon)\nu_2(\varepsilon)=o(1)\quad \text{as}\quad \varepsilon\rightarrow 0,
%\end{equation}
%where $o(1)$ denotes a quantity tending to zero.

As in (\ref{250104}), we have 
\begin{equation}\label{2026100501}
\begin{aligned}
\bar{\partial}_sm ~=~\frac{\tau_f(N_f^2-c_f^2)}{\rho_fN_f(\tau_f^2-\tau^2)}
\bar{\partial}_s\rho_f-
\frac{\tau(N^2-c^2)}{\rho_fN_f(\tau_f^2-\tau^2)}
(t_{+}\bar{\partial}_{+}\rho+t_{-}\bar{\partial}_{-}\rho)\quad \text{along}\quad \widehat{\mathrm{GG}_{\delta}},
\end{aligned}
\end{equation}
where 
$$
t_{+}=\frac{\sin(\phi-\beta)}{\sin(2A)}\quad \text{and}\quad t_{-}=\frac{\sin(\alpha -\phi)}{\sin(2A)}.
$$
Thus, by (\ref{26101002}) we have
\begin{equation}\label{26101004}
(\bar{\partial}_{s}m)(\mathrm{G})=\zeta_2+o(1),
\end{equation}
where $o(1)$ denotes a quantity tending to zero as $\varepsilon\rightarrow 0$.

Recalling (\ref{92001a}) and using (\ref{2026100501}), we have
\begin{equation}\label{2026100502}
\begin{aligned}
  \bar{\partial}_s \phi&~=~\left(-\frac{\sin(2A_f)}{2\rho_f}+\frac{N_{f}\cos^2 A_f}{\rho_{f}L_{f}}\right)\bar{\partial}_{s}\rho_{f}-
\frac{\bar{\partial}_s m}{\rho_fL_f}\\&~=~
\left(-\frac{\sin(2A_f)}{2\rho_f}+\frac{N_{f}\cos^2 A_f}{\rho_{f}L_{f}}-\frac{\tau_f(N_f^2-c_f^2)}{\rho_fN_f(\tau_f^2-\tau^2)}\right)\bar{\partial}_{s}\rho_{f}\\&\qquad\qquad\qquad+
\frac{\tau(N^2-c^2)}{\rho_fN_f(\tau_f^2-\tau^2)}
(t_{+}\bar{\partial}_{+}\rho+t_{-}\bar{\partial}_{-}\rho)\quad \text{along}\quad \widehat{\mathrm{GG}_{\delta}}.
\end{aligned}
\end{equation}

As in (\ref{26101001}), we have
$$
\begin{aligned}
  c \left(\bar{\partial}_s \phi+\frac{\bar{\partial}_s m }{\rho L}\right) &=
\tau c\sin(2A)\bar{\partial}_s \rho+\cos\alpha \bar{\partial}_s u+\sin\alpha \bar{\partial}_s v\\&=
\tau c\sin(2A)(t_{+}\bar{\partial}_{+}\rho+t_{-}\bar{\partial}_{-}\rho)+t_{+}(\cos\alpha \bar{\partial}_+ u+\sin\alpha \bar{\partial}_+ v)
\\&=
\tau c\sin(2A)t_{-}\bar{\partial}_{-}\rho\quad \text{along}\quad \widehat{\mathrm{GG}_{\delta}}.
\end{aligned}
$$
This implies 
\begin{equation}\label{2026100503}
\bar{\partial}_{-}\rho=\frac{\rho}{\sin (\alpha-\phi)}\left(\bar{\partial}_s \phi+\frac{\bar{\partial}_s m }{\rho L}\right)
\quad \text{along}\quad \widehat{\mathrm{GG}_{\delta}}.
\end{equation}
Hence, using (\ref{2026100501}) and (\ref{2026100502}) and recalling (\ref{26101201}), we have
\begin{equation}\label{26101005}
(\bar{\partial}_{-}\rho)(\mathrm{G})=\zeta_1+o(1).
\end{equation}

Using (\ref{7}), (\ref{8}), and (\ref{2026100503}), we have
\begin{equation}\label{2026100504}
\begin{aligned}
&\bar{\partial}_s(\phi-\beta)\\~=~&\bar{\partial}_s\phi-t_{+}\bar{\partial}_{+}\beta-t_{-}\bar{\partial}_{-}\beta
\\~=~&\bar{\partial}_s\phi+\frac{p''\tan A}{2c^2\rho^4}t_{+}\bar{\partial}_{+}\rho+\frac{2t_{+}\sin^2 A}{c}-\left(\frac{\sin(2A)}{\rho}-\frac{p''(\tau)\tan A}{2c^2\rho^4}\right)t_{-}\bar{\partial}_{-}\rho
\\~=~&\left(\frac{p''(\tau)\bar{\partial}_{+}\rho}{4c^2\rho^4\cos^2 A}+\frac{\tan A}{c}\right)\sin(\phi-\beta)+\frac{p''(\tau)\tan A}{2c^2\rho^4}t_{-}\bar{\partial}_{-}\rho-\frac{\bar{\partial}_{s}m}{\rho L}\quad \text{along}\quad \widehat{\mathrm{GG}_{\delta}}.
\end{aligned}
\end{equation}
A direct computation yields
\begin{equation}\label{26101006}
  \left(\frac{p''(\tau)\tan A}{2c^2\rho^4}t_{-}\bar{\partial}_{-}\rho-\frac{\bar{\partial}_{s}m}{\rho L}\right)(\mathrm{G})=\zeta_3+o(1).
\end{equation}

From the second equation of (\ref{81104}), we have
\begin{equation}\label{2026100601}
  \bar{\partial}_{-}\left(\frac{1}{\bar{\partial}_{+}\rho}\right)=-  \frac{p''(\tau)}{4c^2\rho^4\cos^{2}A }-
\left(\frac{\sin2A}{c}+
    \frac{\tau^4p''(\tau)}{4c^2\cos^{2}A }(\mathcal{F}(A, \tau)-1)\bar{\partial}_{-}\rho\right)\frac{1}{\bar{\partial}_{+}\rho}.
\end{equation}
Integrating this along $C_{+}$ characteristic curves emanating from $\widehat{\mathrm{GI}_{\varepsilon}}$ and recalling (\ref{26110301}), we obtain that there exist a small positive constant $\delta$ independent of $\varepsilon$  and a positive constant $\zeta_4$ independent of $\varepsilon$ and $\delta$ such that the solution satisfies 
$$
\frac{1}{\bar{\partial}_{+}\rho}<\zeta_4 \left(\nu_2^{-1}(\varepsilon)-d(\xi, \eta)\right)\quad \mbox{on}\quad \Sigma_{10}^{\varepsilon}(\delta).
$$ 
Thus, 
\begin{equation}\label{26101101}
\frac{1}{\zeta_4 \left(\nu_2^{-1}(\varepsilon)-d(\xi, \eta)\right)}<\bar{\partial}_{+}\rho<0\quad \mbox{on}\quad \Sigma_{10}^{\varepsilon}(\delta).
\end{equation}

Combining (\ref{2026100504}), (\ref{26101006}),  and  (\ref{26101101}), we obtain that there exist a small positive constant $\delta$ independent of $\varepsilon$  and a positive constant $\zeta_5$ independent of $\varepsilon$ and $\delta$ such that the solution satisfies 
\begin{equation}\label{26101102}
\zeta_5<\bar{\partial}_s(\phi-\beta)<2\zeta_5\quad \text{along}\quad \widehat{\mathrm{GG}_{\delta}}.
\end{equation}

By (\ref{26101101}) and (\ref{26101102}), there exists a sufficiently large positive constant $\mathcal{M}$ independent of  $\varepsilon$ and $\delta$ such that the solution satisfies 
$$
\big|(N^2-c^2)t_{+}\bar{\partial}_{+}\rho\big|<\frac{\mathcal{M}\big(\eta+\nu_1(\varepsilon)\big)^2}{\eta-\nu_2^{-1}(\varepsilon)} \quad \text{along}\quad \widehat{\mathrm{GG}_{\delta}}.
$$
Thus, 
using (\ref{2026100501})--(\ref{26101005}) and (\ref{26101101})--(\ref{26101102}), we deduce that when $\delta$ is sufficiently small, 
$$
\frac{3\zeta_1}{4}<\bar{\partial}_{-}\rho<\frac{4\zeta_1}{3}\quad \text{along}\quad \widehat{\mathrm{GG}_{\delta}}.
$$

Integrating the first equation of (\ref{81104}) along the $C_{+}$ characteristic curves emanating from $\widehat{\mathrm{GG}_{\delta}}$, we obtain that  when $\delta$ is sufficiently small,
$$
\frac{\zeta_1}{2}<\bar{\partial}_{-}\rho<2\zeta_1\quad \text{on}\quad 
\overline{\Sigma_{10}^{\varepsilon}(\delta)}.
$$

This completes the proof of the lemma. 
\end{proof}

Using (\ref{11}), (\ref{72804}), (\ref{61513}), and Lemma \ref{lem10.2}, we see that there exists a positive function $\mathcal{C}(r)$, $0<r<\delta$, which independent of $\varepsilon$, such that the solution of the regularized free boundary problem (\ref{PEC}), (\ref{102905a}),  (\ref{102906a}), (\ref{102904an}) satisfies
$$
\max\left\{|\nabla u|, |\nabla v|, |\nabla\rho|\right\}<\mathcal{C}(r)\quad \text{in}\quad \overline{\Sigma_{10}^{\varepsilon}(\delta)}\setminus \mathbf{B}_{r}(\mathrm{G}),
$$
where $\mathbf{B}_{r}(\mathrm{G})=\{(\xi, \eta)\mid (\xi-\xi_{_\mathrm{G}})^2+(\eta-\eta_{_\mathrm{G}})^2<r^2\}$.
Therefore, by the Ascoli-Arzel\`{a} theorem and a standard diagonal argument, we obtain a Lipschitz-continuous solution to the free boundary problem  (\ref{PEC}), (\ref{102904a})--(\ref{102906a}). This completes the proof of Theorem \ref{thm10.1}.

\section{Interaction between a shock–fan composite wave and a fan wave}
\begin{figure}[htbp]
\begin{center}
\includegraphics[scale=0.55]{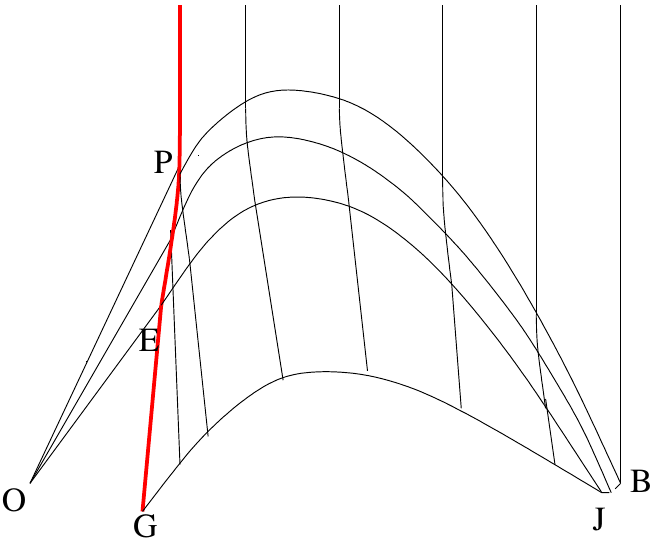}
\caption{\footnotesize Interaction of a shock–fan composite wave with a centered simple wave.}
\label{Fig8}
\end{center}
\end{figure}

We study the 2D Riemann problem (\ref{Potential}), (\ref{initial})--(\ref{boundary}) for $\tau_1^e\leq\tau_0<\tau_2^i$ and $u_0>\frac{\rho_{po}c_{po}}{\rho_0}$.
Then, away from the corner, the oncoming flow and the vacuum
are connected by a planar shock–fan composite wave which can be represented by (\ref{71501c}). In this section, we assume that at the corner the oncoming flow is deflected by a centered simple wave to a constant state $(u_{w}, v_{w}, \tau_{w})$, and the flow near the corner can be represented by (\ref{100901}). See Figure \ref{Fig8}.

%=============
%==========

The two waves start to interact with each other at a point $\mathrm{P}$. The post-sonic shock of the planar shock–fan composite wave will propagate into the centered simple wave region from the point $\mathrm{P}$. 
As in Section 8.1, we can solve a post-sonic shock-free boundary problem to obtain the post-sonic shock emanating from $\mathrm{P}$ and the flow downstream of it.
This post-sonic shock intersects the tail of the centered simple wave at a point $\mathrm{E}$. It will continue to propagate from the point $\mathrm{E}$ with a constant speed until it intersects the ramp $\partial\mathcal{W}_{+}$ or the sonic circle $(\xi-u_d)^2+(\eta-v_d)^2=c_d^2$ at a point $\mathrm{G}$. Then, solving a simple wave problem shows that the interaction of the shock–fan composite wave with the centered simple wave generates a non-planar shock–fan composite wave. Since the discussion of this section is similar to that of the last section, we omit the details.

%%%%%%%%%%%%%%%%%%%%%%%%%%%%%%%%%%%%%%%%%%%%%%%%%%%%%%%%%%%%%%%%%%%%%%%%%%%%%%%%%%%%%%

\section{On vacuum-free initial data}
%Motivated by the work of Li and Zheng \cite{LJ5} on the 2D Riemann problem with four Rarefaction Waves, 
The method in this paper also applies when the initial data do not contain a vacuum region.
Consider (\ref{Potential}) with the initial condition
 \begin{equation}\label{initial3}
(\tau, \Phi)(0,x,y)=\left\{
              \begin{array}{ll}
                (\tau_0, u_0x), & \hbox{$x<0$, $y>0$;} \\[2pt]
               (\tau_1, u_1x), & \hbox{$x>0$, $y\geq 0$,}\\[2pt]
 (\tau_1, u_1x), & \hbox{$y<0$, $x>y\cot\theta_w$,}
              \end{array}
            \right.
 \end{equation}
and the slip condition (\ref{boundary}), where $u_1>u_0>c_0$ and $\tau_1>\tau_0$. 
We assume that the gas is a polytropic gas for convenience.

We assume that $\theta_w>\pi+\alpha_v$, where the constant $\alpha_v$ is defined as in Section 3.3. Then, at the corner, the uniform oncoming flow $(0):=(u_0, 0, \tau_0)$ is deflected by a centered simple wave to a constant state $(u_{w}, v_{w}, \tau_{w})$.

Assume further that $u_1\sin\theta_w>\int_{\tau_1}^{+\infty}\sqrt{-p'(\tau)}{\rm d}\tau$. Then,  against the ramp, the constant state flow $(1):=(u_1, 0, \tau_1)$ expands into a vacuum through a planar rarefaction wave $R_4$, which can be represented by 
\begin{equation}\label{26102001}
\left\{
  \begin{array}{ll}
    u=\bar{v}(\zeta)\sin\theta_w+u_1\cos^2\theta_w,  \\[2pt]
    v=-\bar{v}(\zeta)\cos\theta_w+u_1\cos\theta_w\sin\theta_w,\\[2pt]
    \tau=\bar{\tau}(\zeta),
  \end{array}
\right.
\quad\zeta_1<\zeta<\zeta_2,
\end{equation}
where $\zeta=\xi\sin\theta_w-\eta\cos\theta_w$, $\zeta_1=u_1\sin\theta_w-\int_{\tau_1}^{+\infty}\sqrt{-p'(\tau)}{\rm d}\tau$, $\zeta_2=u_1\sin\theta_w+c_1$, and 
the functions $\bar{u}(\zeta)$ and $\bar{\tau}(\zeta)$ are determined by 
$$
\zeta=\bar{u}(\zeta)+\bar{\tau}(\zeta)\sqrt{-p'(\bar{\tau}(\zeta))}, \quad \bar{u}(\zeta)+\int_{\tau_1}^{\bar{\tau}(\zeta)}\sqrt{-p'(\tau)}{\rm d}\tau=u_1\sin\theta_w.
$$

  %======================
\begin{figure}[htbp]
\begin{center}
\qquad\qquad\includegraphics[scale=0.45]{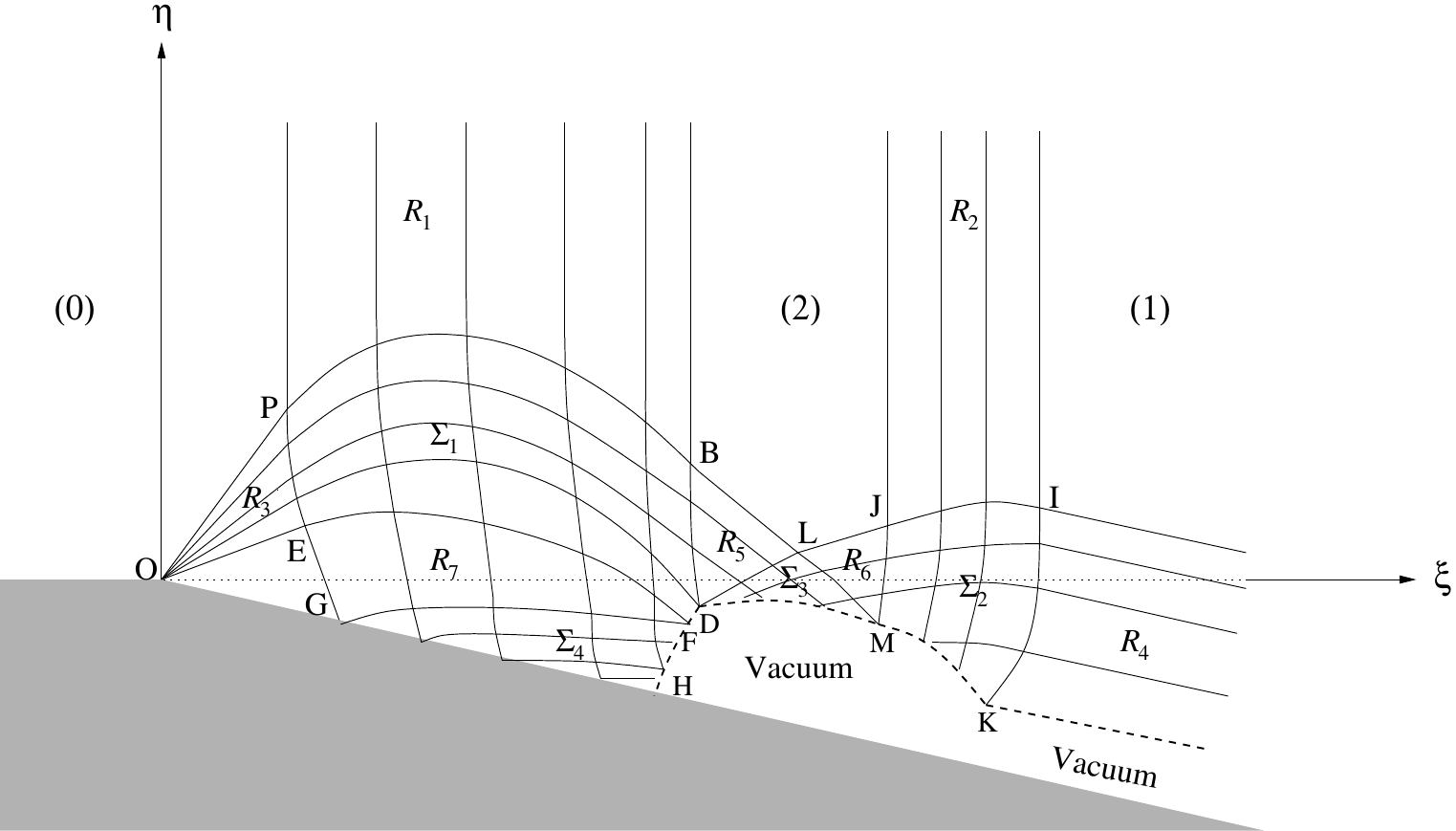}
\caption{\footnotesize Interaction of four rarefaction waves.}
\label{A3}
\end{center}
\end{figure}
%==========

%\subsection{Interaction of four rarefaction waves} 
Assume that 
$$
 u_0+\int_{\tau_0}^{\tau_1}\sqrt{-p'(\tau)}{\rm d}\tau<u_1<u_0+\int_{\tau_0}^{+\infty}\sqrt{-p'(\tau)}{\rm d}\tau+
\int_{\tau_1}^{+\infty}\sqrt{-p'(\tau)}{\rm d}\tau.
$$
Then by solving a 1D Riemann problem we see that away from the wedge, the constant states $(0)$ and $(1)$ are connected by two planar rarefaction waves $R_1$ and $R_2$  with a middle state $(u_2, 0, \tau_2)$, determined by 
$$
u_2=u_0+\int_{\tau_0}^{\tau_2}\sqrt{-p'(\tau)}{\rm d}\tau=u_1-
\int_{\tau_1}^{\tau_2}\sqrt{-p'(\tau)}{\rm d}\tau.
$$

Let us describe, without rigorous proof, a wave structure in the $(\xi, \eta)$-plane, as illustrated in Figure \ref{A3}. The fan waves $R_1$ and $R_3$ interact in a curvilinear quadrilateral domain $\Sigma_1$, which is bounded by the characteristic curves $\widetilde{\mathrm{PE}}$, $\widetilde{\mathrm{PB}}$, $\widetilde{\mathrm{BD}}$, and $\widetilde{\mathrm{EF}}$, together with a vacuum boundary $\wideparen{\mathrm{DF}}$.
A rarefaction simple wave $R_7$ emanates from the characteristic curve $\widetilde{\mathrm{EF}}$ and reflects at the ramp, forming a reflection region $\Sigma_4$ bounded by the characteristic curve $\widetilde{\mathrm{GF}}$, the ramp, and a vacuum boundary $\wideparen{\mathrm{FH}}$.
The fan waves $R_2$ and $R_4$ interact in a curvilinear quadrilateral domain $\Sigma_2$, which is bounded by the characteristic curves $\widetilde{\mathrm{IJ}}$, $\widetilde{\mathrm{IK}}$, and $\widetilde{\mathrm{JM}}$, together with a vacuum boundary $\wideparen{\mathrm{KM}}$.
Two non-planar rarefaction simple waves $R_5$ and $R_6$ emanates from $\widetilde{\mathrm{BD}}$ and $\widetilde{\mathrm{JM}}$, respectively. The simple waves $R_5$ and $R_6$ interact in a curvilinear triangular domain $\Sigma_3$
bounded by characteristic curves $\widetilde{\mathrm{LM}}$ and $\widetilde{\mathrm{LD}}$, together with a vacuum boundary $\wideparen{\mathrm{DM}}$. 
All of the above interactions are obtainable by solving standard Goursat problems.

 %======================
\begin{figure}[htbp]
\begin{center}
\qquad \qquad \includegraphics[scale=0.45]{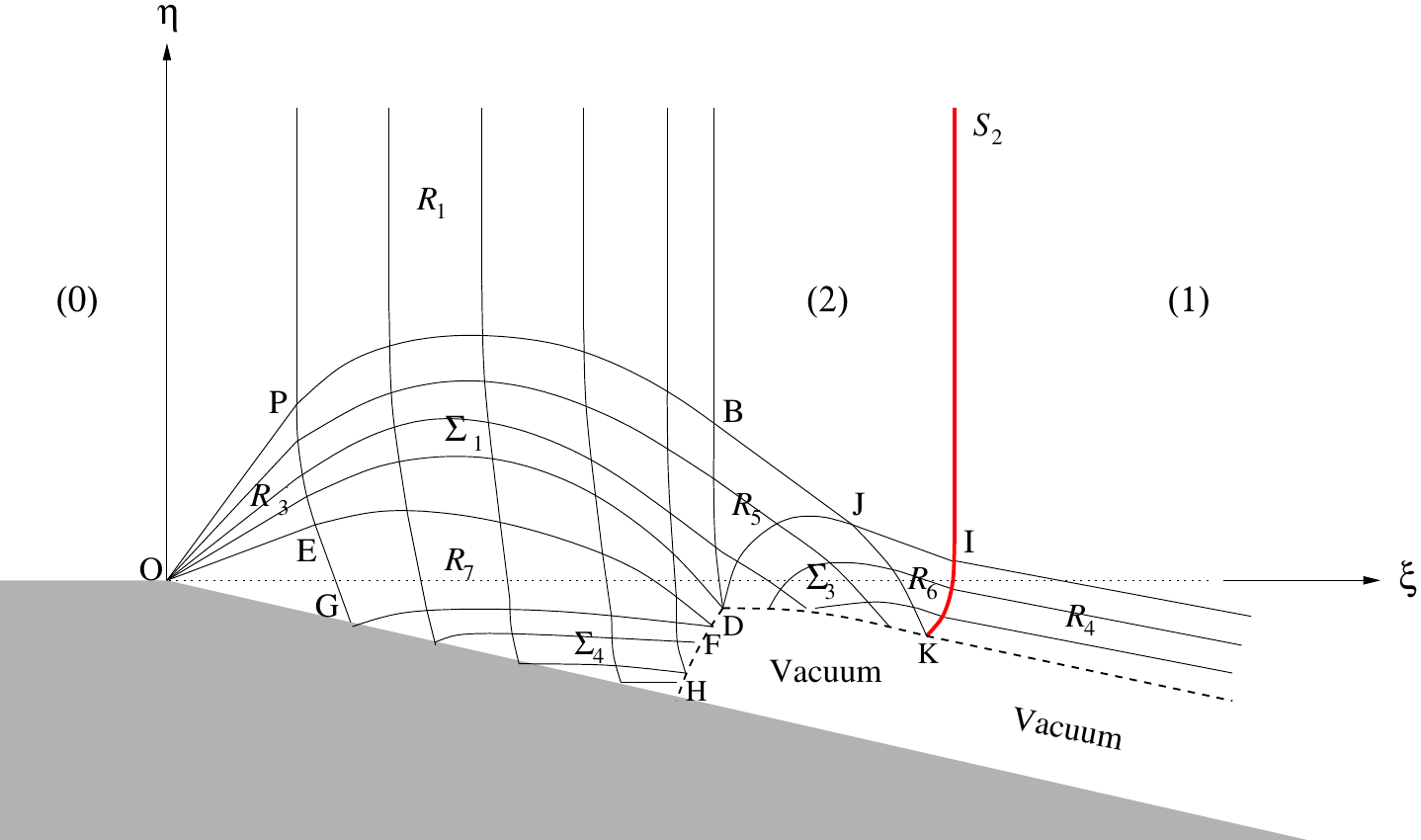}
\caption{\footnotesize Interaction of three rarefaction waves and one shock.}
\label{A2}
\end{center}
\end{figure}
%==========

Assume that 
$$
u_0<u_1<u_0+\int_{\tau_0}^{\tau_1}\sqrt{-p'(\tau)}{\rm d}\tau.
$$
  Then, the constant states $(0)$ and $(1)$ are connected by a planar rarefaction wave $R_1$ and a planar shock $S_2$  with a middle state $(u_2, 0, \tau_2)$, determined by 
$$
u_2=u_0+\int_{\tau_0}^{\tau_2}\sqrt{-p'(\tau)}{\rm d}\tau, \quad \tau_2(u_1-\xi_s)=\tau_1(u_2-\xi_s), 
\quad u_1-\xi_s=\tau_1\sqrt{-\frac{p(\tau_1)-p(\tau_2)}{\tau_1^2-\tau_2^2}},
$$
where $\xi=\xi_s$ is the location of the shock $S_2$; see Figure \ref{A2}.
As depicted in Figure \ref{A2}, the shock wave 
$S_2$ and the fan wave $R_4$ begin to interact at a point $\mathrm{I}$. The shock $S_2$  crosses the entire fan wave $R_4$ and terminates at a vacuum point $\mathrm{K}$. 
A rarefaction simple wave  $R_6$ emanates from the shock curve $\widehat{\mathrm{IK}}$. The simple waves $R_5$ and $R_6$ interact in a curvilinear triangular domain $\Sigma_3$, which is
bounded by the characteristic curves $\widetilde{\mathrm{JD}}$ and $\widetilde{\mathrm{JK}}$, together with a vacuum boundary $\wideparen{\mathrm{DK}}$.

The 2D Riemann initial-boundary value problem \eqref{Potential}, \eqref{boundary}, \eqref{initial3}, which can be viewed as a perturbation of the 1D Riemann problem, also admits other and richer flow structures. These are not discussed in the present paper.

%Let us give a brief description of the wave structure in the $(\xi, \eta)$-plane ().
%As in the last case, The fan waves $R_1$ and $R_3$ interact in the curvilinear quadrilateral domain $\Omega_1$

\section{Conclusion}
This paper presents a systematic study of the two-dimensional Riemann problem for the potential flow equation with a non-convex equation of state, with a focus on ramp-induced wave interaction configurations. We classify the planar elementary waves and the oblique corner wave, identify eight distinct types of wave combinations, and then resolve, case by case, the interactions of two fan waves, two shock–fan composite waves, two fan–shock–fan composite waves, a shock–fan composite wave with an oblique shock, a shock–fan composite wave with a fan–shock composite wave, a fan–shock–fan composite wave with a fan–shock composite wave, and a fan–shock–fan composite wave with a fan wave. Several further structural findings emerge: post-sonic shocks remain post-sonic when they cross a rarefaction simple wave of a different characteristic family; a pre-sonic shock entering an expansion region becomes transonic and may subsequently become post-sonic; and a post-sonic shock interacting with a rarefaction wave of the same family becomes transonic—a transition analysed through a free-boundary problem whose shock direction coincides with the characteristic direction at the corner and whose characteristic data fail to be $C^1$ there.
These results are of direct use in the development and verification of numerical schemes for compressible flow; in particular, the classification of the eight wave combinations provides a library of exact benchmark problems.
Finally, the Mach-like reflection regime, which lies outside the present analytical method, is a concrete target for numerical exploration: simulations can be used to conjecture the reflection configuration and to guide subsequent rigorous analysis.
%The characteristic decompositions and the invariant region for the characteristic angles provide computable indicators of hyperbolicity loss and of the onset of transonic and post-sonic shocks, which can be used to trigger adaptive mesh refinement and to switch between conservative and characteristic-based discretizations near sonic boundaries. shock-capturing, front-tracking and characteristic-based codes. The monotonicity of the hodograph transformation gives a rigorous criterion for monitoring mesh distortion and supports the construction of hodograph-plane solvers that can be inverted back to the physical plane. The singular Cauchy problems of Section 6 and the free-boundary problem of Section 10 are especially valuable as test cases for schemes that must handle shocks which are envelopes of characteristics, since the downstream flow is not $C^1$ smooth up to the shock and standard high-order reconstructions may lose accuracy or stability there. Finally, the Mach-like reflection regime, which lies outside the present analytical method, is a concrete target for numerical exploration: simulations can be used to conjecture the reflection configuration and to guide subsequent rigorous analysis.

%There remain many open questions in this paper.

\vskip 32pt

\section*{Acknowledgement}
This work was partially supported by National Natural Science Foundation of China (No. 12671262) and Natural Science Foundation of Shanghai (No. 23ZR1422100).

\vskip 32pt
%---------------------------------------------------------------------
%%%%%%%%%%%%%%%%%%%%%%%%%%%%%%%%%%%%%%%%%%%%%%%%%%%%%%%%%%%%%%%%%%%%%%%%%%%%%%%%%%%%%%%%
%\small

\end{document}